\documentclass{amsbook}

\usepackage{amsmath,amsfonts,amssymb,amsthm,enumerate,latexsym,mathrsfs,lineno}

\usepackage[dvips]{graphicx}

\newtheorem{thm}{Theorem}

\newtheorem*{thm*}{Theorem}
\newtheorem*{thm**}{Th\'eor\`eme}
\newtheorem{defn}{Definition}
\newtheorem*{defn*}{Definition}
\newtheorem{prop}{Proposition}

\newtheorem*{claim}{Claim}

\newtheorem*{subclaim}{Subclaim}
\newtheorem{lemma}{Lemma}

\newtheorem{cor}{Corollary}

\newcommand{\R}{\mathbb{R}}
\newcommand{\Q}{\mathbb{Q}}

\newcommand{\funct}[2]{#1 \longrightarrow #2}
\newcommand{\arrows}[3]{\longrightarrow {#1}^{#2}_{#3}}
\newcommand{\restrict}[2]{#1 \upharpoonright #2}

\newcommand{\ot}[1]{\textbf{#1},<^{\textbf{#1}}_{lex}}
\newcommand{\om}[1]{\textbf{#1},<^{\textbf{#1}}}
\newcommand{\m}[1]{\textbf{#1}}
\newcommand{\oc}[1]{\widetilde{\textbf{#1}},<^{\widetilde{\textbf{#1}}}}
\newcommand{\mc}[1]{\widetilde{\textbf{#1}}}

\newcommand{\otc}[1]{\widetilde{\textbf{#1}},<^{\widetilde{\textbf{#1}}}_{lex}}

\newcommand{\U}{\mathcal{U}_S}
\newcommand{\Ur}{\textbf{U}}
\newcommand{\uUr}{\textbf{B}_S}
\newcommand{\UU}{\mathcal{U}^{c<}_S}

\newcommand{\cUr}{\widehat{\textbf{B}}_{S}}
\newcommand{\iso}{\mathrm{iso}}
\newcommand{\LO}{\mathrm{LO}}
\newcommand{\cLO}{\mathrm{cLO}}
\newcommand{\nLO}{\mathrm{mLO}} 

\newcommand{\B}{\mathcal{B} _S}
\newcommand{\s}{\textbf{S}}

\newcommand{\M}{\mathcal{M}}
\newcommand{\oM}{\mathcal{M} ^{<}}
\newcommand{\nM}{\mathcal{M} ^{m<}}

\newcommand{\E}{\mathcal{H}}
\newcommand{\ES}{\mathcal{S}}

\newcommand{\dom}{\mbox{$\mathrm{dom}$}}
\newcommand{\ran}{\mbox{$\mathrm{ran}$}}

\makeindex

\begin{document}

\frontmatter

\title{Structural Ramsey theory of metric spaces and topological dynamics of isometry groups}

\author{L. Nguyen Van Th\'e}

\address{Department of Mathematics and Statistics, University of Calgary, 2500 University Drive NW, Calgary, Alberta, Canada, T2N1N4.}

\email{nguyen@math.ucalgary.ca}

\maketitle

\tableofcontents

\chapter*{Abstract}

In 2003, Kechris, Pestov and Todorcevic showed that the structure of certain separable metric spaces - called ultrahomogeneous - is closely related to the combinatorial behavior of the class of their finite metric spaces. The purpose of the present paper is to explore different aspects of this connection.

\newcounter{note1}
\setcounter{note1}{\value{footnote}} 
\def\thefootnote{\*}

\footnote{Received by the editor October 11, 2007

\textsl{2000 Mathematics Subject Classification.} Primary: 03E02. Secondary: 05C55, 05D10, 22A05, 22F05, 51F99

\textsl{Key words and phrases.} Ramsey theory, Metric geometry, Fra\"iss\'e theory, Topological groups actions, Extreme amenability, Universal minimal flows, Oscillation stability, Urysohn metric space.}

\def\thefootnote{\arabic{footnote}} 
\setcounter{footnote}{\value{note1}} 

%

\chapter*{Preface}

This book is based on work carried out between 2003 and 2006 for
the completion of a Ph.D. degree at the 'Equipe de Logique'
(University Paris 7, Denis Diderot), and expanded with recent results
obtained in 2007 thanks to a postdoctoral fellowship at the University
of Calgary. Many people made the realization of such a project
possible, but five of them had a particular influence on it. The first
one is Stevo Todorcevic, who supervised the project from the very
beginning until almost the very end. The second one is Jordi
Lopez-Abad, who also closely followed all of its multiple developments
and whose collaboration led towards the most significant result of the paper. 
The third one is Norbert Sauer, whose
difficult task consisted of verifying the integrity of the whole
construction when submitted as a dissertation. The collaboration which
followed led to the completion of the last step of the main problem of
the thesis. The fourth one is Vladimir Pestov, who also made sure that
all the arguments were robust, and whose ongoing interest has provided
unlimited motivation. The fifth one is the anonymous referee, whose rich and enthusiastic report made the publication of this work as a book possible.  

Several other interactions and discussions helped considerably, in
particular with Gilles Godefroy, Alexander Kechris, Jaroslav
Ne\v{s}et\v{r}il, Maurice Pouzet, Christian Rosendal, and all the
participants of the Set Theory seminar in Paris.

The quality of the paper was substantially improved thanks to
all of these contributions.

And last, this project would not even have existed without the fundamental work of Roland Fra\"iss\'e. This book is dedicated to his memory.

\aufm{Lionel Nguyen Van Th\'e\\
April 9, 2008}


\include{Notations}

\mainmatter





\chapter*{Introduction.}

\section{General notions and motivations.}

The backbone of the present work can be defined as the study of 'Ramsey theoretic properties of finite metric spaces in connection with the structure of separable ultrahomogeneous metric spaces'. Our original motivation comes from the recent paper \cite{KPT} of Kechris, Pestov and Todorcevic connecting various areas of mathematics respectively called 'Fra\"iss\'e theory of amalgamation classes and ultrahomogeneous structures', 'Ramsey theory', and 'topological dynamics of automorphism groups of countable structures'. More precisely, the starting point of our research is a new proof of a theorem by Pestov which provides the computation of a topological invariant attached to the surjective isometry group of a remarkable metric space. This theorem contains two main ingredients. 

The first one is the so-called \emph{universal Urysohn metric space} $\Ur$\index{$\Ur$}. This space, which appeared relatively early in the history of metric geometry (the definition of metric space is given in the thesis of M. Fr\'echet in 1906, \cite{Fre}), was constructed by Paul Urysohn in 1925. Its characterization refers to a property known today as \emph{ultrahomogeneity}: A metric space $\m{X}$ is \emph{ultrahomogeneous} when  
every isometry between finite metric subspaces extends to an isometry of $\m{X}$ onto itself. With this definition in mind, $\Ur$ can be characterized as follows: Up to isometry, it is the unique complete separable ultrahomogeneous metric space which includes all finite metric spaces. As a consequence, it can be proved that $\Ur$ is universal not only for the class of all finite metric spaces, but also for the class of all \emph{separable} metric spaces. This property is essential and is precisely the reason for which Urysohn constructed $\Ur$: Before, it was unknown whether a separable metric space could be universal for the class of all separable metric spaces. However, $\Ur$ virtually disappeared after Banach and Mazur showed that $\mathcal{C}([0,1])$ was also universal and it is only quite recently that it was brought back on the research scene, thanks in particular to the work of Kat\v{e}tov \cite{K} which was quickly followed by several results by Uspenskij \cite{Us1}, \cite{Us2} and later supported by various contributions by Vershik \cite{Ve1}, \cite{Ve2}, Gromov \cite{G}, Pestov \cite{Pe0} and Bogatyi \cite{B1}, \cite{B2}. Today, the study of the space $\Ur$ is a subject of active research and is being carried out by many different authors under many different lights, see \cite{Pr}. It is also worth mentioning that the ideas that were used to construct the space $\Ur$ contain already many of the ingredients that were used twenty-five years later to develop Fra\"iss\'e theory, a theory whose role is nowadays central in model theory and in the present paper. 

Recall now the concept of \emph{extreme amenability} from topological dynamics (Our exposition here follows the introduction of \cite{KPT}). A topological group $G$ is \emph{extremely amenable} or satisfies the \emph{fixed point on compacta property} when every continuous action of $G$ on a compact topological space $X$ admits a fixed point (ie a point $ x \in X$ such that $\forall g \in G \ \ g \cdot x = x$). Extreme amenability of topological groups naturally comes into play in topological dynamics when studying \emph{universal minimal flows}. Given a topological group $G$, a \emph{compact $G$-flow} is a compact topological space $X$ together with a continuous action of $G$ on $X$. A $G$-flow is \emph{minimal} when every orbit is dense. It is easy to show that every $G$-flow includes a minimal subflow. It is less obvious that every topological group $G$ has a \emph{universal minimal flow} $M(G)$, that is a minimal $G$-flow that can be homomorphically mapped onto any other minimal $G$-flow (For a proof, see \cite{A}). Furthermore, it turns out that $M(G)$ is uniquely determined by these properties up to isomorphism (A \emph{homomorphism} between two $G$-flows $X$ and $Y$ is a continuous map $\pi : \funct{X}{Y}$
such that for every $x \in X$ and $g \in G$, $\pi ( g \cdot x) = g \cdot \pi (x)$. An \textit{isomorphism} is a bijective homomorphism). When $G$ is locally compact but non compact, $M(G)$ is an intricate object. However, there are some non-trivial groups $G$ where $M(G)$ trivializes and those are precisely the extremely amenable ones. Pestov theorem provides such an example:

\begin{thm*}[Pestov \cite{Pe0}]

Equipped with the pointwise convergence topology, the group $\iso (\Ur)$ of isometries of $\Ur$\index{$\Ur$} onto itself is extremely amenable. 

\end{thm*}

Most of the techniques used in \cite{Pe0} come from topological group theory. However, a careful analysis of the proof together with another result of Pestov in \cite{Pe-1} according to which the automorphism group $\mathrm{Aut}(\Q , <)$ of all order-preserving bijections of the rationals is also extremely amenable allowed to isolate a substantial combinatorial core. The identification of that core is precisely the content of \cite{KPT} and shows the emergence of two major components: Fra\"iss\'e theory and structural Ramsey theory. 

Developed in the fifties by R. Fra\"iss\'e, Fra\"iss\'e theory provides a general model theoretic and combinatorial analysis of what is called today \emph{countable ultrahomogeneous structures} (Again, our exposition follows here the introduction of \cite{KPT} but a more detailed approach can be found in \cite{Fr} or \cite{H}). Let $L=\{R_i : i \in I\}$ be a fixed relational signature, and $\m{X}$ and $\m{Y}$ be two $L$-structures (that is sets $X$, $Y$ equipped with relations $R_i ^{\m{X}}$ and $R_i ^{\m{Y}}$ for each $i \in I$). An \emph{embedding} from $\m{X}$ to $\m{Y}$
is an injective map $\pi : \funct{X}{Y}$ such that for every $i \in I$ and $x_1, \ldots, x_n \in X$:

\begin{center}
$(x_1, \ldots, x_n) \in R_i^{\m{X}}$ iff $(\pi (x_1), \ldots, \pi (x_n)) \in R_i^{\m{Y}}$.
\end{center}

An \emph{isomorphism} from $\m{X}$ to $\m{Y}$ is a surjective embedding. When there is an isomorphism from
$\m{X}$ to $\m{Y}$, this is written $\m{X} \cong \m{Y}$.
Finally, $\binom{\m{Y}}{\m{X}}$ is defined as:
\[
\binom{\m{Y}}{\m{X}} = \{ \mc{X} \subset \m{Y} : \mc{X} \cong \m{X} \}
\]

When there is an embedding from an $L$-structure $\m{X}$ into another $L$-structure $\m{Y}$, we write $\m{X} \leqslant \m{Y}$.
A class $\mathcal{K}$ of $L$-structures is \emph{hereditary} when for every $L$-structure $\m{X}$ and every $\m{Y} \in \mathcal{K}$:

\begin{center}
$\m{X} \leqslant \m{Y} \rightarrow \m{X} \in \mathcal{K}$.
\end{center}

It satisfies the \emph{joint embedding property} when for every $\m{X}, \m{Y} \in \mathcal{K}$, there is $\m{Z} \in \mathcal{K}$ such that $\m{X}, \m{Y} \leqslant \m{Z}$. It satisfies the \emph{amalgamation property} when for every $\m{X}$, $\m{Y} _0$, $\m{Y} _1 \in \mathcal{K}$ and embeddings $f_0 : \funct{\m{X}}{\m{Y} _0}$ and $f_1 : \funct{\m{X}}{\m{Y}}$, there is $\m{Z} \in \mathcal{K}$ and embeddings $g_0 : \funct{\m{Y} _0}{\m{Z}}$, $g_1 : \funct{\m{Y} _1}{\m{Z}}$ such that $g_0 \circ f_0 = g_1 \circ f_1$. 



Let $\m{F}$ be an $L$-structure. Its \emph{age}, $\mathrm{Age}(\m{F})$, is the collection of all finite $L$-structures that can be embedded into $\m{F}$. $\m{F}$ is \emph{ultrahomogeneous} when every isomorphism between finite substructures of $\m{F}$ can be extended to an automorphism of $\m{F}$. Finally, a class $\mathcal{K}$ of finite $L$-structures is a \emph{Fra\"iss\'e class} when $\mathcal{K}$ contains only countably many structures up to isomorphism, is hereditary, contains structures of arbitrarily high finite size, has the joint embedding property and has the amalgamation property. With these concepts in mind, here is the fundational result in Fra\"iss\'e theory:  

\begin{thm*}[Fra\"iss\'e \cite{Fr0}]

Let $L$ be a relational signature and $\mathcal{K}$ a Fra\"iss\'e class of $L$-structures. Then there is, up to isomorphism, a unique countable ultrahomogeneous $L$-structure $\m{F}$ such that $\mathrm{Age}(\m{F}) = \mathcal{K}$. This structure $\m{F}$ is called the \emph{Fra\"iss\'e limit} of $\mathcal{K}$ and is denoted $\mathrm{Flim}(\mathcal{K})$. 

\end{thm*}  

The fundational result of Ramsey theory is older. It was proved in 1930 by F. P. Ramsey and can be stated as follows. For a set $X$ and an integer $l$, let $[X]^l$ denote the set of subsets of $X$ with $l$ elements:

\begin{thm*}[Ramsey \cite{Ra}]
\index{Ramsey!theorem}
For every $k \in \omega \smallsetminus \{ 0 \}$ and $l, m \in \omega$, there is $p \in \omega$ so that given any set $X$ with $p$ elements, if $[X]^l$ is partitioned into $k$ classes, then there is $Y \subset X$ with $m$ elements such that $[Y]^l$ lies in one of the parts of the partition. 

\end{thm*}

However, it is only in the early seventies thanks to the work of several people, among whom Erd\H{o}s, Graham, Leeb, Rothschild, Ne\v{s}et\v{r}il and R\"odl, that the essential ideas behind this theorem crystallized and expanded to structural Ramsey theory. Here are the related basic concepts: For $k,l \in \omega \smallsetminus \{ 0 \}$ and a
triple $\m{X}, \m{Y}, \m{Z}$ of $L $-structures, $\m{Z} \arrows{(\m{Y})}{\m{X}}{k,l}$ is an abbreviation for the statement: 

\begin{center} 
For any $\chi : \funct{\binom{\m{Z}}{\m{X}}}{k}$
there is $\mc{Y} \in \binom{\m{Z}}{\m{Y}}$ such
that $|\chi ''\binom{\mc{Y}}{\m{X}}| \leqslant l$. 
\end{center} 

When $l = 1$, this is simply written $\m{Z} \arrows{(\m{Y})}{\m{X}}{k}$. Now, given a class $\mathcal{K}$ of finite ordered $L$-structures, say that $\mathcal{K}$ has the \emph{Ramsey property} when for every
$\m{X}$, $\m{Y} \in \mathcal{K} $
and every $k \in \omega \smallsetminus \{ 0 \}$, there is
$\m{Z} \in \mathcal{K} $ such that:

\begin{center} 
$ \m{Z} \arrows{(\m{Y})}{\m{X}}{k}$. 
\end{center} 

The techniques developed in \cite{KPT} show the existence of several bridges between extreme amenability, universal minimal flows, Fra\"iss\'e theory and structural Ramsey theory. For example: Let $L^*$ be a relational signature with a distinguished binary relation symbol $<$. An \emph{order $L^*$-structure} is an $L ^*$-structure $\m{X}$ in which the interpretation $<^{\m{X}}$ of $<$ is a linear ordering. If $\mathcal{K} ^*$ is a class of $L^*$-structures, $\mathcal{K} ^*$ is an \emph{order class} when every element of $\mathcal{K} ^*$ is an order $L ^*$-structure. 

\begin{thm*}[Kechris-Pestov-Todorcevic \cite{KPT}]

Let $L ^* \supset \{ < \}$ be a relational signature, $\mathcal{K} ^*$ a Fra\"iss\'e order class in $L ^*$ and $(\om{F}) = \mathrm{Flim} (\mathcal{K} ^*)$. Then the following are equivalent:

\begin{enumerate}

\item $\mathrm{Aut}(\om{F})$ is extremely amenable. 

\item $\mathcal{K} ^*$ is a Ramsey class.   

\end{enumerate}

\end{thm*}

Together with several similar theorems, this result sets up a general landscape into which the combinatorial attack of extreme amenability can take place. When one is interested in the study of extreme amenability for a group of the form $\mathrm{Aut}(\mathrm{Flim} (\mathcal{K} ^*))$, this theorem can be used directly. However, the range of its applications is not restricted to this particular case. The combinatorial proof of Pestov theorem quoted previously provides a good illustration of that fact. Here are the main ideas. A first step consists in making use of the following Ramsey theorem due to Ne\v{s}et\v{r}il: 

\begin{thm*}[Ne\v{s}et\v{r}il \cite{N1}]

The class $\M _{\Q} ^<$ of all finite ordered metric spaces with rational distances has the Ramsey property. 

\end{thm*}

A second step is to refer to the general theorem. It follows that the group $G : = \mathrm{Aut}(\mathrm{Flim} (\M _{\Q} ^<))$ is extremely amenable. Finally, the last step establishes that $G$ embeds continuously and densely into $\iso (\Ur)$, and that this property is sufficient to transfer extreme amenability from $G$ to $\iso (\Ur)$. 

The success of this strategy led the authors of \cite{KPT} to ask several general questions related to metric Ramsey theory, among which stands the following one:

\

\textbf{Question:} Among the Fra\"iss\'e classes of finite ordered metric spaces, which ones have the Ramsey property? 

\

This general problem can be seen as a metric version of a well-known similar problem for finite ordered graphs out of which originated an impressive quantity of research in the seventies. In our case, it is undoubtedly the main motivation to look for classes of finite ordered metric spaces with the Ramsey property, and several examples will be exposed throughout the present paper. 

Together with Ramsey property, another combinatorial notion related to Fra\"iss\'e classes emerges from \cite{KPT}. It is called \emph{ordering property} and will also receive a particular attention in this article. 

As previously, fix a relational signature $L ^*$ with a distinguished binary relation symbol $<$ and let $L$ be the signature $L^* \smallsetminus \{ < \}$. Given an order class $\mathcal{K} ^*$ of $L ^*$-structures, let $\mathcal{K}$ be the class of $L$-structures defined by:

\begin{center}
$\mathcal{K} = \{ \m{X} : (\om{X}) \in \mathcal{K} ^*\}$.
\end{center}

Say that $\mathcal{K} ^*$ has the \textit{ordering property} when given
$\m{X} \in \mathcal{K}$, there is $\m{Y} \in \mathcal{K}$ such that given any linear orderings $<^{\m{X}}$ and $<^{\m{Y}}$ on $\m{X}$ and $\m{Y}$, if $(\om{X})$ , $(\om{Y}) \in \mathcal{K} ^*$, then $(\om{Y})$ contains an isomorphic
copy of $(\om{X})$. For us, ordering property is relevant because it leads to several interesting notions. 

The first ones are related to topological dynamics and extreme amenability: Still in \cite{KPT}, it is shown that for a certain kind of Fra\"iss\'e order class $\mathcal{K} ^*$, the ordering property provides a direct way to produce minimal $\mathrm{Aut}(\mathrm{Flim} (\mathcal{K}))$-flows. Better: When the Ramsey property and the ordering property are both satisfied, an explicit determination of the universal minimal flow of $\mathrm{Aut}(\mathrm{Flim} (\mathcal{K}))$ becomes available. This fact deserves to be mentioned as before \cite{KPT}, there were only very few cases of non extremely amenable topological groups for which the universal minimal flow was explicitly describable and known to be metrizable. This method allowed to compute the universal minimal flow of the automorphism group of several remarkable Fra\"iss\'e limits like the Rado graph $\mathcal{R}$, the Henson graphs $H_n$, the countable atomless Boolean algebra $\m{B} _{\infty}$ or the $\aleph _0$-dimensional vector space $\m{V} _F$ over a finite field $F$. 

The second kind of notion is purely combinatorial and is called \emph{Ramsey degree}: Given a class $\mathcal{K}$ of $L$-structures and $\m{X} \in \mathcal{K} $, suppose that there is $l \in \omega \smallsetminus \{ 0 \}$ such that for any $\m{Y} \in \mathcal{K}$, and any $k \in \omega \smallsetminus \{ 0 \}$, there exists $\m{Z} \in \mathcal{K} $ such that:

\begin{center}
$\m{Z} \arrows{(\m{Y})}{\m{X}}{k,l}$. 
\end{center}

The \emph{Ramsey degree of $\m{X}$ in $\mathcal{K}$} is then defined as the least such number, and it turns out that its effective computation is possible whenever $\mathcal{K}$ is coming from a $\mathcal{K} ^*$ satisfying both Ramsey and ordering property. 

In fact, the paper \cite{KPT} allows to see the determination of universal minimal flows and the computation of Ramsey degrees as the two sides of a same coin. However, the combinatorial formulation turned out to carry an undeniable advantage: That of allowing a variation which led to a new concept in topological dynamics and which may have appeared much later if not in connection with partition calculus. The variation around the notion of Ramsey degree is called \emph{big Ramsey degree}, while the new concept in topological dynamics is called \emph{oscillation stability for topological groups}. 

A possible way to introduce big Ramsey degrees is to observe that Ramsey degrees can also be introduced as follows: If $\m{F}$ denotes the Fra\"iss\'e limit of a Fra\"iss\'e class $\mathcal{K}$, $\m{X} \in \mathcal{K}$ admits a Ramsey degree in $\mathcal{K}$ when there is $l \in \omega$ such that for any $\m{Y} \in \mathcal{K}$, and any $k \in \omega \smallsetminus \{ 0 \}$,

\begin{center}
$\m{F} \arrows{(\m{Y})}{\m{X}}{k,l}$. 
\end{center} 

The big Ramsey degree corresponds to the exact same notion when this latter result remains valid when $\m{Y}$ is replaced by $\m{F}$. Its value $\mathrm{T}_{\mathcal{K}}(\m{X})$ is the least $l \in \omega$ such that 

\begin{center}
$\m{F} \arrows{(\m{F})}{\m{X}}{k,l}$. 
\end{center} 

Though not with this terminology, Ramsey degrees and big Ramsey degrees have now been studied for a long time in structural Ramsey theory. However, whereas the well-furnished collection of results in finite Ramsey theory very often leads to the determination of the Ramsey degrees, there are only few situations where the analysis of big Ramsey degrees has been completed. Here, we modestly expand those lists with theorems related to classes of finite metric spaces.

Oscillation stability for topological groups is much more recent a notion. Inspired from the Banach-theoretic concept of oscillation stability, it appears for the first time in \cite{KPT} and is more fully explained in the books \cite{Pe1} and \cite{Pe1'} by Pestov. It is important as it captures several deep ideas coming from geometric functional analysis and combinatorics. For a topological group $G$, recall that the \emph{left uniformity} $\mathcal{U} _L (G)$ is the uniformity whose basis is given by the sets of the form $V_L = \{(x,y):x^{-1}y \in V \}$ where $V$ is a neighborhood of the identity. Now, let $\widehat{G}^L$ denote the completion of $(G, \mathcal{U} _L (G))$. The structure $\widehat{G}^L$ may not be a topological group (see \cite{Di}) but is always a topological semigroup (see \cite{RoD}). For a real-valued map $f$ on a set $X$, define the \emph{oscillation} $f$ on $X$ as:

\begin{center}
$\mathrm{osc}(f) = \sup \{ \left| f(y) - f(x) \right| : x, y \in X \}$. 
\end{center} 

Now, let $G$ be a topological group, $ f : \funct{G}{\R}$ be uniformly continuous, and $\hat{f}$ be the unique extension of $f$ to $\widehat{G}^L$ by uniform continuity. Say that $f$ is \emph{oscillation stable} when for every $\varepsilon > 0$, there is a right ideal $\mathcal{I} \subset \widehat{G}^L$ such that 

\begin{center}
$\mathrm{osc}(\restrict{\hat{f}}{\mathcal{I}}) < \varepsilon$. 
\end{center}

Finally, let $G$ be a topological group acting $G$ continuously on a topological space $X$. For $f : \funct{X}{\R}$ and $x \in X$, let $f_{x} : \funct{G}{\R}$ be defined by

\begin{center}
$\forall g \in G \ \ f_{x} (g) = f(gx)$.
\end{center}

Then say that the action is \emph{oscillation stable} when for every $f : \funct{X}{\R}$ bounded and continuous and every $x \in X$, $f_{x}$ is oscillation stable whenever it is uniformly continuous.

The relationship between big Ramsey degrees and oscillation stability can be particularly well understood in the metric context. First, call a metric space $\m{X}$ \emph{indivisible} when for every strictly positive $k \in \omega$ and every $\chi : \funct{\m{X}}{k}$, there is $\mc{X} \subset \m{X}$ isometric to $\m{X}$ on which $\chi$ is constant. It should be clear that when $\m{X}$ is countable and ultrahomogeneous, indivisibility of $\m{X}$ is related to big Ramsey degrees in the Fra\"iss\'e class $\mathrm{Age}(\m{X})$ of all finite metric subspaces of $\m{X}$: The space $\m{X}$ is indivisible iff the $1$-point metric space has a big Ramsey degree in $\mathrm{Age}(\m{X})$ equal to $1$. Observe also that indivisibility can be relaxed in the following sense: If $\m{X} = (X, d^{\m{X}})$ is a metric space, $Y \subset X$ and $\varepsilon > 0$, set 

\begin{center}
$(Y) _{\varepsilon} = \{ x \in X : \exists y \in Y \ \ d ^{\m{X}} (x,y) \leqslant \varepsilon \}$.
\end{center}

Now, say that $\m{X}$ is \emph{$\varepsilon$-indivisible} when for every strictly positive $k \in \omega$, every $\chi : \funct{\m{X}}{k}$ and every $\varepsilon > 0$, there are $i < k $ and $\mc{X} \subset \m{X}$ isometric to $\m{X}$
such that 

\begin{center}
$\mc{X} \subset (\overleftarrow{\chi} \{ i \})_{\varepsilon}$. 
\end{center}

With this concept in mind, here is the promised connection: 

\begin{thm*}[Kechris-Pestov-Todorcevic \cite{KPT}, Pestov \cite{Pe1}, \cite{Pe1'}]
For a complete ultrahomogeneous metric space $\m{X}$, the following are equivalent:

\begin{enumerate}

\item When $\iso (\m{X})$ is equipped with the topology of pointwise convergence, the standard action of $\iso (\m{X})$ on $\m{X}$ is oscillation stable. 

\item For every $\varepsilon > 0$, $\m{X}$ is $\varepsilon$-indivisible. 

\end{enumerate}

\end{thm*}

A consequence of the youth of the notion of oscillation stability for topological groups is that the list of results that can be attached to it is fairly restricted. The most significant result so far in the field was obtained by Hjorth in \cite{Hj}: 
 
\begin{thm*}[Hjorth \cite{Hj}]

Let $G$ be a non-trivial Polish group. Then the action of $G$ on itself by left multiplication is not oscillation stable. 

\end{thm*}

However, some well-known results can also be interpreted in terms oscillation stability. For example, with $\mathbb{S} ^{\infty}$ denoting the unit sphere of the Hilbert space $\ell _2$ (Here, following the standard notation, $\ell_2$ denotes the Banach space of all real sequences $(x_n)_{n \in \omega}$ such that $\sum_{n=0} ^{\infty} |x_n|^2$ is finite), it should be mentioned that a problem equivalent to finding whether the standard action of $\iso (\mathbb{S} ^{\infty})$ on $\mathbb{S} ^{\infty}$ is oscillation stable motivated an impressive amount of research between the late sixties and the early nineties. It is only in 1994 that Odell and Schlumprecht finally presented a solution (cf \cite{OS}), solving the so-called \emph{distortion problem for $\ell _2$}:

\begin{thm*}[Odell-Schlumprecht \cite{OS}]
The standard action of $\iso (\mathbb{S} ^{\infty})$ on $\mathbb{S} ^{\infty}$ is not oscillation stable.
\end{thm*}  

The last part of this work is devoted to the similar problem for another metric space and called the \emph{Urysohn sphere}. From the finite Ramsey theoretic point of view, this space shares many features with the space $\mathbb{S} ^{\infty}$ and for some time, the guess was that this similarity would still hold at the level of oscillation stability. Quite surprisingly, it is not the case, and we will show in section \ref{subsection:Metric oscillation stability of S} that the solution to the distortion problem for the Urysohn sphere goes the opposite direction.

\section{Organization and presentation of the results.}

Chapter 1 is devoted to the presentation of several Fra\"iss\'e classes of finite metric spaces whose role is central in our work. 

One of the most important ones is the class $\M _{\Q}$ of finite metric spaces with rational distances. Its \emph{Urysohn space} (the name given to the Fra\"iss\'e limit in the metric context) is a countable ultrahomogeneous metric space denoted $\Ur _{\Q}$ and called the \emph{rational Urysohn space}. Several variations of $\M _{\Q}$ are also of interest for us: The class $\M _{\Q \cap ]0,1]}$ of finite metric spaces with distances in $\Q \cap ]0,1]$, whose Urysohn space is the \emph{rational Urysohn sphere} $\s _{\Q}$. The class $\M _{\omega}$ of finite metric spaces with distances in $\omega$, leading to the \emph{integral Urysohn space} $\Ur _{\omega}$. And finally the classes $\M _{\omega \cap ]0,m]}$ of finite metric spaces with distances in $\{1,\ldots , m \}$ where $m$ is a strictly positive integer, giving raise to bounded versions of $\Ur _{\omega}$ denoted $\Ur _m$.

Two other kinds of classes appear prominently in our work. The first kind consists of the classes of the form $\U$ of finite ultrametric spaces with distances in a prescribed countable subset $S$ of $]0, + \infty[$. Every $\U$ leads to a so-called \emph{ultrametric Urysohn space} denoted $\uUr$ and which, unlike most of the Urysohn spaces, can be described very explicitly. The second kind consists of the classes $\M _S$ of finite metric spaces with distances in $S$ where $S \subset ]0, + \infty[$ is countable and satisfies the so-called \emph{$4$-values condition}, a condition discovered by Delhomm\'e, Laflamme, Pouzet and Sauer in \cite{DLPS} and which characterizes those subsets $S \subset ]0, +\infty[$ for which the class $\M _S$ of all finite metric spaces with distances in $S$ has the amalgamation property. Every $\M _S$ leads to a space denoted $\Ur _S$ which can also sometimes be described explicitly when $S$ is finite and not too complicated. 

Finally, we finish our list with two classes of finite Euclidean metric spaces, namely the class $\E _S$ of all finite affinely independent metric subspaces of the Hilbert space $\ell _2$ with distances in $S$ where $S$ is a countable dense subset of $]0,+\infty[$, and the class $\ES _S$ of all finite metric spaces $\m{X}$ with distances in $S$ which embed isometrically into the unit sphere $\mathbb{S} ^{\infty}$ of $\ell _2$ with the property that $\{ 0_{\ell _2} \} \cup \m{X}$ is affinely independent ($S$ still being a countable dense subset of $]0,+\infty[$). The corresponding Urysohn spaces are countable metric subspaces of $\ell _2$ and $\mathbb{S} ^{\infty}$ respectively. Unfortunately, because of the combinatorial difficulties which arise when trying to work with those objects, they will only appear anecdotically in our work. 
 
Once those Fra\"iss\'e classes and their related Urysohn spaces are presented, we turn our attention to the interplay between complete separable ultrahomogeneous metric spaces and Urysohn spaces. We start with considerations around the following questions:

\begin{enumerate}  

\item Is the completion of a Urysohn space still ultrahomogeneous ? 

\item Does every complete separable ultrahomogeneous metric space appear as the completion of a Urysohn space ? 

\end{enumerate}

The answer to (1) is negative and is provided by an example taken from an article of Bogatyi \cite{B2}. On the other hand, the answer to (2) turns out to be positive and provides our first substantial theorem, see Theorem \ref{thm:countable dense ultrahomogeneous}:

\begin{thm*}

Every complete separable ultrahomogeneous metric space $\m{Y}$ includes a countable ultrahomogeneous dense metric subspace. 

\end{thm*}  

We then turn to the description of the completion of the Urysohn spaces presented previously. It is the opportunity to present several remarkable spaces, among which the original Urysohn space $\Ur$\index{$\Ur$} (as the completion of $\Ur _{\Q}$), the Urysohn sphere $\s$ (as the completion of $\s _{\Q}$), the Baire space $\mathcal{N}$ (and more generally all the complete separable ultrahomogeneous ultrametric spaces), as well as the Hilbert space $\ell _2$ and its unit sphere $\mathbb{S} ^{\infty}$.  

\

Chapter 2 is devoted to finite metric Ramsey calculus and, as already stressed in the first section of this introduction, is mainly concerned about new proofs along the line of the combinatorial proof of Pestov theorem via Ne\v{s}et\v{r}il theorem and the theory developed in \cite{KPT}. For completeness, we start with a presentation of Ne\v{s}et\v{r}il theorem leading to the following result. For $S \subset ]0, + \infty[$, let $\M ^< _S$\index{$\M ^< _S$} denote the class of all finite ordered metric spaces with distances in $S$. Then (see Theorem \ref{thm:variation RP for MM}): 

\begin{thm*}[Ne\v{s}et\v{r}il \cite{N1}]

Let $T \subset ]0,+ \infty[$ be closed under sums and $S$ be an initial segment of $T$. Then $\M ^< _S$\index{$\M ^< _S$} has the Ramsey property. 

\end{thm*}

Then, we show that similar results hold for other classes of finite ordered metric spaces. The first class is built on the class $\U$: Let $\m{X}$ be an ultrametric space. Call a linear ordering $<$ on $\m{X}$ \textit{convex} when all the metric balls of $\m{X}$ are $<$-convex. For $S \subset ]0, + \infty [$, let $\UU$ denote the class of all finite convexly ordered ultrametric spaces with distances in $S$. Then (see Theorem \ref{thm:RP for UU}):

\begin{thm*}
Let $S \subset ]0, + \infty [$. Then $\UU$ has the Ramsey property.
\end{thm*}

The second kind of class where we can prove Ramsey property is based on the classes $\M _S$. Let $\mathcal{K}$ be a class of metric spaces. Call a distance $s \in ]0 , + \infty[$ \emph{critical for $\mathcal{K}$} when for every $\m{X} \in \mathcal{K}$, one defines an equivalence relation $\approx$ on $\m{X}$ by setting: 

\begin{center}
$\forall x, y \in \m{X} \ \ x \approx y \leftrightarrow d^{\m{X}}(x,y) \leqslant s$. 
\end{center}

The relation $\approx$ is then called a \emph{metric equivalence relation} on $\m{X}$. Now, call a linear ordering $<$ on $\m{X} \in \mathcal{K}$ \emph{metric} if given any metric equivalence relation $\approx$ on $\m{X}$, the $\approx$-equivalence classes are $<$-convex. Given $S \subset ]0, + \infty[$, let $\nM _S$ denote the class of all finite metrically ordered metric spaces with distances in $S$. Then (see Theorem \ref{thm:RP for nM_S}):

\begin{thm*}
Let $S$ be finite subset of $]0, + \infty [$ of size $|S| \leqslant 3$ and satisfying the $4$-values condition. Then $\nM _S$ has the Ramsey property. 
\end{thm*}

After the study of Ramsey property, we turn to ordering property. For $S$ initial segment of $T \subset ]0,+ \infty[$, $T$ closed under sums, ordering property for $\M ^< _S$\index{$\M ^< _S$} can be proved via a probabilistic argument, see \cite{N2}. We present here a proof based on Ramsey property (see Theorem \ref{thm:Ordering Property for metric spaces.}):

\begin{thm*}
Let $T \subset ]0,+ \infty[$ be closed under sums and $S$ be an initial segment of $T$. Then $\M ^< _S$\index{$\M ^< _S$} has the ordering property. 
\end{thm*}

We then follow with the ordering property for $\UU$ and for $\nM _S$, see Theorems \ref{thm:OP for UU} and \ref{thm:OP for nM_S}:

\begin{thm*}

The class $\UU$ has the ordering property. 

\end{thm*}

\begin{thm*}

Let $S$ be finite subset of $]0, + \infty [$ of size $|S| \leqslant 3$ and satisfying the $4$-values condition. Then $\nM _S$ has the ordering property. 
\end{thm*}

As mentioned in the first section of the introduction, Ramsey property together with ordering property allow the computation of Ramsey degrees. In the present situation, we are consequently able to compute the exact value of the Ramsey degrees in the classes $\M _S$ when $S$ is an initial segment of $T$ with $T \subset ]0,+ \infty[$ is closed under sums (see Theorem \ref{thm:variation Rd in M}), $\U$ (see Theorem \ref{thm:Rd in U}) and $\M _S$ where $S$ is a finite subset of $]0, + \infty [$ of size $|S| \leqslant 3$ and satisfying the $4$-values condition (see Theorem \ref{thm:Rd in M_S}). 

Finally, we turn to applications in topological dynamics. We first present the proof of Pestov theorem about the extreme amenability of $\iso (\Ur)$ and then follow with several results about extreme amenability and universal minimal flows. For example (see Theorem \ref{thm:UMF for iso(uUr)}):  

\begin{thm*}
The universal minimal flow of $\mathrm{iso}(\uUr)$ is the set $\cLO (\uUr)$ of
convex linear orderings on $\uUr$ together with the action $\funct{\iso (\uUr) \times \cLO (\uUr)}{\cLO (\uUr)}$, $(g,<)
\longmapsto <^g$ defined
by $x <^g y$ iff $g^{-1}(x) < g^{-1}(y)$.
\end{thm*}

On the other hand, recalling that $\mathcal{N}$ denotes the Baire space (see Theorem \ref{cor:UMF for iso(cUr)}):

\begin{thm*}

The universal minimal flow of $\mathrm{iso}(\mathcal{N})$ is the set $\cLO
(\mathcal{N})$ of convex linear orderings on $\mathcal{N}$ together with the action
$\funct{\iso (\mathcal{N}) \times \cLO (\mathcal{N})}{\cLO (\mathcal{N})}$, $(g,<)
\longmapsto <^g$ defined
by $x <^g y$ iff $g^{-1}(x) < g^{-1}(y)$.
\end{thm*}

As a last example (Theorem \ref{thm:M(iso(Ur_S,<))}):

\begin{thm*}
Let $S$ be finite subset of $]0, + \infty [$ of size $|S| \leqslant 3$ and satisfying the $4$-values condition. Then the universal minimal flow of $\iso (\Ur _S)$ is the set $\nLO (\Ur _S)$ of metric linear orderings on $\Ur _S$ together with the action $\funct{\iso (\Ur _S) \times \nLO (\Ur _S)}{\nLO (\Ur _S)}$, $(g,<)
\longmapsto <^g$ defined
by $x <^g y$ iff $g^{-1}(x) < g^{-1}(y)$.
\end{thm*}

In particular, the underlying spaces of all those universal minimal flow are metrizable.

We finish Chapter 2 with several open questions concerning Ramsey property for the classes $\M _S$ as well as a possible connection between Euclidean Ramsey theory and a theorem by Gromov and Milman.

\

Chapter 3 is devoted to infinite metric Ramsey calculus. We start with a short section on big Ramsey degrees. \emph{Short} cannot be removed from the previous sentence because in most of the cases, the determination of big Ramsey degrees turns out to be too difficult for us to complete. Still, there is one case where we manage to provide a full analysis (see Theorem \ref{thm:Big Ramsey degrees for U, S finite}): 

\begin{thm*}

Let $S$ be a finite subset of $]0, + \infty[$. Then every element of $\U$ has a big Ramsey degree in $\U$.

\end{thm*} 

In fact, we are even able to compute exact the value of the big Ramsey degree. This has to be compared with (see Theorem \ref{thm:no Big Ramsey degrees for U, S infinite}):

\begin{thm*}

Let $S$ be an infinite countable subset of $]0, + \infty[$ and let $\m{X}$ be in $\U$ such that $|\m{X}| \geqslant 2$. Then $\m{X}$ does not have a big Ramsey degree in $\U$.
  
\end{thm*} 

We follow with a section on the indivisibility properties of the Urysohn spaces. Recall that a metric space $\m{X}$ is indivisible when for every strictly positive $k \in \omega$ and every $\chi : \funct{\m{X}}{k}$, there is $\mc{X} \subset \m{X}$ isometric to $\m{X}$ on which $\chi$ is constant. After the presentation of several general results taken from \cite{DLPS}, we provide the details of the proof of the following theorem (see Theorem \ref{thm:s_Q divisible}):

\begin{thm*}[Delhomm\'e-Laflamme-Pouzet-Sauer \cite{DLPS}]
The space $\s _{\Q}$ is not indivisible. 
\end{thm*}

Then, we turn to the study of indivisiblity of simpler Urysohn spaces, namely the spaces $\Ur _m$. We first present the most elementary cases where general theorems such as Milliken theorem or Sauer theorem can be applied. Using techniques inspired from the general theory of indivisibility of countable structures with the so-called \emph{free amalgamation}, we then prove the general case (see Theorem \ref{thm:U_m indiv}):

\begin{thm*}[NVT-Sauer] 

Let $m \in \omega$, $m \geqslant 1$. Then $\Ur _m$ is indivisible. 

\end{thm*} 

We follow with the indivisibility of the ultrametric Urysohn spaces. As for big Ramsey degrees, these cases turn out to be accessible and lead to the following theorem (proved independently of Delhomm\'e, Laflamme, Pouzet and Sauer in \cite{DLPS}), see section \ref{subsection:Indivisibility of Urysohn ultrametric spaces} (\ref{thm:char uUr indivisible, S dually well ordered}):

\begin{thm*}

Let $\m{X}$ be a countable ultrahomogeneous ultrametric space with distance set $S \subset ]0, + \infty[$. Then $\m{X}$ is indivisible iff $\m{X}=\uUr$ and the reverse linear ordering $>$ on $\R$ induces a well-ordering on $S$. 

\end{thm*}

In fact, ultrametric Urysohn spaces behave so nicely that we are even able to establish the following refinement (see Theorem \ref{thm:selectivity for uUr}): 

\begin{thm*}

Let $S$ be an infinite countable subset of $]0, + \infty[$ such that the
reverse linear ordering $>$ on $\R$ induces a well-ordering on $S$. Then given any map $f : \funct{\uUr}{\omega}$, there is an isometric copy $X$ of $\uUr$ inside $\uUr$ such that $f$ is continuous or injective on $X$. 

\end{thm*}

After ultrametric Urysohn spaces, we finish the section on indivisibility with the study of the spaces $\Ur _S$ when  $S$ is finite and satisfies the $4$-values condition. Our proof only covers the case $|S| \leqslant 4$ but even so turns out to be long and tedious (see Theorem \ref{thm:U_S indiv}): 

\begin{thm*}

Let $S$ be finite subset of $]0, + \infty [$ of size $|S| \leqslant 4$ and satisfying the $4$-values condition. 
Then $\Ur _S$ is indivisible.

\end{thm*}

After indivisibility, we turn to oscillation stability. There are some cases where it is easy to study. For example, unsurprisingly in view of the previous results, complete separable ultrahomogeneous ultrametric spaces enter this category (see Theorem \ref{thm:ultrametric oscillation stability}).

\begin{thm*}

Let $\m{X}$ be a complete separable ultrahomogeneous ultrametric space. The following are equivalent: 

\begin{enumerate}
	\item[i)] The standard action of $\iso (\m{X})$ on $\m{X}$ is oscillation stable.
	\item[ii)] $\m{X}=\cUr$ for some $S \subset ]0, + \infty[$ finite or countable on which the reverse linear ordering $>$ on $\R$ induces a well-ordering. 
\end{enumerate}

\end{thm*}

However, in most of the cases, the study of oscillation stability seems to be hard to complete. The case of $\mathbb{S} ^{\infty}$ was already presented in the previous section of this introduction. The last part of this work is devoted to the somehow similar problem for the Urysohn sphere $\s$, namely: Is the standard action of $\iso (\s)$ on $\s$ oscillation stable? We show that the answer is positive (Theorem \ref{thm:s mos}):

\begin{thm*}
The standard action of $\iso (\s)$ on $\s$ is oscillation stable.
\end{thm*}

This result also allows to reach interesting metric partition properties for two remarkable Banach spaces. For example (Theorem \ref{thm:C[0,1] approx indiv}): 

\begin{thm*}

For every $\varepsilon > 0$, the unit sphere of $\mathcal{C}([0,1])$ is $\varepsilon$-indivisible. 

\end{thm*}

On the other hand, Holmes proved in \cite{H} there is a Banach space $\left\langle \Ur \right\rangle$ such that for every isometry $i : \funct{\Ur}{\m{Y}}$ of the Urysohn space $\Ur$ into a Banach space $\m{Y}$ such that $0_{\m{Y}}$ is in the range of $i$, there is an isometric isomorphism between $\left\langle \Ur \right\rangle$ and the closed linear span of $i''\Ur$ in $\m{Y}$. Very little is known about the space $\left\langle \Ur \right\rangle$, but in the present case, Theorem \ref{thm:s mos} directly leads to (see Theorem \ref{thm:<U> approx indiv}):  

\begin{thm*}

For every $\varepsilon > 0$, the unit sphere of the Holmes space is $\varepsilon$-indivisible. 

\end{thm*}

We then close chapter 3 and this work with questions about big Ramsey degrees in the classes $\M _S$, indivisibility of the spaces $\Ur _S$ and the relationship between the oscillation stability problems for the spheres $\mathbb{S}^{\infty}$ and $\s$.

\

Throughout all the present work, we refer as accurately as possible to the original authors and publications for all the results which are not ours. The new results related to finite Ramsey calculus of finite ultrametric spaces and topological dynamics of their Urysohn spaces (Chapter 2) are taken from \cite{NVT1}. Those related to big Ramsey degrees and indivisibility of ultrametric spaces (Chapter 3) are taken from \cite{NVT2}. Finally, those related to the oscillation stability problem for the Urysohn sphere (Chapter 3) were obtained in collaboration with Jordi Lopez-Abad on the one hand and Norbert Sauer on the other hand. They respectively correspond to the papers \cite{LANVT} (volume \cite{Pr} of Topology and its Applications devoted to the universal Urysohn space) and \cite{NVTS}.






\chapter{Fra\"iss\'e classes of finite metric spaces and Urysohn spaces.}

\section{Fundamentals of Fra\"iss\'e theory.}

\label{section:fundamentalsFraisse}

In this section, we introduce the basic concepts related to model theory and Fra\"iss\'e theory. We follow \cite{KPT} but a more detailed approach can be found in \cite{Fr} or \cite{H}. Let $L=\{R_i :~i \in I\}$ be a fixed relational signature (ie a list of symbols to be interpreted later as relations). Let $\m{X}$ and $\m{Y}$ be two $L$-structures (that is, non empty sets $X$, $Y$ equipped with relations $R_i ^{\m{X}}$ and $R_i ^{\m{Y}}$ for each $i \in I$). An \emph{embedding}\index{embedding} from $\m{X}$ to $\m{Y}$
is an injective map $\pi : \funct{X}{Y}$ such that for every $i \in I$ and $x_1, \ldots, x_n \in X$:

\begin{center}
$(x_1, \ldots, x_n) \in R_i^{\m{X}}$ iff $(\pi (x_1), \ldots, \pi (x_n)) \in R_i^{\m{Y}}$.
\end{center}

An \emph{isomorphism}\index{isomorphism} from $\m{X}$ to $\m{Y}$ is a surjective embedding while an \emph{automorphism}\index{automorphism} of $\m{X}$ is an isomorphism from $\m{X}$ onto itself. Of course,
$\m{X}$ and $\m{Y}$ are \emph{isomorphic}\index{isomorphic} when there is an isomorphism from
$\m{X}$ to $\m{Y}$. This is written $\m{X} \cong \m{Y}$.
Finally, $\binom{\m{Y}}{\m{X}}$\index{$\binom{\m{Y}}{\m{X}}$} is defined as:
\[
\binom{\m{Y}}{\m{X}} = \{ \mc{X} \subset \m{Y} : \mc{X} \cong \m{X} \}.
\]

When there is an embedding from an $L$-structure $\m{X}$ into another $L$-structure $\m{Y}$, we write $\m{X} \leqslant \m{Y}$.
A class $\mathcal{K}$ of $L$-structures is \emph{hereditary}\index{hereditary} when for every $L$-structure $\m{X}$ and every $\m{Y} \in \mathcal{K}$:

\begin{center}
$\m{X} \leqslant \m{Y} \rightarrow \m{X} \in \mathcal{K}$.
\end{center}

It satisfies the \emph{joint embedding property}\index{joint embedding property} when for every $\m{X}, \m{Y} \in \mathcal{K}$, there is $\m{Z} \in \mathcal{K}$ such that $\m{X}, \m{Y} \leqslant \m{Z}$. It satisfies the \emph{amalgamation property}\index{amalgamation!amalgamation property} (or is an \emph{amalgamation class}\index{amalgamation!amalgamation class}) when for every $\m{X}$, $\m{Y} _0$, $\m{Y} _1 \in \mathcal{K}$ and embeddings $f_0 : \funct{\m{X}}{\m{Y} _0}$ and $f_1 : \funct{\m{X}}{\m{Y}}$, there is $\m{Z} \in \mathcal{K}$ and embeddings $g_0 : \funct{\m{Y} _0}{\m{Z}}$, $g_1 : \funct{\m{Y} _1}{\m{Z}}$ such that $g_0 \circ f_0 = g_1 \circ f_1$. Finally, $\mathcal{K}$ has the \emph{strong amalgamation property}\index{amalgamation!strong amalgamation property} (or is a \emph{strong amalgamation class}\index{amalgamation!strong amalgamation class}) when one can also fulfill the requirement: 

\begin{center}
$g_0 '' f_0 '' X = g_0 '' Y_0 \cap g_1 '' Y_1 (= g_0 '' f_0 '' X)$.
\end{center}

A structure $\m{F}$ is \emph{ultrahomogeneous}\index{ultrahomogeneous!ultrahomogeneous structure} when every isomorphism between finite substructures of $\m{F}$ can be extended to an automorphism of $\m{F}$. Fra\"iss\'e theory provides a general analysis of countable ultrahomogeneous structures.

Let $\m{F}$ be an $L$-structure. The \emph{age of} $\m{F}$\index{age}, denoted $\mathrm{Age}(\m{F})$\index{$\mathrm{Age}(\m{F})$}, is the collection of all finite $L$-structures that can be embedded into $\m{F}$. Observe also that if $\m{F}$ is countable, then $\mathrm{Age}(\m{F})$ contains only countably many isomorphism types. Abusing language, we will say that $\mathrm{Age}(\m{F})$ is countable. Similarly, a class $\mathcal{K}$ of $L$-structures will be said to be countable if it contains only countably many isomorphism types. 

A class $\mathcal{K}$ of finite $L$-structures is a \emph{Fra\"iss\'e class}\index{Fra\"iss\'e!Fra\"iss\'e class} when $\mathcal{K}$ is countable, hereditary, contains structures of arbitrarily high finite size, has the joint embedding property and the has the amalgamation property (Note that the joint embedding property is not a trivial subcase of the amalgamation property with $\m{X}=\emptyset$ as technically, an $L$-structure is not allowed to be empty). 

It should be clear that if $\m{F}$ is a countable ultrahomogeneous $L$-structure, then $\mathrm{Age}(\m{F})$ is a Fra\"iss\'e class. The following theorem, due to Fra\"iss\'e, establishes a kind of converse: 

\begin{thm}[Fra\"iss\'e \cite{Fr0}]

\label{thm:Fraisse}
\index{Fra\"iss\'e!theorem}

Let $L$ be a relational signature and $\mathcal{K}$ a Fra\"iss\'e class of $L$-structures. Then there is, up to isomorphism, a unique countable ultrahomogeneous $L$-structure $\m{F}$ such that $\mathrm{Age}(\m{F}) = \mathcal{K}$. $\m{F}$ is called the \emph{Fra\"iss\'e limit}\index{Fra\"iss\'e!Fra\"iss\'e limit} of $\mathcal{K}$ and denoted $\mathrm{Flim}(\mathcal{K})$\index{$\mathrm{Flim}(\mathcal{K})$}. 

\end{thm}  

We do not enter the details of the proof here but let us simply mention that uniqueness of the Fra\"iss\'e limit is due to the following fact:

\begin{prop}
Let $\m{F}$ be a countable $L$-structure. Then $\m{F}$ is ultrahomogeneous iff for every finite substructures $\m{X}, \m{Y}$ of $\m{F}$ with $|\m{Y}|=|\m{X}|+1$, every embedding $\m{X} \longrightarrow \m{F}$ can be extended to an embedding $\m{Y} \longrightarrow \m{F}$. 
\end{prop}

Let us now illustrate how these concepts translate in the context of the central objects of this paper: Metric spaces. There are several ways to see a metric space $\m{X} = (X , d^{\m{X}})$\index{$d^{\m{X}}$} as a relational structure. For example, one may consider a binary relation symbol $R _{\delta}$ for every $\delta$ in $\Q \cap ]0, + \infty[$ and set 

\begin{center}
$(x, y) \in R ^{\m{X}}_{\delta} \leftrightarrow d^{\m{X}}(x,y) < \delta$. 
\end{center}

One may also allow $\delta$ to range over $]0, + \infty[$, and define: 

\begin{center}
$(x, y) \in R ^{\m{X}}_{\delta} \leftrightarrow d^{\m{X}}(x,y) = \delta$. 
\end{center}

This latter approach has the disadvantage of requiring the signature to be uncountable if uncountably many distances appear in the metric space we are working with. This is a real issue as Fra\"iss\'e theory really deals with countable signatures, but in the present case, the instances where Fra\"iss\'e theory will be needed will involve only countably many distances so the second way of encoding the distance map by relations will not cause any problem.  

With these facts in mind, substructures in the context of metric spaces really correspond to \emph{metric subspaces} and embeddings are really \emph{isometric embeddings}. It follows that if $\m{X}, \m{Y}$ are metric spaces, then $\binom{\m{Y}}{\m{X}}$ is the set of all isometric copies of $\m{X}$ inside $\m{Y}$. 
 
Other kinds of relational structures will come into play, namely, ordered metric spaces (structures of the form $(\om{X}) = (X , d^{\m{X}}, <^{\m{X}})$ where $\m{X}$ is a metric space and $<^{\m{X}}$ is a linear ordering on $X$), graphs (structures $\m{G}$ in the language $\{R_1\}$ where $R ^{\m{G}} _1$ is binary, symmetric and irreflexive), edge-labelled graphs (structures $\m{G}$ in the language $\{ R_{\delta} : \delta \in ]0, + \infty[ \}$ where each $R _{\delta} ^{\m{G}}$ is binary symmetric and irreflexive), ordered edge-labelled graphs\ldots However, the reader should be aware that in many cases, we will not be too cautious with the notational aspect. In particular, when dealing with a metric space $\m{X}$, we will often use the same notation to denote both the metric space and its underlying set, and we will \emph{almost never} use the relational notation to refer to the distance. Similarly, when dealing with an edge-labelled graph $\m{G}$, we will always work with the \emph{labelling map} $\lambda ^{\m{G}}$\index{$\lambda ^{\m{G}}$} defined on the set $ \bigcup _{\delta \in ]0, + \infty[} R ^{\m{G}} _{\delta}$ by

\begin{center}
$\lambda ^{\m{G}}(x,y) = \delta \leftrightarrow (x,y) \in R ^{\m{G}} _{\delta}$.
\end{center}

A class $\mathcal{K}$ of metric spaces is hereditary when it is closed under isometries and metric subspaces. Next, suppose we want to show that a class $\mathcal{K}$ of finite metric spaces has the strong amalgamation property. We take $\m{X}$, $\m{Y} _0$, $\m{Y} _1 \in \mathcal{K}$, isometric embeddings $f_0 : \funct{\m{X}}{\m{Y} _0}$ and $f_1 : \funct{\m{X}}{\m{Y}}$ and we wish to find $\m{Z} \in \mathcal{K}$ and embeddings $g_0 : \funct{\m{Y} _0}{\m{Z}}$, $g_1 : \funct{\m{Y} _1}{\m{Z}}$ such that $g_0 \circ f_0 = g_1 \circ f_1$. Thanks to the previous comments, we may assume without loss of generality that $\m{X}$ is really a \emph{metric subspace} both of $\m{Y} _0$ and $\m{Y} _1$, and that $\m{Y} _0 \cap \m{Y} _1 = \m{X}$. Hence, the metrics $d^{\m{Y} _0}$ and $d^{\m{Y} _1}$ agree on $X$ and are equal to $d^{\m{X}}$ on $X$. So we will be done if we can prove that $d^{\m{Y} _0} \cup d^{\m{Y} _1} $ can be extended to a metric on $Y _0 \cup Y _1$. As we will see later, the most convenient way to proceed will strongly depend on how $\mathcal{K}$ is defined. 

Let us now examine the meaning of ultrahomogeneity. A metric space $\m{X}$ is ultrahomogeneous\index{ultrahomogeneous!ultrahomogeneous metric space} when any isometry between two finite subspaces can be extended to an isometry of $\m{X}$ onto itself. Throughout this paper, the set of all isometries of a metric space $\m{X}$ onto itself is denoted $\iso (\m{X})$\index{$\iso (\m{X})$}.

In the metric setting, Fra\"iss\'e theorem consequently states:

\begin{thm}[Fra\"iss\'e theorem for metric spaces.]

\index{Fra\"iss\'e!theorem for metric spaces}

Let $\mathcal{K}$ be a Fra\"iss\'e class of metric spaces. Then there is, up to isometry, a unique countable ultrahomogeneous metric space $\m{X}$ whose class of finite metric subspaces is exactly $\mathcal{K}$. This space will be called the \emph{Urysohn space}\index{Urysohn!Urysohn space associated to a class of finite metric spaces} associated to $\mathcal{K}$. 
\end{thm}

As we mentioned when stating the general form of Fra\"iss\'e theorem, uniqueness of the Urysohn space can be shown via a back-and-forth argument after having restated ultrahomogeneity in terms of a certain extension property. The purpose of what follows is to state this extension property, and to show that it is indeed equivalent to ultrahomogeneity. We start with the following important concept:

\begin{defn}

If $\m{X} = (X, d^{\m{X}})$ is a metric space, a map $f : \funct{X}{\R}$ is \emph{Kat\v{e}tov over $\m{X}$}\index{Kat\v{e}tov map} when: 

\begin{center}
$\forall x, y \in X, \ \ |f(x) - f(y)| \leqslant d^{\m{X}} (x,y) \leqslant f(x) + f(y)$. 
\end{center}

If $E(\m{X})$\index{$E(\m{X})$} denotes the set of all Kat\v{e}tov maps over $\m{X}$, $\m{X} \subset \m{Y}$ and $f \in E(\m{X})$, a point $y \in \m{Y}$ \emph{realizes}\index{Kat\v{e}tov map!realization of a Kat\v{e}tov map} $f$ over $\m{X}$ when:

\begin{center}
$\forall x \in \m{X}, \ \ d^{\m{Y}}(x,y) = f(x)$.
\end{center} 

\end{defn}

Equivalently, if $f \in E(\m{X})$, then $f$ can be thought as a potential new point that can be added to the space $\m{X}$. Indeed, if $f$ does not vanish on $\m{X}$, then one can extend the metric $d^{\m{X}}$ on $X \cup \{ f \}$ by defining, for every $x, y$ in X, $ \widehat{d^{\m{X}}} (x, f) = f(x)$ and $\widehat{d^{\m{X}}} (x, y) = d^{\m{X}} (x, y)$. It is not the case when $f$ vanishes at some point $x$ but then, one can check that for every $y \in \m{X}$, $f(y) = d^{\m{X}}(x,y)$ and so $f$ can be identified with $x$. In any case, the corresponding metric space will be denoted $\m{X} \cup \{ f \}$\index{$\m{X} \cup \{ f \}$}.

\begin{prop}

\label{prop:extension}

Let $\m{Y}$ be a countable metric space. Then $\m{Y}$ is ultrahomogeneous iff for every finite subspace $\m{X} \subset \m{Y}$ and every Kat\v{e}tov map $f$ over $\m{X}$, if $\m{X} \cup \{ f \}$ embeds into $\m{Y}$, then there is $y \in \m{Y}$ realizing $f$ over $\m{X}$. The same result holds when $\m{Y}$ is complete separable. 

\end{prop}

\begin{proof}
Assume that $\m{Y}$ is countable (resp. complete separable) and ultrahomogeneous. Consider an embedding $\varphi : \funct{\m{X} \cup \{ f \}}{\m{Y}}$. By ultrahomogeneity of $\m{Y}$, there is an isometry $\psi$ of $\m{Y}$ onto itself such that: 

\begin{center}
$ \forall x \in \m{X}, \ \ \psi (x) = \varphi (x)$. 
\end{center}

Then, $\psi ^{-1} (\varphi (f)) \in \m{Y}$ realizes $f$ over $\m{X}$.  

For the converse, suppose first that $\m{Y}$ is countable. Assume that $\{ x_0 ,\ldots, x_n \}$ and $\{z_0 ,\ldots, z_n \}$ are isometric finite subspaces of $\m{Y}$ and that $\varphi : x_k \mapsto z_k$ is an isometry. We wish to extend $\varphi$ to an isometry of $\m{Y}$ onto itself. We do that thanks to a back and forth method. First, extend $\{ x_0 ,\ldots, x_n \}$ and $\{z_0 ,\ldots, z_n \}$ so that $\{ x_k : k \in \omega \} = \{ z_k : k \in \omega \} = \m{Y}$. For $k\leqslant n$, let $\sigma (k) = \tau (k) = k$. Then, set $\sigma (n+1) = n+1$. Consider the map $f _{n+1}$ defined on $\{\varphi ( x_{\sigma (k)}) : k < n+1 \}$ by:

\begin{center}
$\forall k < n+1, \ \ f_{n+1} (\varphi (x_{\sigma (k)}))  = d^{\m{Y}}(x_{\sigma (n+1)} , x _{\sigma (k)})$.
\end{center}

Observe that $f_{n+1}$ is Kat\v{e}tov over $\{\varphi ( x_{\sigma (k)}) : k < n+1 \}$ and that the space $\{\varphi ( x_{\sigma (k)}) : k < n+1 \} \cup \{ f_{n+1}\}$ is isometric to $\{ x_{\sigma (k)} : k \leqslant n+1 \}$. By hypothesis on $\m{Y}$, we can consequently find $\varphi (x_{\sigma (n+1)})$ realizing $f_{n+1}$ over $\{\varphi ( x_{\sigma (k)}) : k < n+1 \}$. Next, let: 

\begin{center}
$\tau (n+1) = \min \{ k \in \omega : z_k \notin \{ \varphi (x_{\sigma (i)}) : i < n+1 \} \}$ 
\end{center}

Consider the map $g_{n+1}$ defined on $\{ x_{\sigma (k)} : k < n+1 \}$ by:

\begin{center}
$\forall k \leqslant n+1, \ \ g_{n+1} (x_{\sigma (k)}) = d^{\m{Y}}(z_{\tau (n+1)} , \varphi (x_{\sigma (k)}))$. 
\end{center}

Then $g_{n+1}$ is Kat\v{e}tov over the space $\{x_{\sigma (k)} : k < n+1 \}$ and the corresponding union $\{x_{\sigma (k)} : k < n+1 \} \cup \{ g_{n+1}\}$ is isometric to $\{ \varphi (x_{\sigma (k)}) : k < n+1 \} \cup \{ z_{\tau (n+1)} \}$. So again, by hypothesis on $\m{Y}$, we can find $\varphi ^{-1} (z_{\tau (n+1)}) \in \m{Y}$ realizing $g_{n+1}$ over the space $\{x_{\sigma (k)} : k < n+1 \}$. In general, if $\sigma$ and $\tau$ have been defined up to $m$ and $\varphi$ has been defined on $T_m := \{x_{\sigma (0)},\ldots , x_{\sigma (m)} \} \cup \{ \varphi ^{-1} (z_{\sigma (0)}),\ldots , \varphi (z_{\sigma (m)}) \}$, set:

\begin{center}
$\sigma (m+1) = \min \{ k \in \omega : x_k \notin T_m \}$. 
\end{center}

Consider the map $f _{m+1}$ defined on $\varphi '' T_m$ by:

\[
\forall k < m+1, \ \left\{ 
\begin{array}{l} 
f_{m+1} ( \varphi (x _{\sigma (k)})) = d^{\m{Y}}(x_{\sigma (m+1)} , x _{\sigma (k)})\\ 
f_{m+1} ( z _{\tau (k)})) = d^{\m{Y}}(x_{\sigma (m+1)} ,\varphi ^{-1} (z _{\tau (k)}) )
\end{array} \right.
\]

Observe that $f_{m+1}$ is Kat\v{e}tov over $\varphi '' T_m$ and that $\varphi '' T_m \cup \{ f_{m+1}\}$ is isometric to $T_m \cup \{ x_{\sigma (m+1)} \}$. By hypothesis on $\m{Y}$, we can consequently find $\varphi (x_{\sigma (m+1)})$ realizing $f_{m+1}$ over $\varphi '' T_m $. Next, let:

\begin{center}
$\tau (m+1) = \min \{ k \in \omega : z_k \notin \{ \varphi (x_{\sigma (i)}) : i < n+1 \} \}$ 
\end{center}

Consider the map $g_{m+1}$ defined on $T_m$ by:

\[
\forall k < m+1, \ \left\{ 
\begin{array}{l} 
g_{m+1} (x _{\sigma (k)}) = d^{\m{Y}}(z _{\tau (m+1)} , \varphi ( x _{\sigma (k)}))\\ 
g_{m+1} (\varphi ^{-1} (z _{\tau (k)})) = d^{\m{Y}}(z_{\tau (m+1)} , z _{\tau (k)} )
\end{array} \right.
\] 

Then $g_{n+1}$ is Kat\v{e}tov over $T_m$ and the union $T_m \cup \{ g_{m+1}\}$ is isometric to $ \varphi '' T_m \cup \{ z_{\tau (m+1)} \}$. So again, by hypothesis on $\m{Y}$, we can find $\varphi ^{-1} (z_{\tau (m+1)}) \in \m{Y}$ realizing $g_{m+1}$ over $T_m$. After $\omega$ steps, we are left with an isometry $\varphi$ with domain $ \m{Y} = \{ x_k : k \in \omega \}$ and range $\m{Y} = \{ z_k : k \in \omega \}$. This finishes the proof when $\m{Y}$ is countable.

If $\m{Y}$ is complete separable, then the same proof works except that at the very beginning, instead of extending $\{ x_0 ,\ldots, x_n \}$ and $\{z_0 ,\ldots, z_n \}$ so that $\{ x_k : k \in \omega \} = \{ z_k : k \in \omega \} = \m{Y}$, we simply require that $\{ x_k : k \in \omega \}$ and $\{ z_k : k \in \omega \}$ should be dense in $\m{Y}$. At the end of the construction, $\varphi$ is such that $ \{ x_k : k \in \omega \} \dom \varphi$ and $ \{ z_k : k \in \omega \} \subset \ran \varphi$. We can consequently extend it to an isometry of $\m{Y}$ onto itself. \end{proof}

This chapter is organized as follows: In section 2, we present several amalgamation classes of finite metric spaces. In section 3, we present the Urysohn spaces associated to those classes. We finish in section 4, with a section on complete separable ultrahomogeneous metric spaces. 

\section{Amalgamation and Fra\"iss\'e classes of finite metric spaces.}

\subsection{First examples and path distances.} 

\label{subsection:First examples and path distances}
The very first natural example of amalgamation class of finite metric spaces is the class $\M$\index{$\M$} of \emph{all} finite metric spaces. Showing that $\M$ satisfies the amalgamation property (and in fact the strong amalgamation property) is not difficult but the underlying idea will be useful later so we provide a complete proof. \index{amalgamation!strong amalgamation property!for $\M$}

\begin{prop}

\label{prop:pathmetric}

The class $\M$ of all finite metric spaces has the strong amalgamation property. 
\end{prop}

\begin{proof}
Let $\m{X}$, $\m{Y} _0$, $\m{Y} _1 \in \mathcal{M}$ and isometries $f_0 : \funct{\m{X}}{\m{Y} _0}$ and $f_1 : \funct{\m{X}}{\m{Y}}$. We wish to find $\m{Z} \in \mathcal{M}$ and isometries $g_0 : \funct{\m{Y} _0}{\m{Z}}$, $g_1 : \funct{\m{Y} _1}{\m{Z}}$ such that $g_0 \circ f_0 = g_1 \circ f_1$. Equivalently, as mentioned in the previous section, we may assume that $\m{X}$ is a metric subspace both of $\m{Y} _0$ and $\m{Y} _1$, that $\m{Y} _0 \cap \m{Y} _1 = \m{X}$, and that we have to extend $d^{\m{Y} _0} \cup d^{\m{Y} _1} $ to a metric on $Y _0 \cup Y _1$. To achieve that, see $\m{Z} := \m{Y} _0 \cup \m{Y} _1$ as an edge-labelled graph. For $x, y \in Z$, and $n \in \omega$ strictly positive, a define \emph{path from $x$ to $y$ of size $n$}\index{path} as is a finite sequence $\gamma = (z_i)_{i<n}$ such that $z_0 = x$, $z_{n-1} = y$ and for every $i<n-1$, 

\begin{center}
$(z_i, z_{i+1}) \in \dom(\lambda ^{\m{Z}})$. 
\end{center}

The \emph{length}\index{path!length of a path} of $\gamma$ is then defined by\index{$\| \gamma \|$}:

\[ \| \gamma \| = \sum _{i=0} ^{n-1} \lambda ^{\m{Z}}(z_i , z_{i+1} ).\]

Observe that here, the edge-labelled graph $\m{Z}$ is \emph{metric}\index{metric!edge-labelled graph}. This means that for every  $(x,y) \in \dom (\lambda ^{\m{Z}})$ and every path $\gamma$ from $x$ to $y$:

\begin{center} 
$\lambda ^{\m{Z}}(x,y) \leqslant \| \gamma \|$.
\end{center}

This fact allows to define the a metric $d^{\m{Z}}$ as follows: For $x, y$ in $Z$, let $P(x,y)$\index{$P(x,y)$} be the set of all paths from $x$ to $y$. Now, set:

\begin{center}
$d^{\m{Z}}(x,y) = \inf \{ \| \gamma \| : \gamma \in P(x,y)\}$.
\end{center}

Then $d^{\m{Z}}$ is as required. 
\end{proof}

$\M$ is consequently a strong amalgamation class. Not beeing countable, it is not a Fra\"iss\'e class but this can be fixed by restricting the distances to a fixed subset of $]0,+\infty[$ ($0$ is always a distance, so we never mention it as such)\index{amalgamation!strong amalgamation property!for $\M _S$}. The simplest such examples are the classes $\M _{\Q}$\index{$\M _{\Q}$} and $\M _{\omega}$\index{$\M _{\omega}$}, corresponding to the distance sets $\Q \cap ]0, + \infty [$ and $\omega \cap ]0, + \infty [$ respectively. These classes are indeed obviously countable and hereditary. As for the amalgamation property, one can proceed exactly as for $\M$: The fact that the path distance takes its values in $\Q \cap ]0, + \infty [$ or $\omega \cap ]0, + \infty [$ is guaranteed by the fact that these sets are closed under finite sums. Notice also that one may even take bounded subsets of $]0, + \infty[$, say $\Q \cap ]0, r ]$ or $\omega \cap ]0, r ]$ for some strictly positive $r \in \Q$ or $\omega$. In these cases, the previous proof still works provided $\| \gamma \|$ is replaced by $\| \gamma \| _{\leqslant r}$\index{$\| \gamma \| _{\leqslant r}$}:

\[ \| \gamma \| _{\leqslant r} = \min (  \| \gamma \| , r ).\]

\subsection{Ultrametric spaces.}

\label{subsection:ultrametric spaces}

Recall that a metric space $\m{X} = (X, d^\m{X})$ is \emph{ultrametric}\index{ultrametric!ultrametric space} when given any $x, y, z$ in $\m{X}$, 

\begin{center}
$d^\m{X}(x,z) \leqslant \max(d^\m{X}(x,y), d^\m{X}(y,z))$. 
\end{center}

Using the idea of the previous section, one can prove:

\begin{prop}
\index{amalgamation!strong amalgamation property!for $\U$}
Let $S \subset ]0, + \infty [$. Then the class $\U$\index{$\U$} of all finite ultrametric spaces with distances in $S$ has the strong amalgamation property. 
 
\end{prop}

\begin{proof}
Reproduce the proof for $\M$ except that instead of $\| \gamma \| $, use $\| \gamma \| _{\max}$\index{$\| \gamma \| _{\max}$} defined by: 

\[ \| \gamma  \| _{\max} = \max_{0 \leqslant i  \leqslant n-1} \lambda ^{\m{Z}}(z_i , z_{i+1} ). \qedhere \]

\end{proof}

It follows that when $S$ is countable, $\U$ is a Fra\"iss\'e class with strong amalgamation property. In fact, we will see in section \ref{subsection:Ultrametric Urysohn spaces.} that: 

\begin{prop}
Let $\mathcal{K}$ be a Fra\"iss\'e class of finite ultrametric spaces. Assume that $\mathcal{K}$ has the strong amalgamation property. Then there is a countable $S \subset ]0, + \infty [$ such that $\mathcal{K} = \U$. 
\end{prop}

An explicit and detailed study of the classes $\U$ is carried out by Bogatyi in \cite{B1}. Ultrametric spaces are closely related to \emph{trees}\index{tree}. A partially ordered set is a \emph{tree} $\m{T} = (T,<^{\m{T}}) $ when the set $\{ s \in T : s <^{\m{T}} t \}$ is
$<^{\m{T}}$-well-ordered for every element $t
\in T$. When every element of $T$ has finitely many
$<^{\m{T}}$-predecessors, the \emph{height of}\index{height!of an element in a tree} $t
\in \m{T}$ is $\mathrm{ht}(t) = |\{ s \in T : s <^{\m{T}} t \}|$\index{$\mathrm{ht}(t)$}. When $n< \mathrm{ht}(t)$, $t(n)$ denotes the unique
predecessor of $t$ with height $n$. The $m$-th level of $\m{T}$
is $\m{T}(m) = \{t \in T : \mathrm{ht}(t) = m \}$\index{$\m{T}(m)$}. The \emph{height of}\index{height!of a tree} $\m{T}$, $\mathrm{ht}(\m{T})$\index{$\mathrm{ht}(\m{T})$}, is the least $m$ such that $\m{T}(m) = \emptyset$. When $s, t \in \m{T}$, $\Delta (s,t)$ is defined by $\Delta (s,t) = \min \{n < \mathrm{ht}(\m{T}) : s(n) \neq t(n) \}$\index{$\Delta (s,t)$}.

The link between ultrametric spaces and trees is the following: Consider a tree $\m{T}$ of finite height, and where the set $\m{T}^{max}$\index{$\m{T}^{max}$} of 
all $<^{\m{T}}$-maximal elements of $\m{T}$ coincides with the top level set of $\m{T}$ (in other words, all maximal elements have same height). Given such a tree of height $n$ and a finite sequence
$a_0 > a_1> \ldots >a_{n-1}$ of strictly positive real numbers, there is a
natural ultrametric space structure on
$\m{T}^{max}$ if the distance $d$ is defined by:

\begin{center}
$d(s,t) = a_{\Delta(s,t)}$. 
\end{center}

Conversely, given any ultrametric space $\m{X}$ with finitely many distances
given by $a_0 > a_1> \ldots >a_{n-1}$, there is a tree $\m{T}$ of
height $n$ such that $\m{X}$ is the natural ultrametric space associated to $\m{T}$ and $(a_i)_{i<n}$. The
elements of $\m{T}$ are the ordered pairs of the form $( m, b )$ where $0< m < n$ and $b = \{ y \in \m{X} : d^{\m{X}}(y,x) \leqslant a_m \}$ for some $x \in \m{X}$. The structural ordering $<^{\m{T}}$ is
given by:

\begin{center}
$( l, b ) <^{\m{T}} ( m, c )$ iff ($l < m $ and $ b \subset c$).
\end{center}

This connection with trees induces very particular structural properties. For example:

\begin{thm}[Shkarin \cite{S}]

\label{thm:Shkarin}
\index{Shkarin theorem}

Let $\m{X}$ be a finite ultrametric space. Then there is $n \in \omega$ such that $\m{X}$ embeds into any Banach space $\m{Y}$ with $\dim\m{Y} \geqslant n$.   

\end{thm}

This theorem is the last member of a long chain of results concerning isometric embeddings of ultrametric spaces. For example, Vestfrid and Timan proved in \cite{VT1} (see also \cite{VT2}) that any separable ultrametric space is isometric to a subspace of $\ell _2$ (a result also obtained independently by Lemin in \cite{L}). Vestfrid showed later that the result is also true if one replaces $\ell _2$ by $\ell _1$ or $c _0$\index{$\ell_p$, $\ell_{\infty}$}\index{$c_0$} (Recall that $\ell_p$ denotes the Banach space of all real sequences $(x_n)_{n \in \omega}$ such that $\sum_{n=0} ^{\infty} |x_n|^p$ is finite and that $c_0$ is the Banach space of all real sequences converging to $0$ equipped with the supremum norm). Fichet proved that any finite ultrametric space embeds isometrically into $\ell _p$ for every $p \in [1, \infty]$ (Recall also that $\ell_{\infty}$ is the Banach space of all bounded real sequences equipped with the supremum norm), and Vestfrid generalized this fact for a wider class of spaces. For more references, see \cite{S}. Note that it is unknown whether the integer $n$ in Theorem \ref{thm:Shkarin} depends only on the size of $\m{X}$. In other words, is there $n=n(k)$ such that any ultrametric space with size $\leqslant k$ admits an isometric embedding in any $n$-dimensional Banach space? We do not present the proof of Shkarin's theorem here but Fichet's result, which we proved before being aware of the reference, can be obtained easily by combinatorial means: 

\begin{thm}[Fichet \cite{Fi}]

\label{thm:Fichet}
\index{Fichet theorem}

Let $\m{X}$ be a finite ultrametric space. Then there is $n \in \omega$ such that $\m{X}$ embeds into any Banach space $\ell _p ^n$ with $p \in [1 , \infty]$.   

\end{thm}

\begin{proof}

Let $\m{X}$ be a finite ultrametric space with distances
given by $a_0 > a_1> \ldots >a_{n-1}$ and let $\m{T}$ be the finite tree of
height $n$ such that $\m{X}$ is the natural ultrametric space on $\m{T} ^{max}$ associated to $(a_i)_{i<n}$. We show that $n = |\m{T}|$ works. For $p = \infty$, this is a simple consequence of the fact that $\ell ^{|\m{X}|} _{\infty}$ embeds any metric space of size $|\m{X}|$ so we concentrate on the case $p \in [1, \infty[$.  Let $(e_t) _{t \in \m{T}}$ be a subfamily of the canonical basis of $\ell _p$ of size $| \m{T} |$. For $t \in \m{T}$, let

\begin{displaymath} 
\mu(t) = \left \{ \begin{array}{ll}
 (\frac{a_{n-1} ^p}{2})^{\frac{1}{p}} & \textrm{if $\mathrm{ht} (t) = n-1 $} \\
 (\frac{a_i ^p}{2} - \frac{a_{i+1} ^p}{2})^{\frac{1}{p}}  & \textrm{if $\mathrm{ht} (t) = i < n-1 $} 
 \end{array} \right.
\end{displaymath}

Observe then that for every $x, y \in \m{X}$:

\[
d ^{\m{X}} (x, y) = \left( \sum _{\substack{t \leqslant ^{\m{T}} x \\  t \nleqslant ^{\m{T}}  y}} \mu (t) ^p + \sum _{\substack{t \leqslant ^{\m{T}} y \\ t \nleqslant ^{\m{T}} x}} \mu (t) ^p \right) ^{\frac{1}{p}}.
\]

Now, let $\varphi : \funct{\m{X}}{\ell _p}$ be defined by:

\[ \varphi (x) = \sum _{t \leqslant ^{\m{T}} x} \mu (t) e_t .\]

We claim that $\varphi$ is an isometry. Indeed, let $x ,y \in \m{X}$. Then:

\begin{eqnarray*}
\left\| \varphi(y) - \varphi(x) \right\| ^p & = &  \left\| \sum _{t \leqslant ^{\m{T}} y} \mu (t) e_t  - \sum _{t \leqslant ^{\m{T}} x} \mu (t) e_t \right\| ^p \\
& = & \left\| \sum _{\substack{t \leqslant ^{\m{T}} y \\ t \leqslant ^{\m{T}} x}} \mu (t) e_t  + \sum _{\substack{t \leqslant ^{\m{T}}  y \\ t \nleqslant ^{\m{T}} x}} \mu (t) e_t - \sum _{\substack{t \leqslant ^{\m{T}} x \\ t \leqslant ^{\m{T}} y}} \mu (t) e_t - \sum _{\substack{t \leqslant ^{\m{T}} x \\ t \nleqslant ^{\m{T}} y}} \mu (t) e_t \right\|  ^p\\
& = & \left\| \sum _{\substack{t \leqslant ^{\m{T}}  y \\ t \nleqslant ^{\m{T}} x}} \mu (t) e_t - \sum _{\substack{t \leqslant ^{\m{T}} x \\ t \nleqslant ^{\m{T}} y}} \mu (t) e_t \right\|  ^p    \\
& = & \sum _{\substack{t \leqslant ^{\m{T}} x \\ t \nleqslant ^{\m{T}}  y}} \mu (t) ^p + \sum _{\substack{t \leqslant ^{\m{T}} y \\ t \nleqslant ^{\m{T}} x}} \mu (t) ^p \\
& = & d ^{\m{X}} (x, y) ^p . \qedhere
\end{eqnarray*}
\end{proof}

With respect to the comment on Shkarin's theorem mentioned above, note that the previous proof shows that $n$ depends on the size $|\m{X}|$ only. Indeed, notice that if $\m{X}$ is a finite ultrametric space, then the corresponding tree $\m{T}$ associated to $\m{X}$ has the property that each level has strictly less elements than the next level. Therefore, if $\m{X}$ has size $k$, then $\m{T}$ has $k$ maximal elements and at most $k(k+1)/2$ elements. It follows that any ultrametric space $\m{X}$ with size $\leqslant k$ can be embedded into $\ell _p ^n$ where $n = k(k+1)/2$.

\subsection{Amalgamation classes associated to a distance set.}

\label{subsection:Amalgamation classes}

The previous examples are in fact particular instances of a more general case. Indeed, for $S \subset ]0, +\infty[$, let $\M _S$\index{$\M _S$} denote the class of finite metric spaces with distances in $S$. We saw that when $S$ is an initial segment of a set which is closed under finite sums, the path distance allows to prove that $\M _S$ is an amalgamation class. But are there some other cases? For example, can one characterize those subsets $S \subset ]0, +\infty[$ for which $\M _S$ is an amalgamation class? The answer is yes, thanks to a result due to Delhomm\'e, Laflamme, Pouzet and Sauer in \cite{DLPS}.

\begin{defn}

\label{defn:4values}

Let $S \subset ]0, +\infty[$. $S$ satisfies the $4$-\emph{values condition}\index{$4$-values condition} when for every $s_0 , s_1 , s_0 ', s_1 ' \in S$, if there is $t \in S$ such that:
\begin{center} 
$|s_0 - s_1| \leqslant t \leqslant s_0 + s_1$, \ \ $|s_0 '- s_1 '| \leqslant t \leqslant s_0 '+ s_1 '$,
\end{center}
then there is $u \in S$ such that:
\begin{center} 
$|s_0 - s_0 '| \leqslant u \leqslant s_0 + s_0 ' $, \ \ $|s_1 - s_1 '| \leqslant u \leqslant s_1 + s_1 '$.
\end{center}

\end{defn}

In pictures: Assume that the edge-labelled graph $(\{x_0, x_1, y, y'\},\delta)$ described in figure \ref{fig:4value}, where $\delta$ takes values in $S$, is metric. Then $S$ satisfies the $4$-values condition when $\delta$ can be extended to a metric $d$ by setting $d(y,y')=u$ where $u$ is an element of $S$. 

\begin{center}
\begin{figure}[h]
\setlength{\unitlength}{1mm}
\begin{picture}(20,25)(-10,0)

\put(0,0){\circle*{1}}
\put(-10,10){\circle*{1}}
\put(10,10){\circle*{1}}
\put(0,20){\circle*{1}}
\put(0,0){\line(-1,1){10}}
\put(0,0){\line(1,1){10}}
\put(0,0){\line(0,1){20}}
\put(0,20){\line(-1,-1){10}}
\put(0,20){\line(1,-1){10}}
\put(-8,4){$s_1$}
\put(-8,16){$s_0$}
\put(6,4){$s_1 '$}
\put(6,16){$s_0 '$}
\put(1,10){$t$}
\put(-2,-3){$x_1$}
\put(-2,22){$x_0$}
\put(-13,10){$y$}
\put(11,10){$y'$}

\end{picture}
\caption{The edge-labelled graph $(\{x_0, x_1, y, y'\},\delta)$}
\label{fig:4value}
\end{figure}
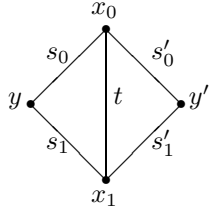
\end{center}

\begin{thm}[Delhomm\'e-Laflamme-Pouzet-Sauer \cite{DLPS}]

\label{thm:4values}
\index{Delhomm\'e-Laflamme-Pouzet-Sauer!theorem on the $4$-values condition}
\index{amalgamation!strong amalgamation property!for $\M _S$}
Let $S \subset ]0, +\infty[$. Then the following are equivalent:

\vspace{0.5em}
\hspace{1em}
i) $\M _S$ has the strong amalgamation property.

\vspace{0.5em}
\hspace{1em}
ii) $\M _S$ has the amalgamation property.

\vspace{0.5em}
\hspace{1em}
iii) $S$ satisfies the $4$-values condition.

\vspace{0.5em}  

\end{thm}
   
\begin{proof}
i) $\rightarrow$ ii) is obvious. For ii)$\rightarrow$ iii), fix $s_0 , s_1 , s_0 ', s_1 ' \in S$ such that there is $t \in S$ with: 

\begin{center} 
$|s_0 - s_1| \leqslant t \leqslant s_0 + s_1$, \ \ $|s_0 '- s_1 '| \leqslant t \leqslant s_0 '+ s_1 '$.
\end{center}

Now, consider $Y := \{ x_0 , x_1 , y \}$ and $Y ' := \{ x_0 , x_1 , y' \}$ and observe that one can define metrics $d^{\m{Y}}$ and $d^{\m{Y}'}$ on $Y$ and $Y'$ by setting:

\[
\left\{
\begin{array}{lll}
d^{\m{Y}} (x_0 , y) = s_0 ,& d^{\m{Y}} (x_1 , y) = s_1 ,& d^{\m{Y}} (x_0 , x_1) = t\\
d^{\m{Y}'} (x_0 , y') = s_0 ' ,& d^{\m{Y}'} (x_1 , y') = s_1 ' ,& d^{\m{Y}'} (x_0 , x_1 ') = t
\end{array} \right.
\]

Therefore, one can obtain a metric space $\m{Z}$ be obtained by amalgamation of $\m{Y}$ and $\m{Y}'$ along $\{ x_0 , x_1 \}$. Then $u=d^{\m{Z}}(y, y')$ is as required. 

For iii) $\rightarrow$ i), consider $\m{Y} _0$ and $\m{Y} _1$ in $\M _S$ such that $d^{\m{Y} _0}$ and $d^{\m{Y} _1} $ agree on $Y_0 \cap Y_1$. We wish to show that $d^{\m{Y} _0} \cup d^{\m{Y} _1} $ can be extended to a metric $d$ on $Y _0 \cup Y _1$. We start with the case where $| Y_0 \smallsetminus Y_1| = | Y_1 \smallsetminus Y_0| = 1$. Set: 

\begin{center}
$Y_0 \smallsetminus Y_1 = \{ y_0 \} , \ \ Y_1 \smallsetminus Y_0 = \{ y_1 \}$.
\end{center}

The only thing we have to do is to define $d$ on $(y_0 , y_1)$. Equivalently, we need to find $u \in S$ such that for every $y \in Y_0 \cap Y_1$:

\begin{center}
$|d^{\m{Y} _0}(y_0 , y) - d^{\m{Y} _1}(y , y_1)| \leqslant u \leqslant d^{\m{Y} _0}(y_0 , y) + d^{\m{Y} _1}(y , y_1)$.
\end{center}

To achieve that, observe that $m \leqslant m'$, where $m$ and $m'$ are defined by: 
\[
\left\{
\begin{array}{l}
m = \max \{|d^{\m{Y} _0}(y_0 , y) - d^{\m{Y} _1}(y , y_1)| : y \in Y_0 \cap Y_1 \} \\
m' = \min \{ d^{\m{Y} _0}(y_0 , y) + d^{\m{Y} _1}(y , y_1) : y \in Y_0 \cap Y_1  \}
\end{array} \right.
\]

Pick witnesses $y$ and $y'$ for $m$ and $m'$ respectively. Then, set:
\[
\left\{
\begin{array}{ll}
s_0 = d^{\m{Y}_0} (y_0 , y) ,&  s_1 = d^{\m{Y}_1} (y_1 , y)\\
s_0 '= d^{\m{Y}_0} (y_0 , y') ,&  s_1 '= d^{\m{Y}_1} (y_1 , y')
\end{array} \right.
\]

Set also:

\begin{center}
$t = d^{\m{Y}_0} (y , y') = d^{\m{Y}_1} (y , y')$.
\end{center}

Then observe that:

\begin{center}
$|s_0 - s_1| \leqslant t \leqslant s_0 + s_1$, \ \ $|s_0 '- s_1 '| \leqslant t \leqslant s_0 '+ s_1 '$.
\end{center} 

So by the $4$-values condition, we obtain the required $u \in S$. We now proceed by induction on the size of the symmetric difference $Y_0 \Delta Y_1$. The previous proof covers the case $|Y_0 \Delta Y_1| \leqslant 2$. For the induction step, let $Y = Y_0 \cup Y_1$. The cases where $Y_0$ and $Y_1$ are $\subset$-comparable are obvious, so we may assume that $Y_0$ and $Y_1$ are $\subset$-incomparable. For $i<2$, pick $y_i \in Y_i \smallsetminus Y_{i-1}$. By induction assumption, obtain a common extension $\m{Z} _0$ of $\m{Y} _0$ and $\m{Y} _1 \smallsetminus \{ y_1\}$ on $Y \smallsetminus \{ y_1 \}$. By induction assumption again, obtain another common extension $\m{Z} _1$ of $\m{Z} _0 \smallsetminus \{y _0 \}$ and $\m{Y} _1$ on $Y \smallsetminus \{ y_0 \}$. Now, observe that $Y = Z_0 \cup Z_1$ and that $|Z_0 \Delta Z_1|=2$, so we can apply the previous case to $\m{Z} _0$ and $\m{Z} _1$ to obtain the required extension. 
\end{proof}

There are some cases where the $4$-values condition is easily seen to hold. For example, if $S \subset [a,2a]$ for some strictly positive $a$, then $S$ satisfies the $4$-values condition. It is also the case when $S$ is closed under sums or absolute value of the difference, which explains why it is possible to restrict distances to $\Q$ or $\omega$. On the other hand, $4$-values condition is also preserved when passing to an initial segment. This allows distance sets of the form $\Q \cap ]0,r]$ or $\omega \cap ]0, r]$. Finally, when $S \subset \{s_n : n \in \mathbb{Z} \}$ with $s_{n} < \frac{1}{2} \ s_{n+1}$, $S$ also satisfes the $4$-values condition as all the elements in $\M _S$ are actually ultrametric. The $4$-values condition consequently covers a wide variety of examples. 

For our purposes, the $4$-values condition is relevant because it allows to produce numerous examples of Fra\"iss\'e classes whose elements can be relatively well handled from a combinatorial point of view. To illustrate that fact, the rest of this section will be devoted to a full classification of the classes $\M _S$ when $|S| \leqslant 3$. This means that we are going to establish a list of classes such that any class $\M _S$ with $|S| \leqslant 3$ will be in some sense isomorphic to some class in the list. More precisely, for finite subsets $S = \{s_0,\ldots , s_m \} _<$, $T = \{t_0,\ldots , t_n \} _<$ of $]0 , + \infty[$, define $S \sim T$\index{$\sim$} when $m=n$ and:

\begin{center}
$\forall i, j, k < m, \ \ s_i \leqslant s_j + s_k \leftrightarrow t_i \leqslant t_j + t_k$.  
\end{center}

Observe that when $S \sim T$, $S$ satisfies the $4$-value condition iff $T$ does and in this case, $S$ and $T$ essentially provide the same amalgamation class of finite metric spaces as any $\m{X} \in \M _S$ is isomorphic to $\m{X} ' = (X , d^{\m{X} '})  \in \M _T$ where:

\begin{center}
$\forall x, y \in X, \ \ d^{\m{X}}(x,y) = s_i \leftrightarrow d^{\m{X} '}(x,y) = t_i$. 
\end{center}

Now, clearly, for a given cardinality $m$ there are only finitely many $\sim$-classes, so we can find a finite collection $\mathcal{S} _m$ of finite subsets of $]0, \infty[$ of size $m$ such that for every $T$ of size $m$ satisfying the $4$-value condition, there is $S \in \mathcal{S} _m$ such that $T \sim S$. Here, we provide such examples of $\mathcal{S} _m$ for $m \leqslant 3$. The reader will find a complete list in Appendix A for $m=4$. This is the largest value we considered as there are already more than 70 $\sim$-equivalence classes on which to test the $4$-values condition. In the sequel, $S = \{s_i : i < |S| \} _<$ is a subset of $]0, +\infty[$. 

The case $|S|=1$ is trivial so we start with $|S|=2$. There are then only $2$ $\sim$-classes corresponding to the following chains of inequalities: 

\begin{center}
(1) $s_0 < s_1 \leqslant 2s_0$. 
	
(2) $s_0 < 2s_0 < s_1$.
\end{center}

(1) is satisfied by the set $\{ 1, 2\}$. The $4$-values condition is satisfied because $\{1, 2 \}$ is an initial segment of $\omega$ which is closed under sums. $\M _{\{ 1, 2\}}$ is consequently a Fra\"iss\'e class. Observe that elements of $\M _{\{ 1, 2\}}$ can be seen as graphs where an edge corresponds to a distance $1$ and a non-edge to a distance $2$. 

(2) is satisfied by the set $\{ 1, 3\}$, which is also a particular case since $1 < \frac{1}{2} \cdot 3$. Thus, elements of $\M _{\{ 1, 3\}}$ are ultrametric and $\M _{\{ 1, 3\}}$ is a Fra\"iss\'e class. 

For $|S|=3$, there are more cases to consider. To list all the relevant chains of inequalities involving elements of $S$, we first write all the relevant inequalities involving $s_0 , s_1$ and their sums. We obtain:

\begin{center}
(1) $s_0 < s_1 \leqslant 2s_0 < s_0 + s_1 < s_1$. 
	
(2) $s_0 < 2s_0 < s_1 < s_0 + s_1 < 2s_1$.
\end{center}

We now look at how $s_2$ may be inserted in these chains. For (1), there are 4 possibilities:

\begin{center}
(1a) $s_0 < s_1 < s_2 \leqslant 2s_0 < s_0 + s_1 < 2s_1$ \hspace{2em} $\{ 2, 3, 4\}$

(1b) $s_0 < s_1 \leqslant 2s_0 < s_2 \leqslant s_0 + s_1 < 2s_1$ \hspace{2em} $\{ 1, 2, 3\}$

(1c) $s_0 < s_1 \leqslant 2s_0 < s_0 + s_1 < s_2 \leqslant 2s_1$ \hspace{2em} $\{ 1, 2, 4\}$

(1d) $s_0 < s_1 \leqslant 2s_0 < s_0 + s_1 < 2s_1 < s_2$ \hspace{2em} $\{ 1, 2 , 5\}$

\end{center} 

We now have to check if the $4$-values condition holds for all the corresponding sets. 

(1a) The set $\{ 2, 3, 4\}$ is an initial segment of $\omega \cap [2, +\infty[$ which is closed under sums. Thus, $\{ 2, 3, 4\}$ satisfies the $4$-values condition. Since there are no non-metric triangles, the elements of $\M _{\{ 2, 3, 4\}}$ can be seen as the edge-labelled graphs with labels in $\{ 2, 3, 4\}$. 

(1b) The set $\{ 1, 2, 3\}$ is also an initial segment of a set which is closed under sums, so it satisfies the $4$-values condition. Note that here, there is a non-metric triangle (corresponding to the distances $1, 1, 3$). 

(1c) The set $\{ 1, 2, 4\}$ does not satisfy the $4$-values condition because of the quadruple $(1, 1, 2, 4)$. $ \M _{\{ 1, 2, 4\}}$ is consequently not a Fra\"iss\'e class. 

\begin{center}
\begin{figure}[h]
\includegraphics[scale=0.6]{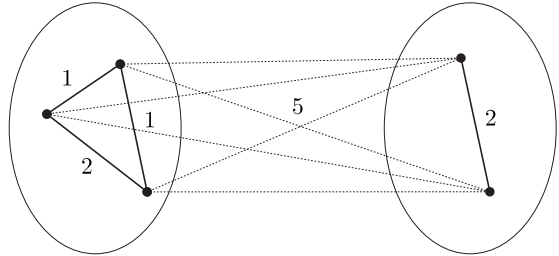}
\caption{An element of $\M _{\{ 1, 2, 5\}}$.}\label{fig125}
\end{figure}
\end{center}

(1d) Finally, the set $\{ 1, 2, 5\}$ satisfies the $4$-values condition but this has to be done by hand (see Appendix A for the details). Simply observe that for $\m{X} \in \M _{\{ 1, 2, 5\}}$, the relation $\approx$ defined by $x \approx y \leftrightarrow d^{\m{X}}(x,y) \leqslant 2$ is an equivalence relation. The $\approx$-classes can be thought as finite graphs with distance $5$ between them. An example is given in Figure \ref{fig125}.

For (2), there are only 3 cases:

\begin{center}

(2a) $s_0 < 2s_0 < s_1 < s_2 \leqslant s_0 + s_1 < 2s_1$ \hspace{2em} $\{ 1, 3 , 4\}$

(2b) $s_0 < 2s_0 < s_1 < s_0 + s_1 < s_2 \leqslant 2s_1$ \hspace{2em} $\{ 1, 3 , 6\}$

(2c) $s_0 < 2s_0 < s_1 < s_0 + s_1 < 2s_1 < s_2$ \hspace{2em} $\{ 1, 3 , 7\}$  

\end{center}

\begin{center}
\begin{figure}[h]
\includegraphics[scale=0.6]{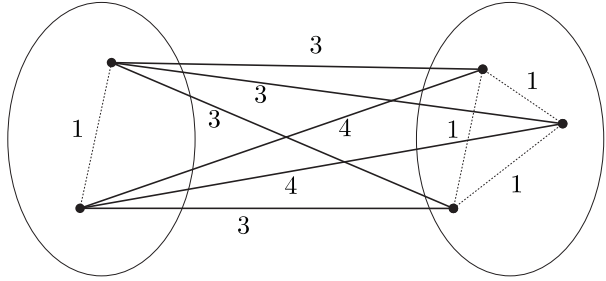}
\caption{An element of $\M _{\{ 1, 3, 4\}}$.}\label{fig134}
\end{figure}
\end{center}

(2a) The $4$-values condition holds for $\{ 1, 3 , 4\}$  but as for $\{1, 2, 5 \}$, this has to be proved by hand. For $\m{X} \in \M _{\{ 1, 3, 4\}}$, the relation $\approx$ defined by $x \approx y \leftrightarrow d^{\m{X}}(x,y) = 1$ is an equivalence relation. Between the elements of two disjoint balls of radius $1$, the distance can be arbitrarily $3$ or $4$. An example is given in Figure \ref{fig134}.

\begin{center}
\begin{figure}[h]
\includegraphics[scale=0.6]{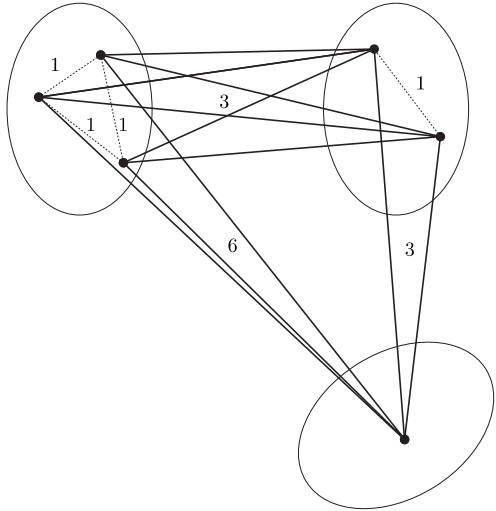}
\caption{An element of $\M _{\{ 1, 3, 6\}}$.}\label{fig136}
\end{figure}
\end{center}

(2b) The set $\{ 1, 3 , 6\}$ also satisfies the $4$-values condition (to be checked by hand). For $\m{X} \in \M _{\{ 1, 3, 6\}}$, the relation $\approx$ defined by $x \approx y \leftrightarrow d^{\m{X}}(x,y) = 1$ is an equivalence relation. Between the elements of two disjoint balls of radius $1$, the distance is either always $3$ or always $6$. An example is provided in figure \ref{fig136}.

(2c) Elements of $\M_{\{ 1, 3 , 7\}}$ are ultrametric. It follows that this class is a Fra\"iss\'e class.

\subsection{Euclidean spaces.}

\label{subsection:Euclidean}

Another way to generate amalgamation classes of finite metric spaces is to fix an ultrahomogeneous metric space and to consider the class of its finite subspaces. For example, if $n \in \omega$ is fixed, the Euclidean space $E_n$\index{$E_n$} of dimension $n$ is ultrahomogeneous (in fact it is even more than ultrahomogeneous as every isometry between \emph{any} two metric subspaces can be extended to an isometry of $E _n$ onto itself). Thus, the class of finite metric subspaces of $E_n$ is an amalgamation class. However, because the dimension is finite, such a class will never have the strong amalgamation property. This requirement being unavoidable when dealing with Ramsey calculus, it will be preferable for us to work with a subclass of the class $\E $\index{$\E $} consisting of all the finite affinely independent metric subspaces of the Hilbert space $\ell _2$. It is easy to see that $\E $ does have the strong amalgamation property. As it is the case for $\M$, $\E $ is not a Fra\"iss\'e class because it is not countable but this can be fixed by restricting the set of distances. For $S$ subset of $]0, +\infty[$, let $\E _S$\index{$\E _S$} denote the class of all elements of $\E $ with distances in $S$. 

\begin{prop}
\index{amalgamation!strong amalgamation property!for $\E _S$}
Let $S$ be dense subset of $]0, +\infty[$. Then $\E _S$ has the strong amalgamation property. 
\end{prop}

\begin{proof}

Following the strategy applied in the previous section, it is enough to show that strong amalgamation holds for $\m{Y} _0$ and $\m{Y} _1$ along $\m{X}$ where 

\begin{center}
$|\m{Y} _0 \smallsetminus \m{Y} _1 | = |\m{Y} _1 \smallsetminus \m{Y} _0| = 1$ and $\m{Y} _i = \m{X} \cup \{ y_i\}$ for each $i<2$. 
\end{center}

Set $n = |\m{X}|$. See $\R ^{n-1}$ as a hyperplane in $\R ^n$ and $\m{X}$ as a metric subspace of $\R ^{n-1}$. Fix $\tilde{y} _0 \in \R ^n$ such that for every $x \in \m{X}$, 

\begin{center}
$\left\| \tilde{y}_0 - x \right\| = d^{\m{Y} _0}(y_0 ,x)$. 
\end{center}

Now, it should be clear that in $\R ^n$ there are exactly two points $y$ such that

\begin{center}
$\forall x \in \m{X}$, $\left\| \tilde{y} - x \right\| = d^{\m{Y} _1}(y_1 , x)$. 
\end{center}

Call them $\tilde{y}_1 ^{min}$ and $\tilde{y}_1 ^{max}$, with $\left\| \tilde{y}_1 ^{min} - \tilde{y}_0 \right\| \leqslant \left\| \tilde{y}_1 ^{max} - \tilde{y}_0 \right\|$. Observe that $\tilde{y}_1 ^{min}$ and $\tilde{y}_1 ^{max}$ are distinct and symmetric with respect to $\R ^{n-1}$. Thus, 

\begin{center}
$\left\| \tilde{y}_1 ^{min} - \tilde{y}_0 \right\| < \left\| \tilde{y}_1 ^{max} - \tilde{y}_0 \right\|$. 
\end{center}

Indeed, if the distances were the same, $\tilde{y} _0$ would be in $\R ^{n-1}$, which is not the case. Now, notice that if we work in $\R ^{n+1}$, we can use rotations to obtain a continuous curve $\varphi : \funct{[0,1]}{\R ^{n+1}}$ such that $\varphi (0) = \tilde{y}_1 ^{min}$, $\varphi (1) = \tilde{y}_1 ^{max}$ and 

\begin{center}
$\forall t \in [0,1] \ \ \forall x \in \m{X} \ \ \left\| \varphi (t) - x  \right\| = d^{\m{Y} _1}(y_1 , x)$. 
\end{center}

Define $\delta : \funct{[0,1]}{\R}$ by: 

\begin{center}
$\delta (t) = \left\| \varphi (t) - \tilde{y} _0  \right\|$
\end{center}

By the intermediate value theorem, $\delta$ takes a value in $S$ for some $t_0 \in ]0,1[$. Then $\m{X} \cup \{ \tilde{y}_0 \} \cup \{ \varphi (t_0) \}$ is the required amalgam. \end{proof}

Observe that a slight modification of the argument allows to show that another class is Fra\"iss\'e and has strong amalgamation: For $\m{X} \in \E$, let $\m{X}^*$\index{$\m{X}^*$} be the edge labelled graph obtained from $\m{X}$ by adjoining an extra point $*$ to $\m{X}$ such that $\lambda ^{\m{X} ^*}(x,*) = 1$ for every $x \in \m{X}$. The class $\ES _S$\index{$\ES _S$} is then defined by the class of all elements $\m{X}$ in $\E _S$ such that $\m{X}^*$ is also in $\E _S$. Equivalently, $\ES _S$ is the class of all elements of $\E _S$ which embed isometrically into the unit sphere $\mathbb{S} ^{\infty}$\index{$\mathbb{S} ^{\infty}$} of $\ell _2$ with the property that $\{ 0_{\ell _2} \} \cup \m{X}$ is affinely independent.

\begin{prop}
\index{amalgamation!strong amalgamation property!for $\ES _S$}
Let $S$ be dense subset of $]0, +\infty[$. Then $\ES _S$ has the strong amalgamation property. 
\end{prop}

\begin{proof}
In the previous proof, simply replace $\m{X}$, $\m{Y} _0$ and $\m{Y} _1$ by $\m{X} ^*$, $\m{Y} ^* _0$ and $\m{Y} ^* _1$ respectively.  
\end{proof}

\textbf{Remark.} It is known that $\ell _2$ is the only separable infinite dimensional ultrahomogeneous Banach space. In fact, much more is known. For example, any separable infinite dimensional Banach space $\m{X}$ where every isometry between finite subsets of size at most $3$ can be extended to an isometry of $\m{X}$ onto itself has to be an inner product space. The problem of whether $3$ can be replaced by $2$ is the content of the famous Banach-Mazur rotation problem. Mazur first proved in \cite{Ma} that the answer is positive in the finite
dimensional case. Pe\l czynski and Rolewicz later showed in \cite{PR} that the
answer is no if one allows $X$ to be non-separable\ldots But in the infinite
dimensional separable case, the problem remains
open, though several partial results seem to suggest that the
answer should be positive (see for example \cite{C2}, \cite{R}, or
\cite{C1} for a survey).

\

We finish this section on Euclidean metric spaces with a further remark about amalgamation property. We saw in section \ref{subsection:First examples and path distances} that when working with metric spaces, an easy way to produce a class of metric spaces with the strong amalgamation property was to start from the class $\M$ of all finite metric spaces and to require that all the distances should be in some $S \subset \R$ that is closed under sums. In particular, we saw that the class $\M _{\omega}$ of all finite metric spaces with distances in $\omega$ has the strong amalgamation property. It turns out that when working with finite Euclidean metric spaces, this is not true anymore: 

\begin{prop}

\label{prop:non amalg E omega}

The class $\E _{\omega}$ does not have the strong amalgamation property. 
\end{prop}

\begin{proof}
Let $\m{X}_n$ denote the finite metric space on $n$ elements where all the distances are equal to $1$. Then $\m{X}_n$ is in $\E _{\omega}$ so one can define $r_n$ the radius of the sphere circumscribed around $\m{X}_n$ in $\R ^{n-1}$. It is easy to show that $(r_n)_{n \in \omega}$ converges to $l = 1/\sqrt{2}$ and since that number is irrational, it follows that for every $\varepsilon >0$, there is $d \in \omega$ such that \[ \left|\left\lceil dl\right\rceil - dl\right|<\varepsilon.\] 

Therefore, for every $\varepsilon >0$, there are $d$ and $n \in \omega$ such that \[ \left|\left\lceil d r_n\right\rceil - d r_n\right|<\varepsilon.\] 

Now, fix $\varepsilon < 1/2$ and consider $d$ and $n$ in $\omega$ as just stated. Let $\m{Y}_n$ denote the finite metric space on $n$ elements where all the distances are equal to $d$. Seeing $\m{Y}_n$ as a subset of $\R^{n-1}$ with isobarycentre $0_{\ell _2}$, let $x \in \R^n$ be orthogonal to $\R^{n-1}$ and such that: \[ \forall y \in \m{Y}_n \ \ \left\|x-y\right\|=\left\lceil d r_n\right\rceil.\] 

Then $\m{Y}_n \cup \{ x\} \in \E _{\omega}$. Note also that \[\left\|x\right\| \leqslant \left|\left\lceil d r_n\right\rceil - d r_n \right| < \varepsilon < 1/2.\] 

As a consequence, one cannot strongly amalgamate two copies of $\m{Y}_n \cup \{ x\}$ by gluing the two copies of $\m{Y}_n$ together while working with distances in $\omega$ only. Indeed, assume that $\m{Y}_n \cup \{ x , x'\}$ is such an amalgam. Then \[\left\|x - x'\right\| \leqslant \left\|x\right\| + \left\|x'\right\| < 1. \qedhere\] 
\end{proof}

The same argument also exhibits a negative amalgamation property for most of the classes $\ES _S$ when $S = \{ k/m : k \in \{ 1,\ldots ,2m\}\}$. Namely, it shows that that there is $M \in \omega$ such that for every integer $m \geqslant M$, the class $\ES _S$ does not have the strong amalgamation property. This fact will be discussed in further detail when we deal with approximations of the spaces $\ell _2$ and $\mathbb{S}^{\infty}$.

\subsection{Other examples.}

\label{subsection:other examples}




There are certainly many more examples of amalgamation classes of finite metric spaces than the ones we mentioned already but as the classification of Fra\"iss\'e classes of finite metric spaces is not known, we will stop our inventory here and refer the interested reader to \cite{B2} by Bogatyi or \cite{W} by Watson. Let us simply mention a very last example, dealing with the class $\mathcal{Q}$\index{$\mathcal{Q}$} of finite metric spaces satisfying the \emph{ultrametric quadrangle inequality}\index{ultrametric!ultrametric quadrangle inequality}. Those are the spaces $\m{X}$ for which given any $x_0 , x_1 , x_2 , x_3 \in \m{X}$, 

\begin{center}
$d^{\m{X}}(x_0, x_1) + d^{\m{X}}(x_2, x_3) \leqslant \max \{ d^{\m{X}}(x_0, x_2) + d^{\m{X}}(x_1, x_3) , d^{\m{X}}(x_0, x_3) + d^{\m{X}}(x_1, x_2) \}$. 
\end{center}

It turns out that $\mathcal{Q}$ is in fact exactly the class of all finite metric spaces which can be embedded into $\R$-\emph{trees}\index{tree!$\R$-tree}. $\R$-trees are defined as follows. For a metric space $\m{Y}$ and $y_0 , y_1 \in \m{Y}$, a \emph{geodesic segment in \m{Y} joining $y_0$ to $y_1$}\index{geodesic segment} is an isometric embedding $g: \funct{[0, d^{\m{Y}}(y_0 , y_1)]}{\m{Y}}$ with $g(0) = y_0$ and $g(d^{\m{Y}}(y_0 , y_1)) = y_1$. Now, a metric space $\m{T}$ is an $\R$-tree if i) For any two distinct points of $\m{T}$, there is a geodesic segment joining them, and ii) If two geodesic segments have exactly one common boundary point, then their union is also a geodesic segment. Using this characterization of $\mathcal{Q}$, one can show that $\mathcal{Q}$ (resp. $\mathcal{Q} _{\Q}$\index{$\mathcal{Q} _{\Q}$}, the class obtained by restricting the distances to $\Q$) is an amalgamation class. $\R$-trees play an important role in so-called asymptotic geometry, but the purpose for which we introduce them here is that they will provide an easy counterexample in section \ref{section:completions} of the present chapter.

\section{Urysohn spaces.}
Recall that according to Fra\"iss\'e theorem, there is a particular countable ultrahomogeneous metric space $\m{X}$ attached to any Fra\"iss\'e class $\mathcal{K}$ of metric spaces: The \emph{Urysohn space}\index{Urysohn!Urysohn space associated to a class of finite metric spaces} associated to $\mathcal{K}$. The purpose of this section is to provide some information about the Urysohn spaces associated to the classes we introduced previously. However, before we start, we should mention that in most of the cases, we will not be able to provide a concrete description of the space. This phenomenon is explained by a general result due to Pouzet and Roux \cite{PouRou} concerning Fra\"iss\'e limits and implying that in some sense, given a countable language $L$ and a Fra\"iss\'e class $\mathcal{K}$ of $L$-structures, the Fra\"iss\'e limit is \emph{generic}\index{generic} among all the countable $L$-structures whose age is included in $\mathcal{K}$. More precisely, when the set of all the countable $L$-structures whose age is included in $\mathcal{K}$ is equipped with the relevant topology, the set of all countable $L$-structures isomorphic to $\mathrm{Flim} (\mathcal{K})$ is a dense $G_{\delta}$ (countable intersection of open sets). This fact is to be compared with the well-known result of Erd\H os and R\'enyi \cite{ER} according to which a random countable graph (obtained by choosing edges independently with probability $1/2$ from a given countable vertex set) is isomorphic to the Rado graph with probability $1$.

\subsection{The spaces $\Ur _{\Q}$ and $\s _{\Q}$.}
The first Urysohn space we present is the space $\Ur _{\Q}$\index{$\Ur _{\Q}$} associated to the class $\M _{\Q}$. This space is called the \emph{rational Urysohn space}\index{Urysohn!rational Urysohn space|see{$\Ur _{\Q}$}} and deserves a particular treatment. It can indeed be seen as the initial step in the construction of Urysohn to provide the very first example of universal separable metric space. The original construction is quite technical but in essence contains the same ideas as the ones that were used some thirty years later in the work of Fra\"iss\'e. The first observation is that to build $\Ur _{\Q}$, it is enough to construct a countable metric space $\m{Y}$ with rational distances such that given any finite subspace $\m{X}$ of $\m{Y}$ and every Kat\v{e}tov map $f$ over $\m{X}$ with rational values, there is $y \in \m{Y}$ realizing $f$ over $\m{X}$. Indeed, for such a $\m{Y}$, ultrahomogeneity is guaranteed by the equivalence provided in proposition \ref{prop:extension}. On the other hand, the set of all finite subspaces is clearly included in $\M _{\Q}$. Consequently, to prove that the finite subspaces of $\m{Y}$ is exactly $\M _{\Q}$, it suffices to show that every element of $\M _{\Q}$ appears as a finite subspace of $\m{Y}$. This is done via the following induction argument: For $\m{X} \in \M _{\Q}$, fix an enumeration $\{ x_n : n < |\m{X}|\}$. Then construct an isometric copy $\mc{X}$ of $\m{X}$ inside $\m{Y}$ by starting with an arbitrary $\tilde{x} _0$ in $\m{Y}$ and by choosing $\tilde{x} _{n+1}$ in the induction step realizing the Kat\v{e}tov map $f_{n+1}$ defined over $\{\tilde{x} _0 ,\ldots \tilde{x} _n \}$ by: 

\begin{center}
$f_{n+1}(\tilde{x} _k) = d^{\m{X}}(x_{n+1} , x_k)$.  
\end{center}

The construction of $\m{Y}$ can be achieved via some kind of exhaustion argument: Start with a singleton $\m{X} _0$. Then, if $\m{X} _n$ is constructed for some $n \in \omega$, $\m{X} _{n+1}$ is build so as to be countable with rational distances, including $\m{X} _n$, and such that given every finite subspace $\m{X} \subset \m{X} _{n}$ and every Kat\v{e}tov map $f$ over $\m{X}$ with rational values, there is $y \in \m{X} _{n+1}$ realizing $f$ over $\m{X}$. Then $\m{Y} = \bigcup _{n \in \omega} \m{X} _n$ is as required. An elegant way to perform the induction step is to follow the method due to Kat\v{e}tov in \cite{K}. It is based on the observation that if $\m{X}$ is a finite subspace of $\m{X} _n$ and $f$ is Kat\v{e}tov over $\m{X}$, then there is a natural way to extend $f$ to a map $k_{\m{X}_n} (f)$ defined on the whole space $\m{X}_n$: Consider the strong amalgam $\m{Z}$ of $\m{X} \cup \{ f\}$ and $\m{X}_n$ along $\m{X}$ obtained using the path metric presented in Proposition \ref{prop:pathmetric}. Then $k_{\m{X}_n} (f)$ is defined by:

\begin{center}
$\forall y \in \m{X}_n, \ \ k_{\m{X}_n} (f) (y) = d^{\m{Z}}(f,y) \ (=\min \{ d^{\m{X}_n} (y,x) + f(x) : x \in \m{X}\})$.
\end{center}

Then, let: \[ X_{n+1} = \bigcup \{ k_{\m{X}_n} (f) : f \in E(\m{X}), \m{X} \subset \m{Y}, \m{X} \  \mathrm{finite} \}.\] 

Equipped with the sup norm, $X_{n+1}$ becomes a metric space $\m{X} _{n+1}$. The map $x \mapsto d^{\m{X}_{n+1}}(x,\cdot)$ then defines an isometric embedding of $\m{X}_n$ into $\m{X}_{n+1}$. The space $\m{X}_n$ can consequently be thought as a subspace of $\m{X}_{n+1}$ and one can check that $\m{X} _{n+1}$ has the required property with respect to $\m{X} _n$. 

A bounded variation of $\Ur _{\Q}$
is obtained by considering the class $\M _{\Q \cap ]0,1]}$. It can be shown that the corresponding Urysohn space, $\s _{\Q}$\index{$\s _{\Q}$}, is isometric to any sphere of radius $1/2$ in the space $\Ur _{\Q}$. For that reason, it is called the \emph{rational Urysohn sphere}\index{Urysohn!rational Urysohn sphere|see{$\s _{\Q}$}}. It will receive a particular interest when we deal with indivisibility.

\subsection{Ultrametric Urysohn spaces.}

\label{subsection:Ultrametric Urysohn spaces.}

We saw that when $S \subset ]0, + \infty [$, the class $\U$\index{$\U$} of finite ultrametric spaces with distances in $S$ is an amalgamation class. So when $S$ is at most countable, the class $\U$ is a Fra\"iss\'e class whose corresponding Urysohn space is denoted $\uUr$\index{$\uUr$}\index{Urysohn!ultrametric Urysohn space|see{$\uUr$}}. A particular feature of this space is that unlike most of the other Urysohn spaces, it admits a very explicit description. Namely, $\uUr$ can be seen as the set of all finitely supported elements of $\Q ^S$ equipped with the distance $d^{\uUr}$ defined by:

\begin{center}
$d^{\uUr}(x,y) = \max \{s \in S : x(s) \neq y(s) \}.$
\end{center}

The spaces $\uUr$ are well-known. They appear together with a study of the classes $\U$ in the article \cite{B1} by Bogatyi but were already studied from a model-theoretic point of view by Delon in \cite{D} and mentioned by Poizat in \cite{P}. More recently, they appeared in \cite{GK} by Gao and Kechris for the study of the isometry relation between ultrahomogeneous discrete Polish ultrametric spaces from a descriptive set-theoretic angle. They are also central in 
\cite{DLPS} where homogeneity in ultrametric spaces is studied in detail. In this paper, these spaces will play a crucial role when we come to the study of big Ramsey degrees as they represent the only case where a complete analysis can be carried out. 

Using the tree representation, one can show that every countable ultrahomogeneous ultrametric space admits a similar description: 

\begin{prop}

\label{prop:ultraultra}
Let $\m{X}$ be a countable ultrahomogeneous ultrametric space. Then there is $S \subset ]0, + \infty[$ at most countable and a family $(A_s)_{s\in S}$ of elements of $\omega \cup \{ \Q \}$ with size at least $2$ such that $\m{X}$ is the set of all finitely supported elements of $\prod _{s\in S} A_s$ equipped with the distance $d$ defined by:  

\begin{center}
$d (x,y) = \max \{s \in S : x(s) \neq y(s) \}$. 
\end{center}
\end{prop}

Note that it is easy to verify that when one of the elements of $(A_s)_{s\in S}$ is finite, the class of its finite subspaces does not have strong amalgamation property. As a consequence, we obtain the following fact mentioned in section \ref{subsection:ultrametric spaces}: The classes $\U$ are the only Fra\"iss\'e classes of finite ultrametric spaces with strong amalgamation property.

\subsection{Urysohn spaces associated to a distance set.}

Similarly, we saw that when $S \subset ]0, + \infty [$ satisfies the $4$-values condition, the class $\M _S$\index{$\M _S$} of finite metric spaces with distances in $S$ is a strong amalgamation class. So when $S$ is at most countable, the class $\M _S$ is a Fra\"iss\'e class whose corresponding Urysohn space is the \emph{Urysohn space with distances in} $S$\index{Urysohn!Urysohn space with distances in $S$|see{$\Ur _S$}}, denoted $\Ur _S$\index{$\Ur _S$}. The space $\Ur _{\Q}$ is a particular case of such space. Similarly, we may simply take $S = \omega \cap ]0,+\infty[$ to obtain the \emph{integral Urysohn space}\index{Urysohn!integral Urysohn space|see{$\Ur _{\omega}$}} $\Ur _{\omega}$\index{$\Ur _{\omega}$}. For $S = \{ 1, 2,\ldots ,m \}$, one obtains a bounded version of $\Ur _{\omega}$ denoted $\Ur _m$\index{$\Ur _m$}. Observe that for $m = 2$, $\Ur _m$ is really the path distance metric space associated to the Rado graph. Finally, the $4$-values condition allows to consider sets $S$ with a more intricate structure than those considered so far. The corresponding Urysohn spaces may then be quite involved combinatorial objects, even when $S$ is finite. In this subsection, we provide a description of $\Ur _S$ when $|S| \leqslant 3$. For $|S|=4$, some cases will be described in the Appendix in order to study their \emph{indivisibility properties}, a notion introduced in the third chapter of this paper. In what follows, the numbering corresponds to the one introduced in subsection \ref{subsection:Amalgamation classes}. 

For $|S|=1$, there is essentially only one Urysohn space: $\Ur _1$, introduced above.  

For $|S|=2$, there are two distances sets, $\{ 1, 2 \}$ and $\{ 1, 3\}$. We just mentioned the case $S = \{ 1, 2 \}$ where the Urysohn space is the Rado graph. As for $S = \{ 1, 3\}$, it was also already presented: $\Ur _{\{1,3\}}$ is ultrametric and is in fact one of the spaces $\uUr$ described in the previous section.  

For $|S|=3$, there are six distances sets. 

(1a) $S = \{ 2, 3, 4\}$. Elements of $\M _{\{ 2, 3, 4\}}$ are essentially edge-labelled graphs with labels in $\{ 2, 3, 4\}$. Consequently, $\Ur _{\{ 2, 3, 4\}}$ can be seen as a complete version of the Rado graph with three kinds of edges. 

(1b) $S = \{ 1, 2, 3\}$. This case was mentioned above, $\Ur _{\{ 1, 2, 3\}}$ is the space we denoted $\Ur _3$. However, like $\Ur _2$ and unlike the other spaces $\Ur _m$ for $m \geqslant 4$, $\Ur _3$ can be described quite simply. This fact, noticed by Sauer, will be important in the third chapter. The main observation is that the only non metric triangle with labels in $\{ 1, 2, 3\}$ corresponds to the labels $1, 1, 3$. It follows that $\Ur _3$ can be encoded by the countable ultrahomogeneous edge-labelled graph with edges in $\{ 1, 3\}$ and forbidding the complete triangle with labels $1, 1, 3$. The distance is then defined as the standard shortest-path distance. Equivalently, the distance between two points connected by an edge is the label of the edge while the distance between two points which are not connected is $2$.      

\begin{center}
\begin{figure}[h]
\includegraphics[scale=0.55]{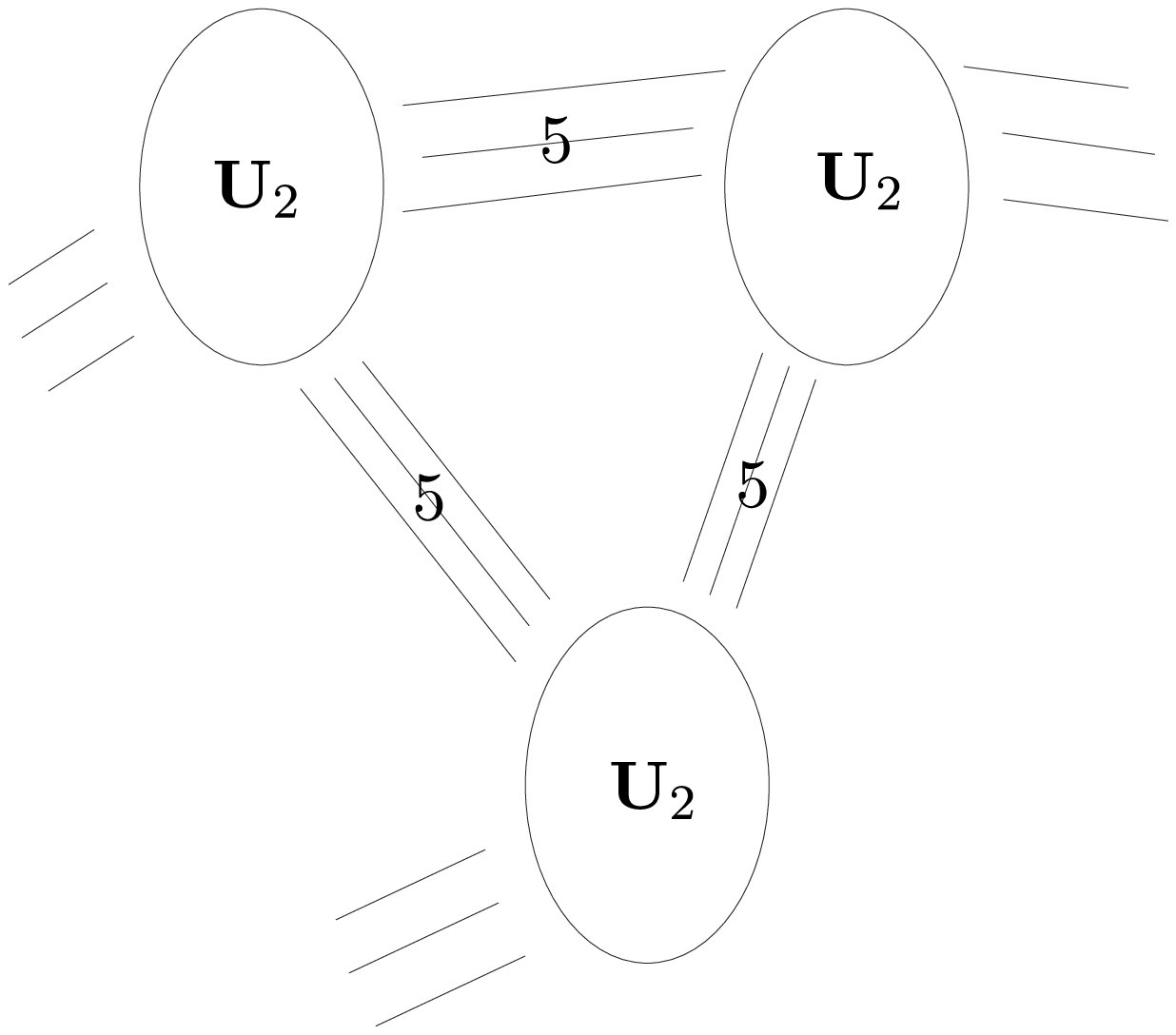}
\caption{$\Ur _{\{ 1, 2, 5\}}$.}\label{figUr125}
\end{figure}
\end{center}

(1d) $S = \{ 1, 2, 5\}$. The structure of the elements of $\M _{\{ 1, 2, 5\}}$ allows to see that $\Ur _{\{ 1, 2, 5\}}$ is composed of countably many disjoint copies of $\Ur _2$, and that the distance between any two points not in the same copy of $\Ur _2$ is always $5$. Figure \ref{figUr125} is an attempt to represent this space.

\begin{center}
\begin{figure}[h]
\includegraphics[scale=0.6]{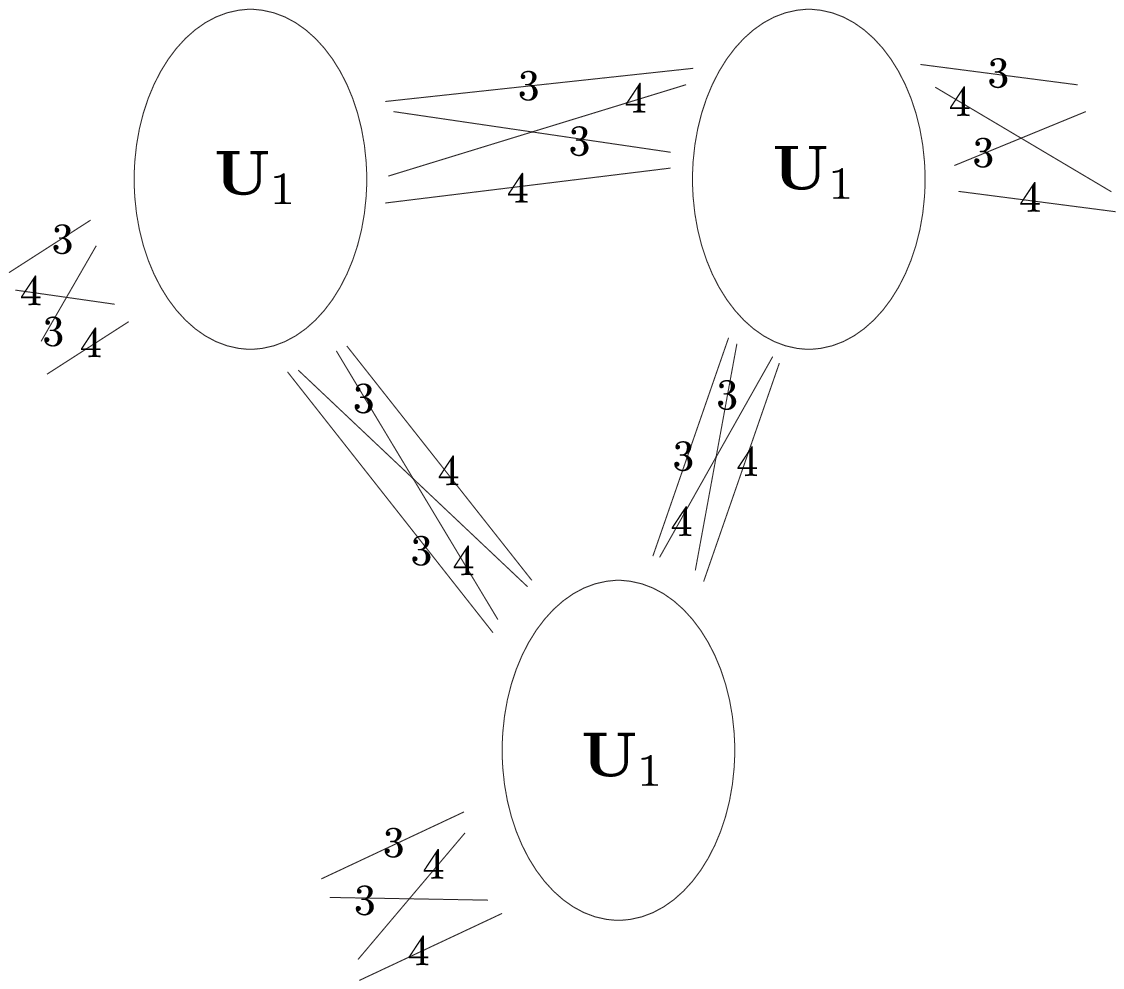}
\caption{$\Ur _{\{ 1, 3, 4\}}$.}\label{figUr134}
\end{figure}
\end{center}

(2a) $S = \{ 1, 3, 4\}$. Here, $\Ur _{ \{ 1, 3, 4\}}$ can be seen as some kind of random partite graph with several kinds of edges. It is composed of countably many disjoint copies of $\Ur _1$ and points belonging to different copies of $\Ur _1$ can be randomly at distance $3$ or distance $4$ apart. Figure \ref{figUr134} is an attempt to represent this space.

\begin{center}
\begin{figure}[h]
\includegraphics[scale=0.5]{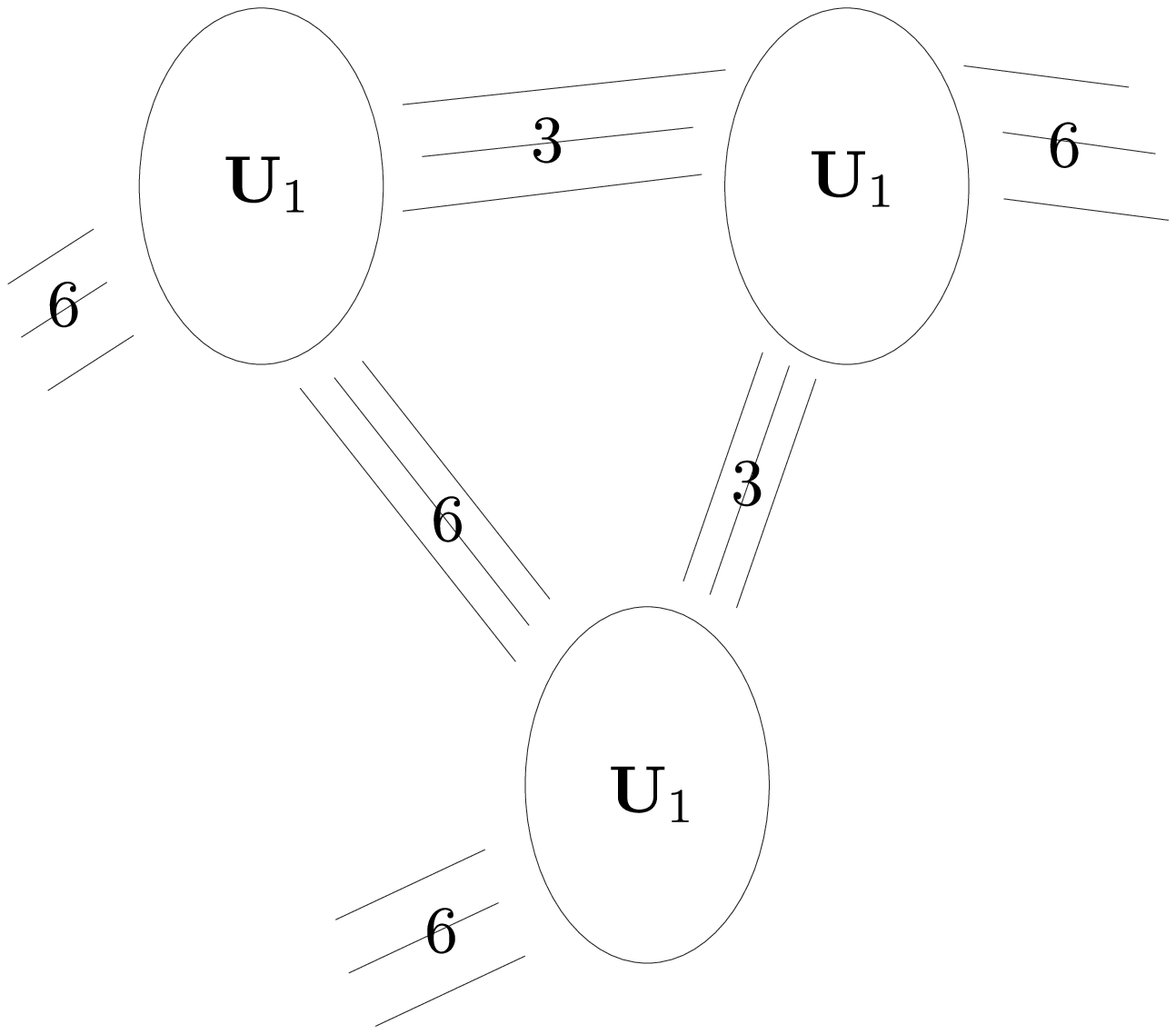}
\caption{$\Ur _{\{ 1, 3, 6\}}$.}\label{figUr136}
\end{figure}
\end{center}

(2b) $S = \{ 1, 3, 6\}$. $\Ur _{\{ 1, 3, 6\}}$ is also composed of countably many disjoint copies of $\Ur _1$ but the distance between points in two fixed disjoint copies of $\Ur _1$ does not vary as in the previous case, and is either $3$ or $6$. A convenient way to construct $\Ur _{\{ 1, 3, 6\}}$ is to obtain it from $\Ur _2$ after having multiplied all the distances by $3$ and blown the points up to copies of $\Ur _1$. Figure \ref{figUr136} is an attempt to represent this space.

(2c) For $S = \{ 1, 3, 7\}$, $\Ur _S$ is again ultrametric, equal to $\uUr$.

\subsection{Countable Hilbertian Urysohn spaces.}

\label{subsection:Hilbertian}

We saw in section \ref{subsection:Euclidean} that when $S$ is a dense subset of $]0, +\infty[$, the class $\E _S$ 
of all finite affinely independent metric subspaces of $\ell _2$ with distances in $S$ is a strong amalgamation class. It follows that the Urysohn space $\textbf{H} _S$\index{$\textbf{H} _S$} associated to $\E _S$ is a countable metric subspace of $\ell _2$ whose elements are all affinely independent. Similarly, the class $\ES _S$ is a strong amalgamation class (recall that $\ES _S$ is the class of all finite metric spaces $\m{X}$ with distances in $S$ and which embed isometrically into the unit sphere $\mathbb{S} ^{\infty}$ of $\ell _2$ with the property that $\{ 0_{\ell _2} \} \cup \m{X}$ is affinely independent). Thus, the associated Urysohn space $\textbf{S} _S$\index{$\textbf{S} _S$} is a countable metric subspace of $\mathbb{S} ^{\infty}$ whose elements are affinely independent. Without being able to go any deeper into the description of those objects, we will see that these spaces have very familiar completions. 

Note on the other hand that still according to results from section \ref{subsection:Euclidean}, the class $\E _{\omega}$ is not a strong amalgamation class. It follows that there is no such a thing as a countable ultrahomogeneous metric space whose class of finite metric subspaces is the class of all affinely independent Euclidean metric spaces with integer distances. Similarly, for $S = \{ k/m : k \in \{ 1,\ldots ,2m\}\}$ with $m$ large enough, the classes $\ES _S$ do not have the strong amalgamation property so there is no countable ultrahomogeneous metric subspace of $\mathbb{S} ^{\infty}$ whose class of finite metric subspaces is $\ES _S$. This comment will be discussed further at the end of Chapter 3.

\section{Complete separable ultrahomogeneous metric spaces.}

\label{section:completions}

It follows from Fra\"iss\'e's theorem that the countable ultrahomogeneous metric spaces are exactly the Fra\"iss\'e limits of the Fra\"iss\'e classes of finite metric spaces. However, many interesting ultrahomogeneous metric are not countable but only separable. We may consequently wonder if there are links between separable ultrahomogeneous metric spaces and countable ones. For example, is the completion of an ultrahomogeneous metric space still ultrahomogeneous? And if so, does every complete separable ultrahomogeneous metric space appear as the completion of a countable ultrahomogeneous metric space? The following theorem provides the answer to the first question.

\begin{prop}[Folklore]

There is an ultrahomogeneous metric space whose completion is not ultrahomogeneous. 

\end{prop}

\begin{proof}

Consider the space $\m{Y}$ defined as follows: Elements of $\m{Y}$ are maps $y : \funct{[0, \rho _y [}{\omega}$ with $\rho _y \in ]0, +\infty[$ and $\{ t \in [0, \rho _y [ : y(t) \neq 0\} \subset \Q$ finite. For $x, y \in \m{Y}$, set:

\begin{center}
$t(x,y) = \min \{ s \in \Q : x(s) \neq y(s) \}$. 
\end{center}

Then, let: 

\begin{center}
$d^{\m{Y}}(x,y) = (\rho _x - t(x,y))+(\rho _y - t(x,y))$.
\end{center}
 
One can check that $\m{Y}$ is complete separable but not ultrahomogeneous. In fact, it is not even point-homogeneous: For $y \in \m{Y}$, if $\rho _y \in \Q$, then $\m{Y} \smallsetminus \{ y \}$ has infininitely many connected components. On the other hand, if $\rho _y \notin \Q$, then $\m{Y} \smallsetminus \{ y \}$ has only two connected components. We now prove the theorem by showing that $\m{Y}$ admits an ultrahomogeneous dense subspace: Consider the subspace $\m{X}$ of $\m{Y}$ corresponding to the elements $x$ of $\m{Y}$ such that $\rho _x \in ]0, +\infty[ \cap \Q$. One can check that $\m{X}$ is countable and dense in $\m{Y}$. But one can also check that $\m{X}$ is ultrahomogeneous by verifying that it is the Fra\"iss\'e limit of the class $\mathcal{Q} _{\Q}$\index{$\mathcal{Q} _{\Q}$} presented in subsection \ref{subsection:other examples}. \end{proof}

The first question above consequently has a negative answer. The purpose of what follows is to show that it is not the case for the second question and that essentially, every complete separable ultrahomogeneous metric space is obtained by completing a countable one. 

\begin{thm}

\label{thm:countable dense ultrahomogeneous}

Every complete separable ultrahomogeneous metric space $\m{Y}$ includes a countable ultrahomogeneous dense metric subspace. 

\end{thm}

\begin{proof}

We provide two proofs. The first one is standard: Let $\m{X} _0 \subset \m{Y}$ be countable and dense. We construct $\m{X}$ countable and ultrahomogeneous such that $\m{X} _0 \subset \m{X} \subset \m{Y}$. We proceed by induction. Assuming that $\m{X} _n \subset \m{Y}$ countable has been constructed, get $\m{X} _{n+1}$ as follows: Consider $\mathcal{F}$ the set of all finite subspaces of $\m{X} _n$. For $\m{F} \in \mathcal{F}$, consider the set $E_n (\m{F})$ of all Kat\v{e}tov maps $f$ over $\m{F}$ with values in the set $\{ d^{\m{Y}} (x,y) : x, y \in \m{X} _n \}$ and such that $\m{F} \cup \{ f \}$ embeds into $\m{Y}$. Observe that $\m{X} _n$ being countable, so are $\{ d^{\m{Y}} (x,y) : x, y \in \m{X} _n \}$ and $E_n (\m{F})$. Then, for $\m{F} \in \mathcal{F}, f \in E_n (\m{F})$, fix $y _{\m{F}} ^f \in \m{Y}$ realizing $f$ over $\m{F}$. Finally, let $\m{X}_{n+1}$ be the subspace of $\m{Y}$ with underlying set $X_n \cup \{ y_F ^f : \m{F} \in \mathcal{F}, f \in E_n (\m{F}) \}$. After $\omega$ steps, set $\m{X} = \bigcup _{n \in \omega} \m{X} _n$. $\m{X}$ is clearly a countable dense subspace of $\m{Y}$. It is ultrahomogeneous thanks to the equivalent formulation of ultrahomogeneity provided in proposition \ref{prop:extension}. Indeed, according to our construction, for every finite subspace $\m{F} \subset \m{X}$ and every Kat\v{e}tov map $f$ over $\m{F}$, if $\m{F} \cup \{ f \}$ embeds into $\m{X}$, then there is $y \in \m{X}$ realizing $f$ over $\m{F}$. This finishes the first proof. 

The second proof was pointed out by Stevo Todorcevic and involves methods from logic. Fix a countable elementary submodel $M \prec H _{\theta}$ for some large enough $\theta$ and such that $Y , d^{\m{Y}} \in M$. Let $\m{X} = M \cap \m{Y}$. We claim that $\m{X}$ has the required property. First, observe that $\m{X}$ is dense inside $\m{Y}$ since by the elementarity of $M$, there is a countable $D \in M$ (and therefore $D \subset M$) which is a dense subset of $\m{Y}$. For ultrahomogeneity, let $\m{F} \subset \m{X}$ be finite and let $f$ be a Kat\v{e}tov map over $\m{F}$ such that $\m{F} \cup \{ f \}$ embeds into $\m{X}$. Observe that $f \in M$. Indeed, $\dom (f) \in M$. On the other hand, let $\mc{F} \cup \{ y \} \subset \m{X}$ be isometric to $\m{F} \cup \{ f \}$ via an isometry $\varphi$. Then for every $x \in \m{F}, d^{\m{Y}}(\varphi (x) ,y) \in M$. But $d^{\m{Y}}(\varphi (x) ,y) = f(x)$. Thus, the range of $f$, $\ran(f)$, is in $M$. It follows that $f$ is an element of $M$. Now, by ultrahomogeneity of $\m{Y}$, there is $y$ in $\m{Y}$ realizing $f$ over $\m{F}$. So by elementarity, there is $x$ in $\m{X}$ realizing $f$ over $\m{F}$. 
\end{proof}

\subsection{The spaces $\Ur$ and $\s$.}
The metric completion $\Ur$\index{$\Ur$}\index{Urysohn!Urysohn space $\Ur$|see{$\Ur$}} of $\Ur _{\Q}$, is known as \emph{the} Urysohn space. It was constructed by Urysohn in 1925 and is, up to isometry, the unique complete separable ultrahomogeneous metric space which contains all finite metric spaces. It follows that $\Ur$ is also universal for the class of all separable metric spaces. This property deserves to be mentioned as historically, $\Ur$ is the first example of separable metric space with this property. However, after Banach and Mazur showed that $\mathcal{C}([0,1])$ was also an example of such a space, the Urysohn space virtually disappeared and it is only after the work of Kat\v{e}tov \cite{K} that $\Ur$ became again subject to research, in particular thanks to the work of Uspenskij, Vershik, Gromov, Bogatyi and Pestov. Today, a complete presentation of the result about the Urysohn space would require much more than what we can provide in the present paper but the reader will find an attempt of survey in the appendix. Let us simply mention the following result due to Pestov \cite{Pe0}: Whenever $\iso (\Ur)$ (equipped with the pointwise convergence topology) acts continuously on a compact space, the action admits a fixed point. We will have the opportunity to come back to this theorem but we would like to mention here once more that its reformulation in terms of structural Ramsey theory by Kechris, Pestov and Todorcevic \cite{KPT} is the starting point of the present paper.   

The metric completion of $\s _{\Q}$\index{$\s _{\Q}$} is the \emph{Urysohn sphere}\index{Urysohn!Urysohn sphere|see{\s}} $\s$\index{$\s$}. Up to isometry, $\s$ is the unique complete separable ultrahomogeneous metric space which contains all finite metric spaces with diameter less or equal to $1$. It is also isometric to any sphere of radius $1/2$ in the Urysohn space $\Ur$, hence the name. The space $\s$ is pretty much as well understood as $\Ur$ is in the sense that most of the proofs working for $\Ur$ can be transposed for $\s$. Later in this paper, we will however study a property called oscillation stability and with respect to which $\Ur$ and $\s$ behave differently.

\subsection{Complete separable ultrahomogeneous ultrametric spaces.} \label{subsection:Complete separable ultrahomogeneous ultrametric spaces}
We now turn to a description of $\cUr$\index{$\cUr$}, the completion of $\uUr$\index{$\uUr$}. Note that if $0$ is not an
accumulation point for $S$, then $\uUr$ is discrete and $\cUr =
\uUr$. Hence, in what follows, we will assume that $0$ is an
accumulation point for $S$. 

\begin{prop}
The completion $\cUr$ of the ultrametric space $\uUr$ is the
ultrametric space with underlying set
the set of all elements $x \in \Q ^S$ for which there is
a strictly decreasing sequence $(s_i)_{i \in \omega}$ of elements of
$S$ converging to $0$ such that $x$ is supported by a subset of $\{s_i : i \in \omega \}$.
The distance is given by 

\begin{center}
$d^{\cUr}(x,y) = \min \{s \in S : \forall t \in
S (s<t \rightarrow x(t) = y(t)) \}$.  
\end{center}

\end{prop}

\begin{proof}
We first check that $\uUr$ is dense in $\cUr$. Let $x \in \cUr$ be
associated to the sequence $(s_i)_{i \in \omega}$. For $n \in \omega$,
let $x_n \in \uUr$ be defined by $x_n (s) = x(s)$ if $s>s_{n+1}$ and by
$x_n (s) = 0$ otherwise. Then $d^{\cUr}(x_n , x) \leqslant s_{n+1}
\longrightarrow 0$, and
the sequence $(x_n)_{n \in \omega}$ converges to $x$.
To prove that $\cUr$ is complete, let $(x_n)_{n \in \omega}$ be a
Cauchy sequence in $\cUr$. Observe first that given any $s \in S$, the
sequence $x_n (s)$ is eventually constant. Call $x(s)$ the
corresponding constant value.

\begin{claim}
$x \in \cUr$.
\end{claim}

The map $x$ is obviously in $\Q ^S$. To show that $x$ is supported by a subset of $\{s_i : i \in \omega \}$ for some strictly decreasing sequence $(s_i)_{i \in \omega}$ of elements of
$S$ converging to $0$, it is enough to show that
given any $s \in S$, there are $t<s<r \in S$ such that $x$ is null 
on $S \cap ]t,s[$ and on $S \cap ]s,r[$. To do that, fix $t' < s $ in
$S$, and take $N \in \omega$ such that $\forall q \geqslant p
\geqslant N$, $d^{\cUr}(x_q , x_p) < t'$. $x_N$ being in $\cUr$, there
are $t$ and $r$ in $S$ such that $t'<t<s<r$ and $x_N$ is null on $S \cap ]t,s[$ and on $S \cap
]s,r[$. We claim that $x$ agrees with $x_N$ on $S \cap ]t', + \infty[$, hence
is null on $S\cap ]t,s[$ and on $S \cap ]s, r[$. Indeed, let $n \geqslant N$. Then $d^{\cUr}(x_n
,x_N) < t' < s$ so $x_n$ and $x_N$ agree on $S \cap ]t', + \infty[$. Hence, for every $u \in S \cap
]t', + \infty [$, the sequence $(x_n (u))_{n \geqslant N}$ is constant and by
definition of $x$ we have $x(u) = x_n (u)$. The claim is proved. 

\begin{claim}
The sequence $(x_n)_{n \in \omega}$ converges to $x$.
\end{claim}

Let $\varepsilon > 0$. Fix $s \in S \cap ]0 , \varepsilon [$ and $N
\in \omega$ such that $\forall q \geqslant p \geqslant N$, $d^{\cUr}(x_q ,x_p) <
\varepsilon$. Then, as in the previous claim, for every $n \geqslant
N$, $x_n$ and $x_N$ (and hence $x$) agree on $S \cap ]s,
+ \infty[$. Thus, $d^{\cUr}(x_n , x) \leqslant s < \varepsilon$. 
\end{proof}

Observe that when $S = \{ 1/(n+1) : n \in \omega \}$, the metric completion of $\uUr$ is the \emph{Baire space}\index{Baire space|see{$\mathcal{N}$}} denoted $\mathcal{N}$\index{$\mathcal{N}$}, a space of particular importance in descriptive set theory.

Note also that the same method can be applied to provide a full description of any complete separable ultrahomogeneous ultrametric space whose distance set admits $0$ as an accumulation point. Indeed, let $\m{X}$ be such a space. According to Theorem \ref{thm:countable dense ultrahomogeneous}, $\m{X}$ admits a countable dense subspace, call it $\m{Y}$. By Proposition \ref{prop:ultraultra}, $\m{Y}$ has a very particular form: It is the space of all finitely supported elements of some product $\prod_{s\in S} A_s$, where each $A_s$ is an integer (seen as a finite set) or $\Q$ and where the distance is defined by 

\[ d^{\m{Y}} (x,y) = \max \{s \in S : x(s) \neq y(s) \}.\] 

Therefore, by the method we just used to describe $\cUr$, the completion of $\m{Y}$ can be described explicitly. Formally: 

\begin{prop}
Let $\m{X}$ be a complete ultrahomogeneous ultrametric space whose distance set $S$ admits $0$ as an accumulation point. Then there is a family $(A_s)_{s\in S}$ of elements of $\omega \cup \{ \Q \}$ with size at least $2$ such that $\m{X}$ is the set of all elements $x \in \prod_{s\in S} A_s$ for which there is
a strictly decreasing sequence $(s_i)_{i \in \omega}$ of elements of
$S$ converging to $0$ such that $x$ is supported by a subset of $\{s_i : i \in \omega \}$. The distance is given by: 

\begin{center}
$d^{\m{X}}(x,y) = \min \{s \in S : \forall t \in S (s<t \rightarrow x(t) = y(t)) \}$.  
\end{center}

\end{prop}

Observe also that in the ultrametric setting, there is no analog of the Urysohn space $\Ur$\index{$\Ur$}: Passing to the completion does not provide a complete separable ultrahomogeneous ultrametric space which is universal for the class of all separable ultrametric spaces. There is a good reason behind this: 

\begin{prop}  
An ultrametric on a separable space takes at most countably many values. 
\end{prop}

\begin{proof}
Let $\m{X}$ be ultrametric and separable with $\m{X}_0 \subset \m{X}$ countable and dense. Then $S := \{ d^{\m{X}}(x,y) : x \neq y \in \m{X}_0 \}$ is countable and $\m{X} _0$ embeds into $\uUr$, so the completion $\widehat{\m{X}} _0$ of $\m{X} _0$ embeds into $\cUr$. But $\m{X} \subset \widehat{\m{X}} _0$. It follows that $\m{X}$ embeds into $\cUr$ and that only countably many distances appear in $\m{X}$. 
\end{proof}

\subsection{$\ell _2$ and $\mathbb{S} ^{\infty}$.}

The purpose of this section is to show how $\ell _2$ or $\mathbb{S} ^{\infty}$ are connected to the spaces introduced in section \ref{subsection:Hilbertian}. We mentioned indeed that for a countable dense $S \subset ]0,+\infty[$, $\E _S$ is a Fra\"iss\'e class whose corresponding Urysohn space $\textbf{H} _S$\index{$\textbf{H} _S$} is a countable metric subspace of $\ell _2$ but that the structure of this space was quite mysterious. The goal of this section is to prove that it is not the case for the completion:

\begin{prop}   
Let $S \subset ]0,+\infty[$ be countable and dense. Then the metric completion of $\textbf{H} _S$ is $\ell _2$.  
\end{prop}

\begin{proof}
It is enough to prove that if $\textbf{H} _S$ is seen as a metric subspace of $\ell _2$ containing $0 _{\ell _2}$, then its closure $\m{X} := \overline{\textbf{H}} _S$ is a vector subspace of $\ell_2$. Indeed, $\m{X}$ will then be an infinite dimensional closed subspace of $\ell_2$, hence isometric to $\ell_2$ itself. 

We first show that if $x \in \m{X}$ and $\lambda \in \R$, then $\lambda x \in \m{X}$. By continuity of $y \mapsto \lambda y$, it suffices to concentrate on the case where $x \in \textbf{H} _S$. Without loss of generality, we may assume $x \neq 0 _{\ell _2}$ and $\lambda \neq 0$. Fix $\varepsilon > 0$. Using the fact that $S$ is dense in $]0,+\infty[$, we can pick $y \in \ell _2$ such that $\{0, x, y \} \in \E _S$ and $\left\| y - \lambda x \right\| < \varepsilon$. By ultrahomogeneity, find $y' \in \textbf{H} _S$ such that $\{ 0_{\ell _2}, x, y' \}$ and $\{ 0_{\ell _2}, x, y \}$ are isometric via the obvious map. Then an easy computation shows that $\left\| y' - \lambda x \right\| < \varepsilon$. Hence, $\lambda x \in \m{X}$. 

Next, we show that $\m{X}$ is closed under sums. As previously, continuity of $+$ allows to restrict ourselves to the case where $x, y \in \textbf{H} _S \smallsetminus \{ 0_{\ell _2} \}$. Fix $\varepsilon > 0$. As previously, find $z \in \ell _2$ be such that $\left\| (x+y) -z \right\| < \varepsilon$ and $\{ 0_{\ell _2} , x, y, z \} \in \E _S$. By ultrahomogeneity, find $z' \in \ell _2$ such that $\{ 0_{\ell _2}, x, y, z' \}$ and $\{ 0_{\ell _2}, x, y, z \}$ are isometric via the obvious map. Then again, an elementary computation shows that $\left\| (x+y) - z' \right\| < \varepsilon$. It follows that $(x+y) \in \m{X}$. \end{proof}

A similar fact holds for $\textbf{S} _S$\index{$\textbf{S} _S$}:

\begin{prop}
Let $S \subset ]0,+\infty[$ be countable and dense. Then the metric completion of $\textbf{S} _S$ is $\mathbb{S} ^{\infty}$.  
\end{prop}

\begin{proof}
See $\textbf{S} _S$ as a metric subspace of $\mathbb{S} ^{\infty}$. Since elements of $\textbf{S} _S \cup \{ 0 _{\ell _2}\}$ are affinely independent, it is enough to prove that $\m{Y} := \overline{\textbf{S}} _S$ is such that the set $\{ \lambda y : \lambda \in \R, \ y \in \m{Y} \}$ is a vector subspace of $\ell_2$. Indeed, $\m{Y}$ will then be the intersection of an infinite dimensional closed subspace of $\ell_2$ with $\mathbb{S} ^{\infty}$, hence isometric to $\mathbb{S} ^{\infty}$ itself. To do that, it suffices to show that $\frac{1}{\left\|x+y\right\|}(x+y) \in \m{Y}$ whenever $x, y \in \m{Y}$ and $x+y \neq 0_{\ell_2}$. By continuity of $\left\|.\right\|$ and of $+$, it is enough to consider the case where $x, y \in \textbf{S} _S$. Fix $\varepsilon > 0$. Find $z \in \mathbb{S} ^{\infty}$ such that $\{ x, y, z \} \in \ES _S$ and $\left\| \frac{1}{\left\|x+y\right\|}(x+y) -z \right\| < \varepsilon$. By ultrahomogeneity, find $z' \in \ell _2$ such that $\{ 0_{\ell _2}, x, y, z' \}$ and $\{ 0_{\ell _2}, x, y, z \}$ are isometric via the obvious map. Then one can check that $\left\| \frac{1}{\left\|x+y\right\|}(x+y) -z' \right\| < \varepsilon$. It follows that $\frac{1}{\left\|x+y\right\|}(x+y) \in \m{Y}$. \end{proof}






\chapter{Ramsey calculus, Ramsey degrees and universal minimal flows.}

\section{Fundamentals of Ramsey theory and topological dynamics.}

In this section, we introduce the basic concepts related to structural Ramsey theory and present the recent results due to Kechris, Pestov and Todorcevic establishing a bridge between structural Ramsey theory and topological dynamics. As for the introductory section in Chapter 1, our main reference here is \cite{KPT}. 

Recall that for $L$-structures $\m{X}, \m{Z}$ in a fixed relational language $L$, $\binom{\m{Z}}{\m{X}}$ denotes the set of all copies of $\m{X}$ inside $\m{Z}$. For $k,l \in \omega \smallsetminus \{ 0 \}$ and a
triple $\m{X}, \m{Y}, \m{Z}$ of $L$-structures, $\m{Z} \arrows{(\m{Y})}{\m{X}}{k,l}$\index{$\m{Z} \arrows{(\m{Y})}{\m{X}}{k,l}$, $\m{Z} \arrows{(\m{Y})}{\m{X}}{k}$} is an abbreviation for the statement: 

\begin{center} 
For any $\chi : \funct{\binom{\m{Z}}{\m{X}}}{k}$
there is $\widetilde{\m{Y}} \in \binom{\m{Z}}{\m{Y}}$ such
that $|\chi ''\binom{\widetilde{\m{Y}}}{\m{X}}| \leqslant l$. 
\end{center} 

When $l = 1$, this is simply
written $\m{Z} \arrows{(\m{Y})}{\m{X}}{k}$. Given a class $\mathcal{K}$ of $L$-structures and $\m{X} \in \mathcal{K} $, suppose that there is $l \in \omega \smallsetminus \{ 0 \}$ such that for any $\m{Y} \in \mathcal{K}$, and any $k \in \omega \smallsetminus \{ 0 \}$, there exists $\m{Z} \in \mathcal{K} $ such that:

\begin{center}
$\m{Z} \arrows{(\m{Y})}{\m{X}}{k,l}$. 
\end{center}

Then we write $\mathrm{t}_{\mathcal{K}}(\m{X})$\index{$\mathrm{t}_{\mathcal{K}}(\m{X})$} for the least such number and call it the \emph{Ramsey degree of} $\m{X}$ \emph{in} $\mathcal{K}$\index{Ramsey!Ramsey degree}. These concepts are closely related to purely Ramsey-theoretic results for classes of \emph{order structures}\index{order!order structure}: Let $L^*$ be a relational signature with a distinguished binary relation symbol $<$. An \emph{order $L^*$-structure} is an $L ^*$-structure $\m{X}$ in which the interpretation $<^{\m{X}}$ of $<$ is a linear ordering. If $\mathcal{K} ^*$ is a class of $L^*$-structures, $\mathcal{K} ^*$ is an \emph{order class}\index{order!order class} when every element of $\mathcal{K} ^*$ is an order $L ^*$-structure. 

Now, given a class $\mathcal{K} ^*$ of finite ordered $L ^*$-structures, say that $\mathcal{K} ^*$ has the \emph{Ramsey property}\index{Ramsey!Ramsey property} (or is a \emph{Ramsey class}\index{Ramsey!Ramsey class}) when for every
$(\om{X}), (\om{Y}) \in \mathcal{K} ^*$
and every $k \in \omega \smallsetminus \{ 0 \}$, there is
$(\om{Z}) \in \mathcal{K} ^*$ such that:

\begin{center} 
$ (\om{Z}) \arrows{(\om{Y})}{(\om{X})}{k}.$ 
\end{center} 

Observe that $k$ can be replaced by $2$ without any loss of generality. On the other hand, given $L ^*$ as above, let $L$ be the signature $L^* \smallsetminus \{ < \}$. Then given an order class $\mathcal{K} ^*$, let $\mathcal{K}$ be the class of $L$-structures defined by:

\begin{center}
$\mathcal{K} = \{ \m{X} : (\om{X}) \in \mathcal{K} ^*\}$.
\end{center}

Say that $\mathcal{K} ^*$ is \emph{reasonable}\index{reasonable} when for every $\m{X}, \m{Y} \in \mathcal{K}$, every embedding $\pi : \funct{\m{X}}{\m{Y}}$ and every linear ordering $\prec$ on $X$ such that $(\m{X}, \prec) \in \mathcal{K} ^*$, there is a linear ordering $\prec '$ on $Y$ such that $\pi$ is also an embedding from $(\m{X}, \prec)$ into $(\m{Y}, \prec ')$. For our purposes, reasonability is relevant because of the following proposition:  

\begin{prop}
Let $L ^* \supset \{ < \}$ be a relational signature, $\mathcal{K} ^*$ a Fra\"iss\'e order class in $L ^*$, $L = L^* \smallsetminus \{ < \}$ and $\mathcal{K} = \{ \m{X} : (\om{X}) \in \mathcal{K} ^*\}$. Let $(\om{F}) = \mathrm{Flim} (\mathcal{K} ^*)$. Then the following are equivalent:

\begin{enumerate}
\item The class $\mathcal{K}$ is a Fra\"iss\'e class and $\m{F} = \mathrm{Flim} (\mathcal{K})$.  
\item The class $\mathcal{K} ^*$ is reasonable. 
\end{enumerate}
\end{prop}

Finally, say that $\mathcal{K} ^*$ has the \emph{ordering property}\index{ordering!ordering property} when given
$\m{X} \in \mathcal{K}$, there is $\m{Y} \in \mathcal{K}$ such that given any linear orderings $<^{\m{X}}$ and $<^{\m{Y}}$ on $\m{X}$ and $\m{Y}$, if $(\om{X})$ , $(\om{Y}) \in \mathcal{K} ^*$, then $(\om{Y})$ contains an isomorphic
copy of $(\om{X})$. Equivalently, for every $(\om{X}) \in \mathcal{K} ^*$, there is $\m{Y} \in \mathcal{K}$ such that for every linear ordering $<^{\m{Y}}$ on $\m{Y}$:

\begin{center}
$(\om{Y}) \in \mathcal{K} ^* \rightarrow \left( (\om{X}) \ \mathrm{embeds \ into} \ (\om{Y}) \right)$. 
\end{center}

Though not exactly stated in the present terminology, the study of the existence and the computation of Ramsey degrees have traditionally been completed for several classes of finite structures such as graphs, hypergraphs and set systems (Ne\v{s}et\v{r}il-R\"odl \cite{NR1}, \cite{NR3}), vector spaces (Graham-Leeb-Rothschild \cite{GLR}), Boolean algebras (Graham-Rothschild \cite{GR}), trees (Fouch\'e \cite{F})\ldots For more information about structural Ramsey theory, the reader should refer to \cite{N0}, to \cite{GRS} or \cite{N2}. As for orderings, it seems that their role was identified quite early (see for example \cite{Lee} or \cite{NR0}). This information, together with many other references about Ramsey and ordering properties, can be found in \cite{N2}. On the other hand, metric spaces do not seem to have attracted much consideration, except maybe when the Ramsey exponent is small (namely, $|\m{X}| = 1$ or $2$, see for example Ne\v{s}et\v{r}il-R\"odl \cite{NR2}). It is only very recently that the first Ramsey class of finite metric spaces was discovered. This result, which is due to Ne\v{s}et\v{r}il and will be presented in the next section, was motivated by the connection we present now between Ramsey theory and topological dynamics.

Let
$G$ be a topological group and $X$ a compact Hausdorff space. A
\emph{$G$-flow}\index{flow} is a continous action $ \funct{G \times X}{X}$. Sometimes,
when the action is understood, the flow is simply referred to as $X$. Given a
$G$-flow $X$, a nonempty compact $G$-invariant subset $Y \subset X$
defines a subflow by restricting the action to $Y$ and $X$ is \emph{minimal}\index{flow!minimal $G$-flow}
when $X$ itself is the only nonempty compact $G$-invariant
set (or equivalently, the orbit of any point of $X$ is dense in
$X$). Using Zorn's lemma, it can be shown that every $G$-flow contains
a minimal subflow. Now, given two $G$-flows $X$ and $Y$, a
\emph{homomorphism}\index{flow!$G$-flow homomorphism} from $X$ to $Y$ is a continuous map $\pi : \funct{X}{Y}$
such that for every $x \in X$ and $g \in G$, $\pi ( g \cdot x) = g
\cdot \pi (x)$. An \emph{isomorphism}\index{flow!$G$-flow isomorphism} from $X$ to $Y$ is a bijective
homomorphism from $X$ to $Y$. The following fact is a standard result
in topological dynamics (a proof can be found in \cite{A}):

\begin{thm}    
Let $G$ be a topological group. Then there is a minimal $G$-flow
$M(G)$ such that for any minimal $G$-flow X there is a surjective
homomorphism $\pi : \funct{M(G)}{X}$. Moreover, up to isomorphism,
$M(G)$ is uniquely determined by these properties.  
\end{thm}

The $G$-flow $M(G)$\index{$M(G)$} is called the \emph{universal minimal flow}\index{flow!universal minimal flow} of $G$. 
When $G$ is locally compact but non compact, $M(G)$ is a highly non-constructive object. Observe also that when
$M(G)$ is reduced to a single point, $G$ has a strong fixed point
property: Whenever $G$ acts continuously on a compact Hausdorff space
$X$, there is a point $x \in X$ such that $ g \cdot x = x$ for every
$g \in G$. $G$ is then said to be \emph{extremely amenable}\index{extreme amenability}.

\begin{thm}[Kechris-Pestov-Todorcevic \cite{KPT}]

\label{thm:KPT ea}
\index{Kechris-Pestov-Todorcevic!theorem on extreme amenability and Ramsey property}

Let $L ^* \supset \{ < \}$ be a relational signature, $\mathcal{K} ^*$ a Fra\"iss\'e order class in $L ^*$ and $(\om{F}) = \mathrm{Flim} (\mathcal{K} ^*)$. Then the following are equivalent:

\begin{enumerate}

\item $\mathrm{Aut}(\om{F})$ is extremely amenable. 

\item $\mathcal{K} ^*$ is a Ramsey class.   

\end{enumerate}

\end{thm}

Let $X _{\mathcal{K} ^*}$ be the set of all \emph{$\mathcal{K} ^*$-admissible orderings}\index{ordering!admissible}, that is linear orderings $\prec$ on $F$ such that for every finite substructure $\m{X}$ of $\m{F}$, $(\m{X} , \restrict{\prec}{\m{X}}) \in \mathcal{K} ^*$. Seen as a subspace of the product $F \times F$ via characteristic functions, the set of all linear orderings on $F$ can be thought as a compact space. As a subspace of that latter space, $X _{\mathcal{K} ^*}$ is consequently compact and acted on continuously by $\mathrm{Aut}(\m{F})$ via the action $\mathrm{Aut}(\m{F}) \times X _{\mathcal{K} ^*} \longrightarrow X _{\mathcal{K} ^*}$, $(g,<) \longmapsto <^g$ defined by $x <^g y$ iff $g^{-1}(x) < g^{-1}(y)$. In other words, $X _{\mathcal{K} ^*}$ can be seen as a compact $\mathrm{Aut}(\m{F})$-flow. The following theorem links minimality of this $\mathrm{Aut}(\m{F})$-flow with the ordering property:

\begin{thm}[Kechris-Pestov-Todorcevic \cite{KPT}]

\label{thm:KPT mf}
\index{Kechris-Pestov-Todorcevic!theorem on minimal flows}

Let $L ^* \supset \{ < \}$ be a relational signature, $L = L^* \smallsetminus \{ < \}$, $\mathcal{K} ^*$ be a reasonable Fra\"iss\'e order class in $L ^*$, and let $\mathcal{K} = \{ \m{X} : (\om{X}) \in \mathcal{K} ^*\}$. Let $(\om{F}) = \mathrm{Flim} (\mathcal{K} ^*)$ and $X _{\mathcal{K} ^*}$ be the set of all $\mathcal{K} ^*$-admissible orderings. Then the following are equivalent: 

\begin{enumerate}

\item $X _{\mathcal{K} ^*}$ is a minimal $\mathrm{Aut}(\m{F})$-flow.  

\item $\mathcal{K} ^*$ satisfies the ordering property.

\end{enumerate}

\end{thm}

Additionally, when Ramsey property and ordering property are satisfied, even more can be said about $X _{\mathcal{K} ^*}$:

\begin{thm}[Kechris-Pestov-Todorcevic \cite{KPT}]

\label{thm:KPT umf}
\index{Kechris-Pestov-Todorcevic!theorem on universal minimal flows}

Let $L ^* \supset \{ < \}$ be a relational signature, $L = L^* \smallsetminus \{ < \}$, $\mathcal{K} ^*$ a reasonable Fra\"iss\'e order class in $L ^*$, and $\mathcal{K}$ defined as $\mathcal{K} = \{ \m{X} : (\om{X}) \in \mathcal{K} ^*\}$. Let $(\om{F}) = \mathrm{Flim} (\mathcal{K} ^*)$ and $X _{\mathcal{K} ^*}$ be the set of all $\mathcal{K} ^*$-admissible orderings. Assume finally that $\mathcal{K} ^*$ has the Ramsey and the ordering properties. Then the universal minimal flow of $\mathrm{Aut}(\m{F})$ is $X _{\mathcal{K} ^*}$. In particular, it is metrizable. 

\end{thm}

Note that this result is not the first one providing a realization of the universal minimal flow of an automorphism group by a space of linear orderings: This approach was first adopted by Glasner and Weiss in \cite{GW} in order to compute the universal minimal flow of the permutation group of the integers. The paper \cite{KPT} continues this trend and provides various other examples. Let us also mention that before \cite{KPT}, the pioneering example by Pestov in \cite{Pe-1} followed by the one by Glasner and Weiss constituted some of the very few known cases of
non extremely amenable topological groups for which the universal
minimal flow was known to be metrizable, a property that
$M(\mathrm{Aut}(\m{F}))$ shares. 

Here, we will be using these theorems to derive results about groups of the form $\iso (\m{X})$ where $\m{X}$ is the Urysohn space or the completion of the Urysohn space attached to a Fra\"iss\'e class of finite metric spaces. 

This chapter is organized as follows: In section 2, we present several Ramsey classes of finite ordered metric spaces. We start with Ne\v{s}et\v{r}il theorem about finite ordered metric spaces, follow with finite convexly ordered ultrametric spaces and finish with results about finite metrically ordered metric spaces. In section 3, we turn to the study of the ordering property and show that all the aforementioned classes satisfy it. We then apply those results to derive several applications. In section 4, we compute Ramsey degrees while in section 5, we use the connection from \cite{KPT} to deduce applications in topological dynamics. We finish in section 6 with some concluding remarks and open problems in metric Ramsey calculus. 

\section{Finite metric Ramsey theorems.}

\subsection{Finite ordered metric spaces and Ne\v{s}et\v{r}il's theorem.}

In what follows, $\M ^<$\index{$\M ^<$} denotes the class of all finite ordered metric spaces. The purpose of this section is to present the proof of the following result, due to Ne\v{s}et\v{r}il. 

\begin{thm}[Ne\v{s}et\v{r}il \cite{N1}]

\label{thm: RP for MM}
\index{Ne\v{s}et\v{r}il!theorem on $\M ^<$}

$\M ^<$ is a Ramsey class. 

\end{thm}

The main idea is to perform a variation of the so-called \emph{partite construction}. This technique is now well-known as its introduction by Ne\v{s}et\v{r}il and R\"odl in the late seventies allowed to solve the long-standing conjecture stating that for every $n \in \omega$, the class of all finite ordered $K_n$-free graphs is a Ramsey class. 

\

\subsubsection{Free amalgamation of edge-labelled graphs.}

The first step is to see finite ordered metric spaces as finite ordered edge-labelled graphs. The result of Ne\v{s}et\v{r}il and R\"odl mentioned above can easily be transposed in the context of edge-labelled graphs (note that the partite construction originally appeared in \cite{NR1}, but the interested reader may refer to \cite{N0} for the details): If one fixes a label set $L$, the class of all finite ordered edge-labelled graphs with labels in $L$ is a Ramsey class. It follows that if $(\om{X})$ and $(\om{Y})$ are finite ordered metric spaces, then there is an edge-labelled graph $(\om{Z})$ with labels in the distance set of $\m{Y}$ such that: \[ (\om{Z}) \arrows{(\om{Y})}{(\om{X})}{2}\] 

The problem here is of course that nothing guarantees that $\m{Z}$ is a metric space. The purpose of what follows is to show that this requirement can be fulfilled. Before going into the details of the proof, observe that ordered edge-labelled graphs satisfy the following version of amalgamation property, called \emph{free amalgamation property}\index{amalgamation!free amalgamation property}: For ordered edge-labelled graphs $(\om{X}), (\m{Y} _0,<^{\m{Y} _0} )$, $(\m{Y} _1,<^{\m{Y} _1} )$ and embeddings $f_0 : \funct{(\om{X})}{(\m{Y} _0,<^{\m{Y} _0} )}$, $f_1 : \funct{(\om{X})}{(\m{Y} _1,<^{\m{Y} _1} )}$, there is a third ordered edge-labelled graph $(\om{Z})$ as well as embeddings $g_0 : \funct{(\m{Y} _0,<^{\m{Y} _0} )}{(\om{Z})}$ and $g_1 : \funct{(\m{Y} _1,<^{\m{Y} _1} )}{(\om{Z})}$ such that:

\vspace{0.5em}
\hspace{1em}
i) $Z = g_0 '' Y_0 \cup g_1 '' Y_1$.

\vspace{0.5em}
\hspace{1em}
ii) $g_0 \circ f_0 = g_1 \circ f_1$, $g_0 '' f_0 '' X = g_0 '' Y_0 \cap g_1 '' Y_1 (= g_0 '' f_0 '' X)$.

\vspace{0.5em}
\hspace{1em}
iii) $\dom(\lambda ^{\m{Z}}) = \bigcup _{i<2} g_i '' \dom(\lambda ^{\m{Y} _i}) = \{ (g_i (x) , g_i (y) ) : (x,y) \in \dom(\lambda ^{\m{Y} _i}) \}$.

\vspace{0.5em}

Such a $(\om{Z})$ is called a \emph{free amalgam}\index{amalgam!free amalgam} of $(\m{Y} _0,<^{\m{Y} _0} )$ and $(\m{Y} _1,<^{\m{Y} _1} )$ over $(\om{X})$. One may think of $(\om{Z})$ as obtained by gluing $(\m{Y} _0,<^{\m{Y} _0} )$ and $(\m{Y} _1,<^{\m{Y} _1} )$ along a prescribed copy of $(\om{X})$. In what follows, free amalgamation will be used to perform the following kind of operation: If an ordered edge-labelled graph $(\om{X})$ embeds into $(\m{Y} _0,<^{\m{Y} _0} )$ and $(\m{Y} _1,<^{\m{Y} _1} )$, then we may obtain a new ordered edge-labelled graph by extending every copy of $(\om{X})$ in $(\m{Y} _1,<^{\m{Y} _1} )$ to a copy of $(\m{Y} _0,<^{\m{Y} _0} )$ and by adding no more connections than necessary.

\

\subsubsection{Hales-Jewett theorem.}

Another ingredient in Ne\v{s}et\v{r}il's proof is the well-known Hales-Jewett theorem coming from combinatorics. A direct combinatorial proof can be found in \cite{GRS}, while a topological proof based on ultrafilters can be found in \cite{T0}. Let $\Gamma$ be a set (the \emph{alphabet}\index{alphabet}), $v \notin \Gamma$ (the \emph{variable}), and $N$ a strictly positive integer. A \emph{word of length $N$ in the alphabet $\Gamma$}\index{word!word in the alphabet $\Gamma$} is a map from $N$ to $\Gamma$. A \emph{variable word in the alphabet $\Gamma$}\index{word!variable word in the alphabet $\Gamma$} is a word in the alphabet $\Gamma \cup \{ v \}$ taking the value $v$ at least once. If $x$ is a variable word and $\gamma \in \Gamma$, $\hat{\gamma} (x)$\index{$\hat{\gamma} (x)$} denotes the word obtained from $x$ by replacing all the occurences of $v$ by $\gamma$ and $\left\langle x \right\rangle$\index{$\left\langle x \right\rangle$} denotes the set defined by 

\begin{center}
$\left\langle x \right\rangle = \{ \hat{\gamma} (x) : \gamma \in \Gamma\}$.   
\end{center}

The set of all words of length $N$ in the alphabet $\Gamma$ is denoted $W(\Gamma , N)$\index{$W(\Gamma , N)$, $V(\Gamma , N)$}, whereas the set of all variable words in the alphabet $\Gamma$ is denoted $V(\Gamma , N)$.  

\begin{thm}[Hales-Jewett \cite{HJ}]

\index{Hales-Jewett theorem}

Let $\Gamma$ be a finite alphabet and $k \in \omega$ strictly positive. Then there exists $N \in \omega$ such that whenever $W(\Gamma , N)$ is partitioned into $k$ many pieces, there is a variable word $x$ of length $N$ in the alphabet $\Gamma$ such that $\left\langle x \right\rangle$ lies in one part of the partition. 
\end{thm}

\subsubsection{Liftings.}

With the previous concepts in mind, we can turn to the first part of Ne\v{s}et\v{r}il's proof. It involves an analog of partite graphs which we will call here $\emph{liftings}$. For an edge-labelled graph $(\om{X})$ and subsets $A$ and $B$ of $X$, write $A < ^{\m{X}} B$\index{$A < ^{\m{X}} B$} when 

\begin{center}
$\forall a \in A \ \forall b \in B \ \ a < ^{\m{X}} b$. 
\end{center}

\begin{defn}
Let $(\om{X})$ with $X = \{x_{\alpha} : \alpha \in |X| \} _{<^{\m{X}}}$ be an ordered edge-labelled graph. A \emph{lifting of}\index{lifting} $(\om{X})$ is an ordered edge-labelled graph $(\om{Y})$ with $Y = \bigcup _{\alpha < |X|} Y_{\alpha}$ such that:

\vspace{0.5em}
\hspace{1em}
i) For every $\alpha < \alpha ' < |X| $, $Y_{\alpha} <^{\m{Y}} Y_{\alpha '}$.

\vspace{0.5em}
\hspace{1em}
ii) For every $\alpha, \alpha ' < |X|$, $y_{\alpha} \in Y_{\alpha}$, $y_{\alpha '} \in Y_{\alpha '}$, 

\begin{displaymath} 
\left \{ \begin{array}{l}
 (y_{\alpha} , y_{\alpha '}) \in \dom (\lambda ^{\m{Y}}) \\
 y_{\alpha} \neq y_{\alpha'}
 \end{array} \right. \rightarrow
\left \{ \begin{array}{l}
\alpha \neq \alpha' \\
(x_{\alpha} , x_{\alpha '}) \in \dom (\lambda ^{\m{X}})\\
\lambda ^{\m{Y}} (y_{\alpha},y_{\alpha '}) = \lambda ^{\m{X}}(x_{\alpha}, x_{\alpha '})
\end{array} \right.
\end{displaymath}

\end{defn}

\begin{lemma}

\label{lem:partite}

Let $(\om{X})$ be a finite ordered metric space and $(\om{Y})$ be a lifting of $(\om{X})$. Then there is a lifting $(\om{Z})$ of $(\om{X})$ such that: \[ (\om{Z}) \arrows{(\om{Y})}{(\om{X})}{2}.\]

\end{lemma}

\begin{proof}

Observe first that since $d ^\m{X}$ is defined everywhere on $X \times X$, $x_{\alpha}
\in Y_{\alpha}$ for every $\alpha < |X|$. More generally, if $(\tilde{x}_{\alpha})_{\alpha < |X|}$ is a strictly increasing
enumeration of some copy $(\oc{X})$ of $(\om{X})$ in $(\om{Y})$, then $\tilde{x}_{\alpha}$ is in $\in Y_{\alpha}$ for every $\alpha < |X|$.  

Moreover, if $\alpha \neq \alpha' < |X|$, then 

\begin{center} $\lambda ^{\m{Y}}(\tilde{x}_{\alpha} , \tilde{x}_{\alpha'}) = \lambda ^{\m{X}} (x_{\alpha} , x_{\alpha'})$.\end{center} 

In other
words, the label of an edge in a copy of $(\om{X})$ in $(\om{Y})$ depends only on the parts where the extremities of this edge live. 
Now, let $N \in \omega$ be large enough so that Hales-Jewett theorem
holds for the colorings of the set $\binom{\om{Y}}{\om{X}} ^N$ with two
colors.

For $\alpha < |X|$, set $Z_{\alpha} = Y_{\alpha}^N$. Now, define $Z = \bigcup_{\alpha < |X|}Z_{\alpha}$. $Z$ is a
subset of $Y^N$ and is consequently linearly ordered by the restriction $<^{\m{Z}}$ of the
lexicographical ordering on $Y^N$. Note that this ordering respects
the parts of the decomposition $Z = \bigcup_{\alpha < |X|}Z_{\alpha}$ ie: \begin{center} $Z_{0} <^{\m{Z}} \ldots <^{\m{Z}} Z_{|X|-1}$. \end{center}

For the edges, proceed as follows: For $\alpha, \alpha ' < |X|$, $z_{\alpha} \in Z_{\alpha}, z_{\alpha '} \in Z_{\alpha '}$, set \begin{center} $(z_{\alpha},z_{\alpha '}) \in \dom (\lambda ^{\m{Z}}) \leftrightarrow \left( \forall n < N \ \ (z_{\alpha}(n),z_{\alpha '}(n)) \in \dom (\lambda ^{\m{Y}}) \right)$. \end{center} 

In this case, set \begin{center} $\lambda ^{\m{Z}}(z_{\alpha},z_{\alpha '}) = \lambda ^{\m{X}} (x_{\alpha},x_{\alpha '})$.\end{center} 

This situation is illustrated in Figure \ref{figNes2}.

\begin{center}
\begin{figure}[h]
\includegraphics[scale=1]{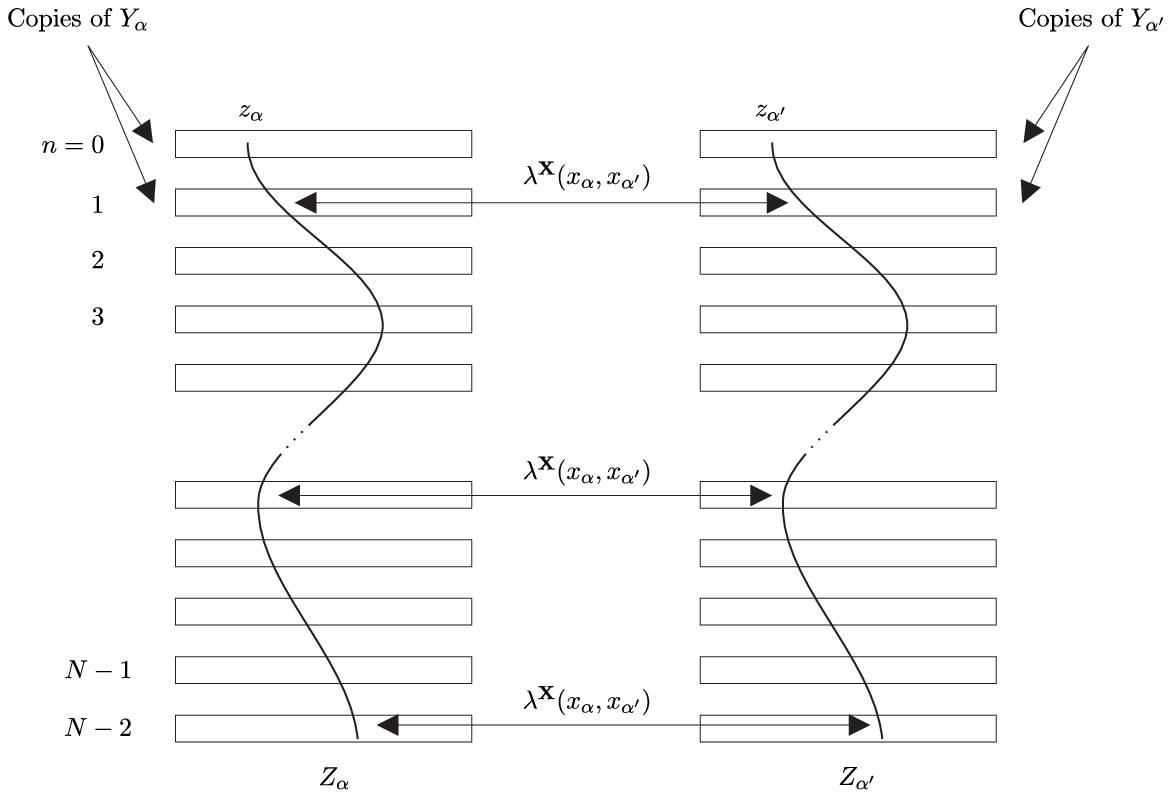}
\caption{An edge $\{ z_{\alpha} , z_{\alpha ' }\}$ with label $\lambda ^{\m{X}}(x_{\alpha} , x_{\alpha '})$.}\label{figNes2}
\end{figure}
\end{center}

It should be clear that the resulting ordered edge-labelled graph $(\om{Z})$ is a lifting of $(\om{X})$. We are now going to show that $(\om{Z}) \arrows{(\om{Y})}{(\om{X})}{2}$. For $n < N$, let $\pi_n$ denote the $n$-th projection from $\m{Z}$ onto
$\m{Y}$, ie: 

\begin{center} $\forall z \in \m{Z} \ \ \pi_{n}(z)=z(n)$. \end{center}

First, observe that copies of $(\om{X})$ are related to their projections. The proof is easy and left to the reader:

\begin{claim}
Let $(\oc{X}) \subset (\om{Z})$. Then: \begin{center} $(\oc{X}) \in \binom{\om{Z}}{\om{X}} \leftrightarrow \left( \forall n < N \ \ \pi_n''(\oc{X}) \in \binom{\om{Y}}{\om{X}} \right)$.\end{center}
\end{claim}

This implies that we can identify $\binom{\om{Z}}{\om{X}}$ with $\binom{\om{Y}}{\om{X}} ^N$, the set of words of length $N$ in the alphabet $\binom{\om{Y}}{\om{X}}$. 


\begin{claim}
Let $U$ be a variable word of length $N$ in the alphabet $\binom{\om{Y}}{\om{X}}$. Then $(\om{Y})$ embeds into $\bigcup \left\langle U \right\rangle$. 
\end{claim}

\begin{proof}
Let $V \subset N$ be the set where the variable lives and let $F = N \smallsetminus V$. For $n \in F$, the $n$th letter of $U$ is a copy $\{ x_{\alpha} ^n : \alpha < |X| \}_{<^{\m{Y}}}$ of $(\om{X})$ in $(\om{Y})$. Now, for $y \in \m{Y}$ with $y \in Y_{\alpha}$, let $e(y)$ be the element of $Z_{\alpha}$ defined by (see Figure \ref{figNes1}): 

\begin{displaymath} 
e(y)(n) = \left \{ \begin{array}{ll}
 x_{\alpha} ^n & \textrm{if $n \in F$,} \\
 y & \textrm{if $n \in V$.} 
 \end{array} \right.
\end{displaymath}

\begin{center}
\begin{figure}[h]
\includegraphics[scale=0.68]{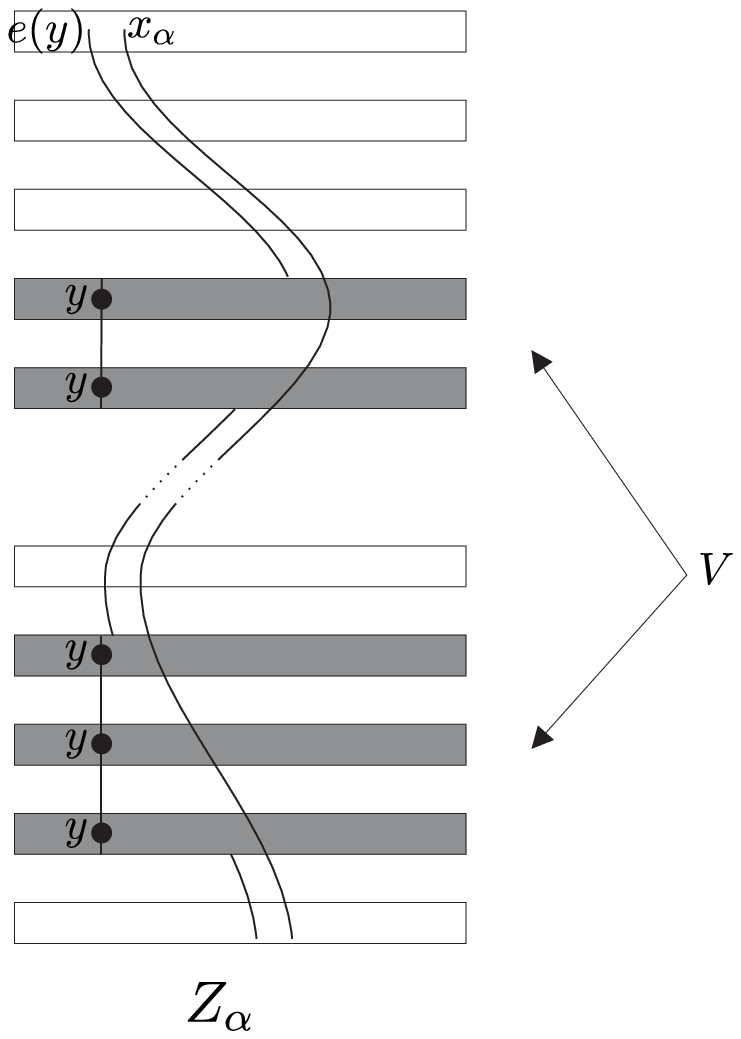}
\caption{$e(y)$ for $y \in Y_{\alpha}$.}\label{figNes1}
\end{figure}
\end{center}

Then $e$ is an embedding from $(\om{Y})$ into $(\om{Z})$ and its direct image $(\oc{Y})$ satisfies: \[\binom{\oc{Y}}{\om{X}} \subset \bigcup \left\langle U \right\rangle. \qedhere\]
\end{proof}

We can now complete the proof of the lemma. Let $\chi : \funct{\binom{\om{Z}}{\om{X}}}{2}$. Thanks to the first claim, $\chi$ transfers to a coloring
$\widehat{\chi} : \funct{\binom{\om{Y}}{\om{X}} ^N}{2}$. Now,
by Hales-Jewett theorem for $\binom{\om{Y}}{\om{X}} ^N$ and two colors,
there is a variable word $U$ of length $N$ in the alphabet
$\binom{\om{Y}}{\om{X}}$ so that $\langle U \rangle$ is monochromatic. This
means that $\binom{\bigcup \langle U \rangle}{\om{X}} $ is
monochromatic. But by the second claim, there is a copy $(\oc{Y})$ of $(\om{Y})$ inside $\bigcup \left\langle U \right\rangle$. Then $\binom{\oc{Y}}{\om{X}}$ is monochromatic. \end{proof}

\subsubsection{Partite construction.}

We start with the following definition, linked to the notion of metric path introduced in Chapter 1. Recall that for an edge-labelled graph $(\om{Z})$, $x, y \in Z$, and $n \in \omega$ strictly positive, a path from $x$ to $y$ of size $n$ as is a finite sequence $\gamma = (z_i)_{i<n}$ such that $z_0 = x$, $z_{n-1} = y$ and for every $i<n-1$, 

\begin{center}
$(z_i, z_{i+1}) \in \dom(\lambda ^{\m{Z}})$. 
\end{center}

For $x, y$ in $Z$, $P(x,y)$ is the set of all paths from $x$ to $y$. If $\gamma = (z_i)_{i<n}$ is in $P(x,y)$, $ \| \gamma \|$ is defined as: \[ \| \gamma \| = \sum _{i=0} ^{n-1} \delta (z_i , z_{i+1} ). \]

On the other hand, for $r \in \R$, $\| \gamma \| _{\leqslant r} $ is defined as: \[ \| \gamma \| _{\leqslant r} = \min (\| \gamma \| , r).\]

\begin{defn}

Let $l \in \omega$ be strictly positive and $\m{X}$ be an edge-labelled graph. $\m{X}$ is \emph{$l$-metric}\index{metric!$l$-metric edge-labelled graph} when for every  $(x,y) \in \dom (\lambda ^{\m{X}})$ and every path $\gamma$ from $x$ to $y$ of size less or equal to $l$:

\begin{center} 
$\lambda ^{\m{X}}(x,y) \leqslant \| \gamma \|$.
\end{center}

\end{defn}

It follows that $\m{X}$ is metric when $\m{X}$ is $l$-metric for every strictly positive $l \in \omega$. Observe that this concept is only relevant when $\lambda ^{\m{X}}$ is not defined everywhere on $X \times X$. 

\begin{prop}
Let $l \in \omega$. Let $\m{Z}$ be a finite $l$-metric edge-labelled graph with label set $L_{\m{Z}}$ such that $l \in \omega$ is such that $\max L_{\m{Z}} \leqslant l \cdot \min L_{\m{Z}}$. Then $\lambda ^{\m{Z}}$ can be extended to a metric on $\m{Z}$.
\end{prop}

\begin{proof}
Using the notation introduced in Chapter 1, simply check that $d^{\m{Z}}$ is as required, where \[ \forall x, y \in Z \ \ d^{\m{Z}}(x,y) = \inf \{ \| \gamma \| _{\leqslant \max L_{\m{Z}} } : \gamma \in P(x,y)\}. \qedhere\]

\end{proof}

Now, let $D _{\m{Y}}$ be the distance set of $\m{Y}$. To show that there is a finite ordered metric space $(\om{Z})$ such that $(\om{Z}) \arrows{(\om{Y})}{(\om{X})}{2}$, it suffices to show that for every strictly positive $l \in \omega$, the statement $\mathcal{H} _l$ holds, where 

\begin{center} 
$\mathcal{H} _l$ : "There is an $l$-metric edge-labelled graph $(\om{Z})$ with $L _{\m{Z}} \subset D _{\m{Y}}$  such that $(\om{Z}) \arrows{(\om{Y})}{(\om{X})}{2}$."
\end{center}

\begin{proof}

We proceed by induction on $l>0$. For $l=1$, there is no restriction on $\m{Z}$, so $\mathcal{H} _1$ is true according to the general theory of Ne\v{s}et\v{r}il and R\"odl.  
Assume now that for a given $l>0$, $\mathcal{H} _l$ holds with witness $(\om{Z}) = \{ z_{\alpha} : \alpha < |\m{Z}|\}$. Let $(\m{P} _0,<^{\m{P} _0} )$ be the lifting of $(\om{Z})$ obtained as follows: The underlying set $P_0$ is obtained by taking a disjoint union of copies of $(\om{Y})$, one for each copy of $(\om{Y})$ in $(\om{Z})$:
\begin{center} $P _0 = \bigcup _{\beta \in \binom{\om{Z}}{\om{Y}}} Y_{\beta}$. \end{center}

For the parts of $\m{P}_0$, given $\beta \in \binom{\om{Z}}{\om{Y}}$, let $\pi _0 ^{\beta}$ be the order preserving isometry from $Y_{\beta}$ onto $\beta$ and let \[ \pi _0  = \bigcup \{ \pi _0 ^{\beta} : {\beta \in \binom{\om{Z}}{\om{Y}}} \}.\]

Then define \[ P _{0 \alpha} = \overleftarrow{\pi _0} \{ z_{\alpha} \}.\]

The construction of $\m{P}_0$ is illustrated in Figure \ref{figNes3}.

\begin{center}
\begin{figure}[h]
\includegraphics[scale=0.75]{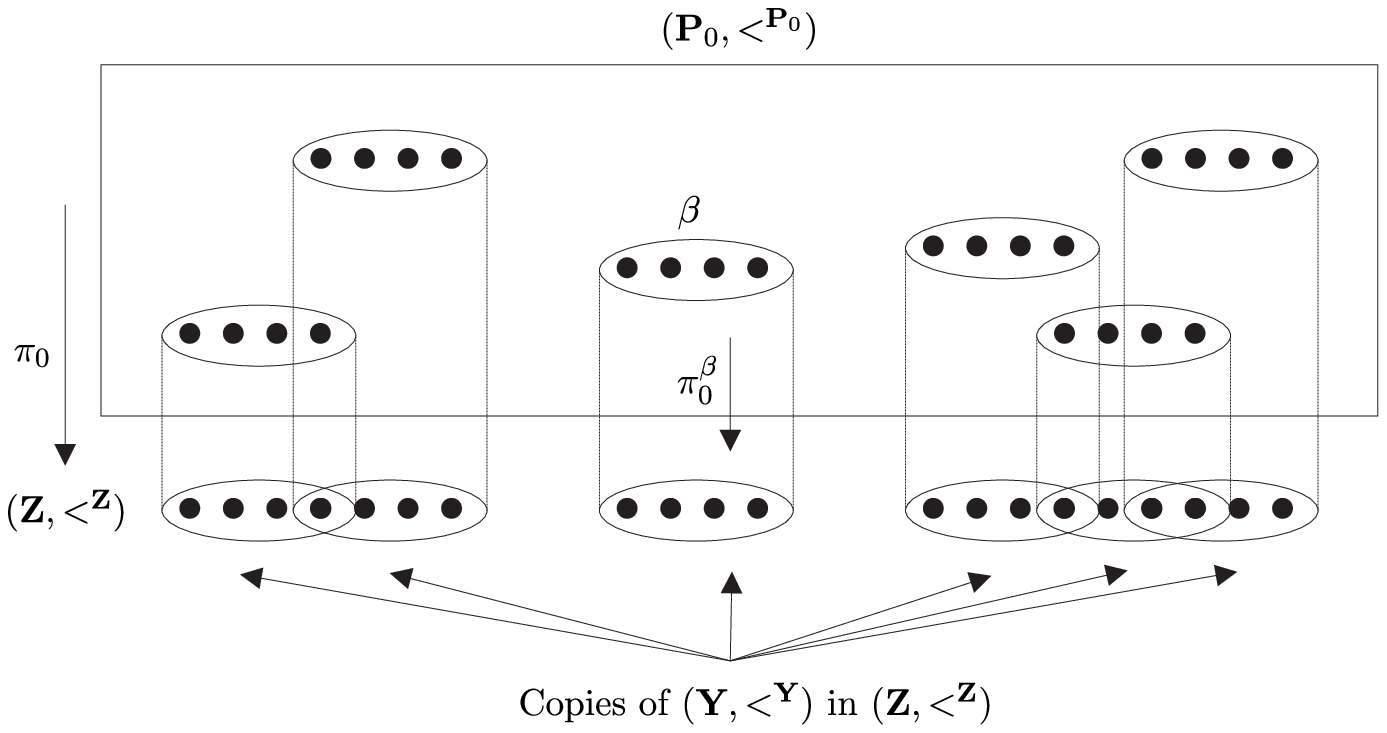}
\caption{Construction of $\m{P} _0$.}\label{figNes3}
\end{figure}
\end{center}


Finally, for the linear ordering $<^{\m{P} _0}$, observe that the linear ordering $<^{\m{Z}}$ already allows to compare points which are not in a same part. By ordering the elements within a same part arbitrarily, one consequently obtains a linear ordering which respects the parts of the decomposition of $P_0$. The resulting lifting of $(\om{Z})$ is $(\m{P} _0,<^{\m{P} _0} )$.

Observe that $\m{P} _0$ is metric, and consequently $(l+1)$-metric. Now, write 

\begin{center}
$\binom{\om{Z}}{\om{X}} = \{\m{X}_1\ldots \m{X}_q \}$.
\end{center}

Inductively, we are now going to construct liftings $(\m{P} _1,<^{\m{P} _1} )$,$\ldots$, $(\m{P} _q,<^{\m{P} _q} )$ of $(\om{Z})$, each of them $(l+1)$-metric, and such that:

\begin{center}    
$(\m{P} _q,<^{\m{P} _q} ) \arrows{(\om{Y})}{(\om{X})}{2} $ 
\end{center}

To construct $(\m{P} _1,<^{\m{P} _1} )$, consider $\overleftarrow{\pi _0}\m{X}_1$. The ordered edge-labelled graph induced on this set, call it $(\m{V} _1,<^{\m{V} _1} )$, is a lifting of $(\om{X})$. Apply lemma \ref{lem:partite} to get a lifting $(\m{W} _1,<^{\m{W} _1} )$ of $(\om{X})$ such that 

\begin{center}
$(\m{W} _1,<^{\m{W} _1}) \arrows{(\m{V} _1,<^{\m{V} _1})}{(\om{X})}{2}$. 
\end{center}

By strong amalgamation property, extend every element of $\binom{\m{W} _1,<^{\m{W} _1}}{\m{V} _1,<^{\m{V} _1}}$ to a copy of $(\m{P} _0,<^{\m{P} _0})$. The resulting finite edge-labelled graph is $\m{P}_1$. Its construction is illustrated in Figure \ref{figNes4}.

\begin{center}
\begin{figure}[h]
\includegraphics[scale=1]{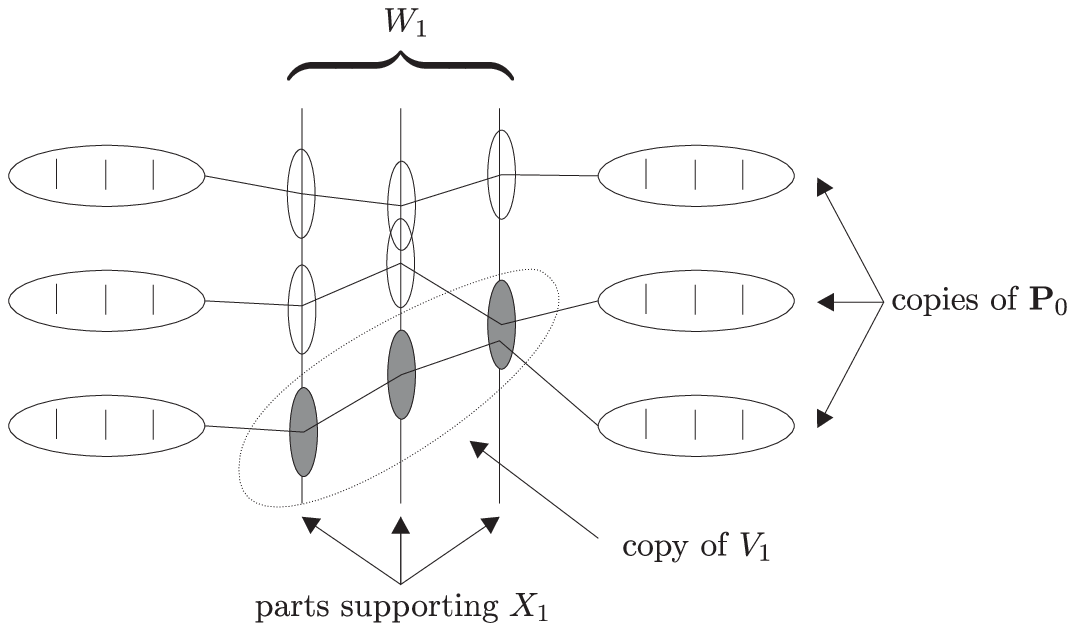}
\caption{Construction of $\m{P}_1$ from $\m{P}_0$.}\label{figNes4}
\end{figure}
\end{center}


It should be clear that associated to $\m{P}_1$ is a natural projection $\pi _1$ from $P_1$ onto $Z$. This allows to define the parts and the ordering on $\m{P}_1$.   

\begin{claim}
$\m{P}_1$ is $(l+1)$-metric. 
\end{claim}

\begin{proof}
Let $x_0,\ldots, x_{l+1}$ be a path in $\m{P} _1$ such that $(x_0 , x_{l+1}) \in \dom (\lambda ^{\m{P} _1})$. We want \[ \lambda ^{\m{P} _1}(x_0 , x_{l+1}) \leqslant \sum _{k=0} ^{l} \lambda ^{\m{P} _1}(x_k , x_{k+1}).\]

Or equivalently \[ \lambda ^{\m{Z}}(\pi_1 (x_0) , \pi_1 (x_{l+1})) \leqslant \sum _{k=0} ^{l} \lambda ^{\m{Z} }(\pi_1(x_k) , \pi_1(x_{k+1})).\]

Since $\m{Z}$ is $l$-metric, the only case to consider is when the only connections occuring between elements of the projection of the path are $(\pi _1 (x_0) , \pi _1 (x_{l+1}))$ and those of the form $(\pi _1 (x_k) , \pi _1 (x_{k+1}))$ where $k \leqslant l$. Since both $\m{W} _1$ and $\m{P} _0$ are $(l+1)$-metric, it is enough to show that the path either stays in $\m{W}_1$, or stays in a fixed copy $\m{P}$ of $\m{P} _0$. So suppose that the path leaves $\m{W}_1$. Using a circular permutation, we may reenumerate the path such that $x_0 \in \m{P} \smallsetminus \m{W}_1$. It follows then that $x_{l+1}$ is also in $\m{P}$. Now, assume now that for some $k$, $x_k \notin \m{P}$. Find $a<k<b$ such that $x_a , x_b \in \m{W} _1$. Observe that because $\pi_1 '' \m{W}_1$ is a copy of $\m{X}$ in $\m{Z}$ (namely $\m{X} _1$), $\pi _1 (x_a)$ and $\pi _1 (x_b)$ are connected. But this is a contradiction: Since $x_0 \notin \m{W}_1$, $\pi _1 (x_0) \notin \{ \pi _1 (x_a) , \pi _1 (x_b) \}$ and so $(\pi _1 (x_a) , \pi _1 (x_b)) \neq (\pi _1 (x_0) , \pi _1 (x_{l+1}))$. On the other hand $a+1 \neq b$. 
\end{proof}

In general, to build $(\m{P} _{i+1},<^{\m{P} _{i+1}} )$ from $(\m{P} _i,<^{\m{P} _i} )$, simply repeat the same procedure: Consider $\overleftarrow{\pi _i}\m{X}_{i+1}$. The ordered edge-labelled graph $(\m{V} _{i+1},<^{\m{V} _{i+1}} )$ induced on this set is a lifting of $(\om{X})$. Apply lemma \ref{lem:partite} to get a lifting $(\m{W} _{i+1},<^{\m{W} _{i+1}} )$ of $(\om{X})$ such that

\begin{center}
$(\m{W} _{i+1},<^{\m{W} _{i+1}}) \arrows{(\m{V} _{i+1},<_{\m{V} _{i+1}})}{(\om{X})}{2}$. 
\end{center}

By strong amalgamation property, extend every element of $\binom{\m{W} _{i+1},<^{\m{W} _{i+1}}}{\m{V} _{i+1},<^{\m{V} _{i+1}}}$ to a copy of $(\m{P} _i,<^{\m{P} _i})$. The resulting finite edge-labelled graph is $\m{P}_{i+1}$. The parts and the ordering on $\m{P}_{i+1}$ are defined according to the natural projection $\pi _{i+1}$ from $\m{P}_{i+1}$ onto $\m{Z}$. $\m{P}_{i+1}$ then becomes a lifting of $\m{Z}$, and one can show that it is $(l+1)$-metric. We now finish the proof by showing that  

\begin{center}
$(\m{P} _q,<^{\m{P} _q} ) \arrows{(\om{Y})}{(\om{X})}{2} $.
\end{center}

For the sake of clarity, we temporarily drop mention of the linear orderings attached to the edge-labelled graphs under consideration. 

Let $\chi : \funct{\binom{\m{P} _q}{\m{X}}}{2}$. We want to find $\mc{Y} \in \binom{\m{P} _q}{\m{Y}} $ such that $\binom{\mc{Y}}{\m{X}}$ is monochromatic.
$\chi$ induces a coloring $\chi : \funct{\binom{\m{W}_q}{\m{X}}}{2}$ and by construction:

\begin{center}
$\m{W}_q \arrows{(\m{V}_q)}{\m{X}}{2}$
\end{center}

Thus, there is a copy $\mc{V}_q$ of $\m{V}_q$ in $\m{W}_q$ so that $\binom{\mc{V}_q}{\m{X}}$ is monochromatic.
Now, when constructing $\m{P} _q $ from $\m{P} _{q-1}$, $\mc{V}_q$ was extended to
$\mc{P} _{q-1} \in \binom{\m{P} _q}{\m{P} _{q-1}}$ for which
$\chi$ induces $\chi : \funct{\binom{\mc{P} _{q-1}}{\m{X}}}{2}$.
Notice that $\mc{V}_q$ is exactly $\mc{P} _{q-1} \cap \overleftarrow{\pi _{q-1}} \m{X}_q$, the subgraph of $\mc{P} _{q-1}$ projecting in $\m{Z}$ onto $\m{X}_q$. $\binom{\mc{V}_q}{\m{X}}$ being monochromatic, every two copies
of $\m{X}$ in $\mc{V} _q$ projecting in $\m{Z}$ onto $\m{X}_q$ have the same color.

Now, consider the natural copy $\mc{W}_{q-1}$ of $\m{W}_{q-1}$ in $\mc{P}_{q-1}$.
$\chi$ induces a $2$-coloring of $\binom{\mc{W}_{q-1}}{\m{X}}$ and $\m{W} _{q-1}$ was chosen so that 

\begin{center}
$\m{W}_{q-1} \arrows{(\m{V}_{q-1})}{\m{X}}{2}$.
\end{center}

Therefore, there is a copy $\mc{V}_{q-1}$ of $\m{V}_{q-1}$ in $\mc{W}_{q-1}$ so that $\binom{\mc{V}_{q-1}}{\m{X}}$ is monochromatic.
Now, knowing how $\m{P} _{q-1}$ is constructed from $\m{P} _{q-2}$, observe that $\mc{V}_{q-1}$ extends to a copy $\mc{P} _{q-2}$ of $\m{P}_{q-2}$ inside $\mc{P}_{q-1}$, with respect to which $\chi$ induces: 

\begin{center}
$\chi : \funct{\binom{\mc{P}_{q-2}}{\m{X}}}{2}$.
\end{center}

As previously, $\mc{V}_{q-1}$ is exactly $\mc{P} _{q-2} \cap \overleftarrow{\pi _{q-2}} \m{X}_{q-1}$, the subgraph of $\mc{P} _{q-2}$ projecting onto $\m{X}_{q-2}$. $\binom{\mc{V}_{q-1}}{\m{X}}$ being monochromatic, every two copies
of $\m{X}$ in $\mc{V} _{q-1}$ projecting in $\m{Z}$ onto $\m{X}_{q-1}$ have the same color. Keep in mind that thanks to the companion result at the previous step, the same holds for those copies
of $\m{X}$ in $\mc{V} _{q-1}$ projecting in $\m{Z}$ onto $\m{X}_q$.  

By repeating this argument $q$ times, we end up with a copy $\mc{P} _0$ of
$\m{P}_0$ in $\m{P} _q$ so that given any $k \in \{1, \ldots ,q \} $,
any two copies of $\m{X}$ in $\mc{P} _0$ projecting in $\m{Z}$ onto $\m{X} _k$ have the same color.
From $\chi$, we can consequently construct a coloring 

\begin{center}
$\widehat{\chi} :
\funct{\{\m{X}_1, \ldots, \m{X}_q \} = \binom{\m{Z}}{\m{X}}}{2} $.
\end{center}

The color $\widehat{\chi}(\m{X}_k)$ is simply the common color of all the copies of $\m{X}$ in $\mc{P}_0$ projecting onto $\m{X} _k$. Now, remember that $\m{Z}$ was chosen so as to satisfy: 

\begin{center}
$\m{Z} \arrows{(\m{Y})}{\m{X}}{2}$.
\end{center}

Thus, there is
$\beta \in \binom{\m{Z}}{\m{Y}}$ such that $\binom{\beta}{\m{X}}$ is $\widehat{\chi}$-monochromatic. At the level of $\mc{P}_0$ and $\chi$, this means that
all the copies of $\m{X}$ in $\mc{P}_0$ projecting in $\m{Z}$ onto a subset of $\beta$ have the same color. But by construction, the subgraph of $\mc{P}_0$ projecting onto $\beta$ includes a copy $\m{Y}$, namely $\m{Y}_{\beta}$. $\m{Y}_{\beta}$ is consequently an
element of $\binom{\m{P} _q}{\m{Y}}$ for which $\binom{ \m{Y} _{\beta}}{\m{X}}$ is
monochromatic. This proves the claim, and finishes the proof of the theorem.
\end{proof}

In fact, the previous proof allows to prove a slightly more general result. For $S \subset ]0, +\infty[$, let $\M ^< _S$\index{$\M ^< _S$} denote the class of all finite ordered metric spaces with distances in $S$.  

\begin{thm}[Ne\v{s}et\v{r}il \cite{N1}]

\label{thm:variation RP for MM}
\index{Ne\v{s}et\v{r}il!theorem on $\M ^< _S$}

Let $T \subset ]0,+ \infty[$ be closed under sums and $S$ be an initial segment of $T$. Then $\M ^< _S$ has the Ramsey property. 

\end{thm}

It follows that in particular, the classes $\M ^< _{\Q}$, $\M ^< _{\Q \cap ]0 , r]}$ with $r > 0$  in $\Q$, $\M ^< _{\omega}$ and $\M ^< _{\omega \cap ]0 , m ]}$ with $m > 0$ in $\omega$ are Ramsey. Let us mention here that the assumption on the behavior of $S$ with respect to sums is not superficial. We will see in the next two subsections that when this requirement is not fulfilled, the situation is pretty different.

\subsection{Finite convexly ordered ultrametric spaces.}

\label{subsection:RP for UU}

The purpose of this subsection is to provide another example of a Ramsey class. Let $\m{X}$ be an ultrametric space. Call a linear ordering $<$ on $\m{X}$ \textit{convex} when all the metric balls of
$\m{X}$ are $<$-convex. For $S \subset ]0, + \infty [$, let $\UU$\index{$\UU$} denote the class of
all finite convexly ordered ultrametric spaces with distances in $S$.

\begin{thm}

\label{thm:RP for UU}

Let $S \subset ]0, + \infty [$. Then $\UU$ has the Ramsey property.
\end{thm}

To prove this result, we first need some notations for the partition
calculus on trees.
Given trees $(\ot{T})$ and $(\ot{S})$ as described in chapter 1, section 2.2, say that
they are \textit{isomorphic} when there is a bijection between them which
preserves both the structural and the lexicographical orderings.
Also, given a tree $(\ot{U})$, set: \[ \binom{\ot{U}}{\ot{T}} = \{ (\otc{T}) : \widetilde{\m{T}} \subset \m{U} \ \mathrm{and} \ (\otc{T}) \cong (\ot{T}) \}.\]

Now, if $(\ot{S}), (\ot{T})$ and $(\ot{U})$ are trees, the symbol \[ (\ot{U}) \arrows{(\ot{T})}{(\ot{S})}{k} \] abbreviates the statement: 

\begin{center} 
For any $\chi : \funct{\binom{\ot{U}}{\ot{S}}}{k}$ there is
$(\otc{T}) \in \binom{\ot{U}}{\ot{T}}$, $i<k$, such that:

$\chi '' \binom{\otc{T}}{\ot{S}} = \{ i \}$.
\end{center}

\begin{lemma}

\label{lemma:2}

Given an integer $k \in \omega \smallsetminus \{ 0 \}$, a finite tree $(\ot{T})$ and a subtree $(\ot{S})$ of $(\ot{T})$ such that
$\mathrm{ht}(\m{T}) = \mathrm{ht}(\m{S})$, there is a finite tree $(\ot{U})$ such that $\mathrm{ht}(\m{U}) = \mathrm{ht}(\m{T})$
and $(\ot{U}) \arrows{(\ot{T})}{(\ot{S})}{k}$.
\end{lemma}

A natural way to proceed is by induction on the height $\mathrm{ht}(\m{T})$ of $\m{T}$. But it is so natural that after having done so, we realized that this method had already been
used in \cite{F} where the exact same
result is obtained. Consequently, we choose to provide a different proof
which uses the notion of ultrafilter-tree.

\begin{proof}

For the sake of clarity, we sometimes do not mention the lexicographical
orderings explicitly. For example, $\m{T}$ stands for $(\ot{\m{T}})$. So let $\m{T} \subset \m{S}$ be some finite trees of height $n$ and set
$\m{U}$ be equal to $\omega^{\leqslant n}$. $\m{U}$ is naturally
lexicographically ordered. To prove the theorem, we only need to prove
that $\m{U} \arrows{(\m{T})}{\m{S}}{k}$. Indeed, even though $\m{U}$
is not finite, a standard compactness argument can take us to the
finite. 

Let $\{s_i : i < |\m{S}| \}_{<^{\m{S}}_{lex}}$ be a strictly $<^{\m{S}} _{lex}$-increasing enumeration
of the elements of $\m{S}$ and define $f : \funct{|\m{S}|}{|\m{S}|}$ such that:

\vspace{0.5em}
\hspace{1em}
i) $f(0) = 0$.

\vspace{0.5em} 
\hspace{1em}
ii) $s_{f(i)}$ is the immediate $<^{\m{S}}$-predecessor of
$s_i$ in $\m{S}$ if $i>0$. 

\vspace{0.5em}
Similarly, define $g : \funct{|\m{T}|}{|\m{T}|}$ for
$\m{T} = \{t_j : j < |\m{T}| \}_{<^{\m{T}}_{lex}}$.
Let also \[\mathscr{S} = \{ X \subset \m{U} : X \sqsubset \m{S} \} \ \ (\mathrm{resp.} \ \mathscr{T} = \{ X \subset \m{U} : X \sqsubset \m{T} \}), \] where $X \sqsubset \m{S}$ means that $X$ is a $<^{\m{U}}
_{lex}$-initial segment of some $\mc{S} \cong \m{S}$. $\mathscr{S}$ (resp. $\mathscr{T}$) has a natural tree structure with respect to
$<^{\m{U}}
_{lex}$-initial segment, has height $|\m{S}|$ (resp. $|\m{T}|$) and \[\mathscr{S} ^{max} = \binom{\m{U}}{\m{S}} \ \  (\mathrm{resp.} \ \mathscr{T} ^{max} = \binom{\m{U}}{\m{T}}).\] 

Now, for $x$ in $\m{U}$, let $\mathrm{IS}_{\m{U}}(x)$ denote the set of immediate
$<^{\m{U}}$-successors of $x$ in $\m{U}$. Then observe that if $X \in \mathscr{S}
\smallsetminus \mathscr{S} ^{max}$ is enumerated as $\{x_i : i < |X| \}_{<^{\m{U}} _{lex}}$ and $u \in \m{U}$ such that $X
<^{\m{U}} _{lex} u$ (that is $x <^{\m{U}} _{lex} u$ for every $x \in X$), then:

\begin{center}
$X \cup \{ u \} \in \mathscr{S}$ iff $u \in \mathrm{IS}_{\m{U}}(x_{f(|X|)})$. 
\end{center}

Consequently, $X, X' \in
\mathscr{S} \smallsetminus \mathscr{S} ^{max}$ can be simultaneously extended
in $\mathscr{S}$ iff: 

\begin{center}
 $x_{f(|X|)} = x'_{f(|X'|)}$. 
\end{center}

Now, for $u \in
\m{U}$, let $\mathcal{W} _u$ be a non-principal ultrafilter on
$\mathrm{IS}_{\m{U}}(u)$ and for every $X \in \mathscr{S}
\smallsetminus \mathscr{S}^{max}$, let $\mathcal{V} _X = \mathcal{W}
_{x_{f(|X|)}}$. Hence, $\mathcal{V} _X$ is an ultrafilter on the set
of all elements $u$ in $\m{U}$ which can be used to extend $X$ in
$\mathscr{S}$. Let $\mathcal{S}$ be a \textit{$\vec{\mathcal{V}}$-subtree} of
  $\mathscr{S}$, that is, a subtree such that for
  every $X \in \mathcal{S} \smallsetminus \mathscr{S} ^{max}$:
  
\begin{center}
$\{ u \in \m{U} : X <^{\m{U}} _{lex} u \ \mathrm{and} \ X \cup \{ u \} \in \mathcal{S}
  \} \in \mathcal{V} _X$.    
\end{center}

\begin{claim}
There is $\mc{T} \in \binom{\m{U}}{\m{T}}$ such that
$\binom{\mc{T}}{S} \subset \mathcal{S} ^{max}$.
\end{claim}

For $X \in \mathcal{S}$, let:

\begin{center}
 $U_X = \{ u
  \in \m{U} : X <^{\m{U}} _{lex} u \ \mathrm{and} \ X \cup \{ u \} \in \mathcal{S}
  \}$.
\end{center} 
  
The tree $\mc{T}$ is constructed inductively. Start with $\tau _0 =
\emptyset$. Generally, suppose that $\tau _0 <^{\m{U}} _{lex} \ldots
<^{\m{U}} _{lex} \tau _j$ were constructed such that:

\begin{center}
$\forall X \subset \{ \tau _0 , \ldots , \tau _j \}, \ \ X \in \mathscr{S} \rightarrow X \in \mathcal{S}$. 
\end{center}

Consider now the family $\mathcal{I}$ defined by:

\begin{center}
$\mathcal{I} = \{ I \subset \{ 0, \ldots , j \}
: \{ t_i : i \in I \} \cup \{ t_{j+1} \} \sqsubset \m{S} \}$
\end{center}

For
$I \in \mathcal{I}$ let: \[X_I = \{ \tau _i : i \in I \}.\] 

The family $(X_I)_{I \in \mathcal{I}}$ is consequently the family of all elements of $\mathscr{S}$ which
need to be extended with $\tau _{j+1}$. In other
words, we have to choose $\tau _{j+1} \in \m{U}$ such that: 

\vspace{0.5em}
\hspace{1em} 
i) $\{ \tau _0 , \ldots , \tau _{j+1} \} \in \mathscr{T}$. 

\vspace{0.5em}
\hspace{1em} 
ii) $X_I \cup \{ \tau _{j+1} \} \in
\mathcal{S}$ for every $I \in \mathcal{I}$. 

\vspace{0.5em}

To do that, notice that
for any $u \in \m{U}$ which satisfies $\tau _j <^{\m{U}} _{lex} u$, we
have:

\begin{center}
$\{ \tau _0 , \ldots , \tau _j , u \} \in \mathscr{T}$ iff $u \in \mathrm{IS}
_{\m{U}} (\tau _{g(j+1)})$.
\end{center}

Now, for any such $u$ and any $I \in
\mathcal{I}$, we have $X_I \cup \{ u \} \in \mathscr{S}$ ie $u$
allows a simultaneous extension of all the elements of $\{ X_I : I \in
\mathcal{I} \}$. Consequently, $\mathcal{V} _{X_I}$ does not depend on
$I \in \mathcal{I}$. Let $\mathcal{V}$ be the corresponding common
value. For every $I \in \mathcal{I}$, we have $U_{X_I} \in
\mathcal{V}$ so one can pick $\tau _{j+1}$ such that:

\[ \tau _j <^{\m{U}} _{lex} \tau _{j+1} \in \bigcap_{I \in
  \mathcal{I}} U_{X_I}. \]
  
Then $\tau _{j+1}$ is as required. Indeed, on the one
hand, because $\tau _{j+1} \in \mathrm{IS} _{\m{U}} (\tau _{g(j+1)})$:

\begin{center}
$\{ \tau _0 , \ldots , \tau _{j+1} \} \in \mathscr{T}$. 
\end{center}

On the other hand, since $\tau _{j+1} \in U_{X_I}$, 

\begin{center}
$X_I \cup \{ \tau _{j+1} \} \in \mathcal{S}$ for every $I \in \mathcal{I}$. 
\end{center}

At the end of the
construction, we are left with $\mc{T} := \{ \tau _j : j \in |\m{T}|
\} \in \mathscr{T}$ such that: \[\binom{\mc{T}}{S} \in \mathcal{S} ^{max}.\]

The claim is proved. The proof of the lemma will be complete if we prove the following claim:

\begin{claim}
Given any $k \in \omega \smallsetminus \{ 0 \}$ and
any $\chi : \funct{\binom{\m{U}}{\m{S}}}{k}$, there is a $\vec{\mathcal{V}}$-subtree
$\mathcal{S}$ of $\mathscr{S}$ such that $\mathcal{S} ^{max}$ is
$\chi$-monochromatic.
\end{claim}

We proceed by induction on the height of $\mathscr{S}$. The case
$\mathrm{ht}(\mathscr{S}) = 0 $ is trivial so suppose that the claim
holds for $\mathrm{ht}(\mathscr{S}) = n$ and consider the case
$\mathrm{ht}(\mathscr{S}) = n+1$. Define a coloring $\Lambda :
\funct{\mathscr{S}(n)}{k}$ by:

\begin{center}
$\Lambda (X) = \varepsilon$ iff $\{ u
\in \m{U} : X \cup \{ u \} \in \mathscr{S}(n+1) \ \mathrm{and} \ \chi (X \cup \{ u
\}) = \varepsilon \} \in \mathcal{V} _X$. 
\end{center}

By induction hypothesis, we can find a
$\vec{\mathcal{V}}$-subtree $\mathcal{S} _n$ of $\restrict{\mathscr{S}}{n}$ (the tree
formed by the $n$ first levels of $\mathscr{S}$) such that
$\mathcal{S} _n ^{max} $ is $\Lambda$-monochromatic with color
$\varepsilon _0$. This means that for every $X \in \mathcal{S} _n$,
the set $V_X$ is in $ \mathcal{V} _X$, where $V_X$ is defined by: \[ V_X := \{ u \in \m{U} : X \cup \{ u \} \in \mathscr{S}(n+1) \ \mathrm{and} \ \chi (X \cup \{ u
\}) = \varepsilon _0 \}.\]

Now, let: \[\mathcal{S} = \mathcal{S} _n \cup \{ X \cup \{ u \} : X \in
\mathcal{S} _n \ \mathrm{and} \ u \in V_X \}.\]

Then $\mathcal{S}$ is a
$\vec{\mathcal{V}}$-subtree of $\mathscr{S}$ and $\mathcal{S} ^{max}$
is $\chi$-monochromatic.
\end{proof}

We now show how to obtain Theorem \ref{thm:RP for UU} from Lemma
\ref{lemma:2}. Fix $S \subset ]0, + \infty [$, let $(\om{X})$, $(\om{Y}) \in \UU$ and consider $(\ot{T})$ associated to
$(\om{Y})$. As presented in section 2, $(\om{Y})$ can be seen as
$(\textbf{T} ^{max}, <^{\textbf{T}} _{lex})$. Now, notice that there
is a subtree $(\ot{S})$ of $(\ot{T})$ such that for every $(\oc{X})
\in \binom{\textbf{T}^{max} , <^{\textbf{T}} _{lex} }{\om{X}}$, the
downward $<^{\textbf{T}}$-closure of $\mc{X}$ is isomorphic to
$(\ot{S})$. Conversely, for any $(\otc{S})$ in $\binom{\ot{T}}{\ot{S}}$, $(\widetilde{\textbf{S}} ^{max} ,
<^{\widetilde{\textbf{S}}} _{lex})$ is in $\binom{\textbf{T}^{max} ,
  <^{\textbf{T}} _{lex} }{\om{X}}$. These facts allow us to build $(\om{Z})$ such that: \[(\om{Z}) \arrows{(\om{Y})}{(\om{X})}{k}.\] 
  
Indeed, apply Lemma \ref{lemma:2} to
  get $(\ot{U})$ of height $\mathrm{ht}(\m{T})$ such that:
  
\begin{center}  
$(\ot{U}) \arrows{(\ot{T})}{(\ot{S})}{k}$. 
\end{center}

Then, simply let $(\om{Z})$ be
  the convexly ordered ultrametric space associated to $(\ot{U})$. To
  check that $(\om{Z})$ works, let: 
  
\begin{center}  
$\chi : \funct{\binom{\om{Z}}{\om{X}}}{k}$.
\end{center}

The map $\chi$ transfers to:

\begin{center}
$\Lambda : \funct{\binom{\ot{U}}{\ot{S}}}{k}$.
\end{center}

Thus, we can find $(\otc{T}) \in
  \binom{\ot{U}}{\ot{T}}$ such that $\binom{\otc{T}}{\ot{S}}$ is
  $\Lambda$-monochromatic. Then the convexly ordered ultrametric space $(\widetilde{\textbf{T}} ^{max}, <^{\widetilde{\textbf{T}}} _{lex})$
 is such that $\binom{\widetilde{\textbf{T}} ^{max},
   <^{\widetilde{\textbf{T}}} _{lex}}{\om{X}}$ is
 $\chi$-monochromatic. But $(\widetilde{\textbf{T}} ^{max},
 <^{\widetilde{\textbf{T}}} _{lex}) \cong (\om{Y})$. Theorem
 \ref{thm:RP for UU} is proved.  

\vspace{1em}
 
\textbf{Remark.} We will see later in this chapter that unlike $\UU$, the class $\mathcal{U} ^< _S$\index{$\mathcal{U} ^< _S$} of all finite ordered ultrametric spaces with distances in $S$ \emph{does not} have the Ramsey property.

\subsection{Finite metrically ordered metric spaces.}

The results of the two previous sections suggest that the metric structure of the spaces under consideration strongly influences the kind of linear orderings to be adjoined in order to get a Ramsey-type result. The present subsection can be seen as an illustration of that fact. Let $\mathcal{K}$ be a class of metric spaces. For $s \in ]0 , + \infty[$ and $\m{X} \in \mathcal{K}$, let $\approx ^{\m{X}} _s$ be the binary relation defined on $\m{X}$ by:

\begin{center}
$\forall x, y \in \m{X} \ \ x \approx ^{\m{X}} _s y \leftrightarrow d^{\m{X}}(x,y) \leqslant s$. 
\end{center}

Say that $s$ is \emph{critical for $\mathcal{K}$}\index{critical distance} when for every $\m{X} \in \mathcal{K}$, $\approx ^{\m{X}} _s$ is an equivalence relation on $\m{X}$. On the other hand, given $\m{X} \in \mathcal{K}$, say that a binary relation $R$ is a \emph{metric equivalence relation}\index{metric!equivalence relation} on $\m{X}$ when there is $s \in ]0 , + \infty[$ critical in $\mathcal{K}$ such that $R = \approx ^{\m{X}} _s$. For example, for the classes $\M _S$\index{$\M _S$}, any $s \in S$ such that $]s , 2s] \cap S = \emptyset$ is critical. Of course, when $S$ is finite, $\max S$ is always critical, but there might be other critical distances. For instance, $2$ is critical for $\M _{ \{ 1, 2, 5\}}$, $1$ is critical for $\M _{ \{ 1, 3, 4\}}$ and for $\M _{ \{ 1, 3, 6\}}$. On the other hand, given $S \subset ]0, + \infty[$, any $s \in S$ is critical for $\U$\index{$\U$}. 

Now, call a linear ordering $<$ on $\m{X} \in \mathcal{K}$ \emph{metric}\index{metric!linear ordering} if given any metric equivalence relation $\approx$ on $\m{X}$, the $\approx$-equivalence classes are $<$-convex. Given $S \subset ]0, + \infty[$, let $\nM _S$\index{$\nM _S$} denote the class of all finite metrically ordered metric spaces with distances in $S$.

\begin{thm}

\label{thm:RP for nM_S}

Let $S$ be finite subset of $]0, + \infty [$ of size $|S| \leqslant 3$ and satisfying the $4$-values condition. Then $\nM _S$ has the Ramsey property. 
\end{thm}

\begin{proof}
The case $|S|=1$ is trivial. 
Recall that for $|S|=2$, there are essentially two cases, namely $S=\{1,2\}$ and $S=\{1,3\}$. When $\m{X} \in \M _{\{1,2\}}$, all the linear orderings on $\m{X}$ are metric so $\nM _{\{1,2\}} = \M ^< _{\{ 1,2\}}$ is a Ramsey class thanks to Theorem \ref{thm:variation RP for MM}. On the other hand, when $\m{X} \in \M _{\{1,3\}}$, $\m{X}$ is ultrametric and the metric linear orderings on $\m{X}$ are the convex ones. Thus, $\nM _{\{1,3\}} = \mathcal{U} ^{c<} _{\{ 1,3\}}$ and has the Ramsey property thanks to Theorem \ref{thm:RP for UU}. 
For $|S|=3$, the cases to consider are:

\begin{center}
(1a) \{2, 3, 4\} \ \ (1b) \{1, 2, 3\} \ \ (1d) \{1, 2, 5\}

(2a) \{1, 3, 4\} \ \ (2b) \{1, 3, 6\} \ \ (2c) \{1, 3, 7\}

\end{center}

(1a) and (1b) are covered by Theorem \ref{thm:variation RP for MM}. (2c) is covered by Theorem \ref{thm:RP for UU}. The remaining cases could be treated one by one but in what follows, we cover them all at once thanks to the following lemma. Let $T :=\{1, 2, 5, 6, 9 \}$. Then: 

\begin{lemma}

\label{lem:RP for nM_T}

$\nM _T$ has the Ramsey property. 
\end{lemma}

\begin{proof}
For $(\om{X}) \in \nM _T$, let $\mathcal{B}_{\m{X}}$ be the set of all balls of $\m{X}$ of radius $2$. Define an ordered graph $(\m{G}_{\m{X}} , <^{\m{G}_{\m{X}}})$ as follows: The set of vertices of $\m{G}_{\m{X}}$ is given by \[ G_{\m{X}} = \bigcup _{b \in \mathcal{B} _{\m{X}}} \{ v ^{\m{X}} _b\} \cup \{ \pi ^{\m{X}} (x) : x \in b\}. \]

The linear ordering $<^{\m{G} _{\m{X}}}$ is such that 

\vspace{0.5em}
\hspace{1em}
i) $v ^{\m{X}} _b <^{\m{G} _{\m{X}}} \{ \pi ^{\m{X}} (x) : x \in b\} <^{\m{G} _{\m{X}}} v ^{\m{X}} _{b'}$ whenever $b < ^{\m{X}} b'$.  

\vspace{0.5em}
\hspace{1em}
ii) $\pi ^{\m{X}}$ is order-preserving.     

\vspace{0.5em} 

The set $E(\m{G}_{\m{X}})$ of edges of $\m{G}_{\m{X}}$ is such that:

\vspace{0.5em}
\hspace{1em}
i) $\{ v ^{\m{X}} _b , v ^{\m{X}} _{b'} \} \in E(\m{G}_{\m{X}})$ iff $(\forall x \in b \ \forall x' \in b' \ d^{\m{X}}(x,x') \in \{ 5, 6\})$. 

\vspace{0.5em}
\hspace{1em} 
ii) For every $b \in \mathcal{B}_{\m{X}}$ and $x \in \m{X}$, $\{ v ^{\m{X}} _b , \pi ^{\m{X}} (x) \} \in E(\m{G}_{\m{X}})$ iff $x \in b$.

\vspace{0.5em}
\hspace{1em}  
iii) $ \{ \pi ^{\m{X}} (x) , \pi ^{\m{X}} (x') \} \in E(\m{G}_{\m{X}})$ iff $d^{\m{X}} (x,x') \in \{ 1, 5\}$.

\vspace{0.5em} 

The construction of $\m{G}_{\m{X}}$ from $\m{X}$ is illustrated in Figure \ref{figRP}.

\begin{center}
\begin{figure}[h]
\includegraphics[scale=0.68]{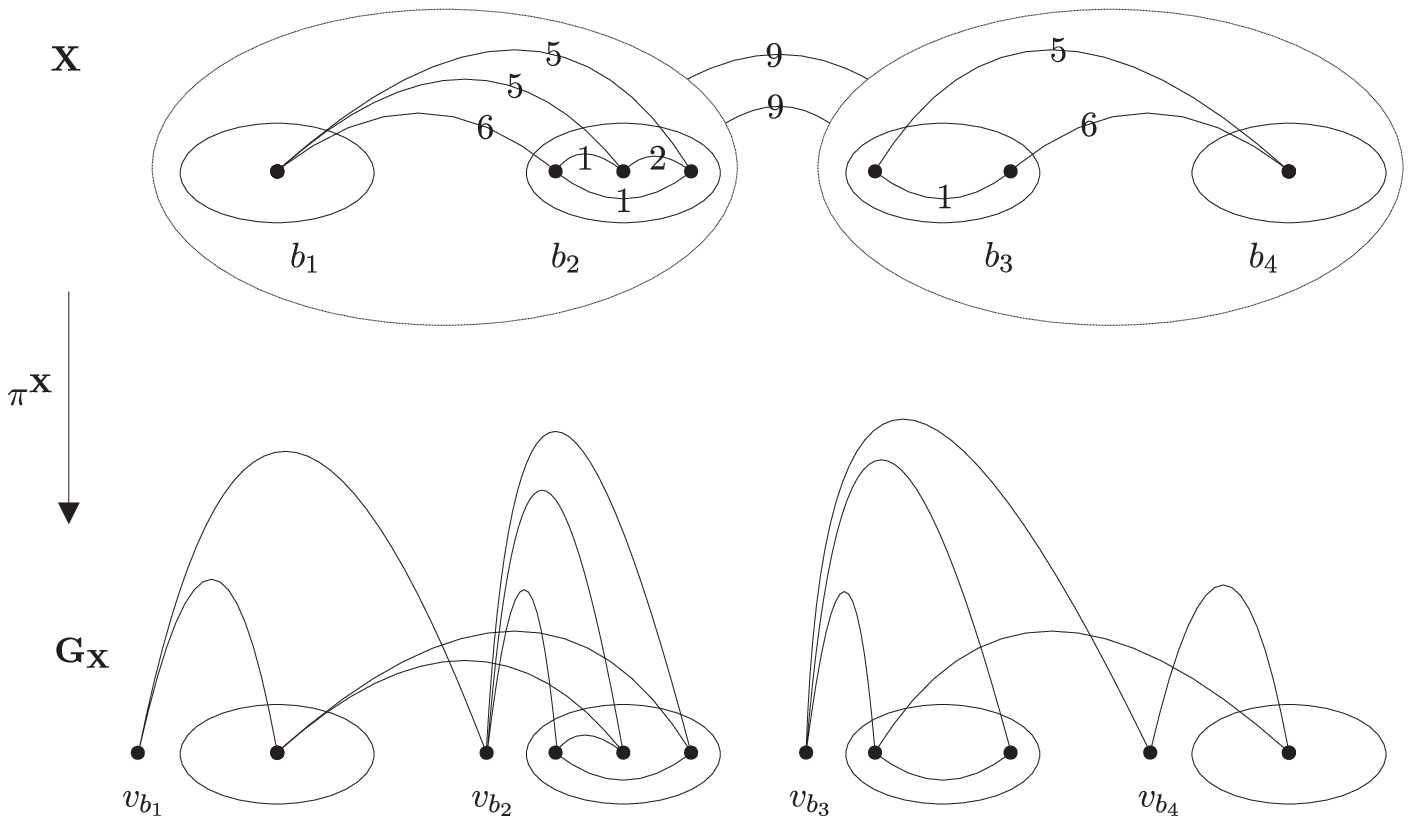}
\caption{Construction of $\m{G}_{\m{X}}.$ from $\m{X}$}\label{figRP}
\end{figure}
\end{center}

Now, define $d^{E(\m{G}_{\m{X}})} (\{ v, w \}_{<^{\m{G} _{\m{X}}}}, \{v', w' \}_{<^{\m{G} _{\m{X}}}})$ by:  

\begin{displaymath}
\left \{ \begin{array}{ll}
 1 & \textrm{if $v=v'$ and $\{ w, w' \} \in E(\m{G}_{\m{X}})$,} \\
 2 & \textrm{if $v=v'$ and $\{ w, w' \} \notin E(\m{G}_{\m{X}})$,} \\
 5 & \textrm{if $v \neq v'$ and $\{ v, v' \} \in E(\m{G}_{\m{X}})$ and $\{ w, w' \} \in E(\m{G}_{\m{X}})$,}\\
 6 & \textrm{if $v \neq v'$ and $\{ v, v' \} \in E(\m{G}_{\m{X}})$ and $\{ w, w' \} \notin E(\m{G}_{\m{X}})$,}\\
 9 & \textrm{if $v \neq v'$ and $\{ v, v' \} \notin E(\m{G}_{\m{X}})$.} 
 \end{array} \right.
\end{displaymath}

\begin{claim}
$d^{E(\m{G}_{\m{X}})}$ is a metric. 
\end{claim}

\begin{proof}
It is enough to show that the triangle inequality is satisfied. Take $\{ v, w \}_{<^{\m{G} _{\m{X}}}}$, $\{ v', w' \}_{<^{\m{G} _{\m{X}}}}$ and $\{ v'', w'' \}_{<^{\m{G} _{\m{X}}}}$ in $E(\m{G}_{\m{X}})$ and set 
\begin{displaymath}
\left \{ \begin{array}{l}
d^{E(\m{G}_{\m{X}})} (\{ v, w \}_{<^{\m{G} _{\m{X}}}}, \{v', w' \}_{<^{\m{G} _{\m{X}}}}) = \alpha \\ d^{E(\m{G}_{\m{X}})} (\{ v', w' \}_{<^{\m{G} _{\m{X}}}}, \{v'', w'' \}_{<^{\m{G} _{\m{X}}}}) = \beta \\ d^{E(\m{G}_{\m{X}})} (\{ v, w \}_{<^{\m{G} _{\m{X}}}}, \{v'', w'' \}_{<^{\m{G} _{\m{X}}}}) = \gamma
\end{array} \right. 
\end{displaymath}

We have to show that we are not in one of the following cases: ($\alpha , \beta \in \{ 1, 2\}$ and $\gamma \geqslant 5$) or ($\alpha \in \{ 1, 2\}$, $\beta \in \{ 5, 6\}$ and $\gamma = 9$). Assume that $\alpha , \beta \in \{ 1, 2\}$. Then $v=v'$ and $v'=v''$. Thus, $v=v''$ and $\gamma < 5$ so the first case is covered. For the second case, assume that $\alpha \in \{ 1, 2\}$ and $\beta \in \{ 5, 6\}$. Then $v=v'$ and $\{ v', v'' \} \in E(\m{G}_{\m{X}})$. It follows that $\{ v, v'' \} \in E(\m{G}_{\m{X}})$ and so $\gamma \neq 9$. 
\end{proof}

For $x \in \m{X}$, let $b(x)$ denote the only element $b$ of $\mathcal{B} _{\m{X}}$ such that $x \in b$ and define a map 
$\varphi _{\m{X}} : \funct{\m{X}}{E(\m{G}_{\m{X}})}$ by $\varphi _{\m{X}}(x) = \{ v ^{\m{X}} _{b(x)} , \pi ^{\m{X}} (x) \}$. Then it is easy to check that when $E(\m{G}_{\m{X}})$ is equipped with the lexicographical ordering:

\begin{claim}
$\varphi _{\m{X}}$ is an order-preserving isometry. 
\end{claim}

The map $(\om{X}) \mapsto (\m{G}_{\m{X}} , <^{\m{G}_{\m{X}}})$ consequently codes the ordered metric space $(\om{X})$ into the ordered graph $(\m{G}_{\m{X}} , <^{\m{G}_{\m{X}}})$. We now prove two essential properties of this coding. Let $(\om{Y})$ be a finite ordered metric space and $(\om{X})$ be a subspace of $(\om{Y})$.

\vspace{0.5em}
\hspace{1em}
1) Every copy of $(\om{X})$ in $(\om{Y})$ gives raise to a copy of $(\m{G}_{\m{X}} , <^{\m{G}_{\m{X}}})$ in 

\hspace{2em} $(\m{G}_{\m{Y}} , <^{\m{G}_{\m{Y}}})$. 
 
\vspace{0.5em}
\hspace{1em} 
2) Conversely, every copy of $(\m{G}_{\m{X}} , <^{\m{G}_{\m{X}}})$ in $(\m{G}_{\m{Y}} , <^{\m{G}_{\m{Y}}})$ codes a copy of 

\hspace{2em} $(\om{X})$ in $(\om{Y})$.

\vspace{0.5em} 

More precisely, for 1), let $(\om{Y}) \in \nM _T$. Thanks to the previous claim, we have: 

\begin{center}
$(\om{Y}) \cong ( \{ \{ v ^{\m{Y}} _{b(y)} , \pi ^{\m{Y}} (y) \} _{<^{\m{G} _{\m{Y}}}} : y \in \m{Y} \} , <_{lex}) =: (\oc{Y})$. 
\end{center}

\begin{claim}
Let $(\oc{X}) \in \binom{\oc{Y}}{\om{X}}$. Then $(\bigcup \mc{X} , \restrict{<^{\m{G} _{\m{Y}}}}{\bigcup \mc{X}}) \cong (\m{G}_{\m{X}} , <^{\m{G}_{\m{X}}})$. 
\end{claim}

\begin{proof}
Since $\varphi _{\m{Y}}$ is an order-preserving isometry, $\overleftarrow{\varphi _{\m{Y}}} \mc{X}$ supports a copy of $(\om{X})$ in $(\om{Y})$. Let $\psi : \funct{\m{X}}{\overleftarrow{\varphi _{\m{Y}}} \mc{X}}$ be the order-preserving isometry witnessing that fact. On the one hand:

\begin{displaymath}
\begin{array}{lll}
\bigcup \mc{X}&=& \{ v ^{\m{Y}} _{b(x)} : x \in \overleftarrow{\varphi _{\m{Y}}} \mc{X} \} \cup \{ \pi ^{\m{Y}}(x) : x \in \overleftarrow{\varphi _{\m{Y}}} \mc{X} \}\\
& = & \{ v ^{\m{Y}} _{b(\psi(x))} : x \in \m{X} \} \cup \{ \pi ^{\m{Y}}(\psi(x)) : x \in \m{X} \}.
\end{array}
\end{displaymath}
 
On the other hand: 

\begin{center}
$\m{G}_{\m{X}} = \{ v ^{\m{X}} _{b(x)} : x \in \m{X} \} \cup \{ \pi ^{\m{X}}(x) : x \in \m{X} \} $. 
\end{center}

Therefore, it is enough to check that the map defined by $v ^{\m{X}} _{b(x)} \mapsto v ^{\m{Y}} _{b(\psi(x))}$ and $\pi ^{\m{X}}(x) \mapsto \pi ^{\m{Y}}(\psi(x))$ for every $x \in \m{X}$ is an ordered graph isomorphism. The fact that the ordering is preserved is obvious. To verify that the edges are also preserved, we have to check that for every $x, x' \in \m{X}$:

\vspace{0.5em}
\hspace{1em}  
i) $\{ v ^{\m{X}} _{b(x)} , v ^{\m{X}} _{b(x')} \} \in E(\m{G}_{\m{X}})$ iff $\{ v ^{\m{Y}} _{b(\psi(x))} , v ^{\m{Y}} _{b(\psi(x'))} \} \in E(\m{G}_{\m{Y}}) $. 

\vspace{0.5em}
\hspace{1em}  
ii) $\{ v ^{\m{X}} _{b(x)} , \pi ^{\m{X}} (x') \} \in E(\m{G}_{\m{X}})$ iff $\{ v ^{\m{Y}} _{b(\psi(x))} , \pi ^{\m{Y}} (\psi(x')) \} \in E(\m{G}_{\m{Y}})$.

\vspace{0.5em}
\hspace{1em}  
iii) $\{ \pi ^{\m{X}} (x) , \pi ^{\m{X}} (x') \} \in E(\m{G}_{\m{X}})$ iff $\{ \pi ^{\m{Y}} (\psi(x)) , \pi ^{\m{Y}} (\psi(x')) \} \in E(\m{G}_{\m{Y}})$.   

\vspace{0.5em}

Let $x \neq x' \in \m{X}$. For i)

\begin{displaymath}
\begin{array}{lll}
\{ v ^{\m{X}} _{b(x)} , v ^{\m{X}} _{b(x')} \} \in E(\m{G}_{\m{X}}) & \leftrightarrow & d^{\m{X}}(x, x') \in \{ 5, 6\} \\
& \leftrightarrow & d^{\m{Y}}(\psi(x), \psi(x')) \in \{ 5, 6\} \\
& \leftrightarrow & \{ v ^{\m{Y}} _{b(\psi(x))} , v ^{\m{Y}} _{b(\psi(x'))} \} \in E(\m{G}_{\m{Y}})
\end{array}
\end{displaymath}

For ii)

\begin{displaymath}
\begin{array}{lll}
\{ v ^{\m{X}} _{b(x)} , \pi ^{\m{X}} (x') \} \in E(\m{G}_{\m{X}}) & \leftrightarrow & d^{\m{X}}(x, x') \in \{ 1, 2\} \\
& \leftrightarrow & d^{\m{Y}}(\psi(x), \psi(x')) \in \{ 1, 2\} \\
& \leftrightarrow & \{ v ^{\m{Y}} _{b(\psi(x))} , \pi ^{\m{Y}} (\psi(x')) \} \in E(\m{G}_{\m{Y}})
\end{array}
\end{displaymath} 

Finally, for iii)

\begin{displaymath}
\begin{array}{lll}
\{ \pi ^{\m{X}} (x), \pi ^{\m{X}} (x') \} \in E(\m{G}_{\m{X}}) & \leftrightarrow & d^{\m{X}}(x, x') \in \{ 1, 5\} \\
& \leftrightarrow & d^{\m{Y}}(\psi(x), \psi(x')) \in \{ 1, 5\} \\
& \leftrightarrow & \{ \pi ^{\m{Y}} (\psi(x))  , \pi ^{\m{Y}} (\psi(x')) \} \in E(\m{G}_{\m{Y}})
\end{array}
\end{displaymath} 
\end{proof}

For 2), we need to show how, given a copy of $(\m{G}_{\m{X}} , <^{\m{G}_{\m{X}}})$, one can reconstruct a 'natural' copy of $(\om{X})$. We proceed as follows: Let $(\om{G})$ be a copy of $(\m{G}_{\m{X}} , <^{\m{G}_{\m{X}}})$ and let $\sigma$ be an order-preserving graph isomorphism from $(\m{G}_{\m{X}} , <^{\m{G}_{\m{X}}})$ onto $(\om{G})$. Then the ordered 
metric subspace of $(E(\m{G} _{\m{X}}) , <_{lex} )$ supported by $\{ \{ \sigma (v ^{\m{X}} _{b(x)}) , \sigma (\pi ^{\m{X}}(x) ) \} : x \in \m{X} \} $ is isomorphic to $(\om{X})$. In the sequel, it will be denoted $\m{X} _{\m{G}}$ and will be called the \emph{natural} copy of $(\om{X})$ inside $(E(\m{G} _{\m{X}}) , <_{lex} )$.

We can now turn to a proof of the lemma. For the sake of clarity, we temporarily drop mention of the linear orderings attached to the graphs and the metric spaces under consideration. Let $\m{X}, \m{Y}$ be in $\nM _T$ and $k > 0$ be in $\omega$. Thanks to Ramsey property for the class of finite ordered graphs, find a finite ordered graph $\m{K}$ such that: \[ \m{K} \arrows{(\m{G}_{\m{Y}})}{\m{G}_{\m{X}}}{k}.\]

Now, let $\m{Z}$ be the ordered metric space $E(\m{K})$ equipped with the metric described previously and ordered lexicographically. We claim that:

\begin{center} 
$\m{Z} \arrows{(\m{Y})}{\m{X}}{k}$.  
\end{center}

Indeed, let $\chi : \funct{\binom{\m{Z}}{\m{X}}}{k}$. The map $\chi$ induces $\Lambda : \funct{\binom{\m{K}}{\m{G}_{\m{X}}}}{k}$ defined by 

\begin{center}
$\Lambda (\m{G}) = \chi (\m{X}_{\m{G}})$. 
\end{center}

Find $\widetilde{\m{G}_{\m{Y}}} \cong \m{G}_{\m{Y}}$ such that $\binom{\widetilde{\m{G}_{\m{Y}}}}{\m{G}_{\m{X}}}$ is $\Lambda$-monochromatic. Call its color $\varepsilon$ and let $\mc{Y}$ be the natural copy of $\m{Y}$ inside $E(\m{G}_{\m{Y}})$. Then $\binom{\mc{Y}}{\m{X}}$ is $\chi$-monochromatic: Indeed, if $\mc{X} \in \binom{\mc{Y}}{\m{X}}$, then by a previous claim $\bigcup \mc{X} \cong \m{G}_{\m{X}}$. It follows that $\chi (\mc{X}) = \Lambda (\bigcup \mc{X}) = \varepsilon$. This finishes the proof of the lemma. \end{proof}

We now deduce Theorem \ref{thm:RP for nM_S} from Lemma \ref{lem:RP for nM_T}. To show that $\nM _{\{1, 2, 5\}}$ has the Ramsey property, let $(\om{X})$, $(\om{Y})$ be in $\nM _{\{1, 2, 5\}}$. Then $(\om{X})$ are also $(\om{Y})$ in $\nM _T$ so we can find $(\om{Z})$ in $\nM _T$ such that \[ (\om{Z}) \arrows{(\om{Y})}{(\om{X})}{2}. \] 

Now, define a new metric $d^{\{1, 2, 5\}}$ on $Z$ by:

\begin{displaymath}
d^{\{1, 2, 5\}}(x,y) = \left \{ \begin{array}{ll}
 1 & \textrm{if \ $d^{\m{Z}}(x,y) = 1$} \\
 2 & \textrm{if \ $d^{\m{Z}}(x,y) = 2$} \\ 
 5 & \textrm{if \ $d^{\m{Z}}(x,y) \geqslant 5$ } 
 \end{array} \right.
\end{displaymath}

Then, observe that $(Z , d' , <^{\m{Z}})$ in $\nM _{\{1, 2, 5\}}$ is such that 

\begin{center}
$(Z , d' , <^{\m{Z}}) \arrows{(\om{Y})}{(\om{X})}{2}$.  
\end{center}

For $\nM _{\{1, 3, 4\}}$, the proof is the same except that $d^{\m{Z}}$ is not replaced by $d^{\{1, 2, 5\}}$ but by $d^{\{1, 3, 4\}}$ defined by:

\begin{displaymath}
d^{\{1, 3, 4\}}(x,y) = \left \{ \begin{array}{ll}
 1 & \textrm{if \ $d^{\m{Z}}(x,y) \in \{1 , 2 \}$} \\
 3 & \textrm{if \ $d^{\m{Z}}(x,y) = 5$} \\ 
 4 & \textrm{if \ $d^{\m{Z}}(x,y) \geqslant 6$} 
 \end{array} \right.
\end{displaymath}  

Finally, for $\nM _{\{1, 3, 6\}}$, replace $d^{\m{Z}}$ by $d^{\{1, 3, 6\}}$ defined by:

\begin{displaymath}
d^{\{1, 3, 6\}}(x,y) = \left \{ \begin{array}{ll}
 1 & \textrm{if \ $d^{\m{Z}}(x,y) \in \{1 , 2 \}$} \\
 3 & \textrm{if \ $d^{\m{Z}}(x,y) \in \{ 5, 6 \}$} \\ 
 6 & \textrm{if \ $d^{\m{Z}}(x,y) = 9$} 
 \end{array} \right. \qedhere
\end{displaymath}

\end{proof}

\section{Ordering properties.}

After Ramsey property, we turn to the study of ordering properties. As we will see, ordering property is usually much easier to prove than Ramsey property. 

\subsection{Finite ordered metric spaces.}

We start with a case for which the ordering property is a consequence of the Ramsey property. 

\begin{thm}

\label{thm:Ordering Property for metric spaces.}

$\M ^<$\index{$\M ^<$} has the ordering property. 

\end{thm}

\begin{proof}

Let $D$ be the largest distance appearing in $\m{X}$. Observe that $(\om{X})$ can be embedded into $(\oc{X})$ such that $(\oc{X})$ and $(\mc{X}, ^{\mc{X}}>)$ are isomorphic. There is consequently no loss of generality if we assume that $(\om{X})$ and $(\m{X}, ^{\m{X}}>)$ are isomorphic. We first construct $(\om{Z})$ including $(\om{X})$ as a subspace and such that given any $x < ^{\m{X}} y \in \m{X}$, there is $z \in \m{Z}$ such that:

\begin{center}
$x < ^{\m{Z}} z < ^{\m{Z}} y$ and $d^{\m{Z}} (x,z) = d^{\m{Z}} (z,y)$.
\end{center}

A way to obtain such an $(\om{Z})$ is to proceed as follows. Seeing $(\om{X})$ as a finite ordered edge-labelled graph, connect any two distinct points by a broken line consisting of two edges with label $D$. Observe that the corresponding edge-labelled graph is $l$-metric for every $l$ so the labelling can be extended using the shortest path distance. Therefore, the corresponding metric space $\m{Z}$ does include $\m{X}$ as a subspace. We now have to order $\m{Z}$. Take $x <^{\m{X}} y \in \m{X}$. When expanding $\m{X}$ to $\m{Z}$, a broken line $\{ x, z, y \}$ was added with $d^{\m{Z}}(x,z) = d^{\m{Z}}(y,z) = D$. Define a linear ordering $<^{ \{ x, y \}}$ on this line by:

\begin{center}
$x <^{ \{ x, y \}} z <^{ \{ x, y \}} y$. 
\end{center}


Now, concatenate all the orderings of the form $<^{ \{ x, y \}}$ according to the lexicographical ordering on the the set of edges $\{ \{ x, y \} _{<^{\m{X}}}: x, y \in X \}$ in order to obtain $<^{\m{Z}}$. Then, the finite ordered metric space $\m{Z}$ is as required. Now, let $(\om{T})$ be the unique ordered metric space with two points and distance $D$ between them, and let $(\om{Y})$ be such that:

\begin{center}
$(\om{Y}) \arrows{(\om{T})}{(\om{Z})}{2}$. 
\end{center}

\begin{claim}
Given any linear ordering $<$ on $\m{Y}$, $(\m{Y} , <)$ includes a copy of $(\om{X})$. 
\end{claim}
  
To prove that claim, let $<$ be a linear ordering on $\m{Y}$ and let $\chi : \funct{\binom{\om{Y}}{\om{T}}}{2}$ be such that:

\begin{center}
$\chi (\{ x, y\}) = 1$ iff $<^{\m{Y}}$ and $<$ agree on $\{ x, y \}$. 
\end{center}

By construction, we can find a copy $(\oc{Z})$ of $(\om{Z})$ in $(\om{Y})$ with $\binom{\om{Z}}{\om{T}}$ monochromatic. Call $\varepsilon$ the correspondong color. Now, let $(\oc{X})$ be a copy of $(\om{X})$ inside $(\oc{Z})$. 

\begin{subclaim}
$(\mc{X} , <) \cong (\om{X})$. 
\end{subclaim}

There are two cases, according to the value of $\varepsilon$. If $\varepsilon = 1$, we prove that given any $x, y \in \mc{X}$, $<$ and $<^{\m{X}}$ agree on $\{ x, y \}$. This will show $(\mc{X} , <) \cong (\oc{X})$. So let $x < ^{\mc{X}} y$. Find $z \in \mc{Z}$ such that $x <^{\mc{Z}} z <^{\mc{Z}} y$ and $d^{\mc{Z}} (x,z) = d^{\mc{Z}} (x,z) = D$. Since $\varepsilon = 1$, $<$ and $<^{\mc{Z}}$ agree on $\{ x, z\}$ and $\{ z, y\}$. Thus, $x<z<y$ and so $x<z$. If $\varepsilon = 0$, we prove that given any $x, y \in \mc{X}$, $<$ and $<^{\m{X}}$ disagree on $\{ x, y \}$. This will show $(\mc{X} , < ) \cong (\mc{X} , ^{\mc{X}} >)$ and since $(\mc{X} , ^{\mc{X}} >) \cong (\oc{X})$, we will get $(\mc{X} , <) \cong (\oc{X})$. Let $x < ^{\mc{X}} y$. Pick $z \in \mc{Z}$ such that $x <^{\mc{Z}} z <^{\mc{Z}} y$ and $d^{\mc{Z}} (x,z) = d^{\mc{Z}} (x,z) = D$. Since $\varepsilon = 0$, $<$ and $<^{\mc{Z}}$ disagree on $\{ x, z\}$ and $\{ z, y\}$. Thus, $x>z>y$ and so $x>z$. This proves the subclaim, finishes the proof of the claim and completes the proof of the lemma.  \end{proof}
  
The proof we presented here makes use of Ramsey property but we should mention here that this is not the only way to proceed. See for example \cite{N2} where the same result is proved thanks to a probabilistic argument.  

Observe also that as for Ramsey property, the previous proof allows to prove ordering property for classes $\M ^< _S$\index{$\M ^< _S$} whenever $S$ is an initial segment of some $T \subset ]0,+ \infty[$ which is closed under sums:  

\begin{thm}

\label{thm:variation OP for MM}

Let $T \subset ]0,+ \infty[$ be closed under sums and $S$ be an initial segment of $T$. Then $\M ^< _S$\index{$\M ^< _S$} has the ordering property. 

\end{thm}

Thus, in particular, all the classes $\M ^< _{\Q}$, $\M ^< _{\Q \cap ]0 , r]}$ with $r > 0$  in $\Q$, $\M ^< _{\omega}$ and $\M ^< _{\omega \cap ]0 , m ]}$ with $m > 0$ in $\omega$ have the ordering property.

\subsection{Finite convexly ordered ultrametric spaces.}

The next case of ordering property shows that ordering property can be proved completely independently of Ramsey property.  

\begin{thm}

\label{thm:OP for UU}

$\UU$ has the ordering property. 

\end{thm}

We begin with a simple observation coming from the tree
representation of elements of $\UU$. 

\begin{lemma}

\label{lemma: 1 OP for UU}

$\UU$ is a reasonable Fra\"iss\'e order class.  
\end{lemma}

\begin{proof}
The proof is left to the reader. Let us simply mention that it suffices to show that given $\m{X} \subset \m{Y}$ in $\U$\index{$\U$} and $<^{\m{X}}$ a convex
linear ordering on $\m{X}$, there is a convex linear ordering $<^{\m{Y}}$ on $\m{Y}$
such that $\restrict{ <^{\m{Y}} }{\m{X} } = <^{\m{X}}$.
\end{proof}

Call an element $\m{Y}$ of $\U$ \textit{convexly order-invariant} when
$(\m{Y},<_1) \cong (\m{Y},<_2)$ whenever $<_1 , <_2$ are convex
linear orderings on $\m{Y}$.
The following result is a direct consequence of the previous lemma: 

\begin{lemma}
Let $(\om{X}) \in \UU$ and assume that $\m{X} \subset \m{Y}$ for some convexly 
order-invariant $\m{Y}$ in $\U$. Then given any convex linear ordering $<$ on $\m{Y}$,
$(\om{X})$ embeds into $(\m{Y}, <)$.
\end{lemma}

\begin{proof}
Let $<^{\m{Y}}$ be as in the previous lemma. Let also $<$ be a convex linear orderings on $\m{Y}$. Then $(\om{X})$ embeds into $(\om{Y}) \cong (\m{Y}, <)$.
\end{proof}

We now show that any element of $\U$ embeds into a convexly order-invariant
one.

\begin{lemma}
Let $\m{X} \in \U$. Then $\m{X}$ embeds into $\m{Y}$ for some convexly order-invariant $\m{Y} \in \U$. 
\end{lemma}

\begin{proof}
Let $a_0 > a_1 > \ldots > a_{n-1}$ enumerate the distances appearing
in $\m{X}$. The tree representation of $\m{X}$ has $n$ levels. Now,
observe that such a tree can be embedded into a tree of height $n$ where all the nodes of a same level have the same number of
immediate successors, and that the ultrametric space associated to that
tree is convexly order-invariant. 
\end{proof}

Theorem \ref{thm:OP for UU} follows then directly. We finish this subsection with the justification of the remark at the end of \ref{subsection:RP for UU} stating that the class $\mathcal{U} ^< _S$\index{$\mathcal{U} ^< _S$} of all finite ordered ultrametric spaces with distances in $S$ does not have the Ramsey property. We start with:

\begin{thm}

\label{thm:non OP for U^<}

$\mathcal{U} ^< _S$ does not have the ordering property. 

\end{thm}

\begin{proof}
Let $(\om{X})$ be in $\mathcal{U} ^< _S$ and such that the ordering $<^{\m{X}}$ is not convex on $\m{X}$. Let $\m{Y}$ be in $\U$. Then there is a linear ordering $<$ on $\m{Y}$ such that $(\om{X})$ does not embed into $(\m{Y} , <)$. Namely, any convex linear ordering $<$ on $\m{Y}$ works. \end{proof}

We now show how this result can be used to prove:

\begin{thm}

\label{thm:non RP for U^<}

$\mathcal{U} ^< _S$ does not have the Ramsey property. 

\end{thm}

\begin{proof}
Assume for a contradiction that $\mathcal{U} ^< _S$ does have the Ramsey property. Then by a proof similar to the proof of Theorem \ref{thm:Ordering Property for metric spaces.}, $\mathcal{U} ^< _S$ would also have the ordering property, which is not the case. 
\end{proof}

\subsection{Finite metrically ordered metric spaces.}

Finally, we show how the methods used in the two previous subsections can be combined to prove that the ordering property holds for other classes of finite ordered metric spaces. 

\begin{thm}

\label{thm:OP for nM_S}

Let $S$ be a finite subset of $]0, + \infty [$ of size $|S| \leqslant 3$ and satisfying the $4$-values condition. Then $\nM _S$\index{$\nM _S$} has the ordering property. 
\end{thm}

\begin{proof}
As usual, the case $|S|=1$ is obvious. For $S = \{ 1, 2\}, \{2, 3, 4 \}$ or $\{ 1, 2, 3\}$, every linear ordering is metric so $\nM _S$ is really $\M ^< _S$\index{$\M ^< _S$} and as for Theorem \ref{thm:Ordering Property for metric spaces.}, ordering property is a consequence of Ramsey property. For $S = \{ 1, 3\}$ or $\{ 1, 3, 7\}$, the metric linear orderings are the convex ones, so ordering property is given by Theorem \ref{thm:OP for UU}. So the only remaining cases are the cases where $S$ is $\{ 1, 2, 5\}, \{1, 3, 6 \}$ and $\{ 1, 3, 4\}$. 

For $\{ 1, 2, 5\}$, ordering property comes from ordering property for finite graphs. To prove that fact, recall that for $\m{X} \in \M _{\{ 1, 2, 5\}}$, balls of radius $\leqslant 2$ are disjoint and can be seen as finite graphs with distance $5$ between them. Observe now that given $(\om{X}) \in \nM _{\{ 1, 2, 5\}}$, we can embed $(\om{X})$ into $(\om{Y}) \in \nM _{\{ 1, 2, 5\}}$ where all the balls of radius $2$ are isomorphic (as ordered graphs) to a same finite ordered graph $(\om{H})$. So $\m{Y} \cong \dot{\bigcup} _{i<k} \m{Y}_i$ for some $k \in \omega$, with $\m{Y}_0 < ^{\m{Y}}\ldots <^{\m{Y}} \m{Y}_{k-1}$ and $(\m{Y}_i , \restrict{<^{\m{Y}}}{Y_i}) \cong (\om{H})$ for every $i<k$. Let $\m{K}$ be a finite graph such that given any linear ordering $<$ on $\m{K}$,
$(\om{H})$ embeds into $(\m{K}, <)$. Then the metric space $\m{Z}$ defined by $\m{Z} \cong \dot{\bigcup} _{i<k} \m{Z}_i$ with $\m{Z}_i \cong \m{K}$ for every $i<k$ is such that for every metric linear ordering $<$ on $\m{Z}$, $(\om{Y})$ and hence $(\om{X})$ embeds into $(\m{Z} , <)$. 

For $\{ 1, 3, 6\}$, ordering property also comes from ordering property about finite graphs. Recall that in that case, balls of radius $1$ can be seen as complete graphs, and that between any two such balls, the distance between any two points is either always $3$ or always $6$. Let $(\om{X})$ be in $\nM _{\{ 1, 3, 6\}}$. Embed $(\om{X})$ into $(\om{Y})$ $\in \nM _{\{ 1, 3, 6\}}$ where all balls of radius $1$ have the same size $m$. Define now a graph $\m{G}_{\m{Y}}$ on the set $G_{\m{Y}}$ of balls of radius $1$ of $\m{Y}$ by connecting two balls iff the distance between any two of their points is equal to $3$. Observe that the ordering $<^{\m{Y}}$ beeing natural, it induces a linear ordering $G_{\m{Y}}$. Observe also that given a linear ordering on $G_{\m{Y}}$, there is a unique metric linear ordering on $\m{Y}$ extending it. Now, let $\m{K}$ be a finite graph such that given any linear ordering on $K$, $(\m{G}_{\m{Y}} , <^{\m{G}_{\m{Y}}})$ embeds into $(\m{K},<)$. Let $\m{Z}$ be the metric space whose space of balls is isomorphic to the graph $\m{K}$ and where every ball of radius $1$ has size $m$. Then given any metric linear ordering $<$ on $\m{Z}$, $(\om{X})$ embeds into $(\m{Z} , <)$.   

For $\{ 1, 3, 4\}$, the proof is a bit more involved. Fix $(\om{X}) \in \nM _{\{ 1, 3, 4\}}$. Recall that the relation $\approx$ defined by $x \approx y \leftrightarrow d^{\m{X}}(x,y) = 1$ is an equivalence relation. However, unlike the previous cases, the distance between the elements of two disjoint balls of radius $1$ can be arbitrarily $3$ or $4$. For 
$(\om{Y}) \in \nM _{\{ 1, 3, 4\}}$, say that a linear ordering $<$ on $Y$ is a \emph{local perturbation of} $< ^{\m{Y}}$\index{ordering!local perturbation} when 

\begin{center}
$\forall x, y \in Y \ \ d^{\m{Y}} (x,y) \geqslant 3 \rightarrow (x < y \leftrightarrow x < ^{\m{Y}} y ) $
\end{center}

\begin{lemma}
There is $(\om{Y}) \in \nM _{\{ 1, 3, 4\}}$ such that for any local perturbation $<$ of $< ^{\m{Y}}$, $(\om{X})$ embeds into $(\m{Y} , <)$. 
\end{lemma}

\begin{proof}

First, define a new linear ordering $<^{\m{X}} _*$ on $X$ by setting

\begin{displaymath}
\forall x, y \in X  \left \{ \begin{array}{l}
 d^{\m{X}} (x,y) = 1 \rightarrow (x <^{\m{X}} _* y \leftrightarrow y < ^{\m{X}} x ) \\
 d^{\m{X}} (x,y) \geqslant 3 \rightarrow (x <^{\m{X}} _* y \leftrightarrow x < ^{\m{X}} y ) 
 \end{array} \right.
\end{displaymath}  

Now, let $(\om{T})$ be the ordered metric space with two points and distance $1$ between them. Let also $(\m{X}_1 , <^{\m{X}_1 })$ be in $\nM _{\{ 1, 3, 4\}}$ and such that $(\om{X})$ and $(\m{X} , <^{\m{X}} _*)$ embed into $(\m{X}_1 , <^{\m{X}_1 })$. By Ramsey property, find $(\om{Y})$ such that 

\begin{center}
$(\om{Y}) \arrows{(\m{X}_1 , <^{\m{X}_1 })}{(\om{T})}{2}$.
\end{center}

We claim that $(\om{Y})$ is as required: Let $<$ be a local perturbation of $<^{\m{Y}}$. Then, define $\chi : \funct{\binom{\om{Y}}{\om{T}}}{2}$ by 

\begin{center}
$\chi (\oc{T}) = 1$ iff $<$ and $<^{\m{Y}}$ agree on $(\oc{T})$. 
\end{center}

By construction, there is a copy  $(\mc{X}_1 , <^{\mc{X}_1 })$ of  $(\m{X}_1 , <^{\m{X}_1 })$ such that $\binom{\mc{X}_1 , <^{\mc{X}_1 }}{\om{T}}$ is $\chi$-monochromatic with color $\varepsilon$. If $\varepsilon = 0$, consider $\mc{X} \subset \mc{X}_1$ such that 

\begin{center}
$(\mc{X} , \restrict{<^{\mc{X} _1}}{\mc{X}}) \cong (\m{X},<^{\m{X}} _*)$. 
\end{center}

Then 

\begin{center}
$(\mc{X} , \restrict{<}{\mc{X}}) \cong (\om{X})$.   
\end{center}

On the other hand, if $\varepsilon = 1$, consider $\mc{X} \subset \mc{X}_1$ such that 

\begin{center}
$(\mc{X} , \restrict{<^{\mc{X} _1}}{\mc{X}}) \cong (\m{X},<^{\m{X}})$. 
\end{center}

Then \[ (\mc{X} , \restrict{<}{\mc{X}}) \cong (\om{X}). \qedhere \]
\end{proof}

\begin{lemma}
There is $\m{Z} \in \M _{\{ 1, 3, 4\}}$ such that for any metric linear ordering $\prec$ on $\m{Z}$, there is a local perturbation $<$ of $< ^{\m{Y}}$ such that $(\m{Y} , <)$ embeds into $(\m{Z} , \prec)$. 
\end{lemma}

\begin{proof}

Define a new linear ordering $<^{\m{Y}} _{**}$ on $Y$ by

\begin{displaymath}
\forall x, y \in Y  \left \{ \begin{array}{l}
 d^{\m{X}} (x,y) = 1 \rightarrow (x <^{\m{Y}} _{**} y \leftrightarrow x < ^{\m{X}} y ) \\
 d^{\m{X}} (x,y) \geqslant 3 \rightarrow (x <^{\m{Y}} _{**} y \leftrightarrow y < ^{\m{X}} x ) 
 \end{array} \right.
\end{displaymath}  

Now, let $(\om{U})$ be the ordered metric space with two points and distance $3$ between them. Let also $(\m{Y}_1 , <^{\m{Y}_1 })$ be in $\nM _{\{ 1, 3, 4\}}$ such that $(\om{Y})$, $(\m{Y} , <^{\m{Y}} _{**})$ embed into $(\m{Y}_1 , <^{\m{Y}_1 })$ and such that any two balls of radius $1$ contain two points with distance $3$ between them. Still by Ramsey property, find $(\om{Z})$ such that 

\begin{center}
$(\om{Z}) \arrows{(\m{Y}_1 , <^{\m{Y}_1 })}{(\om{U})}{2}.$
\end{center}

Then $\m{Z}$ is as required: Let $\prec$ be a metric linear ordering on $\m{Z}$. Define a coloring $\Lambda : \funct{\binom{\om{Z}}{\om{U}}}{2}$ by 

\begin{center}
$\Lambda (\oc{U}) = 1$ iff $\prec$ and $<^{\m{Z}}$ agree on $(\oc{U})$. 
\end{center}

By construction, there is a copy  $(\mc{Y}_1 , <^{\mc{Y}_1 })$ of  $(\m{Y}_1 , <^{\m{Y}_1 })$ such that $\binom{\mc{Y}_1 , <^{\mc{Y}_1 }}{\om{U}}$ is $\Lambda$-monochromatic with color $\varepsilon$. If $\varepsilon = 0$, consider $\mc{Y} \subset \mc{Y}_1$ such that 

\begin{center}
$(\mc{Y} , \restrict{<^{\mc{Y} _1}}{\mc{Y}}) \cong (\m{Y},<^{\m{Y}} _{**})$. 
\end{center}

Otherwise, $\varepsilon = 1$ and choose $\mc{Y} \subset \mc{Y}_1$ such that 

\begin{center}
$(\mc{Y} , \restrict{<^{\mc{Y} _1}}{\mc{Y}}) \cong (\m{Y},<^{\m{Y}})$. 
\end{center}

In both cases, $(\mc{Y} , \restrict{<}{\mc{Y}}) \cong (\m{Y} , <)$ for some local perturbation $<$ of $<^{\m{Y}}$. \end{proof}

To finish the proof of the theorem, it is now enough to observe that given any metric linear ordering $\prec$ on $Z$, $(\om{X})$ embeds into $(\m{Z} , \prec)$.
\end{proof}

\section{Ramsey degrees.}

In this section, we show how the Ramsey property and the ordering property allow to show the existence and to compute the exact values of Ramsey degrees in various contexts. We start with the results about $\M$. For $\m{X} \in
\M$, let $\mathrm{LO}(\m{X})$\index{$\mathrm{LO}(\m{X})$} denote the set of all linear orderings on $\m{X}$. 
Thus, the number $|\mathrm{LO}(\m{X})|/|\mathrm{iso}(\m{X})|$ is essentially the number of all
nonisomorphic structures one can get by adding a linear ordering on
$\m{X}$. Indeed, if $<_1, <_2$ are linear orderings on
$\m{X}$, then $(\m{X},<_1)$ and $(\m{X},<_2)$ are
isomorphic as finite ordered metric spaces if and only if the
unique order preserving bijection from $(\m{X},<_1)$ to
$(\m{X},<_2)$ is an isometry. This defines an equivalence
relation on the set of all finite ordered metric spaces obtained
by adding a linear ordering on $\m{X}$. In what follows, an
\textit{order type for} $\m{X}$ is an equivalence class
corresponding to this relation. 

\begin{thm}

\label{thm:Rd in M}

Every $\m{X} \in \M $ has a
Ramsey degree $\mathrm{t}_{\M }(\m{X})$ in $\M $ and 

\begin{center}
$\mathrm{t}_{\M }(\m{X}) = |\mathrm{LO}(\m{X})|/|\mathrm{iso}(\m{X})|$. 
\end{center}

\end{thm}

\begin{proof}

Let $\tau (\m{X})$ denote the number $|\mathrm{LO}(\m{X})|/|\mathrm{iso}(\m{X})|$. 
We first prove that $\mathrm{t}_{\M }(\m{X}) \leqslant \tau (\m{X})$, ie that for
every $\m{Y} \in \M  $, $k \in \omega \smallsetminus \{ 0 \}$,
there is $\m{Z} \in \M $ such that 

\begin{center}
$\m{Z} \arrows{(\m{Y})}{\m{X}}{k, \tau (\m{X})}$.
\end{center}

Let $\{ <_{\alpha} : \alpha \in A \}$ be a set of linear orderings
on $\m{X}$ such that for every linear ordering $<$ on
\m{X}, there is a unique $\alpha \in A$ such that
$(\m{X},<)$ and $(\m{X},<_{\alpha})$ are isomorphic as
finite ordered metric spaces. Then $A$ has size
$\tau (\m{X})$ so without loss of generality, $A = \{1,
\ldots, \tau (\m{X}) \}$. Now, let $<^\m{Y}$ be any linear ordering on $Y$. By Ramsey property
for $\oM $ we can find $(\m{Z}_1,<^{\m{Z}_1} ) \in \oM $ such that

\begin{center}
$(\m{Z}_1,<^{\m{Z}_1})
\arrows{(\m{Y},<^\m{Y})}{(\m{X},<_1)}{k} $. 
\end{center}

Now,
construct inductively $(\m{Z}_2,<^{\m{Z}_2}), \ldots,
(\m{Z}_{\tau (\m{X})},<^{\m{Z}_{\tau (\m{X})}})
\in \oM _S$ such that for every $n \in \{1, \ldots,
\tau (\m{X})-1 \}$, 

\begin{center}
$(\m{Z}_{n+1},<^{\m{Z}_{n+1}})
\arrows{(\m{Z}_n,<^{\m{Z}_n})}{(\m{X},<_{n+1})}{k}$.
\end{center}

Finally, let $\m{Z}=\m{Z}_{\tau (\m{X})}$. Then one can check
that $\m{Z} \arrows{(\m{Y})}{\m{X}}{k,\tau (\m{X})}$.

To prove the reverse inequality $\mathrm{t}_{\M }(\m{X}) \geqslant \tau (\m{X})$, we need to show that there is $\m{Y} \in \M $ such
that for every $\m{Z} \in \M $, there is $\chi :
\funct{\binom{\m{Z}}{\m{X}}}{\tau (\m{X})}$ with the
property: 

\begin{center}
$\forall \mc{Y} \in \binom{\m{Z}}{\m{Y}} , \ \  \left| \chi '' \binom{\widetilde{\m{Y}}}{\m{X}} \right| = \tau (\m{X}).$
\end{center}

Fix $\m{X} \in \M $. By ordering property for $\oM $, find $\m{Y} \in \M $ such that for any linear ordering $<$ on $\m{Y}$, $(\m{Y},<)$ contains a copy of each order type of $\m{X}$. Now, let $\m{Z} \in \M $ and pick $<^\m{Z}$ any linear ordering on $\m{Z}$. Define a coloring $\chi : \funct{\binom{\m{Z}}{\m{X}}}{\tau (\m{X})}$ which colors any copy $\mc{X}$ of $\m{X}$ according to the order type of $(\mc{X},\restrict{<^\m{Z}}{\mc{X}})$. Now, if
possible, let $\mc{Y} \in \binom{\m{Z}}{\m{Y}} $. Then $(\mc{Y}, \restrict{<^\m{Z}}{\mc{Y}}) $
contains a copy of every order type of $\m{X}$, and \[ \left| \chi '' \binom{\widetilde{\m{Y}}}{\m{X}} \right| = \tau (\m{X}). \qedhere\] \end{proof}

The exact same proof can be used in different contexts. For example, one can replace $\M$ by $\M _S$\index{$\M _S$} where $S$ is an initial segment of a subset of $]0, +\infty[$ which is closed under sums:

\begin{thm}

\label{thm:variation Rd in M}

Let $T \subset ]0,+ \infty[$ be closed under sums and $S$ be an initial segment of $T$. Then every $\m{X} \in \M _S$ has a Ramsey degree $\mathrm{t}_{\M _S}(\m{X})$ in $\M _S$\index{$\M _S$} and \[ \mathrm{t}_{\M _S}(\m{X}) = |\mathrm{LO}(\m{X})|/|\mathrm{iso}(\m{X})|.\] 
\end{thm}

This fact has two consequences. On the one hand, the only objects for which
$\mathrm{t}_{\M _S}(\m{X}) = 1$ are the equilateral ones. On the other
hand, there are objects for which the Ramsey degree is $\mathrm{LO}(\m{X})$ (ie $|\m{X}|!$), those for which there
is no nontrivial isometry.

We now turn to ultrametric spaces: Given $S \subset ]0, + \infty [$, we showed that the class $\UU$\index{$\UU$} has the Ramsey property and the ordering property. Thus, if for $\m{X} \in
\U$, $\mathrm{cLO}(\m{X})$\index{$\mathrm{cLO}(\m{X})$} denotes the set of all convex linear orderings on $\m{X}$, we obtain:  

\begin{thm}

\label{thm:Rd in U}

Let $S \subset ]0, + \infty [$. Then every $\m{X} \in \U$ has a
Ramsey degree $\mathrm{t}_{\U}(\m{X})$ in $\U$\index{$\U$} and \[\mathrm{t}_{\U}(\m{X}) = |\mathrm{cLO}(\m{X})|/|\mathrm{iso}(\m{X})|.\]
\end{thm}

This fact makes the situation for ultrametric spaces a bit different from the metric case: First, 
the ultrametric spaces for which the true Ramsey property holds are those for which the corresponding
tree is uniformly branching on each level. Hence, in the class $\U$,
every element can be embedded into a Ramsey object, a fact which does
not hold in the class of all finite metric spaces. Second, one can
notice that any finite ultrametric space has a nontrivial isometry
(this fact is obvious via the tree representation). Thus, the Ramsey
degree of $\m{X}$ is always strictly less than
$|\mathrm{cLO}(\m{X})|$. In fact, a simple computation shows that the
highest value $\mathrm{t}_{\U}(\m{X})$ can get if the size of
$\m{X}$ is fixed is $2^{|\m{X}|-2}$ and is realized when the tree
associated to $\m{X}$ is a comb, ie when all the branching nodes are
placed on a same branch.

Finally, for $S$ finite subset of $]0, + \infty [$ of size $|S| \leqslant 3$ and satisfying the $4$-values condition, we saw that the class $\nM _S$\index{$\nM _S$} has the Ramsey and the ordering properties. It follows that if for $\m{X} \in \M _S$, $\mathrm{mLO}(\m{X})$\index{$\mathrm{mLO}(\m{X})$} denotes the set of all metric linear orderings on $\m{X}$, one gets:  

\begin{thm}

\label{thm:Rd in M_S}

Let $S$ be finite subset of $]0, + \infty [$ of size $|S| \leqslant 3$ and satisfying the $4$-values condition. Then every $\m{X} \in \M _S$ has a
Ramsey degree $\mathrm{t}_{\M _S}(\m{X})$ in $\M _S$\index{$\M _S$} and \[\mathrm{t}_{\M _S}(\m{X}) = |\mathrm{mLO}(\m{X})|/|\mathrm{iso}(\m{X})|.\]
\end{thm}

\section{Universal minimal flows and extreme amenability.}

After the study of Ramsey and ordering properties, we turn to applications in topological dynamics.   

\subsection{Pestov theorem.}

In this subsection, we present a proof of the following result:

\begin{thm}[Pestov \cite{Pe0}]

\label{thm:Pestov}
\index{Pestov!theorem on $\iso (\Ur)$}
\index{extreme amenability! of $\iso (\Ur)$}

Equipped with the pointwise convergence topology, the group of isometries $\iso (\Ur)$ of the Urysohn space is extremely amenable (has the fixed-point on comptacta property). 

\end{thm}

In the sequel, we present how this result can be deduced from the general theory exposed in the introduction of this chapter. The proof is taken from \cite{KPT}. 

First, the class $\M _{\Q}$ is a reasonable Fra\"iss\'e class. It follows that $\mathrm{Flim} (\M ^< _{\Q}) = (\Ur _{\Q} , <^{\Ur _{\Q}})$ for some linear ordering $<^{\Ur _{\Q}}$ on $\Ur _{\Q}$. Furthermore, we saw that $\M ^< _{\Q}$ has the Ramsey and the ordering properties. Consequently:

\begin{thm}[Kechris-Pestov-Todorcevic \cite{KPT}]

\label{thm:Aut(Ur_Q , <) ea}
\index{Kechris-Pestov-Todorcevic!theorem on $\mathrm{Aut}(\Ur _{\Q} , <^{\Ur _{\Q}})$}
\index{extreme amenability! of $\mathrm{Aut}(\Ur _{\Q} , <^{\Ur _{\Q}})$}

The group $\mathrm{Aut}(\Ur _{\Q} , <^{\Ur _{\Q}})$ is extremely amenable.

\end{thm}

\begin{thm}[Kechris-Pestov-Todorcevic \cite{KPT}]

\label{thm:M(iso(Ur_Q))}
\index{Kechris-Pestov-Todorcevic!theorem on $M(\iso (\Ur _{\Q}))$}
\index{flow!universal minimal flow!of $\iso (\Ur _{\Q})$}

The universal minimal flow of $\iso (\Ur _{\Q})$ is the set $\LO (\Ur _{\Q})$ of linear orderings on $\Ur _{\Q}$ together with the action $\funct{\iso (\Ur _{\Q}) \times \LO (\Ur _{\Q})}{\LO (\Ur _{\Q})}$, $(g,<)
\longmapsto <^g$ defined by \[ x <^g y \ \ \mathrm{iff} \ \ g^{-1}(x) < g^{-1}(y).\]
\end{thm} 

We now show how to deduce Theorem \ref{thm:Pestov} from those results.

\begin{lemma}
Let $G$, $H$ be topological groups and $\pi : \funct{G}{H}$ be a continuous morphism with dense range. Assume that $G$ is extremely amenable. Then so is $H$. 
\end{lemma}

\begin{proof}
Let $X$ be an $H$-flow. Denote by $\alpha : \funct{H \times X}{X}$ the action. Define now $\bar{\alpha} : \funct{G \times X}{X}$ by $\bar{\alpha}(g,x) = \alpha (\pi (g) , x)$. This turns $X$ into a $G$-flow so there is a fixed point $x_0 \in X$. But since $\pi$ has dense range, $x_0$ is also fixed for the $H$-flow. \end{proof}

Now, recall that $\Ur$\index{$\Ur$} is the completion of $\Ur _{\Q}$ so given any $g \in \iso (\Ur _{\Q})$, there is a unique $\bar{g}$ extending $g$ on $\Ur$. Since every $g \in \mathrm{Aut}(\Ur _{\Q} , <^{\Ur _{\Q}})$ is in particular an isometry of $\Ur _{\Q}$, the map $g \mapsto \bar{g}$ is 1-1 from $\mathrm{Aut}(\Ur _{\Q} , <^{\Ur _{\Q}})$ into $\iso (\Ur)$ and it is easy to check that it is continuous. Consequently, according to the previous lemma, it only remains to show that its range is dense in $\iso (\Ur)$.  

\begin{lemma}
Let $D \subset \iso (\Ur)$. Let $d$ denote the metric on $\Ur _{\Q}$. Assume that:
\begin{center}
$\forall \varepsilon > 0 \ \ \forall x_1\ldots x_n \in \Ur \ \ \forall h \in \iso (\Ur) \ \ \exists x' _1\ldots x'_n , y'_1\ldots y'_n \in \Ur \ \ \exists g \in D$ 

$\forall i \leqslant n \ \ d(x_i , x'_i) < \varepsilon , \ d(h(x_i) , y'_i) < \varepsilon , \ g(x'_i) = y'_i$.
\end{center}

Then $D$ is dense in $\iso (\Ur)$. 

\end{lemma}

\begin{proof}
Fix $\varepsilon > 0$, $h \in \iso (\Ur)$ and $x_1\ldots x_n \in U$. Thanks to the hypothesis, find $x'_1\ldots x'_n , y'_1\ldots y'_n \in U$ and $g \in D$ for $\varepsilon /2$. Then for $i \leqslant n$:

\begin{displaymath}
\begin{array}{lll}
d(g(x_i),h(x_i))& \leqslant & d(g(x_i), g(x' _i)) + d(g(x' _i) , h(x_i))\\
& = & d(x_i, x' _i) + d(y' _i , h(x_i))\\
& < & \varepsilon. \ \  \qedhere
\end{array}
\end{displaymath}
\end{proof}

So to check that $\{ \bar{g} : g \in \mathrm{Aut}(\Ur _{\Q} , <^{\Ur _{\Q}})\}$ is dense in $\iso (\Ur)$, it is enough to show:

\begin{lemma}
Given $x_1\ldots x_n , y_1\ldots y_n \in U$ such that $x_i \mapsto y_i$ is an isometry and given $\varepsilon > 0$, there are $x'_1\ldots x'_n , y'_1\ldots y'_n \in U_{\Q}$ so that $x'_i \mapsto y'_i$ is an order-preserving isometry with respect to $<$ and 

\begin{center}
$\forall i \leqslant n \ \ d(x'_i , x_i) < \varepsilon, \ \ d(y'_i , y_i) < \varepsilon $. 
\end{center}

\end{lemma}

\begin{proof}
We proceed by induction on $n$. For $n=1$, simply choose $x'_i , y'_i \in U _{\Q}$ such that $d(x'_i , x_i) < \varepsilon$ and $d(y'_i , y_i) < \varepsilon $. For the induction step, assume that we are at stage $n$ and wish to step up to $n+1$. Suppose that $x_1,\ldots , x_{n+1} , y_1,\ldots y_{n+1} \in U$ are given so that $x_i \mapsto y_i$ is an isometry. By induction hypothesis, find $x'_1\ldots x'_n$ and $y'_1\ldots y'_n \in U_{\Q}$ so that $x'_i \mapsto y'_i$ is an order-preserving isometry and 

\begin{center}
$\forall i \leqslant n \ \ d(x'_i , x_i) < \varepsilon /2, \ \ d(y'_i , y_i) < \varepsilon /2$. 
\end{center} 

Fix $x ^0 _{n+1} , y ^0 _{n+1} \in U_{\Q}$ such that 

\begin{center}
$d(x^0 _{n+1} , x_{n+1}) < \varepsilon /2, \ \ d(y ^0 _{n+1} , y _{n+1}) < \varepsilon /2$.
\end{center}

For $i \leqslant n$, set $d_i := d(x^0 _{n+1} , x'_i )$ and $d'_i := d(y^0 _{n+1} , y'_i)$. Without loss of generality, we may assume that $\varepsilon < d_i , d'_i$. Therefore: 

\begin{center}
$\left| d_i - d(x_{n+1} , x_i) \right| \leqslant \left| d(x^0 _{n+1} , x_{n+1}) + d(x_i , x' _i)   \right| < \varepsilon$. 
\end{center}

Similarly, 

\begin{center}
$\left| d'_i - d(y_{n+1} , y_i) \right| < \varepsilon$. 
\end{center} 

So

\begin{center}
$\left| d_i - d'_i \right| = \left| d_i - d(x_{n+1} , x_i) + d(x_{n+1} , x_i) - d(y_{n+1} , y_i) + d(y_{n+1} , y_i) - d'_i  \right| < \varepsilon$. 
\end{center}

Now, set $e_i := (d_i + d'_i)/2$ and consider the ordered metric space \begin{center} $(\{x'_1 ,\ldots , x'_n , x^0 _{n+1} , u \} , d' , \prec)$ \end{center} where \begin{center} $d'(x'_i , x'_j) = d(x'_i , x'_j)$, $d'(x'_i , x^0 _{n+1}) = d(x'_i , x^0 _{n+1})$, $d'(u, x'_i) = e_i$ \end{center} and $d'(u , x^0 _{n+1})$ is any irrational number satisfying the inequalities:

\begin{center}
$\forall i \leqslant n \ \ \left| d_i - e_i \right| \leqslant d'(u, x^0 _{n+1}) < 2 \varepsilon < d_i + e_i$. 
\end{center}

Observe that the existence of such a number is guaranteed by the inequalities 

\begin{center}
$d_i + e_i = \frac{3d_i + d'_i}{2} > \varepsilon$ 
\end{center}

and 

\begin{center}
$\left| d_i - e_i \right| = \frac{\left|d_i - d'_i\right|}{2} < \varepsilon$.
\end{center}

As for $\prec $, we let it agree with the ordering $<$ of $U_{\Q}$ for $x'_1 ,\ldots , x'_n , x^0 _{n+1}$ and set $x'_i \prec u $ as well as $x^0 _{n+1} \prec u$. Assuming that $d'$ defines a metric, we finish the proof as follows: By the properties of $(\Ur _{\Q} , <^{\Ur _{\Q}})$, we can find a point $x'_{n+1} \in U_{\Q}$ with $x'_i < x'_{n+1}$ for every $i \leqslant n$, $x^0 _{n+1} < x'_{n+1}$ and $d(x'_{n+1} , x'_i) = e_i$, $d(x'_{n+1} , x^0 _{n+1}) = d'(u, x^0 _{n+1}) < 2 \varepsilon$. Similarly, we can find $y'_{n+1} \in U_{\Q}$ with $y'_i < y'_{n+1}$ for every $i \leqslant n$, $y^0 _{n+1} < y'_{n+1}$ and $d(y'_{n+1} , y'_i) = e_i$, $d(y'_{n+1} , y^0 _{n+1}) = d'(u, x^0 _{n+1}) < 2 \varepsilon$. Then, $x'_i \mapsto y'_i$ defines an order preserving map and 

\begin{center}
$d(x'_{n+1}, x_{n+1}) \leqslant d(x' _{n+1} , x^0 _{n+1}) + d(x^0 _{n+1} , x_{n+1}) < 3 \varepsilon$,
\end{center} which completes the proof. It remains to check that $d'$ indeed defines a metric: 

(i) Since $d'(x^0 _{n+1} , x'_i)=d_i$, $d'(u, x'_i) = e_i$, we need to check that \begin{center}
$\left| d_i - e_i \right| \leqslant d'(u, x^0 _{n+1}) \leqslant d_i + e_i$, \end{center} which is given by the definition of $d'(u, x^0 _{n+1})$. 

(ii) Let $\alpha _{ij} = d(x'i , x'_j)$. We need to verify that 

\begin{center}
$\left| e_i - e_j \right| \leqslant \alpha _{ij} \leqslant e_i + e_j$. 
\end{center}

On the one hand: 

\begin{center}
$\left| d_i - d_j \right| \leqslant \alpha _{ij} \leqslant d_i + d_j$. 
\end{center}

On the other hand, $\alpha _{ij} = d(y'i , y'_j)$ so we also have:

\begin{center}
$\left| d'_i - d'_j \right| \leqslant \alpha _{ij} \leqslant d'_i + d'_j$. 
\end{center}

Adding and dividing by $2$, we obtain the required inequality. \end{proof}

As in previous sections, simple adaptations of the proof allow to deduce similar results for other spaces. Fot example, instead of working with $\M ^< _{\Q}$ and the structure $(\Ur _{\Q},<^{\Ur _{\Q}})$, one can work with the reasonable Fra\"iss\'e class $\M ^< _{\Q \cap ]0,1]}$ and its Fra\"iss\'e limit $(\s _{\Q} , <^{\s _{\Q}})$. Here are the results we obtain in this case: 

\begin{thm}[Kechris-Pestov-Todorcevic \cite{KPT}]

\label{thm:Aut(s_Q , <) ea}
\index{Kechris-Pestov-Todorcevic!theorem on $\mathrm{Aut}(\s _{\Q} , <^{\s _{\Q}})$}

The group $\mathrm{Aut}(\s _{\Q} , <^{\s _{\Q}})$ is extremely amenable.

\end{thm}

\begin{thm}[Kechris-Pestov-Todorcevic \cite{KPT}]

\label{thm:M(iso(s_Q))}
\index{Kechris-Pestov-Todorcevic!theorem on $M(\iso (\s _{\Q}))$}
\index{flow!universal minimal flow!of $\iso (\s _{\Q})$}

The universal minimal flow of $\iso (\s _{\Q})$ is the set $\LO (\s _{\Q})$ of linear orderings on $\s _{\Q}$\index{$\s _{\Q}$} together with the action $\funct{\iso (\s _{\Q}) \times \LO (\s _{\Q})}{\LO (\s _{\Q})}$, $(g,<)
\longmapsto <^g$ defined
by $x <^g y$ iff $g^{-1}(x) < g^{-1}(y)$.

\end{thm} 

\begin{thm}[Pestov \cite{Pe0}]

\label{thm:iso(S) ea}
\index{Pestov!theorem on $\iso (\s)$}
\index{extreme amenability! of $\iso (\s)$}

The group $\iso (\s)$ is extremely amenable. 

\end{thm}

Other interesting examples appear when the distance set $\Q$ is replaced by $\omega$ or $\{ 1,\ldots , m\}$ for some strictly positive $m$ in $\omega$. One then deals with the reasonable Fra\"iss\'e classes $\M ^< _{\omega}$ and $\M ^< _m$ and their Fra\"iss\'e limits $(\Ur _{\omega} , <^{\Ur _{\omega}})$ and $(\Ur _m , <^{\Ur _m})$ respectively:

\begin{thm}[Kechris-Pestov-Todorcevic \cite{KPT}]

\label{thm:Aut(U_N , <) ea}
\index{extreme amenability! of $\mathrm{Aut}(\Ur _{\omega}, <^{\Ur _{\omega}})$}
\index{Kechris-Pestov-Todorcevic!theorem on $\mathrm{Aut}(\Ur _{\omega} , <^{\Ur _{\omega}})$}

The group $\mathrm{Aut}(\Ur _{\omega} , <^{\Ur _{\omega}})$ is extremely amenable.

\end{thm}

\begin{thm}[Kechris-Pestov-Todorcevic \cite{KPT}]

\label{thm:M(iso(U_N))}
\index{Kechris-Pestov-Todorcevic!theorem on $M(\iso (\Ur _{\omega}))$}
\index{flow!universal minimal flow!of $\iso (\Ur _{\omega})$}

The universal minimal flow of $\iso (\Ur _{\omega})$ is the set $\LO (\Ur _{\omega})$ of linear orderings on $\Ur _{\omega}$ together with the action $\funct{\iso (\Ur _{\omega}) \times \LO (\Ur _{\omega})}{\LO (\Ur _{\omega})}$, $(g,<)
\longmapsto <^g$ defined by \[x <^g y \ \ \mathrm{iff} \ \ g^{-1}(x) < g^{-1}(y).\]
\end{thm}

\begin{thm}[Kechris-Pestov-Todorcevic \cite{KPT}]

\label{thm:Aut(U_m , <) ea}
\index{Kechris-Pestov-Todorcevic!theorem on $\mathrm{Aut}(\Ur _{m} , <^{\Ur _{m}})$}
\index{extreme amenability! of $\iso (\Ur _m)$}

The group $\mathrm{Aut}(\Ur _{m} , <^{\Ur _{m}})$ is extremely amenable.

\end{thm}

\begin{thm}[Kechris-Pestov-Todorcevic \cite{KPT}]

\label{thm:M(iso(U_m))}
\index{Kechris-Pestov-Todorcevic!theorem on $M(\iso (\Ur _{m}))$}
\index{flow!universal minimal flow!of $\iso (\Ur _{m})$}

The universal minimal flow of $\iso (\Ur _{m})$ is the set $\LO (\Ur _{m})$ of linear orderings on $\Ur _{m}$ together with the action $\funct{\iso (\Ur _{m}) \times \LO (\Ur _{m})}{\LO (\Ur _{m})}$, $(g,<)
\longmapsto <^g$ defined by \[ x <^g y \ \ \mathrm{iff} g^{-1}(x) < g^{-1}(y).\]

\end{thm}

\subsection{Ultrametric Urysohn spaces.}

After Pestov theorem and its variations, the results we present now deal with ultrametric spaces. In chapter 1, we mentioned that the Urysohn space $\uUr$\index{$\uUr$} of the class $\U$\index{$\U$} when $S$ is a countable distance set can be described explicitly. The class $\UU$\index{$\UU$} being a reasonable Fra\"iss\'e class, its Fra\"iss\'e limit is therefore equal to $(\uUr , <^{\uUr})$ for some linear ordering $<^{\uUr}$ on $\uUr$. It turns out that as $\uUr$, the ordering $<^{\uUr}$ is also easy to describe: It is simply the lexicographical ordering $<^{\uUr} _{lex}$ coming from the natural tree associated to $\uUr$.

\begin{prop}
Let $S \subset ]0,+\infty[$ be countable. Then $\mathrm{Flim}(\UU) = (\uUr , <^{\uUr}_{lex})$. 

\end{prop}

\begin{proof}
The only thing we have to check is that $<^{\uUr}_{lex}$ is the relevant linear ordering on $\uUr$, ie that $(\uUr , <^{\uUr}_{lex})$ is ultrahomogeneous. In what follows, we relax the notation and simply write $d$ (resp. $<$)
instead of $d^{\uUr}$ (resp. $<^{\uUr} _{lex}$). We proceed by induction on the size $n$ of the finite substructures.

For $n=1$, if $x$ and $y$ are in $\uUr$, just define $g :
\funct{\uUr}{\uUr}$ by 

\begin{center}
$g(z) = z + y - x$.
\end{center}

For the induction step, assume that the homogeneity of $(\uUr ,
<)$ is proved for finite substructures of size $n$ and
consider two isomorphic substructures of $(\uUr , <)$ of
size $n+1$, namely $x_1 < \ldots < x_{n+1}$ and $y_1 < \ldots <
y_{n+1}$. By induction hypothesis, find $h \in \mathrm{Aut}(\uUr,<)$
such that for every $1 \leqslant i \leqslant n$, $h(x_i) = y_i$. We
now have to take care of $x_{n+1}$ and $y_{n+1}$. Observe first that
thanks to the convexity of $<$, we have 

\begin{center}
$d(x_n,x_{n+1}) = \min \{
d(x_i,x_{n+1}) : 1 \leqslant i \leqslant n \}$.
\end{center}

Similarly, 

\begin{center}
$d(y_n,y_{n+1}) = \min \{
d(y_i,y_{n+1}) : 1 \leqslant i \leqslant n \}$.
\end{center}

Set 

\begin{center}
$s = d(x_n,x_{n+1}) = d(y_n,y_{n+1})$. 
\end{center}

Note that $y_{n+1}$ and $h(x_{n+1})$ agree on $S \cap ]s,\infty[$. Indeed, 

\begin{eqnarray*}
d(y_{n+1},h(x_{n+1})) & \leqslant & \max
(d(y_{n+1},y_n),d(y_n,h(x_{n+1}))) \\
& \leqslant & \max (d(y_{n+1},y_n),d(h(x_n),h(x_{n+1})))\\
& \leqslant & \max (s,s) = s
\end{eqnarray*}

Note also that since $y_n < y_{n+1}$ (resp. $h(x_n) < h(x_{n+1})$), we
have 

\begin{center}
$y_n (s) < y_{n+1} (s)$. 
\end{center}

Similarly, 

\begin{center}
$y_n (s) = h(x_n)(s) < h(x_{n+1})(s)$. 
\end{center}

So $(\R \smallsetminus \Q) \cap ]y_n (s) , \min (y_{n+1}(s), h(x_{n+1})(s)) [$ is non-empty and
has an element $\alpha$. Next, the set $] \alpha , \infty [ \cap \Q$ is order-isomorphic to $\Q$ so we can find a
strictly increasing bijective $\phi : \funct{] \alpha ,
\infty [ \cap \Q}{] \alpha , \infty [ \cap \Q}$ such that \begin{center} $\phi
(h(x_{n+1})(s)) = y_{n+1}(s)$. \end{center} 

Now, define $j : \funct{\uUr}{\uUr}$ by $j(x) = x$ if $d(x,y_{n+1}) > s$. Otherwise (when $d(x,y_{n+1}) \leqslant s$), set 

\begin{displaymath}
j(x)(t) = \left \{ \begin{array}{ll}
 x(t) & \textrm{if $t>s$,} \\
 x(t) & \textrm{if $t=s$ and $x(t)<\alpha$,} \\
 \phi ( x(t)) & \textrm{if $t=s$ and $\alpha < x(t)$,}\\
 x(t) + y_{n+1}(t) - h(x_{n+1})(t) & \textrm{if $t<s$.} 
 \end{array} \right.
\end{displaymath}

One can check that $j \in \mathrm{Aut}(\uUr ,<)$ and that for every $1
\leqslant i \leqslant n$, $j(y_i) = y_i$. Now, let $g = j \circ h$. We
claim that for every $1 \leqslant i \leqslant n+1$, $g(x_i) =
y_i$. Indeed, if $1 \leqslant i \leqslant n$ then $g(x_i) = j(h(x_i))
= j(y_i) = y_i$. Moreover, 

\begin{eqnarray*}
g(x_{n+1})(t) & = & j(h(x_{n+1}))(t) \\
& = & \left \{ \begin{array}{ll}
 h(x_{n+1})(t) & \textrm{if $t>s$,} \\
 \phi (h(x_{n+1})(t)) = y_{n+1}(t) & \textrm{if $t=s$,} \\
 h(x_{n+1})(t) + y_{n+1}(t) - h(x_{n+1})(t) = y_{n+1}(t)& \textrm{if $t<s$.} 
 \end{array} \right.
\end{eqnarray*}

ie $g(x_{n+1}) = y_{n+1}$. \end{proof}

Therefore, Ramsey property together with ordering property for $\UU$ lead to the following result in topological dynamics:     

\begin{thm}

\label{thm:Aut(uUr, <) ea}
\index{extreme amenability! of $\mathrm{Aut}(\uUr , <^{\uUr}_{lex})$}

The group $\mathrm{Aut}(\uUr , <^{\uUr}_{lex})$ is extremely amenable.
\end{thm}

\begin{thm}

\label{thm:UMF for iso(uUr)}
\index{flow!universal minimal flow!of $\iso (\uUr)$}

The universal minimal flow of $\mathrm{iso}(\uUr)$ is the set $\cLO (\uUr)$ of
convex linear orderings on $\uUr$\index{$\uUr$} together with the action $\funct{\iso (\uUr) \times \cLO (\uUr)}{\cLO (\uUr)}$, $(g,<)
\longmapsto <^g$ defined
by $x <^g y$ iff $g^{-1}(x) < g^{-1}(y)$.
\end{thm}

\textbf{Remark.} In \cite{KPT}, Theorem 6.6, it is mentioned that for $S = 2$, Theorem \ref{thm:Aut(uUr, <) ea} can actually be proved directly using preservation of extreme amenability under direct and semi-direct products of topological groups. More recently, we were informed by Christian Rosendal that it is also the case for any countable $S$. Had this result been known to us before Theorem \ref{thm:RP for UU}, the equivalence provided by Theorem \ref{thm:KPT ea} would have allowed to deduce Theorem \ref{thm:RP for UU} from it. 

\

We now use these results to compute the universal minimal flow of the metric completion $\cUr$\index{$\cUr$} of $\uUr$. We follow the scheme adopted in the previous section. Let $<^{\cUr} _{lex}$ be
the natural lexicographical ordering on $\cUr$.

\begin{lemma}

\label{lem:Aut(Ur,<) embeds densely}

There is a continuous group morphism for which
$\mathrm{Aut}(\uUr , <^{\uUr} _{lex})$ embeds densely into
$\mathrm{Aut}(\cUr , <^{\cUr} _{lex})$.

\end{lemma}

\begin{proof}
Every $g \in \mathrm{iso}(\uUr)$ has unique extension $\hat{g} \in
\mathrm{iso}(\cUr)$. Moreover, observe that
$<^{\cUr}_{lex}$ can be reconstituted from
$<^{\uUr} _{lex}$. More precisely, if $\hat{x}, \hat{y} \in \cUr$,
and $x, y \in \uUr$ such that
$d^{\cUr}(x,\hat{x}), d^{\cUr}(y,\hat{y}) < d^{\cUr}(\hat{x},\hat{y})$, then
\begin{center} $\hat{x}<^{\cUr} _{lex} \hat{y}$ iff $x<^{\uUr} _{lex}y$.
\end{center} 

Note that
this is still true when $<^{\cUr} _{lex}$ and $<^{\uUr} _{lex}$ are
replaced by $\prec \in \cLO(\cUr)$ and $\restrict{\prec}{\uUr} \in
\cLO(\uUr)$ respectively. Later, we
will refer to that fact as the \emph{coherence property}\index{coherence property}.
Its first consequence is that the map $g \mapsto \hat{g}$ can actually be seen
as a map from $ \mathrm{Aut}(\uUr , <^{\uUr} _{lex})$ to
$\mathrm{Aut}(\cUr , <^{\cUr} _{lex})$. It is easy to check that it
is a continuous embedding. We now prove that it has dense range. Take $h
\in \mathrm{Aut}(\cUr , <^{\cUr} _{lex})$, $\hat{x}_1 <^{\cUr} _{lex} \ldots
<^{\cUr} _{lex} \hat{x}_n$ in $\cUr$, $\varepsilon > 0$, and consider the
corresponding basic open neighborhood $W$ around $h$. Take $\eta
> 0$ such that $\eta < \varepsilon$ and \[ \forall 1 \leqslant i \neq j \leqslant n, \ \ \eta <
d^{\cUr}(\hat{x}_i , \hat{x}_j).\]

Now, pick $x_1 , \ldots , x_n , y_1 ,
\ldots, y_n \in \uUr$ such that \[ \forall 1 \leqslant i \leqslant n, \ \ d^{\cUr}(\hat{x}_i ,
x_i) < \eta \ \ \mathrm{and} \ \ d^{\cUr}(h(\hat{x}_i) , y_i) < \eta.\]

Then one can
check that the map $x_i \mapsto y_i$ is an isometry
from $\{x_i : 1 \leqslant i \leqslant n \}$ to $\{ y_i : 1 \leqslant i
\leqslant n \}$ (because $\cUr$ is ultrametric) which is also order-preserving
(thanks to the coherence property). By ultrahomogeneity of
$(\uUr , <^{\uUr} _{lex})$, we can extend that
map to $g_0 \in \mathrm{Aut}(\uUr , <^{\uUr} _{lex})$. Finally, consider
the basic open neighborhood $V$ around $g_0$ given by $x_1 , \ldots , x_n$
and $\eta$. Then $\{ \hat{g} : g \in V \} \subset W$. Indeed, let $g
\in V$. Then $d^{\cUr}(\hat{g}(\hat{x}_i),h(\hat{x}_i))$ is less or
equal to \begin{center}$\max \{ d^{\cUr}(\hat{g}(\hat{x}_i), \hat{g}(x_i)),
d^{\cUr}(\hat{g}(x_i), \hat{g_0} (x_i)),
d^{\cUr}(\hat{g}_0 (x_i), h(\hat{x}_i)) \}$.\end{center}

Now, since $\hat{g}$ is an isometry,
$d^{\cUr}(\hat{g}(\hat{x}_i), \hat{g}(x_i)) = d^{\cUr}(\hat{x}_i,
x_i) < \eta < \varepsilon$. Also, since
$g \in V$, $d^{\cUr}(\hat{g}(x_i), \hat{g_0} (x_i)) <
\eta < \varepsilon$. Finally, by construction of
$g_0$, 

\begin{center}
$d^{\cUr}(\hat{g}_0 (x_i), h(\hat{x}_i)) = d^{\uUr}(y_i, h(\hat{x_i}))
< \eta < \varepsilon$. 
\end{center}

Thus $d^{\cUr}(\hat{g}(\hat{x}_i),h(\hat{x}_i))
< \varepsilon$ and $\hat{g} \in W$. \end{proof}

As a direct corollary, we obtain:

\begin{thm}
\label{cor:Aut(cUr,<) ea}
\index{extreme amenability! of $\mathrm{Aut}(\cUr ,<^{\cUr} _{lex})$}

The group $\mathrm{Aut}(\cUr ,<^{\cUr} _{lex})$ is extremely amenable.
\end{thm}

Let us now look at the topological dynamics of the isometry group
$\iso (\cUr)$. Note that $\iso (\cUr)$ is not extremely amenable as
its acts continuously on the space of all convex linear orderings
$\cLO (\cUr)$ on $\cUr$ with no fixed point. The following result
shows that in fact, this is its universal minimal compact action.

\begin{thm}

\label{cor:UMF for iso(cUr)}
\index{flow!universal minimal flow!of $\iso (\cUr)$}

The universal minimal flow of $\mathrm{iso}(\cUr)$ is the set $\cLO
(\cUr)$ together with the action
$\funct{\iso (\cUr) \times \cLO (\cUr)}{\cLO (\cUr)}$, $(g,<)
\longmapsto <^g$ defined by \[ x <^g y \ \ \mathrm{iff} \ \ g^{-1}(x) < g^{-1}(y).\]
\end{thm}

\begin{proof}

Equipped with the topology for which the basic open sets
are those of the form $\{ \prec \in \cLO (\cUr) : \restrict{\prec}{X}
= \restrict{<}{X} \}$ (resp. $\{ \prec \in \cLO (\uUr) : \restrict{\prec}{X}
= \restrict{<}{X} \}$) where $X$ is a finite subset of $\cUr$
(resp. $\uUr$), the space $\cLO (\cUr)$ (resp. $\cLO (\uUr)$) is
compact. To see that the action is continuous,
let $< \in \cLO (\cUr)$, $g \in \iso (\cUr)$ and $W$ a basic open
neighborhood around $< ^g$ given by a finite $X \subset \cUr$.
Now take $\varepsilon >
0$ strictly smaller than any distance in $X$
and consider \begin{center} $U = \{ h \in \iso (\cUr) : \forall x \in X
(d^{\cUr}(g^{-1}(x),h^{-1}(x)) < \varepsilon) \}$.\end{center} 

Let also \begin{center}$V = \{ \prec \in  \cLO (\cUr) :
\restrict{\prec}{\overleftarrow{g}X} =
\restrict{\prec}{\overleftarrow{h}X} \}$. \end{center}

We claim that
for every $(h, \prec) \in U \times V$, we have $\prec ^h \in W$. To
see that, observe first that if
$x,y \in X$, then $h^{-1}(x) \prec
h^{-1}(y)$ iff $g^{-1}(x) \prec g^{-1}(y)$ (this is a consequence of
the coherence property). So if $(h, \prec) \in U \times V$ and $x,y
\in X$ we have  
\begin{eqnarray*}
x \prec ^h y & \textrm{iff} & h^{-1}(x) \prec h^{-1}(y) \quad
\textrm{by definition of} \, \,\prec ^h \\
& \textrm{iff} & g^{-1}(x) \prec g^{-1}(y) \quad \textrm{by the
  observation above} \\
& \textrm{iff} & g^{-1}(x) < g^{-1}(y) \quad
\textrm{since} \, \, h \in U \\
& \textrm{iff} & x <^g y \quad \textrm{by definition of} <^g
\end{eqnarray*}

So $\prec ^h \in W$ and the action is continuous. 

To complete the proof of the theorem, notice that the restriction map $\psi$ defined by $\psi : \funct{\cLO (\cUr)}{\cLO (\uUr)}$ with $\psi (<) = \restrict{<}{\uUr}$ is actually a homeomorphism. The proof of that fact is easy thanks to the coherence property and is left to the reader. It follows that $\cLO(\cUr)$ can be seen as the universal minimal flow of $\iso (\uUr)$ via the action $\alpha : \funct{\iso (\uUr) \times \cLO (\cUr)}{\cLO (\cUr)}$ defined
by 

\begin{center}
$\alpha (g , <) = \psi ^{-1} ( \psi (<)^g)$. 
\end{center}

Now, observe that if $g \in \iso (\uUr)$ and $< \in \cLO (\cUr)$, then 

\begin{center}
$\restrict{<^{\varphi (g)}}{\uUr} = (\restrict{<}{\uUr})^g$. 
\end{center}

It follows that $\psi (<^{\varphi (g) }) = \psi (<)^g$ and thus $\alpha (g, <) = \psi ^{-1} ( \psi (<)^g) = <^{\varphi (g)}$. Observe also that there is a natural dense embedding $\varphi : \funct{\iso (\uUr)}{\iso (\cUr)}$ (recall that $\iso (\uUr)$ is equipped with the
pointwise convergence topology coming from the discrete topology on
$\uUr$ whereas $\iso (\cUr)$ is equipped with the pointwise
convergence topology coming from the metric topology on
$\cUr$).

Now, let $X$ be a minimal $\iso (\cUr)$-flow. Since $\varphi$ is continuous with dense range, the action $\beta : \funct{\iso (\uUr) \times X}{X}$ defined by $\beta (g,x) = \varphi (g) \cdot x$ is continuous with dense orbits and allows to see $X$ as a minimal $\iso (\uUr)$-flow. Now, by one of the previous comments, $\cLO (\cUr)$ is the universal minimal $\iso (\uUr)$-flow so there is a continuous and onto $\pi : \funct{\cLO (\cUr)}{X}$ such that for every $g$ in $\iso (\uUr)$ and every $<$ in $\cLO (\cUr)$, $\pi ( \alpha (g, <)) = \beta (g, \pi (<))$, i.e. $\pi (<^{\varphi (g) }) = \varphi (g) \cdot \pi (<)$. To finish the proof, it suffices to show that this equality remains true when $\varphi (g)$ is replaced by any $h$ in $\iso (\cUr)$. But this is easy since $\varphi$ is continuous with dense range, $\pi$ is continuous, and the actions of $\iso (\cUr)$ on $\cLO (\cUr)$ and $X$ considered here are continuous. \end{proof}

We finish with several remarks. The
first one is a purely topological comment along the lines of the
remark following Theorem \ref{thm:UMF for iso(uUr)}: To show that the underlying
space related to the universal minimal flow of $\iso(\cUr)$ is $\cLO (\cUr)$, we used the fact
that the restriction map $\psi : \funct{\cLO (\cUr)}{\cLO (\uUr)}$
defined by $\psi (<) = \restrict{<}{\uUr}$ is a homeomorphism. The space $\cLO (\uUr)$ being metrizable, we consequently get:

\begin{thm}
The underlying space of the universal minimal flow of $\iso(\cUr)$ is metrizable.
\end{thm}

The second consequence is based on the simple observation that when
the distance set $S$ is $\{1/n : n \in \omega \smallsetminus \{ 0 \} \}$, $\cUr$ is the Baire
space $\mathcal{N}$\index{$\mathcal{N}$}. Hence:

\begin{thm}
When $\mathcal{N}$ is equipped with the product metric, the universal
minimal flow of $\mathrm{iso}(\mathcal{N})$ is the set of all convex
linear orderings on $\mathcal{N}$. 
\end{thm}

\subsection{Urysohn spaces $\Ur _S$.}

We finish this section on topological dynamics with results about the spaces $\Ur _S$ associated to the classes $\M _S$\index{$\M _S$}. 
When $S$ is a subset of $]0, + \infty [$ satisfying the $4$-values condition, the class $\nM _S$\index{$\nM _S$} is a reasonable Fra\"iss\'e class. It follows that $\mathrm{Flim} (\nM _S) = (\Ur _S , <^{\Ur _S})$ for some metric linear ordering $<^{\Ur _S}$ on $\Ur _S$. Furthermore, we saw that $\nM _S$ has the Ramsey and the ordering properties whenever $S$ has size less or equal to $3$. Consequently:

\begin{thm}

\label{thm:Aut(Ur_S,<) ea}
\index{extreme amenability! of $\mathrm{Aut} (\Ur _S , <^{\Ur _S})$}

Let $S$ be finite subset of $]0, + \infty [$ of size $|S| \leqslant 3$ and satisfying the $4$-values condition. 
Then $\mathrm{Aut} (\Ur _S , <^{\Ur _S})$ is extremely amenable. 
\end{thm}

\begin{thm}

\label{thm:M(iso(Ur_S,<))}
\index{flow!universal minimal flow!of $\iso (\Ur _S)$}

Let $S$ be finite subset of $]0, + \infty [$ of size $|S| \leqslant 3$ and satisfying the $4$-values condition. Then the universal minimal flow of $\iso (\Ur _S)$ is the set $\nLO (\Ur _S)$ of metric linear orderings on $\Ur _S$ together with the action $\funct{\iso (\Ur _S) \times \nLO (\Ur _S)}{\nLO (\Ur _S)}$, $(g,<)
\longmapsto <^g$ defined
by $x <^g y$ iff $g^{-1}(x) < g^{-1}(y)$.

\end{thm}

\section{Concluding remarks and open problems.}

The purpose of this section is to present several questions related to the Ramsey calculus of finite metric spaces that we were not able to solve. 

\subsection{Classes $\nM _S$ when $|S|$ is finite. }

The first question we would like to present concerns the generalization of Theorem \ref{thm:RP for nM_S} and Theorem 
\ref{thm:OP for nM_S}. We showed that when $S$ is a finite subset of $]0, + \infty [$ of size $|S| \leqslant 3$ satisfying the $4$-values condition, the class $\nM _S$\index{$\nM _S$} of all finite metrically ordered metric spaces with distances in $S$ has the Ramsey property and the ordering property. For $|S|=4$, the verification is being carried out. So far, all the results provide a positive answer to: 

\

\textbf{Question 0.} Let $S$ be a finite subset of $]0, + \infty [$ satisfying the $4$-values condition. Does the class $\nM _S$ have the Ramsey property and the ordering property? If so, is finiteness of $S$ really necessary?   

\

\textbf{Remark.} We mentioned after Theorem \ref{thm:Aut(uUr, <) ea} that extreme amenability results can sometimes be proved directly via algebraic methods and may allow to deduce new Ramsey theorems. The classes $\nM _S$ where $|S| \leqslant 3$ and $S$ satisfies the $4$-values condition provide other illustrations of that fact. For example, the group $\mathrm{Aut} (\Ur _{\{1, 2, 5\}} , <^{\Ur _{\{1, 2, 5\}}})$ can be seen as a semi-direct product of $\mathrm{Aut} (\Q , <)$ and $\mathrm{Aut} (\mathcal{R} , < ^{\mathcal{R}})^{\Q}$ where $(\mathcal{R} , < ^{\mathcal{R}})$ is the Fra\"iss\'e limit of the class $\mathcal{G} ^<$ of all finite ordered graphs. The group $\mathrm{Aut} (\Q , <)$ is extremely amenable because thanks to the usual finite Ramsey theorem, the class $\mathcal{LO}$ of all the finite linear orderings is a Ramsey class (extreme amenability of $\mathrm{Aut} (\Q , <)$ was originally proved by Pestov\index{Pestov!theorem on $\mathrm{Aut} (\Q , <)$} in \cite{Pe-1} before \cite{KPT} and corresponds to one of the very first examples of non-trivial extremely amenable groups). On the other hand, the group $\mathrm{Aut} (\mathcal{R} , < ^{\mathcal{R}})$ is extremely amenable because $\mathcal{G} ^<$ is a Ramsey class. It follows that $\mathrm{Aut} (\Ur _{\{1, 2, 5\}} , <^{\Ur _{\{1, 2, 5\}}})$ is extremely amenable. The same holds for $\mathrm{Aut} (\Ur _{\{1, 3, 6\}} , <^{\Ur _{\{1, 3, 6\}}})$, which can be seen as a semi-direct product of $\mathrm{Aut} (\mathcal{R} , < ^{\mathcal{R}})$ and $\mathrm{Aut} (\Q , <) ^{\Q}$. Unfortunately there are some cases like $S = \{ 1, 3, 4\}$ where such an analysis does not seems to be possible (it \emph{is} unfortunate because such a generalized phenomenon might have allowed to attack the first part of Question 0 by induction on the size of $S$).

\subsection{Euclidean metric spaces.}

The second question we would like to present is related to a field that we mentioned in chapter 1 but that we did not even touch: Euclidean Ramsey theory. To make the motivation clear, let us start with the following results in topological dynamics: 

\begin{thm}[Gromov-Milman \cite{GM}]

\label{thm:Gromov-Milman}
\index{Gromov-Milman theorem}

Equipped with the pointwise convergence topology, the group $\iso (\mathbb{S} ^{\infty})$ of all surjective isometries of $\mathbb{S} ^{\infty}$\index{$\mathbb{S} ^{\infty}$} is extremely amenable.  

\end{thm}

\begin{thm}[Pestov \cite{Pe0}]

\label{cor:de la Harpe-Valette}
\index{Pestov!theorem on $\iso (\ell _2)$}

Equipped with the pointwise convergence topology, the group $\iso (\ell _2)$ of all surjective isometries of $\ell_2$ is extremely amenable.  

\end{thm} 

In \cite{Pe0}, Theorem \ref{thm:Gromov-Milman} is proved thanks to the same method as the one used to prove Theorem \ref{thm:Pestov}. This latter result being the consequence of the Ramsey property for $\M ^< _{\Q}$, it is therefore conceivable that a Ramsey result is hidden behind Theorem \ref{thm:Gromov-Milman} and and Theorem \ref{cor:de la Harpe-Valette}. Some theorems from Euclidean Ramsey theory seem to suggest that there is some hope: 
Recall that $\E$ is the class consisting of all the finite affinely independent metric subspaces of the Hilbert space $\ell _2$. Let $\m{K}_1$ denote the unique element of $\E$ with only one point.  

\begin{thm}[Frankl-R\"odl \cite{FR}]

\label{thm:Frankl-Rodl}
\index{Frankl-R\"odl theorem}

Let $\m{Y} \in \E $ and $k>0$ be in $\omega$. Then there is a finite metric subspace $\m{Z}$ of $\ell_2$ such that $\m{Z} \arrows{(\m{Y})}{\m{K}_1}{k}$. 
\end{thm}

A result of similar flavor holds for the class of $\ES $ of all elements $\m{X}$ of $\E$ which embed isometrically into $\mathbb{S} ^{\infty}$ with the property that $\{ 0_{\ell _2} \} \cup \m{X}$ is affinely independent.

\begin{thm}[Matou\v{s}ek-R\"odl \cite{MR}]

\label{thm:Matousek-Rodl}
\index{Matou\v{s}ek-R\"odl theorem}

Let $\m{Y} \in \ES $ and $k>0$ be in $\omega$. Then there is a finite metric subspace $\m{Z}$ of $\mathbb{S} ^{\infty}$ such that $\m{Z} \arrows{(\m{Y})}{\m{K}_1}{k}$. 
\end{thm} 

Recall that we proved in the previous chapter that the classes $\E _S$ and $\ES _S$ when $S \subset ]0, + \infty[$ is dense and countable are strong amalgamation classes, and that the metric completions of the corresponding Fra\"iss\'e limits are $\ell _2$ and $\mathbb{S} ^{\infty}$ respectively. Therefore, Theorem \ref{thm:Frankl-Rodl} and Theorem \ref{thm:Matousek-Rodl} may be seen as the first steps towards general Ramsey theorems about Euclidean metric spaces. However, the difficulty posed by the combinatorics of those spaces has so far kept us away from any progress in this direction. This may not be so surprising to the combinatorialist: Euclidean Ramsey theory is a well-known source of difficult problems. For example, following Graham in \cite{G}, say that a finite metric subspace of $\ell_2$ is \emph{spherical} if it can be embedded into a sphere (of finite radius). A known result due to Erd\H os, Graham, Montgomery, Rothschild, Spencer and Straus, asserts that:

\begin{thm}[Erd\H os et al. \cite{EetAl}]
Let $\m{Y}$ be a finite metric subspace of $\ell _2$ such that for every $k >0$ in $\omega$, there is a finite metric subspace $\m{Z}$ of $\ell _2$ such that $\m{Z} \arrows{(\m{Y})}{\m{K}_1}{k}$. Then $\m{Y}$ is spherical. 
\end{thm}

On the other hand, knowing whether the converse of this theorem holds or not is probably the most important open problem in Euclidean Ramsey theory. Following the tradition initiated by Erd\H os, there is even a $\$1000$ reward for the solution! Note that Theorem \ref{thm:Frankl-Rodl} quoted above provides a partial result towards a positive answer. 

Another very similar open problem asks for a characterization of those finite metric subspaces $\m{Y}$ of $\ell _2$ for which for every strictly positive $k \in \omega$ there is a finite spherical $\m{Z}$ such that $\m{Z} \arrows{(\m{Y})}{\m{K}_1}{k}$. A strong version of Theorem \ref{thm:Matousek-Rodl} actually says that every affinely independent $\m{Y}$ has this property, but to our knowledge this is the only known case so far.  

As for the problems we are interested in, they look slightly different, but still may be subject to the same kind of difficulties. In particular, we are able to prove that the metric space $\m{Z}$ from Theorem \ref{thm:Frankl-Rodl} and Theorem \ref{thm:Matousek-Rodl} can be constructed so as to stay in the relevant class (meaning $\E _S$ or $\ES _S$) but cannot show that we can work with ordered metric spaces instead of $\m{Y}$ and $\m{Z}$. The kind of linear orderings to be considered is consequently unclear, even though the results of the previous sections strongly suggest that the class of all linear orderings is the most relevant one. We state all these guesses precisely: 

\

\textbf{Question 1.} Let $S$ be a dense subset of $]0, + \infty[$. Is the class $\E ^< _S$\index{$\E ^< _S$} consisting of all the finite ordered affinely independent metric subspaces of the Hilbert space $\ell _2$ with distances in $S$ a Ramsey class (such a result would be, in some sense, a generalization of Theorem \ref{thm:Frankl-Rodl})? Does it have the ordering property?  

\

\textbf{Question 2.} Same question with the class $\ES ^< _S$\index{$\ES ^< _S$} of all finite ordered $\m{X}$ of $\E$ with distances in $S$ and which embed isometrically into $\mathbb{S} ^{\infty}$ with the property that $\{ 0_{\ell _2} \} \cup \m{X}$ is affinely independent (such a result would, in turn, provide a generalization of Theorem \ref{thm:Matousek-Rodl}).







\chapter{Big Ramsey degrees, indivisibility and oscillation stability.}

\section{Fundamentals of infinite metric Ramsey calculus and oscillation stability.}

\label{section:Definitions and notations}

Recall that given a Fra\"iss\'e class $\mathcal{K}$ of $L$-structures and $\m{X} \in \mathcal{K} $, the Ramsey degree $\mathrm{t}_{\mathcal{K}}(\m{X})$ of $\m{X}$ in $\mathcal{K}$ is defined when there is $l \in \omega $ such that for any $\m{Y} \in \mathcal{K}$, and any $k \in \omega \smallsetminus \{ 0 \}$, there exists $\m{Z} \in \mathcal{K} $ such that:

\begin{center}
$\m{Z} \arrows{(\m{Y})}{\m{X}}{k,l}$. 
\end{center}

In this case, $\mathrm{t}_{\mathcal{K}}(\m{X})$ is simply the least such $l$. Equivalently, if $\m{F}$ denotes the Fra\"iss\'e limit of $\mathcal{K}$, $\m{X}$ admits a Ramsey degree in $\mathcal{K}$ when there is $l \in \omega$ such that for any $\m{Y} \in \mathcal{K}$, and any $k \in \omega \smallsetminus \{ 0 \}$,

\begin{center}
$\m{F} \arrows{(\m{Y})}{\m{X}}{k,l}$. 
\end{center} 

If this latter result remains valid when $\m{Y}$ is replaced by $\m{F}$, we say, following \cite{KPT}, that $\m{X}$ has a \emph{big Ramsey degree in $\mathcal{K}$}\index{Ramsey!big Ramsey degree}. Its value $\mathrm{T}_{\mathcal{K}}(\m{X})$\index{$\mathrm{T}_{\mathcal{K}}(\m{X})$} is the least $l \in \omega$ such that 

\begin{center}
$\m{F} \arrows{(\m{F})}{\m{X}}{k,l}$. 
\end{center} 

The notion of big Ramsey degree can be seen as a generalization of the notion of \emph{indivisibility}\index{indivisibility}. $\m{F}$ is \emph{indivisible}\index{indivisibility!for a structure} when for every strictly positive $k \in \omega$ and every $\chi : \funct{\m{F}}{k}$, there is $\mc{F} \subset \m{F}$ and isomorphic to $\m{F}$ on which $\chi$ is constant. When $\mathcal{K}$ is a class of finite metric spaces, $\m{F}$ is the Urysohn space associated to $\mathcal{K}$ and it is indivisible\index{indivisibility!for metric space} when given every strictly positive $k \in \omega$ and every $\chi : \funct{\m{F}}{k}$, there is an isometric copy $\mc{F}$ of $\m{F}$ included in $\m{F}$ on which $\chi$ is constant. It turns out that as pointed out in \cite{DLPS}, the notion of indivisiblity is too strong a concept to be studied in a general setting. For example, as soon as a complete separable metric space $\m{X}$ is uncountable, there is a partition of $\m{X}$ into two pieces such that none of the pieces includes a copy of the space via a continuous $1-1$ map. This is the reason for which relaxed versions of indivisibility were introduced. If $\m{X} = (X, d^{\m{X}})$ is a metric space, $Y \subset X$ and $\varepsilon > 0$, set \index{$(Y) _{\varepsilon}$} 

\begin{center}
$(Y) _{\varepsilon} = \{ x \in X : \exists y \in Y \ \ d ^{\m{X}} (x,y) \leqslant \varepsilon \}$
\end{center}

Now, say that $\m{X}$ is \emph{$\varepsilon$-indivisible}\index{indivisibility!$\varepsilon$-indivisibity} when for every strictly positive $k \in \omega$ and every $\chi : \funct{\m{X}}{k}$, there is $i<k$ and $\mc{X} \subset \m{X}$ isometric to $\m{X}$ such that 

\begin{center}
$\mc{X} \subset (\overleftarrow{\chi} \{ i \})_{\varepsilon}$. 
\end{center}

Equivalently, $\m{X}$ is $\varepsilon$-indivisible when for every finite cover $\gamma$ of $\m{X}$ there is $A \in \gamma$ and $\mc{X} \subset \m{X}$ isometric to $\m{X}$ such that 

\begin{center}
$\mc{X} \subset (A)_{\varepsilon}$. 
\end{center}

When $\m{X}$ is $\varepsilon$-indivisible for every $\varepsilon > 0$, $\m{X}$ is \emph{approximately indivisible}\index{indivisibility!approximate indivisibility}. When $\m{X}$ is complete and ultrahomogeneous metric space, this notion corresponds to the notion of \emph{oscillation stability}\index{oscillation!oscillation stability} introduced in \cite{KPT}. To present this concept, we start with a short reminder about \emph{uniform spaces}\index{uniform!uniform space}. Given a set $X$, a \emph{uniformity}\index{uniformity} on $X$ is a collection $\mathcal{U}$ of subsets of $X \times X$ called \emph{entourages}\index{entourage} satisfying the following properties: 

\begin{enumerate}

\item  $\mathcal{U}$ is closed under finite intersections and supersets.

\item Every $V \in \mathcal{U}$ includes the diagonal $\Delta = \{ (x,x) : x \in X \}$. 

\item If $V \in \mathcal{U}$, then $V^{-1}:= \{ (y,x) : (x,y) \in V \} \in \mathcal{U}$.

\item If $V \in \mathcal{U}$, there exists $U \in \mathcal{U}$ such that 

\begin{center} $U \circ U := \{ (x,z) : \exists y \in U \ \ ((x,y) \in U \ \mathrm{and} \ (y,z) \in U )\} \subset V$. \end{center}

\end{enumerate}

$(X,\mathcal{U})$ is then called a \emph{uniform space}. A \emph{basis}\index{uniformity!basis for a uniformity} for $\mathcal{U}$ is a family $\mathcal{B} \subset \mathcal{U}$ such that for every $U, V \in \mathcal{U}$, there is $W \in \mathcal{B}$ such that $W \subset U \cap V$.  

Every uniform space $(X , \mathcal{U})$ carries a structure of topological space $(X , T_{\mathcal{U}})$ by declaring a subset $O$ of $X$ to be open if and only if for every $x$ in $O$ there exists an entourage $V$ such that $\{ y \in X : (x, y) \in V \}$ is a subset of $O$. $(X, \mathcal{U})$ is \emph{separated}\index{uniform!separated uniform space} when $(X , T_{\mathcal{U}})$ is, or equivalently when $\bigcap \mathcal{U} = \Delta$. A sequence $(x_n)_{n \in \omega}$ of elements of $X$ is \emph{Cauchy}\index{uniform!Cauchy sequence in a uniform space} when \begin{center} $\forall V \in \mathcal{U} \ \exists N \in \omega \ \forall p,q \in \omega \ \ (q \geqslant p \geqslant N \rightarrow (x_q , x_p) \in V)$ \end{center} and $(X , \mathcal{U})$ is \emph{complete}\index{uniform!complete uniform space} when every Cauchy sequence in $(X , \mathcal{U})$ converges in $(X , T_{\mathcal{U}})$. Uniform spaces constitute the natural setting where \emph{uniform continuity}\index{uniform!uniform continuity} can be defined: Given two uniform spaces $(X,\mathcal{U})$ and $(Y,\mathcal{V})$, a map $f : \funct{X}{Y}$ is \emph{uniformly continuous} when 

\begin{center}
$\forall V \in \mathcal{V} \ \exists U \in \mathcal{U} \ \forall (x,y) \in X\times X  \ \ ((x,y)\in U \rightarrow (f(x),f(y)) \in V)$.
\end{center} 

When additionally $f$ is bijective and $f^{-1}$ is uniformly continuous, $f$ is called a \emph{uniform homeomorphism}\index{uniform!uniform homeomorphism}. Given a separated uniform space $(X , \mathcal{U})$, there is, up to uniform homeomorphism, a unique complete uniform space $(\widehat{X} , \widehat{\mathcal{U}})$ including $(X , \mathcal{U})$ as a dense uniform subspace, called the \emph{completion}\index{uniform!completion of a uniform space} of $(X , \mathcal{U})$. In what follows, we will be particularly interested in uniform structures coming from topological groups. In particular, for a topological group $G$, the \emph{left uniformity}\index{uniformity!left uniformity} $\mathcal{U} _L (G)$\index{$\mathcal{U} _L (G)$} is the uniformity whose basis is given by the sets of the form $V_L = \{(x,y):x^{-1}y \in V \}$ where $V$ is a neighborhood of the identity. Now, let $\widehat{G}^L$\index{$\widehat{G}^L$} denote the completion of $(G, \mathcal{U} _L (G))$. In general,  $\widehat{G}^L$ is not a topological group (see \cite{Di}). However, it is always a topological semigroup (see \cite{RoD}). For a real-valued map $f$ on a set $X$, define the \emph{oscillation}\index{oscillation!oscillation of a map} of $f$ on $X$ as\index{$\mathrm{osc}(f)$}:

\begin{center}
$\mathrm{osc}(f) = \sup \{ \left| f(y) - f(x) \right| : x, y \in X \}$. 
\end{center} 

\begin{defn}

Let $G$ be a topological group, $ f : \funct{G}{\R}$ be uniformly continuous, and $\hat{f}$ be the unique extension of $f$ to $\widehat{G}^L$ by uniform continuity. $f$ is \emph{oscillation stable}\index{oscillation!oscillation stable map} when for every $\varepsilon > 0$, there is a right ideal $\mathcal{I} \subset \widehat{G}^L$ such that 

\begin{center}
$\mathrm{osc}(\restrict{\hat{f}}{\mathcal{I}}) < \varepsilon$. 
\end{center}

\end{defn}

\begin{defn}

Let $G$ be a topological group acting $G$ continuously on a topological space $X$. For $f : \funct{X}{\R}$ and $x \in X$, let $f_{x} : \funct{G}{\R}$ be defined by

\begin{center}
$\forall g \in G \ \ f_{x} (g) = f(gx)$.
\end{center}

Then the action is \emph{oscillation stable}\index{oscillation!oscillation stable action} when for every $f : \funct{X}{\R}$ bounded and continuous and every $x \in X$, $f_{x}$ is oscillation stable whenever it is uniformly continuous. 

\end{defn}

With these concepts in mind, we are now ready to link oscillation stability to the Ramsey-type properties introduced previously: It turns out that when $G$ is the group $\iso (\m{X})$ of all isometries from $\m{X}$ onto itself equipped with the pointwise convergence topology, $\widehat{G}^L$ can be thought as a topological subsemigroup of the topological semigroup $\mathrm{Emb}(\m{X})$\index{$\mathrm{Emb}(\m{X})$} of all isometric embeddings from $\m{X}$ into itself. 

\begin{thm}[Kechris-Pestov-Todorcevic \cite{KPT}, Pestov \cite{Pe1}, \cite{Pe1'}]

\index{Kechris-Pestov-Todorcevic!theorem on oscillation stability}

Let $G$ be a topological group acting continuously and transitively on a complete metric space $\m{X}$ by isometries. Then the following are equivalent:

\begin{enumerate}

\item The action of $G$ on $\m{X}$ is oscillation stable. 

\item Every bounded real-valued $1$-Lipschitz map $f$ on $\m{X}$ is oscillation stable. 

\item For every strictly positive $k \in \omega$, every $\chi : \funct{\m{X}}{k}$ and every $\varepsilon > 0$, there are $g \in \widehat{G}^L$ and $i < k $ such that $g''X \subset (\overleftarrow{\chi} \{ i \})_{\varepsilon}$. 

\end{enumerate}

\end{thm}

When one of those equivalent conditions is fullfilled, $\m{X}$ is \emph{oscillation stable}\index{oscillation!oscillation stable metric space}. In addition, one can check that when the metric space $\m{X}$ is ultrahomogeneous, then $\widehat{G}^L$ is actually equal to $\mathrm{Emb}(\m{X})$. For that reason, in the realm of ultrahomogeneous metric spaces the previous theorem can be stated as follows:

\begin{cor}
For a complete ultrahomogeneous metric space $\m{X}$, the following are equivalent:

\begin{enumerate}

\item When $\iso (\m{X})$ is equipped with the topology of pointwise convergence, the standard action of $\iso (\m{X})$ on $\m{X}$ is oscillation stable. 

\item For every bounded $1$-Lipschitz map $f : \funct{\m{X}}{\R}$ and every $\varepsilon > 0$, there is an isometric copy $\mc{X}$ of $\m{X}$ in $\m{X}$ such that \begin{center} $\mathrm{osc}(\restrict{f}{\mc{X}}) < \varepsilon$. \end{center}

\item $\m{X}$ is approximately indivisible. 

\end{enumerate}

\end{cor}

In particular, for complete ultrahomogeneous metric spaces, oscillation stability and approximate indivisibility coincide. In the more general context of structural Ramsey theory, big Ramsey degrees and oscillation stability for topological groups are also closely linked. For more information about this connection, see \cite{KPT}, section 11(E), or the books \cite{Pe1}, \cite{Pe1'}. 

\

\textbf{Remark.} Though quite close in essence, the concept of oscillation stability presented here is, except in the notable case of the Hilbert space, \emph{not} the same as the classical concept of oscillation stability used in Banach space theory. For more details, see the remark at the end of the introduction of section \ref{section:Almost indivisibility and oscillation stability}.

\

This chapter is organized as follows. In section 2, we cover the only case for which the analysis of the big Ramsey degree can be carried out: Ultrametric spaces. In section 3, we study the indivisibility properties of the countable Urysohn spaces. We finish in section 4 with a solution of the oscillation stability problem (equivalently, of the approximate indivisiblity problem) in two particular cases: The complete separable ultrahomogeneous ultrametric spaces and the Urysohn sphere.

\section{Big Ramsey degrees.}

In this section, we present the only case where we were able to provide a complete analysis for the big Ramsey degree: Ultrametric spaces. 

\label{section:Big Ramsey degrees}

\begin{thm}

\label{thm:Big Ramsey degrees for U, S finite}

Let $S$ be a finite subset of $]0, + \infty[$. Then every element of $\U$\index{$\U$} has a big Ramsey degree in $\U$.

\end{thm}

\begin{thm}

\label{thm:no Big Ramsey degrees for U, S infinite}

Let $S$ be an infinite countable subset of $]0, + \infty[$ and let $\m{X}$ be in $\U$ such that $|\m{X}| \geqslant 2$. Then $\m{X}$ does not have a big Ramsey degree in $\U$.
  
\end{thm}

The ideas we use to reach this goal are not new. The way we met them is through some unpublished work of Galvin, but in \cite{M}, Milner writes that they were also known to and exploited by several other authors, among whom Hajnal (who apparently realized first the equivalent of lemma \ref{lemma:3} and stated it explicitly in \cite{H}), and Haddad and Sabbagh (\cite{HS1}, \cite{HS2} and \cite{HS3}).



Recall that when $S$ is finite and given by elements $s_0 > s_1 \ldots > s_{|S|-1} > 0$, it is convenient to see the space $\uUr$\index{$\uUr$} as the set $\omega^{|S|}$ of maximal nodes of the tree $\omega^{\leqslant |S|} = \bigcup_{i\leqslant |S|} \omega^i$ ordered by set-theoretic inclusion and equipped with the metric defined for $x \neq y$ by \begin{center} $d(x,y) = s_{\Delta (x,y)}$\end{center} where \[\Delta (x,y) = \min \{k < |S|-1 : s(k) \neq t(k) \}\index{$\Delta (s,t)$}.\]

For $A \subset \omega ^{|S|}$, set \index{$ A^\downarrow$} \begin{center} $ A^\downarrow = \{ \restrict{a}{k} : a \in A \ \mathrm{and} \ k \leqslant n \} $.\end{center} 

It should be clear that when $A, B \subset \omega ^{|S|}$, then $A$ and $B$ are isometric iff $A^\downarrow \cong B^\downarrow$. Consequently, when $\m{X} \in \U$, one can define the natural tree associated to $\m{X}$ in $\U$ to be the unique (up to isomorphism) subtree $\m{T}_{\m{X}}$ of $\omega^{\leqslant |S|}$ such that for any copy $\mc{X}$ of $\m{X}$ in $\uUr$, $\mc{X} ^\downarrow \cong \m{T}_{\m{X}}$.

Given a subtree $\m{T}$ of $\omega ^{|S|}$, set \begin{center} $ \binom{\omega ^{\leqslant|S|}}{\m{T}} = \{ \mc{T} : \mc{T} \subset \omega ^{\leqslant|S|} \ \mathrm{and} \ \mc{T} \cong \m{T} \}$.\end{center} 

When $k, l \in \omega \smallsetminus \{ 0 \}$ and for any $\chi : \funct{\binom{\omega ^{\leqslant|S|}}{\m{T}}}{k}$ there is
$\m{U} \in \binom{\omega ^{\leqslant|S|}}{\omega ^{\leqslant|S|}}$ such that $\chi$ takes at most $l$ values on $\binom{\m{U}}{\m{T}}$, we write \begin{center} $ \omega ^{\leqslant|S|}
\arrows{(\omega ^{\leqslant|S|})}{\m{T}}{k,l}$.\end{center} 

If there is $l \in \omega \smallsetminus \{ 0 \}$ such that for any $k \in \omega \smallsetminus \{ 0 \}$, $\omega ^{\leqslant|S|}
\arrows{(\omega ^{\leqslant|S|})}{\m{T}}{k,l}$, the least such $l$ is called the \emph{Ramsey degree} of $\m{T}$ in $\omega ^{\leqslant|S|}$. 

\begin{lemma}

\label{lemma:3}

Let $X \subset \omega ^{|S|}$ and let $\m{T} = X^\downarrow$. Then $\m{T}$ has a Ramsey degree in $\omega ^{\leqslant|S|}$ equal to $|e(\m{T})|$. 
\end{lemma}

\begin{proof}
Say that a subtree $\m{U}$ of $\omega ^{\leqslant|S|}$ is \emph{expanded}\index{expanded subtree of $\omega ^{\leqslant|S|}$} when: 

\begin{enumerate}
 
\item Elements of $\m{U}$ are strictly increasing. 

\item For every $u, v \in \m{U}$ and every $k \in |S|$, \begin{center} $u(k) \neq v(k) \rightarrow (\forall j \geqslant k \ \ u(j) \neq v(j))$.\end{center} 

\end{enumerate}

Note that every expanded $\mc{T} \in \binom{\omega ^{\leqslant|S|}}{\m{T}}$ is linearly ordered by $\prec^{\mc{T}}$ defined by \begin{center} $s \prec^{\mc{T}} t$ iff ($s=\emptyset$ or $s(|s|) < t(|t|)$). \end{center} 

Note also that then $\prec^{\mc{T}}$ is a linear extension of the tree ordering on $\mc{T}$. Now, given $\prec \in e(\m{T})$, let $\binom{\omega ^{\leqslant|S|}}{\m{T},\prec}$ denote the set of all expanded $\mc{T} \in \binom{\omega ^{\leqslant|S|}}{\m{T}}$ \emph{with type} $\prec$, that is, such that the order-preserving bijection between the linear orderings $(\mc{T}, \prec^{\mc{T}})$ and $(\m{T}, \prec)$ induces an isomorphism between the trees $\mc{T}$ and $\m{T}$. Define the map $\psi _{\prec} : \funct{\binom{\omega ^{\leqslant|S|}}{\m{T},\prec}}{[\omega]^{|\m{T}|-1}}$ by 
\begin{center} $ \psi _{\prec} (\mc{T}) = \{ t(|t|) : t \in \mc{T} \smallsetminus \{ \emptyset \} \}.$ \end{center}  

Then $\psi _{\prec}$ is a bijection. Call $\varphi _{\prec}$ its inverse map. Now, let $k \in \omega \smallsetminus \{ 0 \}$ and $\chi : \funct{\binom{\omega ^{\leqslant|S|}}{\m{T}}}{k}$. Define $\Lambda : \funct{[\omega]^{|T|-1}}{k^{e(\m{T})}}$ by 
\begin{center} $ \Lambda (M) = (\chi (\varphi _{\prec}(M)))_{\prec \in e(\m{T})} $.\end{center}   

By Ramsey's theorem, find an infinite $N \subset \omega$ such that $\Lambda$ is constant on $[N]^{|\m{T}|-1}$. Then, on the subtree $N^{\leqslant |S|}$ of $\omega^{\leqslant |S|}$, any two expanded elements of $\binom{\omega ^{\leqslant|S|}}{\m{T}}$ with same type have the same $\chi$-color. Now, let $\m{U}$ be an expanded everywhere infinitely branching subtree of $N^{\leqslant |S|}$. Then $\m{U}$ is isomorphic to $\omega^{\leqslant |S|}$ and $\chi$ does not take more than $|e(\m{T})|$ values on $\binom{\m{U}}{\m{T}}$.

To finish the proof, it remains to show that $|e(\m{T})|$ is the best possible bound. To do that, simply observe that for any $\m{U} \in \binom{\omega^{\leqslant |S|}}{\omega^{\leqslant |S|}}$, every possible type appears on $\binom{\m{U}}{\m{T}}$. \end{proof}

This lemma has two direct consequences concerning the existence of big Ramsey degrees in $\U$. Indeed, it should be clear that when $\m{X} \in \U$, $\m{X}$ has a big Ramsey degree in $\U$ iff $\m{T}_{\m{X}}$ has a Ramsey degree in $\omega ^{\leqslant|S|}$ and that these degrees are equal. Thus, Theorem \ref{thm:Big Ramsey degrees for U, S finite} follows. 

On the other hand, observe that if $S \subsetneq S'$ are finite and $\m{X} \in \U$ has size at least two, then the big Ramsey degree $T_{\mathcal{U} _{S'}}(\m{X})$ of $\m{X}$ in $\mathcal{U} _{S'}$ is strictly larger than the big Ramsey degree of $\m{X}$ in $\U$. In particular, $T_{\mathcal{U} _{S'}}(\m{X})$ tends to infinity when $|S'|$ tends to infinity. That fact can be used to prove Theorem \ref{thm:no Big Ramsey degrees for U, S infinite}. 

\begin{proof}[Proof of Theorem \ref{thm:no Big Ramsey degrees for U, S infinite}.]

It suffices to show that for every $k \in \omega \smallsetminus \{ 0 \}$, there is $k' > k$ and a coloring $\chi : \funct{\binom{\uUr}{\m{X}}}{k'}$ such that for every $B \in \binom{\uUr}{\uUr}$, the restriction of $\chi$ on $\binom{B}{\m{X}}$ has range $k'$. Thanks to the previous remark, we can fix $S' \subset S$ finite such that $X \in \mathcal{U}_{S'}$ and the big Ramsey degree $k'$ of $\m{X}$ in $\mathcal{U}_{S'}$ is larger than $k$. Recall that $\uUr \subset \omega ^S$ so if $\textbf{1} _{S'} : \funct{S}{2}$ is the characteristic function of $S'$, it makes sense to define $f : \funct{\uUr}{\textbf{B}_{S'}}$ by \[ f(x) = \textbf{1}_{S'} x.\] 

Observe that $d(f(x),f(y)) = d(x,y)$ whenever $d(x,y) \in S'$. Thus, given any $B \in \binom{\uUr}{\uUr}$, the direct image $f''B$ of $B$ under $f$ is in $\binom{\textbf{B}_{S'}}{\textbf{B}_{S'}}$. Now, let $\chi ' : \funct{\binom{\textbf{B}_{S'}}{\m{X}}}{k'}$ be such that for every $B' \in \binom{\textbf{B}_{S'}}{\textbf{B}_{S'}}$, the restriction of $\chi '$ to $\binom{B'}{\m{X}}$ has range $k'$. Then $\chi = \chi ' \circ f $ is as required. \end{proof}

\section{Indivisibility.}

As stated in the introduction of this chapter, indivisibility corresponds to the most elementary case in the analysis of the big Ramsey degrees, so one might wonder why the part of this paper devoted to indivisibility is much larger than the one about big Ramsey degrees. Here is the reason: With the exception of ultrametric spaces, the obstacles posed by indivisibility are in most of the cases substantial enough for many problems to remain open. Fortunately, there were also some recent progress, in particular thanks to the paper \cite{DLPS} by Delhomm\'e, Laflamme, Pouzet and Sauer where a detailed analysis of metric indivisibility is carried out. For example, we already mentioned a general observation from \cite{DLPS} in the introduction: No uncountable complete separable metric space is indivisible. Here is another restriction to indivisibility: 

\begin{prop}

Let $\m{X}$ be a metric space whose distance set is unbounded. Then $\m{X}$ is divisible. 

\end{prop}

\begin{proof}
We follow \cite{DLPS}. Observe that inductively, we can construct a sequence of reals $(r_n)_{n \in \omega}$ with $r_0 = 0$ together with a sequence $(x_n)_{n \in \omega}$ of elements of $X$ such that

\begin{center} 
$\forall n < \omega \ \ 2 r_n < d^{\m{X}}(x_0 , x_{n+1}) < r_{n+1} - r_n$. 
\end{center}

Now, define $\chi : \funct{\m{X}}{2}$ by setting:

\[ 
\forall x \in X \ \ \chi(x) = 0 \leftrightarrow \left( d^{\m{X}}(x_0 , x ) \in \bigcup _{n \in \omega}[ r_{2n} , r_{2n+1} [ \right). 
\]

We claim that $\chi$ divides $\m{X}$: Let $\varphi : \funct{\m{X}}{\m{X}}$ be an isometric embedding. Let $n \in \omega$ be such that $ d^{\m{X}}(x_0 , \varphi (x_0) ) \in [ r_{n} , r_{n+1} [$. Then one can check that $d^{\m{X}} (x_0 , \varphi (x_{n+2})) \in [ r_{n+1} , r_{n+2} [ $, and so $\chi (\varphi (x_0)) \neq \chi (\varphi (x_{n+2})) $. \end{proof}

It follows that even if we restrict our attention to the Urysohn spaces associated to the Fra\"iss\'e classes of finite metric spaces, some spaces may have a trivial behaviour as far as indivisibility is concerned. For example, $\Ur _{\Q}$ and $\Ur _{\omega}$ are divisible. However, we will see that when the two obstacles of cardinality and unboundedness are avoided, indivisibility can be substantially more difficult to study. During the past three years, the space whose indivisibility properties attracted most of the attention is $\s _{\Q}$\index{$\s _{\Q}$}. The question of knowing whether $\s _{\Q}$ is indivisible or not is explicitly stated in \cite{N1}, \cite{Pe1} and \cite{Pe1'}. This problem was solved in \cite{DLPS} by Delhomm\'e, Laflamme, Pouzet and Sauer, and we present their result in subsection \ref{subsection:Divisibility of S_Q}. In subsection \ref{subsection: Are the U_m 's indivisible?}, we present the first results concerning indivisibility of the spaces $\Ur _m$ when $m \in \omega$. The general solution is then presented in \ref{subsection:U_m are indiv}. In \ref{subsection:Indivisibility of Urysohn ultrametric spaces}, we consider the case of the countable ultrahomogeneous ultrametric spaces before turning to the study of indivisibility for the spaces $\Ur _S$ with $|S| \leqslant 4$ in subsection \ref{subsection:Indivisibility of U_S}.

\subsection{Divisibility of $\s _{\Q}$.}

\label{subsection:Divisibility of S_Q}

Apart from the intrinsic combinatorial interest, the motivation attached to this problem comes from the problem of the approximate indivisibility for the Urysohn sphere $\s$\index{$\s$}. Indeed, had $\s _{\Q}$ been indivisible, $\s$ would have been approximately indivisible and the standard action of $\iso (\s)$ on $\s$ would have been oscillation stable. We will however see now that the actual answer for the indivisibility problem for $\s _{\Q}$ is not the one that was hoped for. All the concepts and results presented in this subsection come from \cite{DLPS} and are due to Delhomm\'e, Laflamme, Pouzet and Sauer. 

\begin{thm}[Delhomm\'e-Laflamme-Pouzet-Sauer \cite{DLPS}]

\label{thm:s_Q divisible}
\index{Delhomm\'e-Laflamme-Pouzet-Sauer!theorem on the divisibility of $\s _{\Q}$}
\index{indivisibility!for $\s _{\Q}$}

$\s _{\Q}$ is divisible. 

\end{thm}

\begin{proof} 
Call a sequence of elements $x_0 ,\ldots , x_n$ of $\s _{\Q}$ an \emph{$\varepsilon$-chain from $x_0$ to $x_n$}\index{$\varepsilon$-chain} if for every $i<n$, $d^{\s _{\Q}}(x_i , x_{i+1}) \leqslant \varepsilon$. The key idea is the following simple geometrical fact: Let $y \in \s _{\Q}$, $r \in [0,1]$ irrational, $x \in \s _{\Q}$ and $n \in \omega$ strictly positive such that 

\[ 
d^{\s _{\Q}}(y , x) < r\cdot\left(1- \frac{1}{n+1} \right).
\]  
 
Let also $x' \in \s _{\Q}$ be such that 

\[
d^{\s _{\Q}}(x , x') > r.
\]

Finally, let $\varepsilon > 0$ be such that 

\[
\varepsilon < \frac{1}{(n+1)(n+2)}.
\] 

Then for every $\varepsilon$-chain $(x_i)_{i\leqslant n}$ from $x$ to $x'$, there is $i \leqslant n$ such that 

\[
r\cdot\left(1 - \frac{1}{n+1} \right) \leqslant d^{\s _{\Q}}(y , x_i) < r\cdot\left(1 - \frac{1}{n+2} \right).
\]

With this fact in mind, we now prove that $\s _{\Q}$ is divisible. First, construct inductively a subset $Y$ of $\s _{\Q}$ together with a family $(r_y)_{y \in Y}$ of irrationals in $]0,1/2[$ such that 

\begin{center}
$\forall x \in \s _{\Q} \ \exists ! y_x \in Y \ \ d^{\s _{\Q}} (y_x , x) < r_x$.
\end{center}

Now, let $\chi : \funct{\s _{\Q}}{2}$ be defined by 

\[
\chi (x) = 0 \leftrightarrow \left(  \exists n > 0 \ \ r_{y_x} . \left(1 - \frac{1}{2n} \right) \leqslant d^{\s _{\Q}} (y_x , x) < r_{y_x} . \left(1 - \frac{1}{2n+1} \right) \right). 
\]

We claim that $\chi$ divides $\s _{\Q}$. Indeed, let $\mc{S} _{\Q}$ be an isometric copy of $\s _{\Q}$ in $\s _{\Q}$. Fix $x \in \mc{S} _{\Q}$, and consider $n > 0$ such that 

\[
r_{y_x}\cdot\left(1 - \frac{1}{n} \right) \leqslant d^{\s _{\Q}} (y_x , x) < r_{y_x}\cdot\left(1 - \frac{1}{n+1} \right).
\]

In $\mc{S} _{\Q}$, there is $x'$ such that $d^{\s _{\Q}} (x , x') > r_{y_x}$. Fix $\varepsilon > 0$ with 

\[
\varepsilon < \frac{1}{(n+1)(n+2)}.
\]   

Then in $\mc{S} _{\Q}$, there is an $\varepsilon$-chain $(x_i)_{i \leqslant n}$ from $x$ to $x'$. But by the previous property, there is $i \leqslant n$ such that 

\[
r\cdot\left(1 - \frac{1}{n+1} \right) \leqslant d^{\m{X}}(y , x_i) < r\cdot\left(1 - \frac{1}{n+2} \right).
\]

Then $\chi (x) \neq \chi (x_i)$. \end{proof}

Theorem \ref{thm:s_Q divisible} is actually only a particular case of a more general result which can be proved using the same idea. For a metric space $\m{X}$, $x \in \m{X}$, and $\varepsilon > 0$, let $\lambda _{\varepsilon} (x)$\index{$\lambda _{\varepsilon} (x)$} be the supremum of all reals $l \leqslant 1$ such that there is an $\epsilon$-chain $(x_i)_{i \leqslant n}$ containing $x$ and such that $d^{\m{X}}(x_0 , x_n) \geqslant l$. Then, define \index{$\lambda (x)$}

\begin{center}
$\lambda (x) = \sup \{ l \in \R : \forall \varepsilon > 0 \ \ \lambda _{\varepsilon} (x) \geqslant l \}$.
\end{center}

\begin{thm}[Delhomm\'e-Laflamme-Pouzet-Sauer \cite{DLPS}]

\label{thm:lambda divisible}
\index{Delhomm\'e-Laflamme-Pouzet-Sauer!theorem on divisibility of countable metric spaces}

Let $\m{X}$ be a countable metric space. Assume that there is $x_0 \in \m{X}$ such that $\lambda (x_0) > 0$. Then $\m{X}$ is divisible. 

\end{thm}

Theorem \ref{thm:s_Q divisible} then follows since in $\s _{\Q}$ every $x$ is such that $\lambda (x) = 1$.

\subsection{Are the $\Ur _m$'s indivisible?}

\label{subsection: Are the U_m 's indivisible?}
\index{indivisibility!for $\Ur _m$}

We mentioned earlier that $\Ur _{\Q}$ is divisible because its distance set is unbounded. We also saw in the previous subsection that unboundedness is not the only reason for this phenomenon as the bounded counterpart $\s _{\Q}$ of $\Ur _{\Q}$ is not indivisible either. In this subsection, we try to answer the same question when $\Ur _{\Q}$ is replaced by $\Ur _{\omega}$. This latter space is divisible because its distance set is unbounded. However, what if one works with one of its bounded versions, namely a space of the form $\Ur _m$ when $m \in \omega$? Of course, when $m = 1$, the space $\Ur _m$ is indivisible. The first non-trivial case is consequently for $m=2$. However, we mentioned in chapter 1 that $\Ur _2$ is really the Rado graph $\mathcal{R}$ where the distance is $1$ between connected points and $2$ between non-connected distinct points. Therefore, indivisiblity for $\Ur _2$\index{indivisibility!for $\Ur _2$} is equivalent to indivisibility of $\mathcal{R}$, a problem whose solution is well-known:

\begin{prop}

\label{prop:R indivisible}

The Rado graph $\mathcal{R}$ is indivisible. 
\end{prop}

\begin{proof}
Let $k \in \omega$ be strictly positive and  $\chi : \funct{\mathcal{R}}{k}$. Let $\{ x_n : n \in \omega\}$ be an enumeration of the vertices of $\mathcal{R}$. If all vertices have color $0$, we are done. Otherwise, choose $\tilde{x} _0$ such that $\chi(\tilde{x} _0) = 0$. In general, assume that $\tilde{x}_0 ,\ldots , \tilde{x}_n$ were constructed with $\chi$-color $0$ and such that

\begin{center}
$\forall i, j \leqslant n \ \ \{ \tilde{x}_i , \tilde{x}_j \} \in E^{\mathcal{R}} \leftrightarrow \{ x_i , x_j \} \in E^{\mathcal{R}}$. 
\end{center}

Now, consider the set $E$ defined by 

\begin{center}
$E = \{ x \in \mathcal{R} : \forall i \leqslant n \ \ \left( \{ \tilde{x}_i , x \} \in E^{\mathcal{R}} \leftrightarrow \{ x_i , x_{n+1} \} \in E^{\mathcal{R}} \right) \} \smallsetminus \{x_0 ,\ldots , x_n \}$.
\end{center}

If $\chi$ does not take the value $0$ on $E$, observe that the subgraph of $\mathcal{R}$ supported by $E$ is ultrahomogeneous and includes an isomorphic copy of every finite graph. Therefore, this subgraph is isomorphic to $\mathcal{R}$ itself and $\chi$ is constant on it with value $1$, so we are done. Otherwise, $\chi$ takes the value $0$ on $E$ and we choose $x_{n+1}$ in $E$ and such that $\chi(x_{n+1}) = 0$. Thus, if the construction stops at some stage, then we are left with a copy of $\mathcal{R}$ with $\chi$-color $1$. Otherwise, after $\omega$ steps, we are left with $\{ \tilde{x} _n : n \in \omega \}$ isomorphic to $\mathcal{R}$ and with $\chi$-color $0$. \end{proof}

Another possible proof for the indivisibility of $\mathcal{R}$ uses a Ramsey-type theorem known as Milliken's theorem. 
This result will be useful later in this paper to prove that Urysohn spaces more sophisticated than $\Ur _2$ are indivisible, so we present it now. The main concept attached to Milliken's theorem is the concept of \emph{strong subtree}:\index{strong subtree} Fix a downwards closed finitely branching subtree $\m{T}$ of the tree $\omega ^{< \omega}$ with height $\omega$. Thus, $\m{T}$ has a root (a smallest element), namely, the empty sequence, and the height of a node $t \in \m{T}$ is the integer $|t|$ such that $t : \funct{|t|}{\omega}$. Say that a subtree $\m{S}$ of $\m{T}$ is \emph{strong} when 

\vspace{0.5em}
\hspace{1em} 
i) $\m{S}$ is rooted,

\vspace{0.5em}
\hspace{1em}
ii) Every level of $\m{S}$ is included in a level of $\m{T}$, 

\vspace{0.5em}
\hspace{1em}
iii) For every $s \in S$ not of maximal height in $\m{S}$ and every immediate successor 

\ \hspace{1.9em} $t$ of $s$ in $\m{T}$ there is exactly one immediate successor of $s$ in $\m{S}$ extending $t$.   

\vspace{1em}

For $s , t \in T$, set \index{$s \wedge t$}

\begin{center}
$s \wedge t = \max \{ u \in T : u \subset s, \ u \subset t \}$.
\end{center}

Now, for $A \subset T$, set \index{$A^{\wedge}$}

\begin{center}
$A^{\wedge} = \{ s \wedge t : s, t \in A\}$. 
\end{center}

Note that $A \subset A^{\wedge}$ and that $A^{\wedge}$ is a rooted subtree of $\m{T}$. For $A, B \subset T$, write $A \mathrm{Em} B$\index{$A \mathrm{Em} B$} when there is a bijection $f : \funct{A^{\wedge}}{B^{\wedge}}$ such that for every $s, t \in A^{\wedge}$:

\vspace{0.5em}
\hspace{1em} 
i) $s \subset t \leftrightarrow f(s) \subset f(t)$,

\vspace{0.5em}
\hspace{1em} 
ii) $|s|<|t| \leftrightarrow |f(s)|<|f(t)|$,

\vspace{0.5em}
\hspace{1em} 
iii) $s \in A \leftrightarrow f(s) \in B$,

\vspace{0.5em}
\hspace{1em} 
iv) $t(|s|)=f(t)(|f(s)|)$ whenever $|s|<|t|$. 

\vspace{0.5em}

It should be clear that the relation $\mathrm{Em}$\index{$\mathrm{Em}$} is an equivalence relation. Given $A \subset T$, the $\mathrm{Em}$-equivalence class of $A$ is written $[A]_{\mathrm{Em}}$\index{$[A]_{\mathrm{Em}}$}. Finally, for a strong subtree $\m{S}$ of $\m{T}$, let $\restrict{[A]_{\mathrm{Em}}}{S}$\index{$\restrict{[A]_{\mathrm{Em}}}{S}$} denote the set of all elements of $[A]_{\mathrm{Em}}$ included in $S$. With these notions in mind, the version of Milliken's theorem we need can be stated as follows:

\begin{thm}[Milliken \cite{Mi}] 

\label{thm:Milliken}
\index{Milliken theorem}

Let $\m{T}$ be a nonempty downwards closed finitely branching subtree $\m{T}$ of $\omega ^{< \omega}$ with height $\omega$. Let $A$ be a finite subset of $\m{T}$. Then for every strictly positive $k \in \omega$ and every $k$-coloring of $[A]_{\mathrm{Em}}$, there is a strong subtree $\m{S}$ of $\m{T}$ with height $\omega$ such that $\restrict{[A]_{\mathrm{Em}}}{S}$ is monochromatic. 

\end{thm}

For more on this theorem and its numerous applications, the reader is referred to \cite{T1}. We now show how to deduce proposition \ref{prop:R indivisible} from Theorem \ref{thm:Milliken}. 

\begin{proof}
Let $\m{T}$ be the complete binary tree $2 ^{< \omega}$. On $\m{T}$, define the following graph structure (sometimes called the \emph{standard graph structure on $2 ^{< \omega}$}\index{standard graph structure on $2 ^{< \omega}$}) by:

\begin{center}
$\forall s < t \in 2^{< \omega} \ \ \{ s , t \} \in E \leftrightarrow \left( |s| < |t| , \ t(|s|) = 1 \right)$. 
\end{center}

Now, observe that $\mathcal{R}$ embeds into the corresponding resulting graph. Indeed, let $\{ x_n : n \in \omega\}$ be an enumeration of the vertices of $\mathcal{R}$. Set $t_0 = \emptyset$. In general, assume that $t_0 ,\ldots , t_n$ were constructed such that $|t_i| = i $ for every $i$ and 

\begin{center}
$\forall i, j \leqslant n \ \ \left( \{ t_i , t_j \} \in E \leftrightarrow \{ x_i , x_j \} \in E^{\mathcal{R}} \right)$. 
\end{center}

Choose $t_{n+1} \in 2^{< \omega}$ with height $n+1$ and such that 

\begin{center}
$\forall k \leqslant n \ \ t_{n+1}(i) = 1 \leftrightarrow \{ x_k , x_{n+1} \} \in E^{\mathcal{R}}$.
\end{center}

Then after $\omega$ steps, we are left with $\{ t_n : n \in \omega \}$ isomorphic to $\mathcal{R}$. In fact, observe that this construction can be carried out inside any strong subtree $\m{S}$ of $\m{T}$. On the other hand, it follows that $\mathcal{R}$ is indivisible iff $(2^{< \omega} , E)$ is. But now, indivisibility of $(2^{< \omega} , E)$ is guaranteed by Milliken's theorem: Let $A$ denote any $1$-point subset of $2^{< \omega}$. Then $[A]_{\mathrm{Em}}$ is simply $2^{< \omega}$ itself. So given $k \in \omega$ strictly positive and a coloring $\chi : \funct{2^{< \omega}}{k}$, one can find a $\chi$-monochromatic strong subtree $\m{S}$ of $2^{< \omega}$. The subgraph of $(2^{< \omega}, E)$ supported by $S$ being isomorphic to $(2^{< \omega}, E)$ itself, $S$ provides the required $\chi$-monochromatic copy of $(2^{< \omega}, E)$. \end{proof}

The following case to consider is $\Ur _3$, which turns out to be another particular case. As mentioned already in chapter 1, $\Ur _3$ can be encoded by the countable ultrahomogeneous edge-labelled graph with edges in $\{ 1, 3\}$ and forbidding the complete triangle with labels $1, 1, 3$. The distance between two points connected by an edge is the label of the edge while the distance between two points which are not connected is $2$. This fact allows to show: 

\begin{thm}[Delhomm\'e-Laflamme-Pouzet-Sauer \cite{DLPS}]

\label{thm:U_3 indivisible}
\index{Delhomm\'e-Laflamme-Pouzet-Sauer!theorem on $\Ur _3$}
\index{indivisibility!for $\Ur _3$}

$\Ur _3$ is indivisible. 

\end{thm}

The proof of this theorem can be deduced from the proof of the indivisibility of the $\m{K}_n$-free ultrahomogeneous graph by El-Zahar and Sauer in \cite{EZS1}. We do not provide the details here but mention few facts which will be useful for us later in subsection \ref{subsection:Indivisibility of U_S}. The presentation we adopt follows \cite{DLPS}. Fix a relational signature $L$ and consider an $L$-structure $\m{H}$. A nonempty subset $O$ of $H$ is an \emph{orbit}\index{orbit} if it is an orbit for the action of the automorphism group $\mathrm{Aut}(\m{H})$ on $\m{H}$ which pointwise fixes a finite subset of $H$. Now, given two $L$-structures $\m{R}$ and $\m{S}$, write $\m{R} \prec \m{S}$ when there is a partition of $R$ into finitely many parts $R_0 , \ldots , R_n$ such that for every $i \leqslant n$, $\m{R}_i$ embeds into $\m{S}$.  The following theorem follows from results in \cite{EZS2} and \cite{Sa1}. For the definition of free amalgamation see chapter 2 of the present article, subsection on Ne\v{s}et\v{r}il's theorem. 

\begin{thm}[El-Zahar - Sauer \cite{EZS2}, Sauer \cite{Sa1}]

\label{thm:Sauer ultrahomogeneous}
\index{El-Zahar - Sauer theorem}
\index{Sauer theorem}

Let $L$ be a finite binary signature and $\m{H}$ a countable ultrahomogeneous $L$-structure whose age has free amalgamation. Then $\m{H}$ is indivisible iff any two orbits of $\m{H}$ are related under $\prec$. 

\end{thm}

It follows that to prove that $\Ur _3$ is indivisible, it suffices to show that the countable ultrahomogeneous edge-labelled graph with edges in $\{ 1, 3\}$ and forbidding the complete triangle with labels $1, 1, 3$ satisfies those conditions, which in the present case is easy to check. We will see later that this method is actually useful in many cases. However, it does not allow to solve all the indivisibility problems that we are interested in. In particular, the indivisibility problem for $\Ur _m$ when $m\geqslant4$ is still, at that stage, left open. The purpose of the following section is to fill that gap. 

\subsection{The $\Ur _m$'s are indivisible.}
\index{indivisibility!for $\Ur _m$}

\label{subsection:U_m are indiv}

In the present section, we show that: 

\begin{thm}[NVT-Sauer \cite{NVTS}] 

\label{thm:U_m indiv}

Let $m \in \omega$, $m \geqslant 1$. Then $\Ur _m$ is indivisible. 

\end{thm}

The basic methods used in the proof have been developed in the sequence of papers \cite{EZS1}, \cite{EZS1'}, \cite{S3}, \cite{EZS2}, \cite{Sa1} dealing with partition results of countable ultrahomogeneous structures with free amalgamation. However, because the spaces $\Ur _m$ do not enter the framework provided by free amalgamation, substantial modifications were needed to prove Theorem \ref{thm:U_m indiv}. 

The proof is organized as follows. In section \ref{subsubsection:orbits} and \ref{subsubsection:largeness}, the essential ingredients, the main technical results (Lemma \ref{thm:orbit} and Lemma \ref{lem:psi}) as well as the general outline of the proof of Theorem \ref{thm:U_m indiv} are presented. Finally, the proof of Lemma \ref{thm:orbit} is presented in \ref{section:indiv}-\ref{subsection:red}.

\

\subsubsection{Kat\v{e}tov maps and orbits.}

\label{subsubsection:orbits}

We start with a reminder about Kat\v{e}tov maps\index{Kat\v{e}tov map}. Those objects were already defined in Chapter 1, section \ref{section:fundamentalsFraisse} but because of their omnipresence in the following pages, a bit of repetition will not harm. Recall that given a metric space $\m{X} = (X, d^{\m{X}})$, a map $f~:~\funct{X}{]0,+\infty[}$ is \emph{Kat\v{e}tov over $\m{X}$} when \[\forall x, y \in X, \ \ |f(x) - f(y)| \leqslant d^{\m{X}} (x,y) \leqslant f(x) + f(y).\] 

Equivalently, one can extend the metric $d^{\m{X}}$ to $X \cup \{ f \}$ by defining, for every $x, y$ in X, $ \widehat{d^{\m{X}}} (x, f) = f(x)$ and $\widehat{d^{\m{X}}} (x, y) = d^{\m{X}} (x, y)$. The corresponding metric space is then written $\m{X} \cup \{ f
\}$. 

The set of all Kat\v{e}tov maps over $\m{X}$ is written $E(\m{X})$. For a metric subspace $\m{X}$ of $\m{Y}$, a Kat\v{e}tov map $f \in E(\m{X})$ and a point $y \in \m{Y}$, then $y$ \emph{realizes $f$ over $\m{X}$} if \[ \forall x \in \m{X} \ \ d^{\m{Y}}(y,x) = f(x). \] 

The set of all $y \in \m{Y}$ realizing $f$ over $\m{X}$ is then written $O(f,\m{Y})$\index{$O(f,\m{Y})$, $O(f)$} and is called the \emph{orbit of $f$ in $\m{Y}$}\index{orbit}. When $\m{Y}$ is clear from the context, the set $O(f,\m{Y})$ is simply written $O(f)$. Again, the concepts of Kat\v{e}tov map and orbit are relevant because of the following standard reformulation of the notion of ultrahomogeneity, which will be used extensively in the sequel: 

\begin{lemma}
\label{prop:extension1} Let $\m{X}$ be a countable metric space. Then $\m{X}$ is ultrahomogeneous
iff for every finite subspace  $\m{F} \subset \m{X}$ and every Kat\v{e}tov map $f$ over $\m{F}$, if
$\m{F} \cup \{ f \}$ embeds into $\m{X}$, then $O(f, \m{X}) \neq \emptyset$.
\end{lemma}

\subsubsection{A notion of largeness.} 

\label{subsubsection:largeness}

In this section, $p$ is a fixed strictly positive integer. 

\begin{defn}

\label{defn:P}

The set $\mathbb{P}$\index{$\mathbb{P}$} is the set of all ordered pairs of the form $s = (f_s,\m{C}_s)$ where 
\begin{enumerate}
\item $\m{C}_s \in \binom{\Ur _p}{\Ur _p}$.
\item $f_s$ is a map with finite domain $\dom f_s \subset \m{C}_s$ and with values in $\{ 1,\ldots , p\}$.
\item $f_s \in E(\dom f_s)$, ie $f_s$ is Kat\v{e}tov on its domain.
\end{enumerate}

The set $\mathbb{P}$ is partially ordered by the relation $\leqslant$\index{$\leqslant$, $\leqslant_k$} defined by \[ \forall s, t \in \mathbb{P} \ \ t \leqslant s \leftrightarrow \left(\dom f_s \subset \dom f_t \subset \m{C}_t \subset \m{C}_s \ \ \mathrm{and}  \ \ \restrict{f_t}{\dom f_s} = f_s\right). \] 

Finally, if $k \in \omega$, then $t \leqslant _k s $ stands for \[ t \leqslant s \ \ \mathrm{and} \ \ \min f_t = \left \{ \begin{array}{cl} 
 \min f_s - k & \textrm{if $\min f_s > k$,} \\
 1 & \textrm{otherwise.}
 \end{array} \right.\]

\end{defn}

Observe that if $s \in \mathbb{P}$, then the ultrahomogeneity of $\Ur _p$ ensures that the set $O(f_s , \m{C} _s)$ is not empty and isometric to $\Ur _n$ where $n = \min (2\min f_s , p)$ (indeed, $O(f_s , \m{C} _s)$ is countable ultrahomogeneous with distances in $\{1,\ldots ,n \}$ and embeds every countable metric space with distances in $\{1,\ldots ,n \}$). Observe also that there is always a $t \in \mathbb{P}$ such that $t \leqslant _1 s$. Observe finally that unlike the relations $\leqslant$ and $\leqslant _0$, the relation $\leqslant _k$ is not transitive when $k>0$.

\begin{defn}

\label{defn:psi}

Let $s \in \mathbb{P}$ and $\Gamma \subset \Ur _p$. The notion of \emph{largeness of $\Gamma$ relative to $s$}\index{largeness} is defined recursively as follows:

If $\min f_s = 1$, then $\Gamma$ is large relative to $s$ iff \[ \forall t \leqslant _0 s  \ \left(O(f_t , \m{C}_t) \cap  \Gamma \ \textrm{is infinite}\right).\]

If $\min f_s > 1$, then $\Gamma$ is large relative to $s$ iff \[ \forall t \leqslant _0 s  \ \ \exists u \leqslant _1 t \ \ \left(\textrm{$\Gamma$ is large relative to $u$}\right) .\]

\end{defn}

The idea behind the definition of largeness is that if $\Gamma$ is large relative to $s$, then inside $\m{C}_s$ the set $\Gamma$ should represent a substantial part of the orbit of $f_s$. This intuition is made precise by the following Lemma \ref{thm:orbit}: 

\begin{lemma}

\label{thm:orbit}

Let $s \in \mathbb{P}$. Assume that $\Gamma$ is large relative to $s$. Then there exists an isometric copy $\m{C}$ of $\Ur _p$ inside $\m{C}_s$ such that:
\begin{enumerate} 
\item $\dom f_s \subset \m{C}$.
\item $O(f_s,\m{C}) \subset \Gamma$.  
\end{enumerate}  
\end{lemma}

In words, Lemma \ref{thm:orbit} means that by thinning up $\m{C}_s$, it is possible to ensure that the whole orbit of $f_s$ is included in $\Gamma$. The requirement $\dom f_s \subset \m{C}$ guarantees that the orbit of $f_s$ in the new space has the same metric structure as the orbit of $f_s$ in the original space. The proof of Lemma \ref{thm:orbit} represents the core of the proof of Theorem \ref{thm:U_m indiv} and is detailed in section \ref{section:indiv}. The second crucial fact about $\mathbb{P}$ and largeness lies in:  

\begin{lemma}

\label{lem:psi}

Let $s \in \mathbb{P}$ be such that $\Gamma$ is not large relative to $s$. Then there is $t \leqslant _0 s$ such that $\Ur _p \smallsetminus \Gamma$ is large relative to $t$.  
\end{lemma}

\begin{proof}

We proceed by induction on $\min f_s$. If $\min f_s = 1$, then there is $t \leqslant _0 s $ such that \[ O(f_t , \m{C}_t) \cap \Gamma \ \ \textrm{is finite}.\] 

It is then clear that $\Ur _p \smallsetminus \Gamma$ is large relative to $t$. On the other hand, if $\min f_s > 1$, then there is $t \leqslant _0 s$ such that \[ \forall w \leqslant _1 t \ \ \textrm{$\Gamma$ is not large relative to $w$}.\] 

We claim that $\Ur _p \smallsetminus \Gamma$ is large relative to $t$: let $u \leqslant _0 t$. We want to find $v \leqslant _1 u $ such that $\Ur _p \smallsetminus \Gamma$ is large relative to $v$. Let $w \leqslant _1 u$. Then $w \leqslant _1 t$ and it follows that $\Gamma$ is not large relative to $w$. By induction hypothesis, since $\min f_w < \min f_u = \min f_t$ there is $v \leqslant _0 w $ such that $\Ur _p \smallsetminus \Gamma$ is large relative to $v$. Additionally $v \leqslant _1 u$. Thus $v$ is as required. \end{proof}

When combined, Lemma \ref{thm:orbit} and Lemma \ref{lem:psi} lead to Theorem \ref{thm:U_m indiv} as follows: Take $p=m$ and consider a finite partition $\gamma$ of $\Ur_m$. Without loss of generality, $\gamma$ has only two parts, namely $\Pi$ and $\Omega$. Fix $t \in \mathbb{P}$ such that $\min f _t = m$. According to Lemma \ref{lem:psi}, either $\Pi$ is large relative to $t$ or there is $u \leqslant _0 s$ such that $\Omega$ is large relative to $u$. In any case, there are $s \in \{t, u\}$ and $\Gamma \in \{\Pi,\Omega\}$ such that $\min f_s = m$ and $\Gamma$ is large relative to $s$. Applying Lemma \ref{thm:orbit} to $s$, we obtain a copy $\m{C}$ of $\Ur _m$ inside $\m{C}_s$ such that $\dom f_s \subset \m{C}$ and $O(f_s , \m{C}) \subset \Gamma$. Observe that $O(f_s , \m{C})$ is isometric to $\Ur _m$. $\qed$ 

\

The remaining part of the proof is therefore devoted to a proof of Lemma \ref{thm:orbit}. 

\

\subsubsection{Proof of Lemma \ref{thm:orbit}, general strategy.}

\label{section:indiv}

From now on, the integer $p > 0$ is fixed together with $\Gamma \subset \Ur _p$. We proceed by induction and prove that for every strictly positive $m \in \omega$ with $m \leqslant p$ the following statement $\mathcal{J}_m$ holds: 

\
$\mathcal{J} _m$\index{$\mathcal{J} _m$} : ''For every $s \in \mathbb{P}$ such that $\min f_s = m$, if $\Gamma$ is large relative to $s$, then there exists an isometric copy $\m{C}$ of $\Ur _p$ inside $\m{C}_s$ such that:
\begin{enumerate} 
\item $\dom f_s \subset \m{C}$.
\item $O(f_s , \m{C}) \subset \Gamma$.''  
\end{enumerate} 

The proof is organized as follows. In subsection \ref{subsection:reformulation}, we show that the statement $\mathcal{J}_m$ is equivalent to a stronger statement denoted $\mathcal{H}_m$. This is achieved thanks to a technical lemma (Lemma \ref{lem:red}) about the structure of the orbits in $\Ur _p$ and whose proof is postponed to subsection \ref{subsection:red}. In subsection \ref{subsection:J_1}, we initiate the proof by induction and show that the statement $\mathcal{J}_1$ holds. We then show that if $\mathcal{H}_j$ holds for every $j<m$, then $\mathcal{J}_m$ holds. The general strategy of the induction step is presented in subsection \ref{subsection:induction}, while \ref{subsection:sequences} provides the details for the most technical aspects.  

\

\subsubsection{Reformulation of $\mathcal{J}_m$.}

\label{subsection:reformulation}

As mentioned previously, we start by reformulating the statement $\mathcal{J}_m$ under a form which will be useful when performing the induction step. Consider the following statement, denoted $\mathcal{H}_m$: 

$\mathcal{H} _m$\index{$\mathcal{H} _m$} : ''For every $s \in \mathbb{P}$ and every $F \subset \dom f_s$ such that $\min \restrict{f_s}{F} = \min f_s = m$, if $\Gamma$ is large relative to $s$, then there exists an isometric copy $\m{C}$ of $\Ur _p$ inside $\m{C}_s$ such that:
\begin{enumerate} 
\item $\dom f_s \cap \m{C} = F$.
\item $O(\restrict{f_s}{F} , \m{C}) \subset \Gamma$.''  
\end{enumerate}

The statement $\mathcal{J}_m$ is clearly implied by $\mathcal{H}_m$: Simply take $F = \dom f_s$. The purpose of the following lemma is to show that the converse is also true. 

\begin{lemma}

\label{lem:JmHm}

The statement $\mathcal{J}_m$ implies the statement $\mathcal{H}_m$.  

\end{lemma}

\begin{proof}
Our main tool here is the following technical result, whose proof is postponed to section \ref{subsection:red}. 

\begin{lemma}

\label{lem:red}

Let $G_0 \subset G$ be finite subsets of $\Ur_p$, $\mathcal{G}$ a family of Kat\v{e}tov maps with domain $G$ and such that for every $g, g' \in \mathcal{G}$: \[ \max(\restrict{|g - g'|}{G_0}) = \max | g - g'|, \] \[ \min(\restrict{(g+g')}{G_0}) = \min(g + g').\] 

Then there exists an isometric copy $\m{C}$ of $\Ur _p$ inside $\Ur _p$ such that:
\begin{enumerate} 
\item $G \cap \m{C} = G_0$.
\item $\forall g \in \mathcal{G} \ \ O(\restrict{g}{G_0} , \m{C}) \subset O(g, \Ur _p).$  
\end{enumerate} 

\end{lemma} 

Assuming Lemma \ref{lem:red}, here is how $\mathcal{J}_m$ implies $\mathcal{H}_m$: Let $s$ and $F$ be as in the hypothesis of $\mathcal{H}_m$. Apply $\mathcal{J}_m$ to $s$ to get an isometric copy $\mc{C}$ of $\Ur _p$ inside $\m{C}_s$ such that $\dom f_s \subset \mc{C}$ and $O(f_s , \mc{C}) \subset \Gamma$. Apply then Lemma \ref{lem:red} inside $\mc{C}$ to $F \subset \dom f_s$ and the family $\{ f_s\}$ to get an isometric copy $\m{C}$ of $\Ur _p$ inside $\mc{C}$ such that $\dom f_s \cap \m{C} = F$ and $O(\restrict{f_s}{F} , \m{C}) \subset O(f_s, \mc{C})$. Then $\m{C}$ is as required. \end{proof}

\subsubsection{Proof of $\mathcal{J}_1$.} 

\label{subsection:J_1}

Consider an enumeration $\{ x_n : n \in \omega\}$ of $\m{C}_s$ admitting $\dom f_s$ as an initial segment. Assume that the points $\varphi (x_0),\ldots , \varphi(x_n)$ are constructed so that: \begin{itemize}
	\item The map $\varphi$ is an isometry.
	\item $\restrict{\varphi}{\dom f_s} = id_{\dom f_s}$.
	\item $\varphi (x_k) \in \Gamma$ whenever $\varphi(x_k)$ realizes $f_s$ over $\dom f_s$. 
\end{itemize} 

We want to construct $\varphi (x _{n+1})$. Consider $h$ defined on $\{ \varphi (x_k) : k \leqslant n\} $ by: 
\[ \forall k \leqslant n \ \ h(\varphi(x_k)) = d^{\m{C}_s} (x_k , x_{n+1}).\] 

Observe that the metric subspace of $\m{C}_s$ given by $\{x_k : k \leqslant n+1\}$ witnesses that $h$ is Kat\v{e}tov. It follows that the set of all $y \in \m{C}_s \smallsetminus \dom f_s$ realizing $h$ over $\{ \varphi (x_k) : k \leqslant n\}$ is not empty and $\varphi (x_{n+1})$ can be chosen in that set. Additionally, observe that if $\restrict{h}{\dom f_s} = f_s$, then the fact that $\min f_s = 1$ and $\Gamma$ is large relative to $s$ then guarantees that $h$ can be realized by a point in $\Gamma$. We can therefore choose $\varphi (x_{n+1})$ to be one of those points. After $\omega$ steps, the subspace $\m{C}$ of $\m{C}_s$ supported by $\{ \varphi (x_n) : n \in \omega\}$ is as required. $\qed$

\

\subsubsection{Induction step.}

\label{subsection:induction}

Assume that the statements $\mathcal{J}_1\ldots\mathcal{J}_{m-1}$, and therefore the statements $\mathcal{H}_1\ldots\mathcal{H}_{m-1}$ hold. We are going to show that $\mathcal{J}_m$ holds. So let $s \in \mathbb{P}$ such that $\min f_s = m$ and $\Gamma$ is large relative to $s$. To make the notation easier, we assume that $s$ is of the form $(f,\Ur _p)$ and we write $F$ instead of $\dom f$. We need to produce an isometric copy $\m{C}$ of $\Ur _p$ inside $\Ur _p$ such that $F \subset \m{C}$ and $O(f , \m{C}) \subset \Gamma$. This is achieved inductively thanks to the following lemma. Recall that for metric subspaces $\m{X}$ and $\m{Y}$ of $\Ur _p$ and $\varepsilon > 0$, the sets $(\m{X}) _{\varepsilon}$ and $\binom{\m{Y}}{\Ur _p}$ are defined by: 
\[(\m{X}) _{\varepsilon} = \{ y \in \Ur _p : \exists x \in \m{X} \ \ d ^{\Ur _p} (y,x) \leqslant \varepsilon \},\]
\[ \binom{\m{Y}}{\Ur _p} = \{ \mc{U} \subset \m{Y} : \mc{U} \cong \Ur _p \}. \]

\begin{lemma}

\label{lem:ind1}

Let $\m{X}$ be a finite subspace of $\Ur _p$ and $\m{A} \in \binom{\Ur _p}{ \Ur _p}$ such that: 

\vspace{0.5em}
\hspace{1em}
(i) $F \subset \m{X} \subset \m{A}$.

\vspace{0.5em}
\hspace{1em}
(ii) $\left( \m{X} \right)_{m-1} \cap O(f,\m{A}) \subset \Gamma$.

\vspace{0.5em}
\hspace{1em}
(iii) $\forall g \in E(\m{X}) \ \ \restrict{g}{F} = \restrict{f}{F} \rightarrow \left(\textrm{$\Gamma$ is large relative to $(g, \m{A})$}\right)$.

\vspace{0.5em}

Then for every $h \in E(\m{X})$, there are $\m{B} \in \binom{\m{A}}{\Ur _p}$ and $x^* \in \m{B}$ realizing $h$ over $\m{X}$ such that: 

\vspace{0.5em}
\hspace{1em}
(i') $F \subset (\m{X} \cup \{x^*\}) \subset \m{B}$.

\vspace{0.5em}
\hspace{1em}
(ii') $ \left(\m{X} \cup \{ x^*\} \right)_{m-1} \cap O(f, \m{B}) \subset \Gamma$.

\vspace{0.5em}
\hspace{1em}
(iii') $\forall g \in E(\m{X} \cup \{ x^*\}) \ \ \restrict{g}{F} = \restrict{f}{F} \rightarrow \left(\textrm{$\Gamma$ is large relative to $(g, \m{B})$}\right)$.
\end{lemma}

\begin{claim}
Lemma \ref{lem:ind1} implies $\mathcal{J} _m$. 
\end{claim}

\begin{proof}
The required copy of $\m{C}$ can be constructed inductively. We start by fixing an enumeration $\{ x_n : n \in \omega \}$ of $\Ur _p$ such that $F = \{x_0 ,\ldots , x_k \}$ and by setting $\tilde{x}_i = x_i$ for every $i \leqslant k$. Next, we proceed as follows: Set $\m{A}_k = \Ur _p$. Then the subspace of $\Ur _p$ supported by $\{\tilde{x}_0 ,\ldots , \tilde{x}_k\}$ and the copy $\m{A}_k$ satisfy the requirements (i)-(iii) of Lemma \ref{lem:ind1}. Consider then $h_{k+1}$ defined on $\{\tilde{x}_0 ,\ldots , \tilde{x}_k\}$ by: \[ \forall i \leqslant k \ \ h_{k+1}(\tilde{x}_i) = d^{\Ur _p}(x_{k+1} , x_i).\] 

Then $h_{k+1}$ is Kat\v{e}tov over $\{\tilde{x}_0 ,\ldots , \tilde{x}_k\}$ and Lemma \ref{lem:ind1} can be applied to the subspace of $\Ur _p$ supported by $\{\tilde{x}_0 ,\ldots , \tilde{x}_k\}$, the copy $\m{A}_k$ and the Kat\v{e}tov map $h_{k+1}$. It produces $x^*$ and $\m{B}$, and we set $\tilde{x} _{k+1} = x^*$ and $\m{A}_{k+1} =\m{B}$. In general, assume that $\tilde{x}_0 ,\ldots , \tilde{x}_l$ and $\m{A}_k ,\ldots , \m{A}_l$ are constructed so that $\m{A}_l$ and the subspace of $\Ur _p$ supported by $\{\tilde{x}_0 ,\ldots , \tilde{x}_l\}$ satisfy the hypotheses of Lemma \ref{lem:ind1}. Consider $h_{l+1}$ defined on $\{\tilde{x}_0 ,\ldots , \tilde{x}_l\}$ by: \[ \forall i \leqslant l \ \ h_{l+1}(\tilde{x}_i) = d^{\Ur _p}(x_{l+1} , x_i).\] 

Then $h_{l+1}$ is Kat\v{e}tov over $\{\tilde{x}_0 ,\ldots , \tilde{x}_l\}$, Lemma \ref{lem:ind1} can be applied to produce $x^*$ and $\m{B}$, and we set $\tilde{x} _{l+1} = x^*$ and $\m{A}_{l+1} = \m{B}$. After $\omega$ steps, we are left with $\m{C} = \{ \tilde{x}_n : n \in \omega\}$ isometric to $\Ur _p$, as required. \end{proof} 

The remaining part of this section is consequently devoted to a proof of Lemma~\ref{lem:ind1} where $\m{X}$, $\m{A}$ and $h$ are fixed according to the requirements (i)-(iii) of Lemma \ref{lem:ind1}. 

\begin{claim}
If $x^*$ and $\m{B}$ satisfy (i') and (ii') of Lemma~\ref{lem:ind1}, then (iii') is also satisfied. 
\end{claim}

\begin{proof}
Let $g \in E( \m{X} \cup \{ x^*\})$ be such that $\restrict{g}{F} = \restrict{f}{F}$. We need to show that $\Gamma$ is large relative to $(g,\m{B})$. If $\min g \geqslant m$, then $(g,\m{B}) \leqslant _0 (f,\Ur _p)$. Since $\Gamma$ is large relative to $(f, \Ur _p)$, it follows that $\Gamma$ is also large relative to $(g,\m{B})$ and we are done. On the other hand, if $\min g \leqslant m-1$, then \[ O(g, \m{B}) \subset \left(\left(\m{X} \cup \{ x^*\} \right)_{m-1} \cap O(f, \m{B}) \right) \subset \Gamma.\] 

So $\Gamma$ is large relative to $(g,\m{B})$.\end{proof}

With this fact in mind, we define \[ K = \{ \phi \in E(\m{X} \cup \{ h \}) : \restrict{\phi}{F} = \restrict{f}{F} \ \ \mathrm{and} \ \ \phi(h) \leqslant m-1 \}.\] 

The reason for which $K$ is relevant here lies in:

\begin{claim}
Assume that $\m{B} \in \binom{\m{A}}{\Ur _p}$ and $x^* \in \m{B}$ are such that: 

\begin{enumerate}
\item $\m{X} \subset \m{B}$.
\item $x^*$ realizes $h$ over $\m{X}$.
\item For every $\phi \in K$, every point in $\m{B}$ realizing $\phi$ over $\m{X} \cup \{ x ^*\} \cong \m{X} \cup \{ h\}$ is in $\Gamma$. 
\item $x^* \in \Gamma$ if $\restrict{h}{F} = \restrict{f}{F}$ (that is if $x^* \in O(f, \m{B})$). 
\end{enumerate} 

Then $x^*$ and $\m{B}$ satisfy (i') and (ii') Lemma \ref{lem:ind1}.
\end{claim}

\begin{proof}
The requirement (i') is obviously satisfied so we concentrate on (ii'). Let $y \in \left( \m{X} \cup \{ x^*\} \right)_{m-1} \cap O(f,\m{B})$. We need to prove that $y \in \Gamma$. If $y \in \left(\m{X} \right)_{m-1}$, then $y$ is actually in $\left(\m{X}\right)_{m-1} \cap O(f, \m{A}) \subset \Gamma$ and we are done. Otherwise, $y \in \left(\{x^* \}\right)_{m-1}$. If $y=x^*$, there is nothing to do: Since $y$ is in $O(f, \m{B})$, so is $x^*$. Thus, by (iv), $x^* \in \Gamma$, that is $ y \in \Gamma$. Otherwise, let $\phi$ be the Kat\v{e}tov map realized by $y$ over $\m{X} \cup \{ x ^*\} \cong \m{X} \cup \{ h\}$. According to (iii), it suffices to show that $\phi \in K$. This is what we do now. First, the metric space $\m{X}\cup \{ x^*, y\}$ witnesses that $\phi$ is Kat\v{e}tov over $\m{X}\cup \{ h\}$. Next, $y \in O(f,\m{B})$ hence $\restrict{\phi}{F} = \restrict{f}{F}$. Finally, $\phi(h) = d^{\Ur _p} (x^*,y) \leqslant m-1$ since $y \in \left(\{x^* \}\right)_{m-1}$. \end{proof}

The strategy to construct $\m{B}$ and $x^*$ is the following one. Let $\{ \phi _{\alpha} : \alpha < \left|K\right|\}$ be an enumeration of $K$. We first construct a sequence of points $(x_{\alpha})_{\alpha < \left|K\right|}$ and a decreasing sequence $(\m{D}_{\alpha}) _{\alpha < \left|K\right|}$ of copies of $\Ur _p$ so that $x _{\alpha} \in \m{D}_{\alpha}$ and for every $\beta \leqslant \alpha < \left|K\right|$:

\begin{enumerate}
\item $\m{X} \subset \m{D} _{\alpha}$. 
\item $x _{\alpha}$ realizes $h$ over $\m{X}$. 
\item Every point in $\m{D}_{\alpha}$ realizing $\phi _{\beta}$ over $\m{X} \cup \{ x_{\alpha}\} \cong \m{X} \cup \{ h\}$ is in $\Gamma$. 
\end{enumerate}

The details of this construction are provided in section \ref{subsection:sequences}. Once this is done, call $x'=x_{\left|K\right|-1}$, $\m{B}' = \m{D}_{\left|K\right|-1}$. The point $x'$ and the copy $\m{B}'$ are almost as required except that $x'$ may not be in $\Gamma$. If $\restrict{h}{F} \neq \restrict{f}{F}$, this is not a problem and setting $x^* = x'$ and $\m{B} = \m{B}'$ works. On the other hand, if $\restrict{h}{F} = \restrict{f}{F}$, then some extra work is required and we proceed as follows. 

Pick $x^* \in \m{B}'$ realizing $h$ over $\m{X}$ and such that $d^{\Ur _p}(x^*, x')=1$. We will be done if we construct $\m{B} \in \binom{\m{B}'}{\Ur _p}$ so that $\left(\m{X}  \cup \{ x^* , x' \}\right) \cap \m{B} = \m{X}  \cup \{ x^* \}$ and for every $\phi \in K$, every point in $\m{B}$ realizing $\phi $ over $\m{X} ^* \cup \{ x^* \}$ realizes $\phi $ over $\m{X} ^* \cup \{ x' \}$. Here is how this is achieved thanks to Lemma \ref{lem:red}. For $\phi \in K$, define the map $\hat{\phi}$ on $\m{X}  \cup \{ x^* , x' \}$ by
\begin{displaymath}
\left \{ \begin{array}{l}
 \restrict{\hat{\phi}}{\m{X}} = \restrict{\phi}{\m{X}}, \\
 \hat{\phi}(x^*) = \hat{\phi}(x') = \phi(h).
 \end{array} \right.
\end{displaymath} 

Using the fact that $\phi$ is Kat\v{e}tov over $\m{X}\cup \{h\}$ and $\m{X} \cup \{ x^*\} \cong \m{X} \cup \{ x'\} \cong \m{X} \cup \{ h\}$, it is easy to check that $\hat{\phi}$ is Kat\v{e}tov over $\m{X}  \cup \{ x^* , x' \}$ and that for every $\phi, \phi' \in K$: \[ \max(\restrict{|\hat{\phi} - \hat{\phi}'|}{\m{X} \cup \{ x^*\}}) = \max | \hat{\phi} - \hat{\phi}'|, \] \[ \min(\restrict{(\hat{\phi}+\hat{\phi}')}{\m{X} \cup \{ x^*\}}) = \min(\hat{\phi} + \hat{\phi}').\] 

Working inside $\m{B}'$, we can therefore apply Lemma \ref{lem:red} to $\m{X} \cup \{ x^*\} \subset \m{X}  \cup \{ x^* , x' \}$ and the family $(\hat{\phi})_{\phi \in K}$ to obtain $\m{B}$ as required. $\qed$

\

\subsubsection{Construction of the sequences $(x_{\alpha})_{\alpha < \left|K\right|}$ and $(\m{D}_{\alpha}) _{\alpha < \left|K\right|}$.}

\label{subsection:sequences}

The construction of the sequences $(x_{\alpha})_{\alpha < \left|K\right|}$ and $(\m{D}_{\alpha}) _{\alpha < \left|K\right|}$ is carried out thanks to a repeated application of the following lemma:  

\begin{lemma}

\label{lem:ind2}

Let $\mathcal{F} \subset K$ and $\m{D} \in \binom{\m{A}}{\Ur _p}$ be such that $\m{X} \subset \m{D}$. Assume that $u \in \m{D}$ realizes $h$ over $\m{X}$ and is such that for every $\phi \in \mathcal{F}$, every point in $\m{D}$ realizing $\phi$ over $\m{X} \cup \{ u \} \cong \m{X} \cup \{ h \}$ is in $\Gamma$. Let $s \in K \smallsetminus \mathcal{F}$ be such that \[ \forall \phi \in K \ \ \phi (h) > s(h) \rightarrow \phi \in \mathcal{F} \ \ \textrm{and} \ \ \phi (h) < s(h) \rightarrow \phi \notin \mathcal{F}. \ \ \ (*)\] 

Then there are $\m{E} \in \binom{\m{D}}{\Ur _p}$ and $v \in \m{E}$ realizing $h$ over $\m{X} $ such that $\m{X} \subset \m{E}$ and for every $\phi \in \mathcal{F} \cup \{ s \}$, every point in $\m{E}$ realizing $\phi$ over $\m{X} \cup \{ v \} \cong \m{X} \cup \{ h \}$ is in $\Gamma$.
\end{lemma}

Once Lemma \ref{lem:ind2} is proven, here is how the sequences $(x_{\alpha})_{\alpha < \left|K\right|}$ and $(\m{D}_{\alpha}) _{\alpha < \left|K\right|}$ are constructed: Choose the enumeration $\{ \phi _{\alpha} : \alpha < \left|K\right|\}$ of $K$ so that the sequence $(\phi _{\alpha} (h))_{\alpha < \left|K\right|}$ is nondecreasing. Apply Lemma \ref{lem:ind2} to $\mathcal{F} = \emptyset$, $\m{D} = \m{A}$ and $s = \phi_0$ to produce $x_0$ and $\m{D}_0$. In general, apply Lemma \ref{lem:ind2} to $\mathcal{F} = \{ \phi_0\ldots \phi_{\alpha} \}$, $\m{D} = \m{D}_{\alpha}$ and $s=\phi_{\alpha+1}$ to produce $x_{\alpha+1}$ and $\m{D}_{\alpha+1}$. After $|K|$ steps, the sequences $(x_{\alpha})_{\alpha < \left|K\right|}$ and $(\m{D}_{\alpha}) _{\alpha < \left|K\right|}$ are as required.

\begin{proof}[Proof of Lemma \ref{lem:ind2}]
We start with the case where $s(h) \geqslant \min \restrict{s}{\m{X}}$. The map $s$ being in $K$, $s(h) \leqslant m-1$ and so $\min \restrict{s}{\m{X}} \leqslant m-1$. Then, \[ O(\restrict{s}{\m{X}}, \m{D}) \subset \left( \left(\m{X}\right)_{m-1} \cap O(f, \m{D})\right).\] 

But from the requirement (ii) of Lemma \ref{lem:ind1}, \[ \left( \left(\m{X}\right)_{m-1} \cap O(f, \m{D})\right) \subset \Gamma.\] 

Observe now that every point in $\m{D}$ realizing $s$ over $\m{X} \cup \{ u \}$ is in $O(\restrict{s}{\m{X}}, \m{D})$. Thus, according to the previous inclusions, any such point is also in $\Gamma$. So in fact, there is nothing to do: $v = u$ and $\m{E} = \m{D}$ works. 

From now on, we consequently suppose that $s(h) < \min \restrict{s}{\m{X}}$. Let $s_1$ be defined on $\m{X} \cup \{ u \}$ by 
\begin{displaymath}
s_1(x) = \left \{ \begin{array}{cl}
 s(x) & \textrm{if $x \in \m{X}$,} \\
 s(h) + 1 & \textrm{if $x = u$.}
 \end{array} \right.
\end{displaymath}  

\begin{claim}
The map $s_1$ is Kat\v{e}tov. 
\end{claim}

\begin{proof}

The map $s$ is Kat\v{e}tov over $\m{X}$. Hence, it is enough to prove that for every $x \in \m{X}$ \[ \left| s_1(u) - s_1(x) \right| \leqslant d^{\Ur _p}(x,u) \leqslant s_1(u) + s_1(x). \] 

That is \[ \left| s(h) + 1 - s(x) \right| \leqslant h(x) \leqslant s(h) + 1 + s(x).\] 

Because $s$ is Kat\v{e}tov over $\m{X} \cup \{ h \}$, it is enough to prove that \[ s(h) + 1 - s(x) \leqslant h(x). \]

But this holds since $s(h) < \min \restrict{s}{\m{X}}$. \end{proof}

\begin{claim}
$\Gamma$ is large relative to $(s_1, \m{D})$. 
\end{claim}

\begin{proof}
If $s(h) = m-1$, then $\min s_1 = m = \min f$ and so $(s_1 , \m{D}) \leqslant _ 0 (f , \Ur _p)$. Since $\Gamma$ is large relative to $(f ,\Ur _p)$, it is also large relative to $(s_1 , \m{D})$ and we are done. On the other hand, if $s(h) < m-1$, then $s_1 \in K$ and it follows from the hypothesis $(*)$ on $\mathcal{F}$ that $s_1 \in \mathcal{F}$. In particular, every point in $\m{D}$ realizing $s_1 $ over $\m{X} \cup \{ u \}$ is in $\Gamma$, and it follows that $\Gamma$ is large relative to $(s_1 , \m{D})$. \end{proof}

Consequently, there is $(s_2 , \m{D} _2) \leqslant _1 (s_1 , \m{D})$ such that $\Gamma$ is large relative to $(s_2 , \m{D} _2)$. We are now going to construct $v$ and a Kat\v{e}tov extension $s_3$ of $s_2$ such that $v$ realizes $h$ over $\m{X}$, $s_3 (v) = s(h)$ and $(s_3,\m{D}_2) \leqslant_0 (s_2 , \m{D}_2)$. This last requirement will make sure that $\Gamma$ is large relative to $(s_3,\m{D}_2)$. We will then apply Lemma \ref{lem:red} to obtain the copy $\m{E}$ as required. 
Here is how we proceed formally: Fix $w \in O(s_2, \m{D}_2)$ and consider the map $h_1$ defined on $\m{X} \cup \{ u, w \}$ by
\begin{displaymath}
h_1(x) = \left \{ \begin{array}{cl}
 h(x) & \textrm{if $x \in \m{X}$.} \\
 1 & \textrm{if $x = u$.} \\
 s(h) & \textrm{if $x = w$.}
 \end{array} \right.
\end{displaymath}

\begin{claim}
The map $h _1$ is Kat\v{e}tov. 
\end{claim}

\begin{proof}
The metric space $\left(\m{X} \cup \{ h\}\right) \cup \{ s\}$ witnesses that $\restrict{h_1}{\m{X} \cup \{ w \}}$ is Kat\v{e}tov. Next, $\restrict{h_1}{\m{X} \cup \{ u \}}$ is also Kat\v{e}tov: Let $x \in \m{X}$. Then \[ \left|h_1(x) - h_1(u) \right| = h(x) - 1 \leqslant h(x) = d^{\Ur _p} (x,u) \leqslant h(x) + 1 = h_1(x) + h_1(u).\] 

The only thing we still need to show is therefore \[ \left|h_1(u) - h_1(w) \right| \leqslant d^{\Ur _p} (u,w) \leqslant h_1(u) + h_1(w).\] 

But these inequalities hold as they are equivalent to \[ \left| 1 - s(h) \right| \leqslant s(h) + 1 \leqslant 1 + s(h).\qedhere\]  \end{proof}

Let $v \in \m{D} _3$ realizing $h_1$ over $\m{X} \cup \{ u, w\}$. As announced previously, define an extension $s _3$ of $s_2$ on $\dom s_2 \cup \{v\}$ by setting $s_3 (v) = s(h)$. 

\begin{claim}
The map $s_3$ is Kat\v{e}tov and $\Gamma$ is large relative to $(s_3, \m{D}_2)$. 
\end{claim}

\begin{proof}
The point $w$ realizes $s_3$ over $\dom s_2 \cup \{v\}$ and therefore witnesses that $s_3$ is Kat\v{e}tov. As for $\Gamma$, it is large relative to $(s_3, \m{D}_2)$ because it is large relative to $(s_2 , \m{D} _2)$ and $(s_3,\m{D}_2) \leqslant_0 (s_2 , \m{D}_2)$. 
\end{proof}

Observe now that $\min s_3 = s(h) = \min \restrict{s_3}{\m{X} \cup \{ u, v\}} = \min s \leqslant m-1$. Thus, one can apply $\mathcal{H}_{\min s}$ inside $\m{D}_2$ to $s_3$ and $\m{X} \cup \{ u, v\}$ to obtain $\m{D}_3 \in \binom{\m{D}_2}{\Ur _p}$ such that $\dom s_3 \cap \m{D}_3 = \m{X} \cup \{ u, v\}$ and $O(\restrict{s_3}{\m{X} \cup \{ u, v\}}, \m{D}_3) \subset \Gamma$. At that point, both $u$ and $v$ realize $h$ over $\m{X}$ and if $\phi \in \mathcal{F}$, then every point in $\m{D} _3$ realizing $\phi$ over $\m{X} \cup \{ u\}$ is in $\Gamma$. Thus, we will be done if we can construct $\m{E} \in \binom{\m{D}_3}{\Ur _p}$ such that:

\begin{itemize}
\item $\left(\m{X} \cup \{u, v \}\right) \cap \m{E} = \m{X} \cup \{ v \}$.
\item For every $\phi \in \mathcal{F}$, every point in $\m{E}$ realizing $\phi$ over $\m{X} \cup \{ v \}$ realizes $\phi$ over $\m{X} \cup \{ u \}$. 
\item Every point in $\m{E}$ realizing $s$ over $\m{X} \cup \{ v \}$ realizes $s_3$ over $\m{X} \cup \{ u, v\}$. 
\end{itemize} 

Once again, this is achieved thanks to Lemma \ref{lem:red}: For $\phi \in \mathcal{F}$, define the map $\hat{\phi}$ on $\m{X}  \cup \{ u , v \}$ by: 
\begin{displaymath}
\left \{ \begin{array}{l}
 \restrict{\hat{\phi}}{\m{X}} = \restrict{\phi}{\m{X}}, \\
 \hat{\phi}(u) = \hat{\phi}(v) = \phi(h).
 \end{array} \right.
\end{displaymath} 

Using the fact that $\phi$ is Kat\v{e}tov over $\m{X}\cup \{h\}$ and $\m{X} \cup \{ u\} \cong \m{X} \cup \{ v\} \cong \m{X} \cup \{ h\}$, it is easy to check that $\hat{\phi}$ is Kat\v{e}tov over $\m{X}  \cup \{ u, v \}$. Let $\widehat{\mathcal{F}} = (\hat{\phi})_{\phi \in \mathcal{F}}$. Working inside $\m{D}_3$, we would like to apply Lemma \ref{lem:red} to $\m{X} \cup \{ v \} \subset \m{X}  \cup \{ u, v \}$ and the family $\{ s_3\} \cup \widehat{\mathcal{F}}$ to obtain $\m{E}$ as required. It is therefore enough to check: 

\begin{claim}
For every $g, g' \in \{ s_3\} \cup \widehat{\mathcal{F}}$: \[ \max(\restrict{|g - g'|}{\m{X} \cup \{ v \}}) = \max | g - g'|, \] \[ \min(\restrict{(g+g')}{\m{X} \cup \{ v \}}) = \min(g + g').\] 
\end{claim}

\begin{proof}
When $g, g' \in \widehat{\mathcal{F}}$, this is easily done. We therefore concentrate on the case where $g = \hat{\phi}$ for $\phi \in \mathcal{F}$ and $g' = s_3$. What we have to do is to show that: \[ |\hat{\phi}(u) - s_3(u)| \leqslant \max(\restrict{|\hat{\phi} - s_3|}{\m{X} \cup \{ v \}}) \ \ (1)\] \[ \hat{\phi}(u) + s_3(u) \geqslant \min(\restrict{(\hat{\phi}+s_3)}{\m{X} \cup \{ v \}}) \ \ (2)\] 

Recall first that $s_3(u) = s(h) + 1$ and that $s_3(v) = s(h)$. Remember also that according to the properties of $\mathcal{F}$, $s(h)\leqslant \phi(h)$. For $(1)$, if $s(h) < \phi(h)$, then we are done since 
\begin{align*} |\hat{\phi}(u) - s_3(u)| & = |\phi(h) - (s(h) + 1)| \\
& = \phi(h) - (s(h) + 1) \\
& \leqslant \phi(h) - s(h) \\
& = \phi(v) - s_3(v) \\
& \leqslant |\hat{\phi}(v) - s_3(v)|.\end{align*} 

On the other hand, if $\phi(h) = s(h)$, then $|\hat{\phi}(u) - s_3(u)| = 1$ but then this less or equal to $\max(\restrict{|\hat{\phi} - s_3|}{\m{X} \cup \{ v \}})$ as this latter quantity is equal to $\max |\phi - s|$, which is at least $1$ since $\phi \in \mathcal{F}$ and $s \notin \mathcal{F}$. Thus, the inequality $(1)$ holds. As for $(2)$, simply observe that \[ \hat{\phi}(u) + s_3 (u) \geqslant \hat{\phi}(v) + s_3 (v). \qedhere\] \end{proof}
This finishes the proof of Lemma \ref{lem:ind2}. \end{proof}

\subsubsection{Proof of Lemma \ref{lem:red}}

\label{subsection:red}

The purpose of this section is to provide a proof of Lemma \ref{lem:red} which was used extensively in the previous proofs. Let $G_0 \subset G$ be finite subsets of $\Ur_p$, $\mathcal{G}$ a family of Kat\v{e}tov maps with domain $G$ and such that for every $g, g' \in \mathcal{G}$: \[ \max(\restrict{|g - g'|}{G_0}) = \max | g - g'|, \] \[ \min(\restrict{(g+g')}{G_0}) = \min(g + g').\] 

We need to produce an isometric copy $\m{C}$ of $\Ur _p$ inside $\Ur _p$ such that:
\begin{enumerate} 
\item $G \cap \m{C} = G_0$.
\item $\forall g \in \mathcal{G} \ \ O(\restrict{g}{G_0} , \m{C}) \subset O(g, \Ur _p).$  
\end{enumerate} 

First, observe that it suffices to provide the proof assuming that $G$ is of the form $G_0\cup\{ z\}$. The general case is then handled by repeating the procedure. 
  
\begin{lemma}

\label{lem:red1}

Let $\m{X}$ be a finite subspace of $\bigcup \{ O(\restrict{g}{G_0}) : g \in \mathcal{G}\}$. Then there is an isometry $\varphi$ on $\Ur _p$ fixing $G_0 \cup (\m{X} \cap \bigcup \{ O(g) : g \in \mathcal{G}\})$ and such that: \[ \forall g \in \mathcal{G} \ \ \varphi '' \m{X} \cap O(\restrict{g}{G_0}) \subset O(g).\]

\end{lemma}

\begin{proof}
For $x \in \m{X}$, there is a unique element $g_x \in \mathcal{G}$ such that $x \in O(\restrict{g_x}{G_0})$. Let $k$ be the map defined on $G_0 \cup \m{X}$ by
\begin{displaymath}
k(x) = \left \{ \begin{array}{cl}
 d^{\Ur _p}(x,z) & \textrm{if $x \in G_0$,} \\
 g_x(z) & \textrm{if $x \in \m{X}$.}
 \end{array} \right.
\end{displaymath}

\begin{claim}
The map $k$ is Kat\v{e}tov. 
\end{claim}

\begin{proof}
The metric space $G_0\cup \{z\}$ witnesses that $k$ is Kat\v{e}tov over $G_0$. Hence, it suffices to check that for every $x \in \m{X}$ and $y \in G_0\cup \m{X}$, \[ |k(x) - k(y)| \leqslant d^{\Ur _p} (x,y) \leqslant k(x) + k(y).\]

Consider first the case $y \in G_0$. Then $d^{\Ur}(x,y) = g_x(y)$ and we need to check that \[ |g_x(z) - d^{\Ur _p}(y,z)| \leqslant g_x(y) \leqslant g_x(z) + d^{\Ur _p}(y,z).\] 

Or equivalently, \[ |g_x(z) - g_x(y)| \leqslant d^{\Ur _p}(y,z) \leqslant g_x(z) + g_x(y).\] 

But this is true since $g_x$ is Kat\v{e}tov over $G_0\cup \{z\}$. Consider now the case $y\in \m{X}$. Then $k(y) = g_y(z)$ and we need to check \[ |g_x(z) - g_y(z)| \leqslant d^{\Ur _p}(x,y) \leqslant g_x(z) + g_y(z).\] 

But since $\m{X}$ is a subspace of $\bigcup \{ O(\restrict{g}{G_0}) : g \in \mathcal{G}\}$, we have, for every $u \in G_0$, \[ |d^{\Ur _p}(x,u) - d^{\Ur _p}(u,y)| \leqslant d^{\Ur _p}(x,y) \leqslant d^{\Ur _p}(x,u) + d^{\Ur _p}(x,u).\]

Since $x \in O(\restrict{g_x}{G_0})$ and $y \in O(\restrict{g_y}{G_0})$, this is equivalent to \[ |g_x(u) - g_y(u)| \leqslant d^{\Ur _p}(x,y) \leqslant g_x(u) + g_y(u).\] 

Therefore, \[ \max (\restrict{|g_x - g_y|}{G_0}) \leqslant d^{\Ur _p}(x,y) \leqslant \min (\restrict{(g_x + g_y)}{G_0}) .\]

Now, by hypothesis on $\mathcal{G}$, this latter inequality remains valid if $G_0$ is replaced by $G_0\cup \{ z\}$. The required inequality follows. \end{proof}

By ultrahomogeneity of $\Ur _p$, we can consequently realize the map $k$ over $G_0\cup \m{X}$ by a point $z' \in \Ur _p$. The metric space $G_0\cup(\m{X}\cap\bigcup \{ O(g) : g \in \mathcal{G}\})\cup \{k\}$ being isometric to the subspace of $\Ur _p$ supported by $G_0\cup(\m{X}\cap\bigcup \{ O(g) : g \in \mathcal{G}\})\cup \{z\}$, so is the subspace of $\Ur _p$ supported by $G_0\cup(\m{X}\cap\bigcup \{ O(g) : g \in \mathcal{G}\})\cup \{z'\}$. By ultrahomogeneity again, we can therefore find a surjective isometry $\varphi$ of $\Ur _p$ fixing $G_0\cup(\m{X}\cap\bigcup \{ O(g) : g \in \mathcal{G}\})$ and such that $\varphi(z')=z$. Then $\varphi$ is as required: Let $g \in \mathcal{G}$ and $x \in O(\restrict{g}{G_0})$. Then: \[d^{\Ur _p}(\varphi(x),z) = d^{\Ur _p}(\varphi(x), \varphi(z')) = d^{\Ur _p}(x,z') = k(x) = g(z) \] 

That is, $\varphi(x) \in O(g)$. \end{proof}

\begin{lemma}

\label{lem:red2}

There is an isometric embedding $\psi$ of $G_0 \cup \bigcup \{ O(\restrict{g}{G_0}) : g \in \mathcal{G})\}$ into $G_0 \cup \bigcup \{ O(g) : g \in \mathcal{G})\}$ fixing $G_0$ such that: \[ \forall g \in \mathcal{G} \ \ \psi '' O(\restrict{g}{G_0}) \subset O(g).\] 

\end{lemma}

\begin{proof}
Let $\{ x_n : n \in \omega\}$ enumerate $\bigcup \{ O(\restrict{g}{G_0}) : g \in \mathcal{G})\}$. For $n \in \omega$, let $g_n$ be the only $g \in \mathcal{G}$ such that $x_n \in O(\restrict{g_n}{G_0})$. Apply Lemma \ref{lem:red1} inductively to construct a sequence $(\psi _n)_{n \in \omega}$ of surjective isometries of $\Ur _p$ such that for every $n\in \omega$, $\psi_n$ fixes $G_0 \cup \psi_{n-1}''\{x_k : k < n \}$ and $\psi_n(x_n) \in O(g_n)$. Then $\psi$ defined on $G_0 \cup \{ x_n : n \in \omega\}$ by $\restrict{\psi}{G_0} = id_{G_0}$ and $\psi(x_n) = \psi_n(x_n)$ is as required. \end{proof}

\begin{figure}[h]
\begin{center}
\hskip-10pt
\setlength{\unitlength}{1mm}
\begin{picture}(130,100)(0,0)
\includegraphics[width=132.00mm]{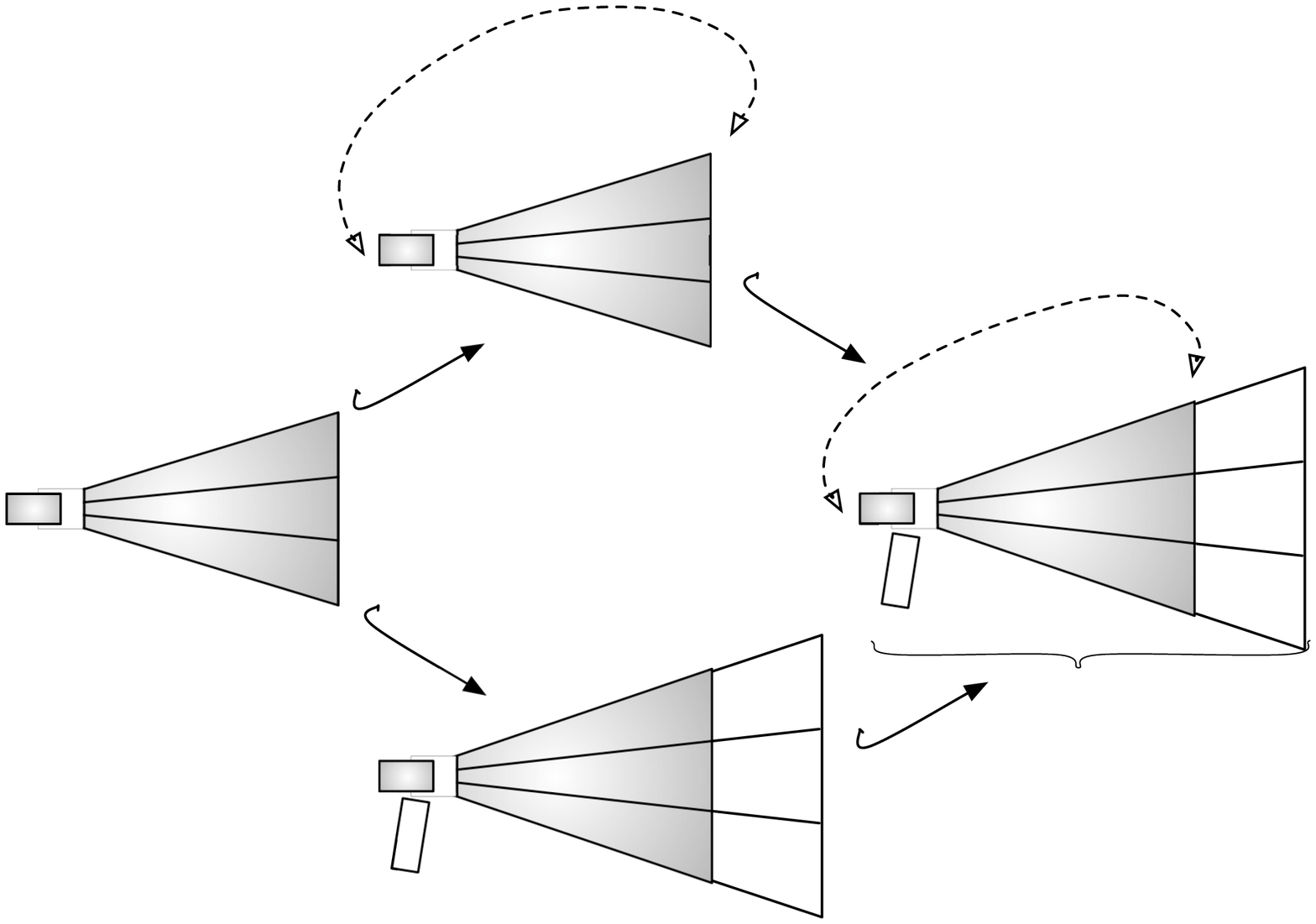}
\put(-75.26,93.04){$\Ur _p$}
\put(-128.91,44.55){$G_0$}
\put(-112.80,45.80){$O(\restrict{g_1}{G_0})$}
\put(-112.80,41){$O(\restrict{g_2}{G_0})$}
\put(-112.80,37){$O(\restrict{g_3}{G_0})$}
\put(-93.35,26.26){$j_0$}
\put(-93,19){$G_0$}
\put(-95.19,10){$G$}
\put(-59.50,22.41){$O(g_1)$}
\put(-59.50,14){$O(g_2)$}
\put(-59.50,7){$O(g_3)$}
\put(-54.64,57.79){$i_1$}
\put(-91.30,50.69){$i_0$}
\put(-24.49,64.92){$ i_1 '' \Ur _p$}
\put(-27.5,23){$\m{W}$}
\put(-40.80,19){$j_1$}
\put(-50.80,36){$j_1''G$}
\end{picture}
\end{center}
\caption{The space $\m{W}$}\label{fig:figW}
\end{figure}

We now turn to the proof of Lemma \ref{lem:red}. Let $\m{Y}$ and $\m{Z}$ be the metric subspaces of $\Ur _p$ supported by $G \cup \bigcup \{ O(g) : g \in \mathcal{G})\}$ and $G_0 \cup \bigcup \{ O(\restrict{g}{G_0}) : g \in \mathcal{G})\}$ respectively. Let $i_0 : \funct{\m{Z}}{\Ur _p}$ be the isometric embedding provided by the identity. By Lemma \ref{lem:red2}, the space $\m{Z}$ embeds isometrically into $\m{Y}$ via an isometry $j_0$ that fixes $G_0$. We can therefore consider the metric space $\m{W}$ obtained by gluing $\Ur _p$ and $\m{Y}$ via an identification of $\m{Z} \subset \Ur _p$ and $j_0''\m{Z} \subset \m{Y}$. The space $\m{W}$ is described in Figure \ref{fig:figW}. 

Formally, the space $\m{W}$ can be constructed thanks to a property of the countable metric spaces with distances in $\{ 1, \ldots , p\}$ known as \emph{strong amalgamation}: We can find a countable metric space $\m{W}$ with distances in $\{ 1, \ldots , p\}$ and isometric embeddings $i_1 : \funct{\Ur _p}{\m{W}}$ and $j_1 : \funct{\m{Y}}{\m{W}}$ such that: 
\begin{itemize}
	\item $i_1 \circ i_0 = j_1 \circ j_0$.
	\item $\m{W} = i_1''\Ur _p \cup j_1''\m{Y}$.
	\item $i_1''\Ur _p \cap j_1''\m{Y} = (i_1 \circ i_0)''\m{Z} = (j_1 \circ j_0)'' \m{Z}$. 
	\item For every $x \in \Ur _p$ and $y \in \m{Y}$: \begin{align*}d^{\m{W}}(i_1(x),j_1(y)) & = \min \{ d^{\m{W}}(i_1(x),i_1 \circ i_0 (z)) + d^{\m{W}}(j_1 \circ j_0 (z),j_1(y)) : z \in \m{Z}\} \\
& = \min \{ d^{\Ur _p}(x,i_0 (z)) + d^{\m{Y}}(j_0 (z),y) : z \in \m{Z}\} \\
& = \min \{ d^{\Ur _p}(x,z) + d^{\m{Y}}(j_0 (z),y) : z \in \m{Z}\}.
\end{align*}
\end{itemize}
\vskip-5pt

The crucial point here is that in $\m{W}$, every $x \in i_1 '' \Ur _p$ realizing some $\restrict{g}{G_0}$ over $i_1 '' G_0$ also realizes $g$ over $j_1 '' G$.  

Using $\m{W}$, we show how $\m{C}$ can be constructed inductively: Consider an enumeration $\{ x_n : n \in \omega\}$ of $i_1'' \Ur _p$ admitting $i_1'' G_0$ as an initial segment. Assume that the points $\varphi (x_0),\ldots , \varphi(x_n)$ are constructed so that: 
\begin{itemize}
	\item The map $\varphi$ is an isometry.
	\item $\dom \varphi \subset i_1'' \Ur _p$.
	\item $\ran \varphi \subset \Ur _p$.
	\item $\varphi(i_1(x)) = x$ whenever $x \in G_0$.
	\item $d^{\Ur _p} (\varphi (x_k) , z) = d^{\m{W}} (x_k , j_1(z))$ whenever $z \in G$ and $k \leqslant n$.
\end{itemize}

We want to construct $\varphi (x _{n+1})$. Consider $e$ defined on $\{ \varphi (x_k) : k \leqslant n\} \cup G $ by: 
\begin{displaymath}
\left \{ \begin{array}{l}
 \forall k \leqslant n \ \ e(\varphi(x_k)) = d^{\m{W}} (x_k , x_{n+1}), \\
 \forall z \in G \ \ e(z) = d^{\m{W}} (j_1(z) , x_{n+1}).
 \end{array} \right.
\end{displaymath}

Observe that the metric subspace of $\m{W}$ given by $\{x_k : k \leqslant n+1\} \cup j_1'' G $ witnesses that $e$ is Kat\v{e}tov. It follows that the set $E$ of all $y \in \Ur _p$ realizing $e$ over the set $\{ \varphi (x_k) : k \leqslant n\} \cup G $ is not empty and $\varphi (x_{n+1})$ can be chosen in $E$. $\qed$

\subsection{Indivisibility of ultrametric Urysohn spaces.}

\label{subsection:Indivisibility of Urysohn ultrametric spaces}
\index{indivisibility!for $\uUr$}

We saw in section \ref{section:Big Ramsey degrees} that the classes of ultrametric spaces $\U$\index{$\U$} were the only case where
we were able to compute the big Ramsey degree explicitly. However, Theorem \ref{thm:Big Ramsey degrees for U, S finite} and Theorem \ref{thm:no Big Ramsey degrees for U, S infinite} leave an open case: Nothing is said about the big Ramsey degree of the $1$-point ultrametric space when the set $S$ is infinite. In other words, Theorems \ref{thm:Big Ramsey degrees for U, S finite} and \ref{thm:no Big Ramsey degrees for U, S infinite} do not solve the indivisibility problem for $\uUr$\index{$\uUr$} when $S$ is infinite. The purpose of this subsection is to fix that flaw.

\begin{thm}

\label{thm:uUr divisible, S not dually well ordered}

Let $S \subset ]0, + \infty[$ countable. Assume that the reverse linear ordering $>$ on $\R$ does not induce a well-ordering on $S$. Then there is a map $\chi : \funct{\uUr}{\omega}$ whose restriction on any isometric copy $X$ of $\uUr$ inside $\uUr$ has range $\omega$. 
\end{thm} 

In particular, in this case, $\uUr$ is divisible. This result should be compared with the following one:

\begin{thm}

\label{thm:uUr indivisible, S dually well ordered}

Let $S\ \subset ]0, + \infty[$ be finite or countable. Assume that the reverse linear ordering $>$ on $\R$ induces a well-ordering on $S$. Then $\uUr$ is indivisible.
\end{thm}

Two remarks before entering the technical parts: First, Theorem \ref{thm:uUr divisible, S not dually well ordered} and Theorem \ref{thm:uUr indivisible, S dually well ordered} were first obtained completely independently of our work by Delhomm\'e, Laflamme, Pouzet and Sauer in \cite{DLPS}. The proofs presented here are ours but the reader should be aware of the fact that for Theorem \ref{thm:uUr indivisible, S dually well ordered}, though the ideas are essentially the same, the proof presented in \cite{DLPS} is considerably shorter. Second remark: It is easy to show that a necessary condition on a countable ultrahomogeneous ultrametric space $\m{X}$ to be indivisible is to be of the form $\uUr$ for some at most countable $S \subset ]0, + \infty[$. Indeed, otherwise, according to Proposition \ref{prop:ultraultra} (Chapter 1), $\m{X}$ is the set of finitely supported elements of $\prod_{s\in S} A_s$ where at least one of the elements of $(A_s)_{s\in S}$, say $A_{s_0}$, is finite. But then, the coloring $\chi : \funct{\m{X}}{A_{s_0}}$ defined by $\chi(x)=x(s_0)$ divides $\m{X}$. Therefore:  

\begin{thm}

\label{thm:char uUr indivisible, S dually well ordered}

Let $\m{X}$ be a countable ultrahomogeneous ultrametric space with distance set $S \subset ]0, + \infty[$. Then $\m{X}$ is indivisible iff $\m{X}=\uUr$ and the reverse linear ordering $>$ on $\R$ induces a well-ordering on $S$. 

\end{thm}

This subsection is organized as follows. Theorem \ref{thm:uUr divisible, S not dually well ordered} is proved in \ref{subsubsection:Proof of theorem uUr divisible, S not dually well ordered.}. Theorem \ref{thm:uUr indivisible, S dually well ordered} is proved in \ref{subsubsection:Proof of theorem uUr indivisible, S dually well ordered.}. Finally, in \ref{subsubsection:An application}, we present an application of Theorem \ref{thm:uUr indivisible, S dually well ordered} dealing with restrictions of maps $f : \funct{\uUr}{\omega}$. 

\

\subsubsection{Proof of Theorem \ref{thm:uUr divisible, S not dually well ordered}.}

\label{subsubsection:Proof of theorem uUr divisible, S not dually well ordered.}

Fix a countable subset $S$ of $]0, + \infty[$ such that the reverse linear ordering $>$ on $\R$ does not induce a well-ordering on $S$. The idea to prove that $\uUr$ is divisible is to use a coloring which is constant on some particular spheres. More precisely, observe that $(S, >)$ not being well-ordered, there is a strictly increasing sequence $(s_i)_{i \in \omega}$ of reals such that $s_0 = 0$ and $s_i \in S$ for every $i > 0$. Observe that we can construct a subset $E$ of $\uUr$ such that given any $y \in \uUr$, there is exactly one $x$ in $E$ such that for some $i < \omega$, $d^{\uUr}(x,y) < s_i$. Indeed, if $\sup _{i < \omega}s_i = \infty$, simply take $E$ to be any singleton. Otherwise, let $\rho = \sup _{i < \omega}s_i$ and choose $E \subset \uUr$ maximal such that 
\[ \forall x, y \in E \ \ d^{\uUr}(x,y) \geqslant \rho. \]

To define $\chi : \funct{\uUr}{\omega}$, let $(A_j)_{j \in \omega}$ be a family of infinite pairwise disjoint subsets of $\omega$ whose union is $\omega$. Then, for $y \in \uUr$, let $e(y)$ and $i(y)$ be the unique elements of $E$ and $\omega$ respectively such that $d^{\uUr}(e(y),y) \in [s_{i(y)} , s_{i(y)+1}[$, and set \begin{center} $\chi (y) = j$ iff $i(y) \in A_j$. \end{center} 

\begin{claim}
$\chi$ is as required.
\end{claim}

\begin{proof}
Let $Y \subset \uUr$ be isometric to $\uUr$. Fix $y \in Y$. For every $j \in \omega$, pick $i_j > i(y) + 1$ such that $i_j \in A_j$. Since $Y$ is isometric to $\uUr$, we can find an element $y_j$ in $Y$ such that $d^{\uUr}(y, y_j) = s_{i_j}$. We claim that $\chi (y_j) = j$, or equivalently $i(y_j) \in A_j$. Indeed, consider the triangle $\{e(y), y, y_j \}$. Observe that in an ultrametric space every triangle is isosceles with short base and that here, 
\[ d^{\uUr}(e(y),y) < s_{i_j} = d(y,y_j). \]

Thus, 
\[ d^{\uUr}(e(y),y_j) = d^{\uUr}(y,y_j) \in [s_{i_j},s_{i_j +1}[. \]

And therefore $e(y_j) = e(y)$ and $i(y_j) = i_j \in A_j$. \end{proof}

\

\subsubsection{Proof of Theorem \ref{thm:uUr indivisible, S dually well ordered}.}

\label{subsubsection:Proof of theorem uUr indivisible, S dually well ordered.}

When $S \subset ]0 , + \infty[$ is finite, it follows from the proof of section \ref{section:Big Ramsey degrees} that the $1$-point ultrametric space has a big Ramsey degree equal to $1$. Thus, $\uUr$ is indivisible. From now on, we consequently concentrate on the case where $S$ is infinite. Fix an infinite countable subset $S$ of $]0, + \infty[$ such that the reverse linear ordering $>$ on $\R$ induces a well-ordering on $S$. Our goal here is to show that the space $\uUr$ is indivisible. For convenience, we will simply write $d$ instead of $d^{\uUr}$. 

Observe first that the collection $\B$ of metric balls of $\uUr$ is a
tree when ordered by reverse set-theoretic inclusion. When $x \in \uUr$
and $r \in S$, $B(x,r)$ denotes the set $\{ y \in \uUr : d^{\uUr}(x,y)
\leqslant r \}$. $x$ is called a \emph{center} of the ball and $r$ a
\emph{radius}. Note that in $\uUr$, non empty balls have a unique radius but
admit all of their elements as centers. Note also that when $s > 0 $
is in $S$, the fact that $(S,>)$ is well ordered  allows to define \index{$s^-$}
\[ s^- = \max \{t \in S : t < s \}. \]

The main ingredients are contained in the following definition and lemma. 

\begin{defn}

\label{defn:1}

Let $A \subset \uUr$ and $b \in \B$ with radius $r \in S \cup \{ 0 \}$. Say that \emph{$A$ is
small in $b$} when $r=0$ and $A \cap b = \emptyset$, or $r > 0$
and $A \cap b$ can be covered by finitely many balls of radius $r^-$.
\end{defn}

We start with an observation. Assume that $\{ x_n : n \in \omega \}$ is an enumeration of $\uUr$, and that we are trying to build inductively a copy $\{ a_n : n \in \omega \}$ of $\uUr$ in $A$ such that for every $n, m \in \omega$, $d(a_n , a_m) = d(x_n , x_m)$. Then the fact that we may be blocked at some finite stage exactly means that at that stage, a particular metric ball $b$ with $A \cap b \neq \emptyset$ is such that $A$ is small in $b$. This idea is expressed in the following lemma.

\begin{lemma}

\label{lemma:1}

Let $X \subset \uUr$. The following are equivalent:

i) $\binom{X}{\uUr} \neq \emptyset$.

ii) There is $Y \subset X$ such that $Y$ is not small in $b$ whenever $b \in \B$ and $Y \cap b \neq \emptyset$.

\end{lemma}

\begin{proof}
Assume that i) holds and let $Y$ be a copy of $\uUr$ in $X$. Fix $b \in \B$ with radius $r$ and such that $Y \cap b \neq \emptyset$. Pick $x \in Y \cap b $ and let $E \subset \uUr$ be an infinite subset where all the distances are equal to $r$. Since $Y$ is isometric to $\uUr$, $Y$ includes a copy $\tilde{E}$ of $E$ such that $x \in \tilde{E}$. Then $\tilde{E} \subset Y \cap b$ and cannot be covered by finitely many balls of radius $r^-$, so ii) holds.

Conversely, assume that ii) holds. Let $\{x_n : n \in \omega \}$ be an enumeration of the elements of $\uUr$. We are going to construct inductively a sequence $(y_n)_{n \in \omega}$ of elements of $Y$ such that \[ \forall m, n \in \omega \ \ d(y_m , y_n) = d(x_m , x_n).\]

For $y_0$, take any element in $Y$. In general, if $(y_n)_{n \leqslant k}$ is built, construct $y_{k+1}$ as follows. Consider the set $E$ defined as 
\[ E = \{ y \in \uUr : \forall \ n \leqslant k \ \ d(y,y_n) = d(x _{k+1} , x_n) \}. \]

Let also 
\[ r = \min \{ d(x_{k+1},x_n) : n \leqslant k \}. \]

and
\[ M = \{ n \leqslant k : d(x_{k+1},x_n) = r \}. \] 

We want to show that $E \cap Y \neq \emptyset$. Observe first that for every $m, n \in M$, $d(y_m , y_n) \leqslant r$. Indeed,
\[d(y_m , y_n ) =  d(x_m , x_n) \leqslant \max (d(x_m , x_{k+1}) , d(x_{k+1} , x_n)) = r.\]

So in particular, all the elements of $\{ y_m : m \in M \}$ are contained in the same ball $b$ of radius $r$.

\begin{claim}
$E = b \smallsetminus \bigcup _{m \in M} B(y_m, r^-)$.
\end{claim}

\begin{proof}
It should be clear that \[ E \subset b \smallsetminus \bigcup _{m \in M} B(y_m, r^-). \] 

On the other hand, let $y \in b \smallsetminus \bigcup _{m \in M} B(y_m, r^-)$. Then for every $m \in M$, \[ d(y,y_m) = r = d(x_{k+1}, x_m). \] 

Therefore, it remains to show that $d(y,y_n) = d(x_{k+1}, x_n)$ whenever $n \notin M$. To do that, we use again the fact that every triangle is isosceles with short base. Let $m \in M$. In the triangle $\{x_m, x_n, x_{k+1} \}$, we have $d(x_{k+1}, x_n) > r$ so \[ d(x_m, x_{k+1}) = r < d(x_n , x_m) = d(x_n , x_{k+1}). \]

Now, in the triangle $\{ y_m, y_n, y\}$, $d(y, y_m) = r$ and  $d(y_m , y_n) = d(x_m, x_n) > r$. Therefore, \[ d(y, y_n) = d(y_m, y_n) = d(x_m, x_n) = d(x_{k+1}, x_n). \qedhere \] \end{proof}

We consequently need to show that $(b \smallsetminus \bigcup _{m \in M} B(y_m, r^-)) \cap Y \neq \emptyset$. To achieve that, simply observe that when $m \in M$, we have $y_m \in Y \cap b$. Thus, $Y \cap b \neq \emptyset$ and by property ii), $Y$ is not small in $b$. In particular, $Y \cap b $ is not included in $\bigcup _{m \in M} B(y_m, r^-)$. \end{proof}

We are now ready to prove Theorem \ref{thm:uUr indivisible, S dually well ordered}. However, before we do so, let us make another observation concerning the notion smallness. Let $\uUr = A \cup B$. 

Note that if $A$ is small in $b \in \B$, then 1) $A \cap b$ cannot contribute to build a copy of $\uUr$ in $A$ and 2) $B \cap b$ is isometric to $b$. So intuitively, everything happens as if $b$ were completely included in $B$. So the idea is to remove from $A$ all those parts which are not essential and to see what is left at the end. More precisely, define a sequence $(A _{\alpha})_{\alpha \in \omega _1}$ recursively as follows: 

\begin{itemize}

\item $A_0 = A$.

\item $A_{\alpha + 1} = A_{\alpha} \smallsetminus \bigcup \{ b : A_{\alpha} \ \mathrm{is \ small \ in \ b} \}$. 

\item For $\alpha < \omega _1$ limit, $A_{\alpha} = \bigcap _{\eta < \alpha} A_{\eta}$. 

\end{itemize}

Since $\uUr$ is countable, the sequence is eventually constant. Set \begin{center} $\beta = \min \{ \alpha < \omega _1 : A_{\alpha +1} = A_{\alpha} \}$. \end{center}

Observe that if $A_{\beta}$ is non-empty, then $A_{\beta}$ is not small in any metric ball it intersects. Indeed, suppose that $b \in \B$ is such that $A_{\beta}$ is small in $b$. Then $A_{\beta + 1} \cap b = \emptyset$. But $A_{\beta + 1} = A_{\beta}$ so $A_{\beta} \cap b = \emptyset$. Therefore, since $A_{\beta} \subset A$, $A$ satisfies condition ii) of lemma \ref{lemma:1} and $\binom{A}{\uUr} \neq \emptyset$.

It remains to consider the case where $A_{\beta} =
\emptyset$. According to our second observation, the intuition is that
$A$ is then unable to carry any copy of $\uUr$ and is only composed of
parts which do not affect the metric structure of $B$. Thus, $B$
should include an isometric copy of $\uUr$. For $\alpha < \omega _1$,
let $\mathcal{C}_{\alpha}$ be the set of all minimal elements (in the
sense of the tree structure on $\B$) of the collection $\{ b \in \B :
A_{\alpha} \ \mathrm{is \ small \ in \ b} \}$. Equivalently, $\mathcal{C}_{\alpha}$ is the set of elements of $\{ b \in \B : A_{\alpha} \ \mathrm{is \ small \ in \ b} \}$ with largest radius. Note that since all
points of $B$ can be seen as balls of radius $0$ in which $A$ is
small, we have $B \subset \bigcup \mathcal{C}_0$. Note also that
$(\bigcup \mathcal{C}_{\alpha})_{\alpha < \omega _1}$ is increasing. By induction on $\alpha > 0$, it
follows that
\[ \forall \ 0 < \alpha < \omega _1 \ \ A_{\alpha} = \uUr \smallsetminus
\bigcup_{\eta < \alpha} \bigcup \mathcal{C}_{\eta} \ \ \ \ (*)\] 

\begin{claim}
Let $\alpha < \omega _1$, $b \in \mathcal{C}_{\alpha}$ with radius $r \in S$. Then $b \smallsetminus \bigcup _{\eta < \alpha} \bigcup \{ c \in \mathcal{C}_{\eta} : c \subset b \}$ is small in $b$. 

\end{claim}

\begin{proof}
$A_{\alpha}$ is small in b so find $c_0 \ldots c_{n-1} \in \B$ with radius $r^-$ and included in $b$ such that \[ A_{\alpha} \cap b \subset \bigcup_{i<n}c_i .\]

Then thanks to $(*)$ \[ b \smallsetminus \bigcup_{i<n}c_i \subset \bigcup_{\eta < \alpha} \bigcup \mathcal{C}_{\eta}. \]

Note that by minimality of $b$, if $\eta < \alpha$, then $b \subsetneq c $ cannot happen for any element of $\mathcal{C}_{\eta}$. It follows that either $c \cap b = \emptyset$ or $c \subset b$. Therefore, \[ b \smallsetminus \bigcup_{i<n}c_i \subset \bigcup_{\eta < \alpha} \bigcup \{ c \in \mathcal{C}_{\eta} : c \subset b\}. \qedhere \]  

\end{proof}

\begin{claim}
Let $\alpha < \omega _1$ and $b \in \mathcal{C}_{\alpha}$. Then $\binom{B \cap b}{b} \neq \emptyset$. 
\end{claim}

\begin{proof}
We proceed by induction on $\alpha < \omega _1$. 

For $\alpha = 0$, let $b \in \mathcal{C}_0$. If the radius $r$ of $b$ is $0$, there is nothing to do. If $r>0$, then $r\in S$. $A_0 = A$ is small in $b$ so find $c_0, \ldots , c_{n-1}$ with
radius $r^-$ such that $A \cap b \subset \bigcup_{i<n}c_i$. Then $b
\smallsetminus \bigcup_{i<n}c_i$ is isometric to $b$ and is included
in $B \cap b$.

Suppose now that the claim is true for every $\eta < \alpha$. Let $b \in \mathcal{C}_{\alpha}$ with radius $r \in S$. Thanks to the previous claim, we can find $c_0 \ldots c_{n-1} \in \B$ with radius $r^-$ and included in $b$ such that \[
b = \bigcup _{i<n} c_i \cup \bigcup _{\eta < \alpha} \bigcup \{ c \in \mathcal{C}_{\eta} : c \subset b \}.\]

Observe that \[ \bigcup _{\eta < \alpha} \bigcup \{ c \in \mathcal{C}_{\eta} : c \subset b \} = \bigcup \{ c \in \bigcup _{\eta < \alpha} \mathcal{C}_{\eta} : c \subset b \}. \]

Define $\mathcal{D}_{\alpha}$ as the set of all minimal
elements (still in the sense of the tree structure on $\B$) of the
collection \[ \{ c \in \bigcup _{\eta < \alpha} \mathcal{C}_{\eta} : c \subset b \ \mathrm{and} \ \forall i < n \ \ c \cap c_i = \emptyset \}. \]

Then $\{c_i : i < n \} \cup \mathcal{D}_{\alpha}$ is a collection of pairwise disjoint balls and $\bigcup \mathcal{D}_{\alpha}$ is isometric to $b$. By induction hypothesis, $\binom{B \cap c}{c} \neq \emptyset$ whenever $c \in \mathcal{D}_{\alpha}$ and there is an isometry $\varphi _c : \funct{c}{B \cap c}$. Now, let $\varphi : \funct{\bigcup \mathcal{D}_{\alpha}}{B \cap b}$ be defined as \[ \varphi = \bigcup_{c \in \mathcal{D}_{\alpha}} \varphi _c .\] 

We claim that $\varphi$ is an isometry. Indeed, let $x, x' \in \bigcup \mathcal{D}_{\alpha}$. If there is $c \in \mathcal{D}_{\alpha}$ such that $x, x' \in c$ then \[ d(\varphi (x) , \varphi (x')) = d(\varphi _c (x) , \varphi _c (x')) = d(x,x'). \]

Otherwise, find $c \neq c' \in \mathcal{D}_{\alpha}$ with $x \in c$ and $x' \in c'$. Observe that since we are in an ultrametric space, we have \[ \forall y, z \in c \ \ \forall y', z' \in c' \ \ d(y,y') = d(z,z'). \] 

Thus, since $x, \varphi (x) \in c$ and $x', \varphi(x') \in c'$, we get \[ d(\varphi (x) , \varphi (x')) = d(x,x'). \qedhere \] \end{proof}

To finish the proof of Theorem \ref{thm:uUr indivisible, S dually well ordered}, it suffices to notice that as a metric ball (the unique ball of radius $\max S$), $\uUr$ is in $\mathcal{C}_{\beta}$. So according to the previous claim, $\binom{B}{\uUr} \neq \emptyset$ and we are done. 

\vspace{1em}

\subsubsection{An application of Theorem \ref{thm:uUr indivisible, S dually well ordered}.}

\label{subsubsection:An application}

Let $S \subset ]0, + \infty[$ be infinite and countable such that the
reverse linear ordering $>$ on $\R$ induces a well-ordering on $S$. We
saw that $\uUr$ is then indivisible but that there is no big
Ramsey degree for any $\m{X} \in \U$ as soon as $|\m{X}| \geqslant
2$. In other words, in the present context, the analogue of infinite
Ramsey's theorem holds in dimension $1$ but fails for higher
dimensions. Still, one may ask if some partition result fitting in
between holds. For example, given any $f : \funct{\uUr}{\omega}$, is
there an isometric copy of $\uUr$ inside $\uUr$ on which $f$ is
constant or injective? Such a property is sometimes refered to as \emph{selectivity}\index{selectivity}. Selectivity can be thought as an intermediate Ramsey-type result between dimension $1$ and $2$. It is indeed clearly stronger than the $1$-dimensional result, but is in turn implied by the $2$ dimensional one if one considers the $2$-coloring $\chi$ defined by $\chi (\{ x,y\}) = 1$ iff $f(x) = f(y)$. It turns out that in the present case, selectivity does not hold. To see
that, consider a family $(b_n)_{n \in \omega}$ of disjoint balls
covering $\uUr$ whose sequence of corresponding radii $(r_n)_{n \in
  \omega}$ decreases strictly to $0$ and define $f :
\funct{\uUr}{\omega}$ by $f(x) = n$ iff $x \in b_n$. Then $f$ is not constant or injective on any isometric copy of $\uUr$. Observe in fact that $f$ is neither uniformly continuous nor injective on any isometric copy of
$\uUr$. However, if ``uniformly continuous'' is replaced by
``continuous'', then the result becomes true:

\begin{thm}

\label{thm:selectivity for uUr}

Let $S$ be an infinite countable subset of $]0, + \infty[$ such that the
reverse linear ordering $>$ on $\R$ induces a well-ordering on $S$. Then given any map $f : \funct{\uUr}{\omega}$, there is an isometric copy $X$ of $\uUr$ inside $\uUr$ such that $f$ is continuous or injective on $X$. 

\end{thm}

The purpose of what follows is to provide a proof of that fact. The reader will notice the similarities with the proof of Theorem \ref{thm:uUr indivisible, S dually well ordered}.

\begin{defn}

Let $f : \funct{\uUr}{\omega}$, $Y \subset \uUr$ and $b \in \B$ with radius $r > 0$. Say that \emph{$f$ has almost finite range on $b$ with respect to $Y$}\index{almost finite range} when there is a finite family $(c_i)_{i < n}$ of elements of $\B$ with radius $r^-$ such that $f$ has finite range on $Y \cap (b \smallsetminus \bigcup _{i<n} c_i)$.  

\end{defn}

\begin{lemma}

\label{lemma:4}

Let $f : \funct{\uUr}{\omega}$ and $Y \subset \uUr$ such that for every $b \in \B$ meeting $Y$, $f$ does not have almost finite range on $b$ with respect to $Y$. Then there is an isometric copy of $\uUr$ included in $Y$ on which $f$ is injective. 

\end{lemma}

\begin{proof}

Let $\{x_n : n \in \omega \}$ be an enumeration of the elements of $\uUr$. Our goal is to construct inductively a sequence $(y_n)_{n \in \omega}$ of elements of $Y$ on which $f$ is injective and such that \[ \forall m, n \in \omega \ \ d(y_m , y_n) = d(x_m , x_n). \] 

For $y_0$, take any element in $Y$. In general, if $(y_n)_{n \leqslant k}$ is built, construct $y_{k+1}$ as follows. Consider the set $E$ defined as 
\[ E = \{ y \in \uUr : \forall \ n \leqslant k \ \ d(y,y_n) = d(x _{k+1} , x_n) \}. \]

As in lemma \ref{lemma:1}, there is $b \in \B$ with radius $r > 0$ intersecting $Y$ and a set $M$ such that \[ E = b \smallsetminus \bigcup _{m \in M} B(y_m, r^-).\] 

Since $f$ does not have almost finite range on $b$ with respect to $Y$, $f$ takes infinitely many values on $E$ and we can choose $y_{k+1} \in E$ such that \[ \forall n \leqslant k \ \ f(y_n) \neq f(y_{k+1}). \qedhere  \] \end{proof}

We now turn to a proof of Theorem \ref{thm:selectivity for uUr}. Here, our strategy is
to define recursively a sequence $(Q _{\alpha})_{\alpha \in \omega
  _1}$ whose purpose is to get rid of all those parts of $\uUr$ on which $f$ is essentially of finite range: 

\begin{itemize}

\item $Q_0 = \uUr$.

\item $Q_{\alpha + 1} = Q_{\alpha} \smallsetminus \bigcup \{ b : \mathrm{f \ has \ almost \ finite \ range \ on} \ b \ \mathrm{with \ respect \ to} \ Q_{\alpha} \}$. 

\item For $\alpha < \omega _1$ limit, $Q_{\alpha} = \bigcap _{\eta < \alpha} Q_{\eta}$. 

\end{itemize}

$\uUr$ being countable, the sequence is eventually constant. Set \begin{center} $\beta = \min \{ \alpha < \omega _1 : Q_{\alpha +1} = Q_{\alpha} \}$. \end{center}

If $Q_{\beta}$ is non-empty, then $f$ and $Q_{\beta}$ satisfy the hypotheses of lemma \ref{lemma:4}. Indeed, suppose that $b \in \B$ is such that $f$ has almost finite range on $b$ with respect to $Q_{\beta}$. Then $Q_{\beta + 1} \cap b = \emptyset$. But $Q_{\beta + 1} = Q_{\beta}$ so $Q_{\beta} \cap b = \emptyset$. 

Consequently, suppose that $Q_{\beta} = \emptyset$. The intuition is that on any ball $b$, $f$ is essentially of finite range. Consequently, we should be able to show that there is $X \in \binom{\uUr}{\uUr}$ on which $f$ is continuous. 

For $\alpha < \omega _1$,
let $\mathcal{C}_{\alpha}$ be the set of all minimal elements of the collection \[ \{ b : \mathrm{f \ has \ almost \ finite \ range \ on} \ b \ \mathrm{with \ respect \ to} \ Q_{\alpha} \}.\]

Then  

\[ \forall \ 0 < \alpha < \omega _1 \ \ Q_{\alpha} = \uUr \smallsetminus
\bigcup_{\eta < \alpha} \bigcup \mathcal{C}_{\eta} \ \ \ \ (**)\] 

\begin{claim}
Let $\alpha < \omega _1$ and $b \in \mathcal{C}_{\alpha}$. Then there is $\tilde{b} \in \binom{b}{b}$ on which $f$ is continuous.  
\end{claim}

\begin{proof}
We proceed by induction on $\alpha < \omega _1$. 

For $\alpha = 0$, let $b \in \mathcal{C}_0$. $f$ has almost finite range on $b$ with respect to $Q_0 = \uUr$ so find $c_0, \ldots , c_{n-1}$ with radius $r^-$ such that $f$ has finite range on $b \smallsetminus \bigcup_{i<n}c_i$. Then $b
\smallsetminus \bigcup_{i<n}c_i$ is isometric to $b$. Now, by Theorem \ref{thm:uUr indivisible, S dually well ordered}, $b$ is indivisible. Therefore, there is $\tilde{b} \in \binom{b}{b}$ on which $f$ is constant, hence continuous.

Suppose now that the claim is true for every $\eta < \alpha$. Let $b \in \mathcal{C}_{\alpha}$ with radius $r \in S$. 
Find $c_0 \ldots c_{n-1} \in \B$ with radius $r^-$ and included in $b$ such that $f$ has finite range on $Q_{\alpha} \cap (b \smallsetminus \bigcup _{i<n} c_i) $. Then $b' := b \smallsetminus \bigcup _{i<n} c_i$ is isometric to $b$ and thanks to $(**)$, \[ b' = (b'\cap Q_{\alpha}) \cup (b' \cap \bigcup _{\eta < \alpha} \bigcup \mathcal{C}_{\eta}).
\]

Now, let $\mathcal{D}_{\alpha}$ be defined as the set of all minimal elements of the collection \[ \{ c \in \bigcup _{\eta < \alpha} \mathcal{C}_{\eta} : c \subset b \ \mathrm{and} \ \forall i < n \ \ c \cap c_i = \emptyset \}. \]

Then, for the same reason as in section 3, we have \[ b' = (b' \cap Q_{\alpha}) \cup \bigcup \mathcal{D}_{\alpha}. \]

Thanks to Theorem \ref{thm:uUr indivisible, S dually well ordered}, $b' \cap Q_{\alpha}$ or $\bigcup \mathcal{D}_{\alpha}$ includes an isometric copy $\tilde{b}$
of $b$. If $b' \cap Q_{\alpha}$ does, then for every $i<n$, $c_i \cap
\tilde{b}$ is a metric ball of $\tilde{b}$ of same radius as
$c_i$. Thus, $\tilde{b} \smallsetminus \bigcup _{i<n} c_i$ is an
isometric copy of $b$ on which $f$ takes only finitely many values and
Theorem \ref{thm:uUr indivisible, S dually well ordered} allows to conclude. Otherwise, suppose
that $\bigcup \mathcal{D}_{\alpha}$ includes an isometric copy of $b$. Note
that $\bigcup \mathcal{D} _{\alpha}$ includes an isometric
copy of itself on which $f$ is continuous. Indeed, by induction
hypothesis, for every $c \in \mathcal{D} _{\alpha}$, there is an
isometry $\varphi _c : \funct{c}{c}$ such that $f$ is
continuous on the range $\varphi _c ''c$ of $\varphi _c$. As in the
previous section, one obtains an isometry by setting $\varphi := \funct{\bigcup \mathcal{D}_{\alpha}}{\bigcup \mathcal{D}_{\alpha}}$ defined as \[ \varphi = \bigcup_{c \in
  \mathcal{D}_{\alpha}} \varphi _c .\] 

Thus, its range $\varphi '' \bigcup \mathcal{D}
_{\alpha} $ is an isometric copy of $\bigcup
\mathcal{D}_{\alpha}$ on which $f$ is continuous. Now, since
$\bigcup \mathcal{D}_{\alpha}$ includes an isometric copy of $b$, so does $\varphi '' \bigcup \mathcal{D}
_{\alpha} $ and we are done. \end{proof}

We conclude with the same argument we used at the end of Theorem \ref{thm:uUr indivisible, S dually well ordered}: As a metric ball, $\uUr$ is in $\mathcal{C}_{\beta}$. Thus, there is an isometric copy $X$ of $\uUr$ inside $\uUr$ on which $f$ is continuous.

\subsection{Indivisibility of $\Ur _S$.}

\label{subsection:Indivisibility of U_S}
\index{indivisibility!for $\Ur _S$}

The last spaces we will be studying in this section on indivisibility are the spaces $\Ur _S$ where $S$ is a finite set satisfying the $4$-values condition. We saw already that they provided a wide variety of combinatorial objects and that the classes $\M _S$\index{$\M _S$} to which they are attached seemingly behave quite well from a Ramsey-theoretic point of view. The purpose of this subsection is to show that to some extend, this apparent good behaviour of the $\M _S$'s also appears at the level of their Urysohn spaces. The first result here reads as follows: 

\begin{thm}

\label{thm:U_S indiv ext}

Let $S = \{s_0,\ldots, s_m \}$ be finite subset of $]0, + \infty [$ satisfying the $4$-values condition and such that for every $i<m$, $s_{i+1} \leqslant 2s_i$. Then $\Ur _S$ is indivisible. 

\end{thm}

The proof of this theorem comes from a direct adaptation of the proof of Theorem \ref{thm:U_m indiv}, noting that the method used for the spaces $\Ur_m$ actually applies for $\Ur_S$ provided that $S$ does not have any large gap. However, if one tries to get rid of that requirement, serious obstacles appear and the result we obtain is at the price of a serious restriction on the size of $S$: 

\begin{thm}

\label{thm:U_S indiv}

Let $S$ be finite subset of $]0, + \infty [$ of size $|S| \leqslant 4$ and satisfying the $4$-values condition. Then $\Ur _S$ is indivisible. 

\end{thm}

\begin{proof}
When the proofs are not elementary, their main ingredients are Milliken's theorem (Theorem \ref{thm:Milliken}), Sauer's theorem (Theorem \ref{thm:Sauer ultrahomogeneous}) or Theorem \ref{thm:U_m indiv} stated in \ref{subsection: Are the U_m 's indivisible?} and \ref{subsection:U_m are indiv}. As mentioned in chapter 1, there are many classes $\M _S$, and hence many spaces $\Ur _S$ when $S$ has size 4 and satisfies the $4$-values condition. Thus, we only cover here the cases where $|S|\leqslant 3$. The cases where $|S|=4$ are treated in appendix.   

For $|S|=1$, the result is trivial. 

For $|S|=2$: When $S = \{ 1, 2 \}$, the Urysohn space is the Rado graph equipped with the path metric. The Rado graph being indivisible, so is $\Ur _{\{1,2\}}$. When $S = \{ 1, 3\}$, $\Ur _{\{1,3\}}$ is ultrametric and is indivisible thanks to Theorem \ref{thm:uUr indivisible, S dually well ordered}.   

For $|S|=3$: 

(1a) $S = \{ 2, 3, 4\}$. The space $\Ur _{\{ 2, 3, 4\}}$ can be seen as a complete version of the Rado graph with three kinds of edges. An easy variation of the proof working for the Rado graph shows that $\Ur _{\{ 2, 3, 4\}}$ is indivisible. 

(1b) $S = \{ 1, 2, 3\}$. The space $\Ur _{\{ 1, 2, 3\}}$ is the space we denoted $\Ur _3$ and we saw in Theorem \ref{thm:U_3 indivisible} that it is indivisible.  

(1d) $S = \{ 1, 2, 5\}$. The space $\Ur _{\{ 1, 2, 5\}}$ is composed of countably many disjoint copies of $\Ur _2$, and the distance between any two points not in the same copy of $\Ur _2$ is always $5$. The indivisibility of $\Ur _2$ consequently implies that $\Ur _{\{ 1, 2, 5\}}$ is indivisible. 

(2a) $S = \{ 1, 3, 4\}$. The space $\Ur _{ \{ 1, 3, 4\}}$ is composed of countably many disjoint copies of $\Ur _1$ and points belonging to different copies of $\Ur _1$ can be randomly at distance $3$ or distance $4$ apart. As for $\Ur _2$, its indivisibility can be proved via Milliken theorem: Fix an $\omega$-linear ordering $<$ on $2 ^{< \omega}$ extending the tree ordering and consider the standard graph structure on $2^{< \omega}$: 

\begin{center}
$\forall s < t \in 2^{< \omega} \ \ \{ s , t \} \in E \leftrightarrow \left( |s| < |t| , t(|s|) = 1 \right)$. 
\end{center}

Now, define a map $d$ on the set $[2 ^{< \omega}]^2$ of pairs of $2 ^{< \omega}$ as follows: Let $\{ s, t \}_<$, $\{ s' , t' \}_< $ be in $[2 ^{< \omega}]^2$. Then define $d (\{ s, t \}_< , \{ s' , t' \}_< )$ as:

\begin{displaymath}
\left \{ \begin{array}{ll}
 1 & \textrm{if $s=s'$}\\
 3 & \textrm{if $s \neq s'$ and $\{ t, t'\} \in E$.}\\
 4 & \textrm{if $s \neq s'$ and $\{ t, t' \} \notin E$.}
 \end{array} \right.
\end{displaymath}

It is easy to check that $d$ is a metric. Since $d$ takes its values in $\{ 1, 3, 4\}$, $([2 ^{< \omega}]^2 , d)$ embeds into $\Ur _{ \{ 1, 3, 4\}}$. We now claim that the space $\Ur _{ \{ 1, 3, 4\}}$ embeds into $([2 ^{< \omega}]^2 , d)$. To do that, we actually show that $\Ur _{ \{ 1, 3, 4\}}$ embeds into the subspace $\m{X}$ of $([2 ^{< \omega}]^2 , d)$ supported by the set 

\begin{center}
$ X = \{ \{ s , t \}_< \in [2 ^{< \omega}]^2 : |s| < |t|, \ s <_{lex} t, \ t(|s|) = 0 \}$.
\end{center} 

The embedding is constructed inductively. Let $\{ x_n : n \in \omega \}$ be an enumeration of $\Ur _{ \{ 1, 3, 4\}}$. We are going to construct a sequence $(\{ s_n , t_n \})_{n \in \omega}$ of elements in $X$ such that 

\begin{center} 
$\forall m, n \in \omega \ \ d(\{ s, t \}_< , \{ s' , t' \} _<) = d^{\Ur _{ \{ 1, 3, 4\}}} (x_m , x_n)$. 
\end{center}

For $\{ s_0 , t_0 \} _<$, take $s_0 = \emptyset$ and $t_0 = 0$. Assume now that $\{ s_0 , t_0 \} _< ,\ldots , \{ s_n , t_n \} _<$ are constructed such that all the elements of $\{ s_0 ,\ldots , s_n\} \cup \{ t_0 ,\ldots , t_n \}$ have different heights and all the $s_i$'s are strings of $0$'s. Set 

\begin{center}
$M = \{ m \leqslant n : d^{\Ur _{\{ 1, 3, 4\}}} (x_m , x_{n+1}) = 1 \}$. 
\end{center}

If $M = \emptyset $, choose $s_{n+1}$ to be a string of $0$'s longer that all the elements constructed so far. Otherwise, there is $s \in 2^{< \omega}$ such that 

\begin{center} 
$\forall m \in M \ \ s_m = s$. 
\end{center}

Set $s_{n+1} = s$. Now, choose $t_{n+1}$ above all the elements constructed so far and such that 

\vspace{0.5em}
\hspace{1em}
i) $\forall m \notin M \ \ (t_{n+1} (|t_m|) = 1) \leftrightarrow (d^{\Ur _{ \{ 1, 3, 4\}}}(x_{n+1} , x_m) = 3)$.

\vspace{0.5em}
\hspace{1em}
ii) $\{ s_{n+1} , t_{n+1}\}_< \in X$. 

\vspace{0.5em} 

The requirement i) is easy to satisfy because all the $t_m$'s have different heights. As for ii), $|s_{n+1}| < |t_{n+1}|$ and $ t_{n+1} (|s_{n+1}|) = 0$ are also easy (again because all heights are different) while $s_{n+1} <_{lex} t_{n+1}$ is satisfied because $s_{n+1}$ being a $0$ string, $|s_{n+1}| < |t_{n+1}|$ implies $s_{n+1} <_{lex} t_{n+1}$. After $\omega$ steps, we are left with a set $\{ \{ s_n , t_n \} : n \in \omega \} \subset \m{X}$ isometric to $\Ur _{ \{ 1, 3, 4\}}$. Observe that actually, this construction shows that $\Ur _{ \{ 1, 3, 4\}}$ embeds into any subspace of $([2 ^{< \omega}]^2 , d)$ supported by a strong subtree of $2^{< \omega}$. 

Now, to prove that $\Ur _{ \{ 1, 3, 4\}}$ is indivisible, it suffices to prove that given any $\chi : \funct{([2 ^{< \omega}]^2 , d)}{k}$ where $k \in \omega$ is strictly positive, there is a strong subtree $\m{T}$ of $2^{< \omega}$ such that $\chi$ is constant on $[T]^2 \cap X$. But this is guaranteed by Milliken theorem: Indeed, consider the subset $A := \{ 0, 01\}$. Then using the notation introduced for Theorem \ref{thm:Milliken}, $[A]_{\mathrm{Em}} = X$. So the restriction $\restrict{\chi}{[A]_{\mathrm{Em}}}$ is really a coloring of $X$, and there is a strong subtree $\m{T}$ of height $\omega$ such that $\restrict{[A]_{\mathrm{Em}}}{T} = [T]^2 \cap X$ is $\chi$-monochromatic. 

(2b) $S = \{ 1, 3, 6\}$. The space $\Ur _{\{ 1, 3, 6\}}$ is obtained from $\Ur _2$ after having multiplied all the distances by $3$ and blown the points up to copies of $\Ur _1$. Its indivisibility is a direct consequence of the basic infinite pigeonhole principle and of the indivisibility of $\Ur _2$.  

(2c) $S = \{ 1, 3, 7\}$. The space $\Ur _S$ is indivisible because ultrametric. \end{proof}

At that point, a comment can be made about the general problem of indivisibility of the spaces $\Ur_S$: Theorem \ref{thm:U_S indiv} is proved thanks to a case by case analysis. There is therefore very little hope that this method will lead to the proof of the general case. Still, our feeling is that Theorem \ref{thm:U_S indiv} should be thought as a good intermediate result towards a general solution. Indeed, even though $|S|\leqslant 4$ is a severe restriction, the large panel of combinatorial situations it provides seems to us of a reasonable variety. Our guess is therefore that given every $S$, the space $\Ur_S$ is indivisible.



\section{Approximate indivisibility and oscillation stability.}

\label{section:Almost indivisibility and oscillation stability}

After the study of indivisibility of countable Urysohn spaces, we now turn to the study of approximate indivisibility of complete separable metric spaces. As presented in section \ref{section:Definitions and notations}, in the realm of ultrahomogeneous metric spaces, approximate indivisibility corresponds to oscillation stability whose formulation brings topological groups into the picture. This fact is worth being mentioned as one of the most significant metric Ramsey-type theorems, namely Milman's theorem, appeared in close connection with topological groups dynamics. For $N \in \omega$ strictly positive, let $\mathbb{S}^N$ denote the unit sphere of the $(N+1)$-dimensional Euclidean space. Recall also $\mathbb{S} ^{\infty}$\index{$\mathbb{S} ^{\infty}$} denotes the unit sphere of the Hilbert space. Milman's theorem can then be stated as follows:

\begin{thm}[Milman \cite{Mil}]

\label{thm:Milman}
\index{Milman theorem}

Let $f : \funct{\mathbb{S}^{\infty}}{\R}$ be uniformly continuous. Then for every $\varepsilon > 0$ and every $N \in \omega$, there is a vector subspace $V$ of $\ell _2$ with $\dim V = N$ such that

\begin{center}
$\mathrm{osc}(\restrict{f}{V \cap \mathbb{S} ^{\infty}}) < \varepsilon$. 
\end{center}

\end{thm}

Equivalently:

\begin{thm}[Milman \cite{Mil}]

\label{thm:Milman'}

Let $\gamma$ be a finite cover of $\mathbb{S} ^{\infty}$. Then for every $\varepsilon > 0$ and every $N \in \omega$, there is $A \in \gamma$ and an isometric copy $\widetilde{\mathbb{S}}^N$ of $\mathbb{S}^N$ in $\mathbb{S} ^{\infty}$ such that $\widetilde{\mathbb{S}}^N \subset (A)_{\varepsilon}$. 

\end{thm}

Milman's theorem is at the heart of the recent books \cite{Pe1} and \cite{Pe1'}, where the interested reader will find a wide variety of  its developments in geometric functional analysis, topological group theory and combinatorics. One of the most famous questions raised after the discovery of Milman's theorem is known as the \emph{distortion problem for $\ell _2$}\index{distortion problem for $\ell _2$} and asks the following: Does Milman's theorem still hold when $N$ is replaced by $\infty$? In other words, if $f : \funct{\mathbb{S}^{\infty}}{\R}$ is uniformly continuous and $\varepsilon > 0$, is there an infinite-dimensional subspace $V$ of $\ell _2$ such that $\mathrm{osc}(\restrict{f}{V \cap \mathbb{S} ^{\infty} }) < \varepsilon$? Or, with the terminology introduced in section \ref{section:Definitions and notations}: Is $\mathbb{S} ^{\infty}$ approximately indivisible? This problem remained opened for about 30 years, until the solution of Odell and Schlumprecht in \cite{OS}:

\begin{thm}[Odell-Schlumprecht \cite{OS}]

\label{thm:Odell-Schlumprecht}
\index{Odell-Schlumprecht theorem}
\index{indivisibility!approximate indivisibility!for $\mathbb{S} ^{\infty}$}

$\mathbb{S} ^{\infty}$ is not approximately indivisible.  

\end{thm}

However, quite surprisingly, this solution is not based on an analysis of the intrinsic geometry of $\ell _2$. For that reason, it is sometimes felt that something essential is still to be discovered about the metric structure of $\mathbb{S} ^{\infty}$. This impression is certainly one of the motivations for the introduction of the concept of oscillation stability as presented in section \ref{section:Definitions and notations}. From this point of view, the approximate indivisibility problem for the Urysohn sphere $\s$ inherits a special status: Behind a solution based on the geometry of $\s$, a better understanding of $\mathbb{S} ^{\infty}$ might be hidden\ldots But at the present moment, it is unclear whether such a belief is justified or not. What \emph{is} clear is that very little is currently known about approximate indivisibility of ultrahomogeneous complete separable metric spaces or even about oscillation stability for topological groups in general. With the exception of Theorem \ref{thm:Odell-Schlumprecht}, the most significant result so far in the field was obtained by Hjorth in \cite{Hj}: 
 
\begin{thm}[Hjorth \cite{Hj}]

\label{thm:Hjorth}
\index{Hjorth!theorem on oscillation stability}

Let $G$ be a non-trivial Polish group. Then the action of $G$ on itself by left multiplication is not oscillation stable. 

\end{thm}

This section is organized as follows: In \ref{subsection:Metric oscillation stability for ultrametric Urysohn spaces.}, we solve the approximate indivisibility problem for the ultrametric Urysohn spaces. We then turn in \ref{subsection:Metric oscillation stability of S} to the approximate indivisibility problem for the Urysohn sphere. 

\

\textbf{Remark.} Before the concept of oscillation stability for topological groups was introduced by Kechris, Pestov and Todorcevic, Milman's work led to a notion which we will call here \emph{classical oscillation stability}\index{oscillation!classical oscillation stability}. This concept has now been central in geometric functional analysis for several decades and is already visible in the formulation of Theorem \ref{thm:Milman}: Given a Banach space $E$, a function $f : \funct{\mathbb{S}_E}{\R}$ defined on the unit sphere $\mathbb{S}_E$ of $E$ is \emph{oscillation stable in the classical sense} if for every infinite-dimensional closed subspace $Y$ of $E$, and every $\varepsilon > 0$, there is a infinite-dimensional closed subspace $Z$ of $Y$ such that  

\begin{center}
$\mathrm{osc}(\restrict{f}{Z \cap \mathbb{S}_E}) < \varepsilon$. 
\end{center}

Now, say that $E$ is \emph{oscillation stable in the classical sense} if every uniformly continuous $f : \funct{\mathbb{S}_E}{\R}$ is oscillation stable in the classical sense. In spirit, classical oscillation stability and oscillation stability for topological groups are consequently closely related. In some cases, they even coincide: When $\mathbb{S}_E$ is ultrahomogeneous as a metric space, classical oscillation stability for a Banach space $E$ is equivalent to oscillation stability of its unit sphere in the sense of \cite{KPT}. However, this case is quite exceptional: When $\mathbb{S}_E$ is \emph{not} ultrahomogeneous (which actually holds as soon as $E$ is not a Hilbert space), this equivalence does not hold anymore and there is no direct connection between classical oscillation stability and oscillation stability for topological groups.


\subsection{Approximate indivisibility for complete separable ultrametric spaces.}

\label{subsection:Metric oscillation stability for ultrametric Urysohn spaces.}
\index{indivisibility!approximate indivisibility!for complete ultrahomogeneous ultrametric spaces}

We saw in \ref{subsection:Indivisibility of Urysohn ultrametric spaces} that the indivisibility problem was completely solved for ultrametric Urysohn spaces. When passing to the metric completion, this allows to solve the approximate indivisibility problem for the complete separable ultrahomogeneous ultrametric spaces:  

\begin{thm}

\label{thm:ultrametric oscillation stability}

Let $\m{X}$ be a complete separable ultrahomogeneous ultrametric space. The following are equivalent: 

\begin{enumerate}
	\item[i)] $\m{X}$ is approximately indivisible. 
	\item[ii)] $\m{X}=\cUr$ for some $S \subset ]0, + \infty[$ finite or countable on which the reverse linear ordering $>$ on $\R$ induces a well-ordering. 
\end{enumerate}

\end{thm}

\begin{proof}

The implication $ii) \rightarrow i)$ is a consequence of Theorem \ref{thm:uUr indivisible, S dually well ordered}, which specifies that if $S \subset ]0, + \infty[$ is finite or countable such that the reverse linear ordering $>$ on $\R$ induces a well-ordering, then $\uUr$ is indivisible. For $i) \rightarrow ii)$, let $S$ denote the distance set of $\m{X}$.

We start by considering the case where $0$ is not an accumulation point of $S$. Then $\m{X}$ is discrete and therefore countable. Take $\varepsilon > 0$ such that $\varepsilon < \min S$. Since $\m{X}$ is approximately indivisible, it is in particular $\varepsilon$-indivisible, which in the present case truly means indivisible. So by Theorem \ref{thm:char uUr indivisible, S dually well ordered}, $\m{X} = \uUr \ (= \cUr)$ and the reverse linear ordering $>$ on $\R$ induces a well-ordering on $S$.  

Assume now that that $0$ is an accumulation point of $S$. Then thanks to a result in Chapter 1, Section \ref{subsection:Complete separable ultrahomogeneous ultrametric spaces}, there is a sequence $(A_s)_{s\in S}$ of elements of $\omega \cup \{ \Q \}$ with size at least $2$ such that $\m{X}$ is the set of all elements $x \in \prod_{s\in S} A_s$ whose support is a subset of $\{s_i : i \in \omega \}$ for some strictly decreasing sequence $(s_i)_{i \in \omega}$ of elements of $S$ converging to $0$. The distance is given by: 

\begin{center}
$d^{\m{X}}(x,y) = \min \{s \in S : \forall t \in S (s<t \rightarrow x(t) = y(t)) \}$.  
\end{center}

Because $\m{X}$ is approximately indivisible, no element of $(A_s)_{s\in S}$ is finite: If, say, $A_{s_0}$, were finite, then the coloring $\chi : \funct{\m{X}}{A_{s_0}}$ defined by $\chi(x)=x(s_0)$ would contradict $\varepsilon$-indivisibility for any $\varepsilon < t_0$. Hence, $A_s = \Q$ for every $s\in T$ and so $\m{X} = \cUr$. It remains to show that the reverse linear ordering $>$ on $\R$ induces a well-ordering on $S$. Assume not. Observe then that the extension $\widehat{\chi}$ to $\m{X}$ of the coloring $\chi$ used in the proof of Theorem \ref{thm:uUr divisible, S not dually well ordered} to divide $\uUr$ contradicts the fact that $\m{X}$ is approximately indivisible. \qedhere

\end{proof}

\subsection{Approximate indivisibility of $\s$.}

\label{subsection:Metric oscillation stability of S}
\index{indivisibility!approximate indivisibility!for $\s$}

As already mentioned in section \ref{subsection:Divisibility of S_Q}, the first attempt towards the approximate indivisibility for $\s$\index{$\s$} corresponds to the study of the indivisiblity problem for $\s _{\Q}$\index{$\s _{\Q}$}: Had $\s _{\Q}$ been indivisible, $\s$ would have been approximately indivisible. However, we saw with Theorem \ref{thm:s_Q divisible} that $\s _{\Q}$ is not indivisible. Worse: The proof of that fact does provide any information about $\s$, so the approximate indivisibility problem for $\s$ has to be attacked from another direction. The purpose of this subsection is to provide such an alternative. In essence, the idea remains the same: Approximate indivisibility for $\s$ should be attacked via the study of the exact indivisibility of simpler spaces. $\s _{\Q}$ was the first natural candidate because it is a very good countable approximation of $\s$. But this good approximation is paradoxically responsible for the divisibility of $\s _{\Q}$: The distance set of $\s _{\Q}$ is too rich and allows to create a dividing coloring. A natural attempt at that point is consequently to replace $\s _{\Q}$ by another space with a simpler distance set but still allowing to approximate $\s$ in a reasonable sense. There are natural candidates for this position, namely, the spaces obtained from the $\Ur _m$'s after having rescaled the distances in $[0,1]$. In the sequel, these spaces will be denoted $\s _m$'s. Formally, for $m \in \omega$ strictly positive, if $\Ur _m = (U_m , d^{\Ur _m})$, then \index{$\s _m$}

\[
\s _m = (U_m , \frac{d^{\Ur _m}}{m}). 
\]

This subsection is organized as follows: In \ref{subsubsection:From indivisibility of S_m to  oscillation stability of S.}, we show how to derive approximate indivisiblity of $\s$ from indivisibility of the $\s _m$'s. This proves: 

\begin{thm}  

\label{thm:s mos}

The Urysohn sphere $\s$ is approximately indivisible (equivalently, the standard action of $\iso(\s)$ on $\s$ is oscillation stable). 

\end{thm}

We then show (see \ref{subsubsection:From  oscillation stability of S to approximate indivisibility of S_Q.}): 

\begin{thm}

\label{cor:TFAE i-ii}

The rational Urysohn sphere $\s _{\Q}$ is approximately indivisible. 

\end{thm}

Theorem \ref{thm:s mos} exhibits an essential Ramsey-theoretic distinction between $\mathbb{S}^{\infty}$ and $\s$. At the level of $\iso (\mathbb{S}^{\infty})$ and $\iso (\s)$, it answers a question mentioned by Kechris, Pestov and Todorcevic in \cite{KPT}, Hjorth in \cite{Hj} and Pestov in \cite{Pe1}, \cite{Pe1'}, and highlights a deep topological difference which, for the reasons mentioned previously, was not at all apparent until now. 

Before going deeper into the technical details, let us mention here that part of our hope towards the discretization strategy came from the proof of a famous result in Banach space theory, namely Gowers' stabilization theorem for $c_0$. Recall that $c_0$ is the space of all real sequences converging to $0$ equipped with the $\left\| \cdot \right\|_{\infty}$ norm. Let $\mathbb{S}_{c_0}$ denote its unit sphere and $\mathbb{S}_{c_0} ^+$ denote the set of all those elements of $\mathbb{S}_{c_0}$ taking only positive values. In \cite{Gow}, Gowers studied the indivisibility properties of the spaces $\mathrm{FIN}_m$ (resp. $\mathrm{FIN}_m ^+$) of all the elements of $\mathbb{S}_{c_0}$ taking only values in $\{ k/m : k \in [-m,m] \cap \mathbb{Z} \}$ (resp. $\{ k/m : k \in \{0, 1,\ldots , m \} \}$) where $m$ ranges over the strictly positive integers:

\begin{thm}[Gowers \cite{Gow}] 

Let $m \in \omega$, $m \geqslant 1$. Then $\mathrm{FIN}_m$ (resp. $\mathrm{FIN}_m ^+$) is $1$-indivisible (resp. indivisible).

\end{thm}

A strong form of these results (see \cite{Gow} for the precise statement) then led to:

\begin{thm}[Gowers \cite{Gow}] 
The sphere $\mathbb{S}_{c_0}$ (resp. $\mathbb{S}_{c_0} ^+$) is approximately indivisible. 
\end{thm}

In the present case, Theorem \ref{thm:s mos} actually provides several other results of a similar flavor. For example, it allows to reach the following generalization:

\begin{thm}

\label{thm:univ indiv}

Let $\m{X}$ be a separable metric space with finite diameter $\delta$. Assume that every separable metric space with diameter less or equal to $\delta$ embeds isometrically into $\m{X}$. Then $\m{X}$ is approximately indivisible.  

\end{thm} 

Then, notice that when applied to the unit sphere of certain remarkable Banach spaces, this theorem yields interesting consequences. For example, it is known that every separable metric spaces with diameter less or equal to $2$ embeds isometrically into the unit sphere $\mathbb{S}_{\mathcal{C}([0,1])}$ of the Banach space $\mathcal{C}([0,1])$. It follows that:

\begin{thm}

\label{thm:C[0,1] approx indiv}

The unit sphere of $\mathcal{C}([0,1])$ is approximately indivisible. 

\end{thm}

On the other hand, it is also known that $\mathcal{C}([0,1])$ is not the only space having a unit sphere satisfying the hypotheses of Theorem \ref{thm:univ indiv}. For example, Holmes proved in \cite{H} there is a Banach space $\left\langle \Ur \right\rangle$ such that for every isometry $i : \funct{\Ur}{\m{Y}}$ of the Urysohn space $\Ur$ into a Banach space $\m{Y}$ such that $0_{\m{Y}}$ is in the range of $i$, there is an isometric isomorphism between $\left\langle \Ur \right\rangle$ and the closed linear span of $i''\Ur$ in $\m{Y}$. Very little is known about the space $\left\langle \Ur \right\rangle$, but it is easy to see that its unit sphere embeds isometrically every separable metric space with diameter less or equal to $2$. Therefore:  

\begin{thm}

\label{thm:<U> approx indiv}

The unit sphere of the Holmes space is approximately indivisible. 

\end{thm}

Observe that these result do \emph{not} say that for $\m{X} = \mathcal{C}([0,1])$ or $\left\langle \Ur \right\rangle$, every finite partition $\gamma$ of the unit sphere $\mathbb{S}_{\m{X}}$ of $\m{X}$ and every $\varepsilon > 0$, there is $\Gamma \in \gamma$ and a closed infinite dimensional subspace $\m{Y}$ of $\m{X}$ such that $\mathbb{S}_{\m{X}} \cap \m{Y} \subset (\Gamma)_{\varepsilon}$: According to the classical results about oscillation stability in Banach spaces, this latter fact is false for those Banach spaces into which every separable Banach space embeds linearly, and it is known that both $\mathcal{C}([0,1])$ and $\left\langle \Ur \right\rangle$ have this property. 

On the other hand, these results do not say either that for $\m{X} = \mathcal{C}([0,1])$ or $\left\langle \Ur \right\rangle$ the standard action of the surjective isometry group of the unit sphere of $\m{X}$ on the unit sphere of $\m{X}$ is oscillation stable. Indeed, since the unit sphere of $\m{X}$ is not ultrahomogeneous, the left completion of its surjective isometry group is not the entire semigroup of all isometric embeddings. Therefore, it might very well be that when a finite coloring of those spheres is given, the embedding which provides an almost monochromatic copy is not in the left completion of the surjective isometry group. To draw a parallel with Gowers' theorems mentioned previously, this is exactly what happens in the case of the unit sphere of $c_0$.    

\

\subsubsection{From indivisibility of $\s _m$ to oscillation stability of $\s$.}

\label{subsubsection:From indivisibility of S_m to oscillation stability of S.}

In this section, we show how oscillation stability of $\s$ follows from indivisibility of the spaces $\s _m$. This proof was obtained in collaboration with Jordi Lopez-Abad, and follows the lines of \cite{LANVT}.

\begin{prop}
\label{cor:if s_m indiv then s 1/2m-indiv} Let $m \in \omega$ be strictly positive. Then $\s$ is $1/m$-indivisible.
\end{prop}

\begin{proof}

This is obtained by showing that for every strictly positive  $m \in \omega$, there is an isometric
copy $\s _m ^*$ of $\s _m$ inside $\s$ such that for every $\mc{S} _m \subset \s _m ^* $ isometric
to $\s _m$, $(\mc{S} _m)_{1/m}$ includes an isometric copy of $\s$. This property indeed suffices
to prove Proposition \ref{cor:if s_m indiv then s 1/2m-indiv}: Let $\chi : \funct{\s}{k}$ for some
strictly positive $k \in \omega$. $\chi$ induces a $k$-coloring of the copy $\s _m ^*$. By indivisibility of $\s _m$, find $i<k$ and $\mc{S} _m \subset \s _m ^*$
such that $\chi$ is constant on $\mc{S} _m$ with value $i$. But then, in $\s$, $(\mc{S} _m)_{1/m}$
includes a copy of $\s$. So $(\overleftarrow{\chi} \{ i \})_{1/m}$ includes a copy of $\s$.

We now turn to the construction of $\s _m ^*$. The core of the proof is contained in Lemma
\ref{lem:hedgehog} which we present now. For $m \in \omega$ strictly positive, set \index{$[0,1] _m$}

\begin{center}
$[0,1] _m := \{ k/m : k \in \{ 0, \ldots , m \} \}$.
\end{center}

On the other hand, for $\alpha \in [0,1]$, set \index{$\left\lceil \cdot \right\rceil _m$}

\begin{center}
$\left\lceil \alpha \right\rceil _m = \min [\alpha ,1] \cap [0,1] _m$.
\end{center}

Fix an enumeration $\{y_n : n \in \omega \}$ of $\s
_{\Q}$. Also, let $\m{X} _m$ be the metric space $(\s _{\Q} , \left\lceil d^{\s _{\Q}} \right\rceil _m)$.
The underlying set of $\m{X} _m$ is really $\{y_n : n \in \omega \}$ but to avoid confusion, we
write it $\{x_n : n \in \omega \}$, being understood that for every $n \in \omega$, $x_n = y_n$. On
the other hand, observe that $\s _m$ and $\m{X} _m$ embed isometrically into each other.

\begin{lemma}
\label{lem:hedgehog} There is a countable metric space $\m{Z}$ with distances in $[0,1]$ and
including $\m{X} _m$ such that for every strictly increasing $\sigma : \funct{\omega}{\omega}$
 such that  $x_n \mapsto x_{\sigma (n)}$ is an isometry,   $ (\{ x_{\sigma (n)} : n \in \omega \})_{1/m}$
 includes an isometric copy of $\s _{\Q}$.
\end{lemma}
Assuming Lemma \ref{lem:hedgehog}, we now show how we can construct $\s ^* _m$. $\m{Z}$ is
countable with distances in $[0,1]$ so we may assume that it is a subspace of $\s$. Now, take $\s
_m ^*$ a subspace of $\m{X}_m$ and isometric to $\s _m$. We claim that $\s _m ^* $ works: Let
$\mc{S} _m \subset \s _m ^* $ be isometric to $\s _m$. We first show that $(\mc{S}
_m)_{1/m}$ includes a copy of $\s _{\Q}$. The enumeration $\{x_n : n \in \omega \}$
induces a linear ordering $<$ of $\mc{S} _m$ in type $\omega$. According to lemma \ref{lem:hedgehog}, it suffices to
show that $(\mc{S} _m,<)$ includes a copy of $\{x_n : n \in \omega \}_<$ seen as an ordered metric
space. To do that, observe that since $\m{X} _m$ embeds isometrically into $\s _m$, there is a
linear ordering $<^*$ of $\s _m$ in type $\omega$ such that $\{x_n : n \in \omega \}_<$ embeds into
$(\s _m , <^*)$ as ordered metric space. Therefore, it is enough to show:
\begin{claim}
$(\mc{S} _m,<)$ includes a copy of $(\s _m , <^*)$.
\end{claim}
\begin{proof}
Write
\begin{align*}
(\s _m , <^*)& = \{s_n : n \in \omega \}_{<^*}\\
(\mc{S} _m,<) &  = \{t_n : n \in \omega \}_{<}.
\end{align*}

Let $\sigma (0) = 0$. If $\sigma (0) < \dots < \sigma (n)$ are chosen such that $s_k \mapsto
t_{\sigma (k)}$ is a finite isometry, observe that the following set is infinite
\[\{ i \in \omega : \forall k \leqslant n \ \ d^{\s _m}(t_{\sigma (k)} , t_i) = d^{\s _m}(s_k ,
s_{n+1})\}.\]

Therefore, simply take $\sigma (n+1) $ in that set and larger than $\sigma (n)$.
\end{proof} 

Observe that since the metric completion of $\s _{\Q}$ is $\s$,
the closure of $(\mc{S} _m)_{1/m}$ in $\s$ includes a  copy of $\s$. But $(\mc{S} _m)_{1/m}$ is closed in $\s$, so $(\mc{S} _m)_{1/m}$ includes a copy of $\s$, and we are done. \end{proof}

We now turn to the proof of lemma \ref{lem:hedgehog}. Intuitively, here is the idea: First, construct a metric space $\m{Y}_m$ defined on the set $\s _{\Q}\times \{0,1 \}$ and where the metric $d^{\m{Y}_m}$ satisfies, for every $x, y \in \s _{\Q}$:

\begin{enumerate}
	\item[i)] $d^{\m{Y}_m}((x,1),(y,1))=d^{\s _{\Q}}(x,y)$,
	\item[ii)] $d^{\m{Y}_m}((x,0),(y,0)) = \left\lceil d^{\s _{\Q}}(x,y)\right\rceil_m$,
	\item[iii)] $d^{\m{Y}_m}((x,0),(x,1))=1/m$.
\end{enumerate}

The space $\m{Y}_m$ is really a two-level metric space with a lower level isometric to $\m{X}_m$. Note that in $\m{Y}_m$, $(\m{X}_m)_{1/m}$ includes a copy of $\s _{\Q}$. So the basic idea to construct $\m{Z}$ is to start from $\m{X}_m$ and to use some kind of gluing technique to glue a copy of $\m{Y}_m$ on $\m{X}_m$ along $\mc{X}_m$ whenever $\mc{X}_m$ is a copy of $\m{X}_m$ inside $\m{X}_m$. This process adds a copy of $\s _{\Q}$ inside $(\mc{X}_m)_{1/m}$ whenever $\mc{X}_m \subset \m{X}_m$ is isometric to $\m{X}_m$. There is, however, a delicate part. Namely, the gluing process has to be performed in such a way that $\m{Z}$ is separable. For example, this restriction forbids the brutal use of strong amalgamation, because then, we would go from $\m{X}_m$ to $\m{Z}$ by adding continuum many copies of $\s _{\Q}$ that are pairwise disjoint and at least $1/m$ apart. In spirit, the way this issue is solved is by allowing the different copies of $\s _{\Q}$ we are adding to intersect using some kind of tree-like pattern on the set of copies $\mc{X}_m$ inside $\m{X}_m$. The purpose of what follows is to describe precisely how this can be achieved. We first construct the set $Z$ on which the metric space $\m{Z}$ is supposed to be based, and then argue that the distance $d^{\m{Z}}$ can be obtained as required (Lemmas \ref{lem:3} to \ref{lem:6}). To construct $Z$, proceed as follows: For $t \subset \omega$, write $t$ as the strictly increasing enumeration of its
elements:

\begin{center}
$t = \{t_i : i \in |t| \}_<$.
\end{center}

Now, let $T$ be the set of all finite nonempty subsets $t$ of $\omega$ such that $x_n \mapsto
x_{t_n}$ is an isometry between $\{ x_n : n \in |t| \}$ and $\{x_{t_n} : n \in |t|\}$. This
set $T$ is a tree when ordered by end-extension. Let

\begin{center}
$Z = X_m \overset{.}{\cup} T$.
\end{center}

For $z \in Z$, define
\begin{displaymath}
\pi(z) = \left \{ \begin{array}{cl}
 z & \textrm{if $z \in X_m$,} \\
 x_{\max z } & \textrm{if $z \in T$.}
 \end{array} \right.
\end{displaymath}

Now, consider an edge-labelled graph structure on $Z$ by defining   $\delta$ with domain
$\mathrm{dom} (\delta) \subset Z \times Z $  and range included in $[0,1]$ as follows:
\begin{itemize}
\item If $s, t \in T$, then $(s,t) \in \mathrm{dom}(\delta)$ iff $s$ and $t$ are $<_T$ comparable. In this case (recall that $\{y_n : n \in \omega \}$ is an enumeration of $\s_{\Q}$),
\[\delta (s,t) = d^{\s _{\Q}} (y_{|s|-1 }, y_{|t|-1}).\]
\item If $x, y \in X_m$, then $(x,y)$ is always in $\mathrm{dom}(\delta)$ and
\[\delta (x,y) = d^{\m{X}_m} (x, y).\]
\item If $t \in T$ and $x \in X_m$, then $(x,t)$ and $(t,x)$ are in $\mathrm{dom}(\delta)$ iff $x = \pi (t)$. In this case
\[\delta (x,t) = \delta (t,x) = \frac1{m}.\]
\end{itemize}

For a branch $b$ of $T$ and $i \in \omega$, let $b(i)$ be the unique element of $b$ with height $i$
in $T$. Observe that $b(i)$ is an $(i+1)$-element subset of $\omega$. So: 

\begin{enumerate}
	 \item[i)] $\delta (b(i), b(j)) = d^{\s _{\Q}} (y_{|b(i)|-1} , y_{|b(j)|-1}) = d^{\s _{\Q}} (y_{i+1-1} , y_{j+1-1}) = d^{\s _{\Q}} (y_i , y_j)$. 
\end{enumerate}
 
Observe also that:
\begin{enumerate}
\item[ii)] $\delta(\pi(b(i)),\pi(b(j)))$ is equal to any of the following quantities: 

$d^{\m{X}_m}(x_{\max b(i)},x_{\max b(j)})= d^{\m{X}_m}(x_i,x_j)=  \lceil d^{\s_Q}(y_i,y_j) \rceil_m$,

\item[iii)] $\delta (b(i), \pi (b(i)) = 1/{m}$.

\end{enumerate}

The subspace $b\cup \pi''b$ will really play the role of the space $\m{Y}_m$ we mentioned previously. In particular, if $b$ is a branch of $T$, then $\delta$ induces a metric on $b$ and the map from
$\s _{\Q}$ to $b$ mapping  $y_i$ to $b(i)$ is a surjective isometry. We claim that if we can show
that $\delta$ can be extended to a metric $d^{\m{Z}}$ on $Z$ with distances in $[0,1]$, then Lemma
\ref{lem:hedgehog} will be proved. Indeed, let
\[\mc{X} _m = \{ x_{\sigma (n)} : n
\in \omega \} \subset \m{X} _m,\] with $\sigma : \funct{\omega}{\omega}$ strictly increasing and
$x_n \mapsto x_{\sigma (n)}$ distance preserving. See the range of $\sigma$ as a branch $b$ of $T$. Then
$(b, d^{\m{Z}}) = (b, \delta)$ is isometric to $\s _{\Q}$ and
\[b \subset (\pi '' b )_{ 1/{m}} = (\mc{X} _m)_{1/{m}} .\]

Our goal now is consequently to show that $\delta$ can be extended to a metric on $Z$ with values
in $[0,1]$. Recall that for $x, y \in Z$, and $n \in \omega$ strictly positive, a path from $x$ to
$y$ of size $n$ as is a finite sequence $\gamma = (z_i)_{i<n}$ such that $z_0 = x$, $z_{n-1} = y$
and for every $i<n-1$,
\[(z_i, z_{i+1}) \in \dom(\delta).\]

For $x, y$ in $Z$, $P(x,y)$ is the set of all paths from $x$ to $y$. If $\gamma = (z_i)_{i<n}$ is
in $P(x,y)$, $ \| \gamma \|$ is defined as:
\[ \| \gamma \| = \sum _{i=0} ^{n-1} \delta (z_i , z_{i+1} ).\]

We are going to show that for every  $(x,y) \in \dom (\delta)$, every path $\gamma$
from $x$ to $y$ is metric, that is:
\begin{equation}
\label{hojthurhgr}\delta (x,y) \leqslant \| \gamma \|
\end{equation}

This will prove that the required metric can be obtained by setting \[ d^{\m{Z}}(x,y) = \min (1, \inf \{ \| \gamma \| _{\leqslant 1} : \gamma \in P(x,y)\}).\]

Let $x, y \in Z$. Call a path $\gamma$ from $x$ to $y$ \emph{trivial} when $\gamma = (x,y)$ and
\emph{irreducible} when no proper subsequence of $\gamma$ is a non-trivial path from $x$ to $y$.
Finally, say that $\gamma$ is a \emph{cycle} when $(x,y) \in \dom (\delta)$. It should be clear
that to prove that $d^{\m{Z}}$ works, it is enough to show that the previous inequality
\eqref{hojthurhgr} is true for every irreducible cycle. Note that even though $\delta$ takes only
rational values, it might not be the case for $d^{\m{Z}}$. We now turn to the study of the
irreducible cycles in $Z$.

\begin{lemma}

\label{lem:3}

Let $x, y \in T$. Assume that $x$ and $y$ are not $<_T$-comparable. Let $\gamma$ be an irreducible
path from $x$ to $y$ in $T$. Then there is $z \in T$ such that $z <_T x$, $z <_T y$ and $\gamma =
(x,z,y)$.

\end{lemma}

\begin{proof}

Write $\gamma = (z_i)_{i<n+1}$. $z_1$ is connected to $x$ so $z_1$ is $<_T$-comparable with $x$.
We claim that $z_1 <_T x$ : Otherwise, $x <_T z_1$ and every element of $T$ which is
$<_T$-comparable with $z_1$ is also $<_T$-comparable with $x$. In particular, $z_2$ is
$<_T$-comparable with $x$, a contradiction since $z_2$ and $x$ are not connected. We now claim that
$z_1 <_T y$. Indeed, observe that $z_1 <_T z_2$ : Otherwise, $z_2 <_T z_1 <_T x$ so $z_2 <_T x$
contradicting irreducibility. Now, every element of $T$ which is $<_T$-comparable with $z_2$ is
also $<_T$-comparable with $z_1$, so no further element can be added to the path. Hence $z_2 = y$
and we can take $z_1 = z$. \end{proof}

\begin{lemma}

\label{lem:4}

Every non-trivial irreducible cycle in $X_m$ has size $3$.

\end{lemma}

\begin{proof}

Obvious since $\delta $ induces the metric $d^{\m{X}_m}$ on $X_m$. \end{proof}

\begin{lemma}

\label{lem:4'}

Every non-trivial irreducible cycle in $T$ has size $3$ and is included in a branch.

\end{lemma}

\begin{proof}

Let $c = (z_i)_{i<n}$ be a non-trivial irreducible cycle in $T$. We may assume that $z_0 <_T
z_{n-1}$. Now,  observe that every element of $T$ comparable with $z_0$ is also comparable with
$z_{n-1}$. In particular, $z_1$ is such an element. It follows that $n = 3$ and that $z_0, z_1,
z_2$ are in a same branch. \end{proof}

\begin{lemma}

\label{lem:5}

Every irreducible cycle in $Z$ intersecting both $X_m$ and $T$ is supported by a set whose form is one of the following ones:
\begin{center}
\begin{figure}[h]
\includegraphics[scale=0.73]{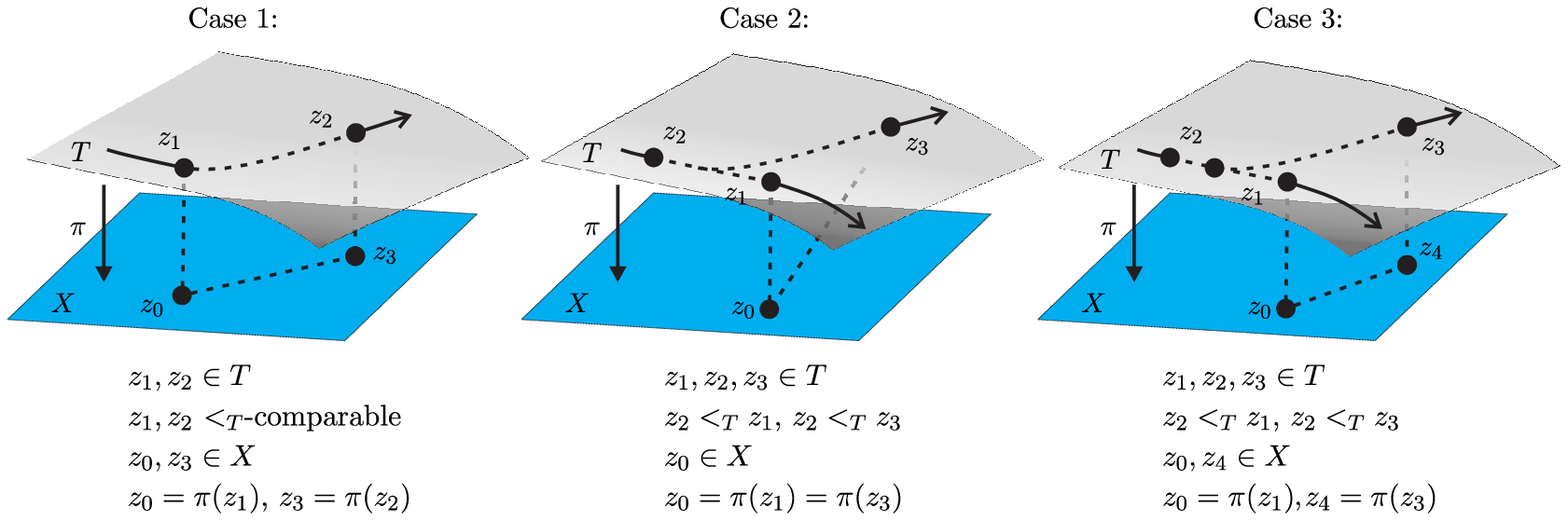}
\caption{Irreducible cycles}\label{figu1}
\end{figure}
\end{center}
%

\end{lemma}

\begin{proof}

Let $C$ be a set supporting an irreducible cycle $c$ intersecting both $X_m$ and $T$. It should be
clear that $|C \cap X_m| \leqslant 2$: Otherwise since any two points in $X_m$ are connected, $c$ would admit a strict subcycle, contradicting irreducibility.

If $C \cap X_m$ has size $1$, let $z_0$ be its unique element. In $c$, $z_0$ is connected to two
elements which we denote $z_1$ and $z_3$. Note that $z_1, z_3 \in T$ so $\pi (z_1) = \pi (z_3) =
z_0$. Since elements in $T$ which are connected never project on a same point, it follows that
$z_1, z_3$ are $<_T$-incomparable. Now, $c$ induces an irreducible path from $z_1$ to $z_3$ in $T$
so from lemma \ref{lem:3}, there is $z_2 \in C$ such that $z_2 <_T z_1$, $z_2 <_T z_3$, and we are
in case 2.

Assume now that $C \cap X_m = \{ z_0 , z_4 \}$. Then there are $z_1, z_3 \in C \cap T$ such that
$\pi(z_1) = z_0$  and $\pi(z_3) = z_4$. Note that since $z_0 \neq z_4$, we must have $z_1 \neq
z_3$. Now, $C \cap T$ induces an irreducible path from $z_1$ to $z_3$ in $T$. By lemma \ref{lem:3},
either $z_1$ and $z_3$ are compatible and in this case, we are in case 1, or $z_1$ and $z_3$ are
$<_T$-incomparable and there is $z_2$ in $C\cap T$ such that $z_2 <_T z_1$, $z_2 <_T z_3$ and we
are in case 3. \end{proof}

\begin{lemma}

\label{lem:6}

Every non-trivial irreducible cycle in $Z$ is metric.

\end{lemma}

\begin{proof}
Let $c$ be an irreducible cycle in $Z$. If $c$ is supported by $X_m$, then by lemma \ref{lem:4} $c$
has size $3$ and  is metric since $\delta$ induces a metric on $X_m$. If $c$ is supported by $T$,
then by lemma \ref{lem:4'} $c$ also has size $3$ and is included in a branch $b$ of $T$. Since
$\delta$ induces a metric on $b$, $c$ is metric. We consequently assume that $c$ intersects both
$X_m$ and $T$. According to lemma \ref{lem:5}, $c$ is supported by a set $C$ whose form is covered
by one of the cases 1, 2 or 3. So to prove the present lemma, it is enough to show every cycle
obtained from a re-indexing of the cycles described in those cases is metric.

Case 1: The required inequalities are obvious after having observed that \[\delta (z_0 , z_3) = \left\lceil \delta (z_1, z_2) \right\rceil _m  \text{ and }  \delta
(z_0 , z_1) = \delta (z_2 , z_3) = \frac1{m}.\]

Case 2: Notice that $\delta (z_0 , z_1) = \delta (z_0 , z_3) = 1/m$. So the inequalities we need to prove are
\begin{align}
\delta (z_1 , z_2) & \leqslant   \delta (z_2, z_3) + \frac2m, \label{lulu1}\\
\delta (z_2 , z_3) & \leqslant  \delta (z_1, z_2) + \frac2m. \label{lulu2}
\end{align}

By symmetry, it suffices to verify that \eqref{lulu1} holds. Observe that since $\pi (z_1) = \pi (z_3) = z_0$, we must have $\left\lceil \delta (z_1 , z_2) \right\rceil_m=\left\lceil \delta(z_2 , z_3)\right\rceil_m$. So: 

\[
\delta (z_1 , z_2) \leqslant  \left\lceil \delta (z_1 , z_2) \right\rceil _m = \left\lceil \delta (z_2 , z_3) \right\rceil _m \leqslant  \delta (z_2 , z_3) + \frac2m.\]

Case 3: Observe that $\delta (z_0 , z_1) = \delta (z_3 , z_4) = 1/m$, so the inequalities we need
to prove are
\begin{align}
\delta (z_1 , z_2) & \leqslant   \delta (z_2, z_3) + \delta (z_0 , z_4) + \frac2m, \label{ohrtjuer1}\\
\delta (z_0 , z_4) & \leqslant  \delta (z_1, z_2) + \delta (z_2 , z_3) + \frac2m. \label{ohrtjuer2}
\end{align}

For \eqref{ohrtjuer1}:
\begin{align*}
\delta (z_1 , z_2) & \leqslant  \left\lceil \delta (z_1 , z_2) \right\rceil _m \\
& =  \delta (\pi (z_1) , \pi (z_2) ) \\
& =  \delta (z_0 , \pi (z_2) ) \\
& \leqslant  \delta (z_0 , z_4) + \delta (z_4 , \pi (z_2)) \\
& =  \delta (z_0 , z_4) + \left\lceil \delta (z_3 , z_2) \right\rceil _m \\
& \leqslant  \delta (z_0 , z_4) + \delta (z_2 , z_3) + \frac2m.
\end{align*}

For \eqref{ohrtjuer2}: Write $z_1 = b(j)$, $z_3 = b'(k)$, $z_2 = b(i) = b'(i)$. Then $z_0 = \pi
(z_1) = x_{\max b(j)}$ and $z_4 = \pi (z_3) = x_{\max b'(k)}$. Observe also that $\delta (z_1 ,
z_2) = d^{\s _{\Q}}(y_j , y_i)$ and that $\delta (z_2 , z_3) = d^{\s _{\Q}}(y_i , y_k)$. So:
\begin{align*}
\delta (z_0 , z_4) & =  d^{\m{X}_m}(x_{\max b(j)}, x_{\max b'(k)})\\
& \leqslant  d^{\m{X}_m}(x_{\max b(j)}, x_{\max b(i)}) + d^{\m{X}_m}(x_{\max b'(i)}, x_{\max b'(k)}) \\
& = d^{\m{X}_m}(x_j, x_i) + d^{\m{X}_m}(x_i, x_k) \\
& =  \left\lceil d^{\s _{\Q}}(y_j, y_i)\right\rceil _m + \left\lceil d^{\s _{\Q}}(y_i, y_k) \right\rceil _m \\
& = \left\lceil \delta (z_1, z_2) \right\rceil _m + \left\lceil \delta (z_2 , z_3) \right\rceil _m \\
& \leqslant  \delta (z_1, z_2) + \frac1m + \delta (z_2 , z_3) + \frac1m \\
& = \delta (z_1, z_2) + \delta (z_2 , z_3) + \frac2m. \qedhere
\end{align*}
\end{proof} 

\

\subsubsection{From oscillation stability of $\s$ to approximate indivisibility of $\s _{\Q}$.}
\label{subsubsection:From  oscillation stability of S to approximate indivisibility of S_Q.}

The purpose of what follows is to prove that the rational Urysohn sphere is approximately indivisible (Theorem \ref{cor:TFAE i-ii}). We start with the following proposition. 

\begin{prop}

\label{thm: s_Q in s}

Suppose that $\s_{\Q}^0$ and $\s_{\Q}^1$ are two copies of $\s_{\Q}$ in $\s$ such that $\s_{\Q}^0$ is dense in $\s$. Then for every $\varepsilon>0$ the  subspace $\s_{\Q}^0\cap (\s_{\Q}^1)_{\varepsilon}$ includes a copy of $\s_\Q$.
\end{prop}

\begin{proof}
We construct the required copy of $\s _{\Q}$ inductively. Let $\{ y_n : n \in \omega \}$ enumerate
 $\s _{\Q}^1$.  For $k \in \omega$, set \[\delta _k = \frac{\varepsilon}{2} \sum
_{i = 0} ^k \frac{1}{2^i}.\] 

Set also
 \[\eta _k = \frac{\varepsilon}{3} \frac{1}{2^{k+1}}.\]
 
$\s_{\Q}^0$ being dense in $\s$, choose $z_0 \in \s _{\Q}^0$ such that $d^{\s}(y_0 , z_0) < \delta _0$.
Assume now that $z_0,\dots, z_n \in \s _{\Q}^0$ were constructed such that for every $k, l \leqslant
n$
\begin{displaymath}
\left \{ \begin{array}{l}
 d^{\s}(z_k , z_l)=d^{\s} (y_k , y_l), \\
 d^{\s} (z_k , y_k) < \delta _k.
 \end{array} \right.
\end{displaymath}

Again by denseness of $\s_{\Q}^0$ in $\s$, fix $z \in \s _{\Q}^0$ such that

\begin{center}
$d^{\s} (z,y_{n+1}) < \eta _{n+1}$.
\end{center}

Then for every $k \leqslant n$,
\begin{align*}
\left| d^{\s} (z , z_k) - d^{\s} (y_{n+1} , y_k) \right| & =  \left| d^{\s} (z , z_k) - d^{\s}(z_k
, y_{n+1})  + d^{\s}(z_k , y_{n+1})
 - d^{\s} (y_{n+1} , y_k) \right| \\
 & \leqslant  d^{\s} (z , y_{n+1}) + d^{\s} (z_k , y_k) \\
 & <  \eta _{n+1} + \delta _k \\
 & <  \eta _{n+1} + \delta_n.
\end{align*}

It follows that there is $z_{n+1} \in \s _{\Q}^0$ such that
\begin{displaymath}
\left \{ \begin{array}{l}
 \forall k \leqslant n \ \ d^{\s} (z_{n+1} , z_{k}) = d^{\s} (y_{n+1} , y_k) \\
 d^{\s} (z_{n+1} , z) < \eta _{n+1} + \delta _n.
 \end{array} \right.
\end{displaymath}

Indeed, consider the map $f$ defined on $\{ z_k : k \leqslant n \} \cup \{ z \}$ by:
\begin{displaymath}
\left \{ \begin{array}{l}
 \forall k \leqslant n \ \ f(z_k) = d^{\s} (y_{n+1} , y_k), \\
 f(z) =  \max\{\left| d^{\s} (z , z_k) - d^{\s} (y_{n+1} , y_k) \right| : k \leqslant n \}.
 \end{array} \right.
\end{displaymath}

\begin{claim}
$f$ is Kat\v{e}tov.
\end{claim}

\begin{proof}

The metric space $\{y_k: k \leqslant n+1 \}$ witnesses that $f$ is Kat\v{e}tov over $\{ z_k : k \leqslant n \}$ so it suffices to prove that for every $k\leqslant n$, \[ \left| f(z) - f(z_k) \right| \leqslant d^{\s}(z,z_k) \leqslant f(z) + f(z_k).\] 

Equivalently, for every $k\leqslant n$, \[ \left| d^{\s}(z,z_k) - f(z_k) \right| \leqslant f(z) \leqslant d^{\s}(z,z_k) + f(z_k).\] 

There is nothing to do for the left-hand side because by definition of $f$, we have \[ f(z) =  \max\{\left| d^{\s} (z , z_k) - f(z_k) \right| : k \leqslant n \}.\] 

For right-hand side, what we need to show is that for every $k,l\leqslant n$, \[ \left|d^{\s} (z , z_l) - d^{\s} (y_{n+1} , y_l)\right| \leqslant d^{\s} (z , z_k) + d^{\s} (y_{n+1} , y_k).\]

Equivalently, 

\begin{displaymath}
\left \{ \begin{array}{l}
 d^{\s} (z , z_l) - d^{\s} (y_{n+1} , y_l) \leqslant d^{\s} (z , z_k) + d^{\s} (y_{n+1} , y_k), \\
 d^{\s} (y_{n+1} , y_l) - d^{\s} (z , z_l) \leqslant d^{\s} (z , z_k) + d^{\s} (y_{n+1} , y_k).
 \end{array} \right.
\end{displaymath}

The first inequality is equivalent to \[ d^{\s} (z , z_l) - d^{\s} (z , z_k) \leqslant d^{\s} (y_{n+1} , y_k) + d^{\s} (y_{n+1} , y_l).\]

But this is satisfied because \[ d^{\s} (z , z_l) - d^{\s} (z , z_k) \leqslant d^{\s}(z_l , z_k) = d^{\s}(y_k,y_l) \leqslant d^{\s} (y_k , y_{n+1}) + d^{\s} (y_{n+1} , y_l).\]

Similarly, the second inequality is equivalent to \[ d^{\s} (y_{n+1} , y_l) - d^{\s} (y_{n+1} , y_k) \leqslant d^{\s} (z , z_k) + d^{\s} (z , z_l).\]

This holds because \[ d^{\s} (y_{n+1} , y_l) - d^{\s} (y_{n+1} , y_k) \leqslant d^{\s}(y_k,y_l) = d^{\s}(z_k,z_l) \leqslant d^{\s} (z , z_k) + d^{\s} (z , z_l). \qedhere\]

\end{proof}

The map $f$ being Kat\v{e}tov, consider a point $z_{n+1} \in \s _{\Q}^0$ realizing $f$ over the set $\{ z_k : k \leqslant n \} \cup \{ z \}$. Observe then that
\begin{eqnarray*}
d^{\s}(z_{n+1} , y_{n+1}) & \leqslant & d^{\s}(z_{n+1},z) + d^{\s}(z,y_{n+1})\\
& < & \eta _{n+1} + \delta _n + \eta _{n+1}\\
& < & \delta _{n+1}.
\end{eqnarray*}

After $\omega$ steps, we are left with $\{ z_n : n \in \omega \} \subset \s _{\Q}^0 \cap
(\s_\Q^1)_{\varepsilon}$ isometric to $\s _{\Q}$. \end{proof}

We now show how to deduce of Theorem \ref{cor:TFAE i-ii} from Proposition
\ref{thm: s_Q in s}:  Let $\varepsilon > 0$, $k \in \omega$ strictly positive and $\chi : \funct{\s
_{\Q}}{k}$. Then in $\s$, seeing $\s_{\Q}$ as a dense subspace:
\[ \s = \bigcup _{i<k} (\overleftarrow{\chi} \{ i \})_{\varepsilon / 2}.\]

By oscillation stability of $\s$, there is $i<k$ and a copy $\mc{S}$ of $\s$ included in $\s$ such
that
\begin{center}
$\mc{S} \subset ((\overleftarrow{\chi} \{ i \})_{\varepsilon/2})_{\varepsilon/4}$.
\end{center}

Since $\widetilde{\s}$ includes copies of $\s_\Q$, and since $\s_\Q$ is dense in $\s$, it follows by Proposition \ref{thm: s_Q in s} that there is a copy $\mc{S} _{\Q}$ of $\s _{\Q}$ in
$\s _{\Q}\cap(\mc{S})_{\varepsilon /4}$. Then in $\s _{\Q}$
\begin{center}
\mbox{ }\hfill$\mbox{ }\hfill \mc{S} _{\Q} \subset (\overleftarrow{\chi} \{ i \})_{\varepsilon}$.
\hspace{\stretch{1}} \qed
\end{center}

\section{Concluding remarks and open problems.}

We mentioned several times in this chapter that for the moment, not much is known as far as big Ramsey degrees are concerned, so this direction already provides a first axis of future research. In fact, this is not particular to metric spaces: Even at the more general level of structural Ramsey theory, very little is known. To our knowledge, apart from ultrametric spaces, the only cases where a complete analysis was carried out correspond essentially to finite linear orderings (Devlin, see section 11 of \cite{KPT} or \cite{T1}) and finite graphs (Laflamme-Sauer-Vuksanovic \cite{LSV}). We should also mention at that stage another recent general result, which is closely linked to Theorem \ref{thm:Hjorth}. In \cite{Hj}, Hjorth proved the following: Let $\mathcal{K}$ be a Fra\"iss\'e class with Fra\"iss\'e limit $\m{F}$ whose automorphism group is non-trivial. Let also $\m{X}$ be a finite substructure of $\m{F}$ with $|\m{X}|\geqslant 2$. Then the action of $\mathrm{Aut}(\m{F})$ on $\mathrm{Aut}(\m{F})/St_{\m{X}}$ (where $St_{\m{X}}$ denotes the pointwise stabilizer of $\m{X}$ in $\mathrm{Aut}(\m{F})$) is not oscillation stable. With respect to big Ramsey degrees, this result is relevant because it implies:

\begin{thm}[Hjorth \cite{Hj}]

\label{thm:Hjorth'}
\index{Hjorth!theorem on Big Ramsey degrees}

Let $\mathcal{K}$ be a Fra\"iss\'e class and $\m{X} \in \mathcal{K}$. Assume that $|\m{X}|\geqslant 2$ and that $\m{X}$ is rigid (ie has a trivial automorphism group). Then the big Ramsey degree of $\m{X}$ in $\mathcal{K}$ is, when defined, at least 2.   

\end{thm}
 
The rigidity hypothesis is really necessary here: If it is dropped, the usual infinite Ramsey theorem provides a counterexample. Note also that when $\mathcal{K}$ is a Fra\"iss\'e order class, every $\m{X}$ in $\mathcal{K}$ is rigid and therefore has a big Ramsey degree at least 2 whenever $|\m{X}|\geqslant 2$. No similar general result is known for upper bounds (or even existence) of big Ramsey degrees. Furthermore, even when big Ramsey degrees are determined, their explicit computation is not always easy. Ultrametric spaces are a good illustration of this phenomenon: For $\m{X} \in \U$, we proved that $T_{\U}(\m{X})$ is equal to the number of linear extensions of the tree associated to $\m{X}$ in $\U$\index{$\U$} but we did not touch the question of how this number can be computed in practice. For graphs, the problem is similar, and it turns out that even in the most simple cases, highly non-trivial combinatorial problems appear (see for example \cite{La}). For more about big Ramsey degrees in structural Ramsey theory, see \cite{KPT}, section 11, or \cite{T1}. Back to the metric context, here is the question which looks like the most reasonable to us:   

\

\textbf{Question 3.} Let $m \in \omega$ be strictly positive. Does every $\m{X}$ in $\M _{\omega \cap ]0,m]}$ have a big Ramsey degree in $\M _{\omega \cap ]0,m]}$? More generally, if $S \subset ]0,+\infty[$ is finite and satisfies the $4$-values condition, does every $\m{X}$ in $\M _S$\index{$\M _S$} have a big Ramsey degree in $\M _S$? 

\

When $\m{X}$ is the $1$-point metric space $\m{K}_1$, this question is closely related to indivisibility. However, as mentioned several times already in the body of this paper, our belief is not only that $\m{K}_1$ has a big Ramsey degree in the class $\M _S$ but that the related Urysohn spaces $\Ur _S$ are indivisible. We also saw that this belief is already supported by several results when extra assumptions are made about $S$ (see Theorem \ref{thm:U_m indiv} and Theorem \ref{thm:U_S indiv ext}), but that the general case remains open. Here is therefore the next question:

\

\textbf{Question 4.} If $S \subset ]0,+\infty[$ is finite and satisfies the $4$-values condition, is $\Ur _S$ indivisible?  

\

Our last question is related to the connection between the approximate indivisibility problems for the sphere $\mathbb{S}^{\infty}$ and the Urysohn sphere $\s$. We saw indeed that the numerous Ramsey-theoretic properties that those two spaces share potentially indicated that solving the approximate indivisibility problem for $\s$ would lead to a better understanding of the result of Odell and Schlumprecht according to which $\mathbb{S}^{\infty}$ is not approximately indivisible. However, we showed with Theorem \ref{thm:s mos} that at the level of approximate indivisibility, $\mathbb{S}^{\infty}$ and $\s$ behave differently. This comment leads to: 

\

\textbf{Question 5.} From a metric point of view, which distinction between $\mathbb{S}^{\infty}$ and $\s$ is responsible for the different behaviors regarding approximate indivisibility?  

\

In particular, where is it that the techniques involved in the proof of approximate indivisibility for $\s$ fail for $\mathbb{S}^{\infty}$? We are currently unable to fully answer that question but a first analysis of the problem suggests that whereas the space $\s$ is easily approximated by a sequence of countable ultrahomogeneous metric spaces with finitely distances (namely, the sequence $(\s _m)_{m \in \omega}$), it may not be the case for $\mathbb{S}^{\infty}$. Indeed, we saw in Chapter 1 as a consequence of Proposition \ref{prop:non amalg E omega} that the class $\ES _S$ (recall that $\ES _S$ is the class of all finite metric spaces $\m{X}$ with distances in $S$ and which embed isometrically into the unit sphere $\mathbb{S} ^{\infty}$ of $\ell _2$ with the property that $\{ 0_{\ell _2} \} \cup \m{X}$ is affinely independent) does not have the strong amalgamation property when $S = \{ k/m : k \in \{ 1,\ldots ,m\}\}$ with $m$ large enough. Therefore, there is no countable ultrahomogeneous metric subspace of $\mathbb{S} ^{\infty}$ whose class of finite metric subspaces is $\ES _S$. In other words, unlike what we did for $\s$, we cannot use the most obvious discretization method to approximate $\mathbb{S} ^{\infty}$ with a sequence of spaces with only finitely many distances and whose indivisibility behaviors reflect the behavior of $\mathbb{S} ^{\infty}$. But is there a deeper reason behind that fact? Could it be that there is no countable ultrahomogeneous metric space with finitely many distances whose divisibility captures the non approximate indivisibility of $\mathbb{S} ^{\infty}$? The exercise is left to the reader.



\backmatter





\chapter*{Appendix A. \ \  Amalgamation classes $\M _S$ when $|S| \leqslant 4$.}

The purpose of this appendix is to provide a list of all the amalgamation classes $\M _S$ when $|S| \leqslant 4$. Thanks to \cite{DLPS}, it is known that $\M _S$\index{$\M _S$} is an amalgamation class iff $S$ satisfies the $4$-values condition. Recall that $S$ satisfies the $4$-\emph{values condition} when for every $s_0 , s_1 , s_0 ', s_1 ' \in S$, if there is $t \in S$ such that:
\begin{center} 
$|s_0 - s_1| \leqslant t \leqslant s_0 + s_1$, \ \ $|s_0 '- s_1 '| \leqslant t \leqslant s_0 '+ s_1 '$,
\end{center}
then there is $u \in S$ such that:
\begin{center} 
$|s_0 - s_0 '| \leqslant u \leqslant s_0 + s_0 ' $, \ \ $|s_1 - s_1 '| \leqslant u \leqslant s_1 + s_1 '$.
\end{center}

\section{$|S|=3$.}

\subsection{$s_0 < s_1 \leqslant 2s_0 < s_0 + s_1 < 2s_1 < s_2 \ \ \{ 1, 2, 5\}$.}

\

For a quadruple $(u_0 , u_1 , u_2, u_3)$ of elements of $S$, let $I(u_0 , u_1 , u_2, u_3)$ be defined as the interval:

\begin{center}
$I(u_0 , u_1 , u_2, u_3) := [\max(|u_0 - u_1|, |u_2 - u_3|) , \min(u_0 + u_1, u_2 + u_3)] $
\end{center}

Call $(u_0 , u_1 , u_2, u_3)$ \emph{good} if $I(u_0 , u_1 , u_2, u_3) \cap S \neq \emptyset$. Otherwise, call it \emph{bad}. Define also $(u_0 , u_1 , u_2, u_3)^* :=(u_0 , u_2 , u_1 , u_3)$. So $S$ satisfies the $4$-values condition iff for every $(u_0 , u_1 , u_2, u_3) \in S^4$, $(u_0 , u_1 , u_2, u_3)$ is good iff $(u_0 , u_1 , u_2, u_3)^*$ is good. Also, call a permutation $\sigma$ of $ \{ 0, 1, 2, 3 \}$ \emph{trivial} if:

\begin{center}

$\forall (u_0 , u_1 , u_2, u_3) \in S^4, I(u_{\sigma (0)} , u_{\sigma (1)} , u_{\sigma (2)}, u_{\sigma (3)}) = I(u_0 , u_1 , u_2, u_3)$. 
\end{center}

Equivalently, $\sigma$ is trivial when $\sigma '' \{ 0, 1\} \in  \{ \{ 0, 1\} , \{ 2, 3\} \}$. Now, set:

\begin{center}
$A := \{ |s - s'| : s, s' \in S\}$ \ \ $B := \{ s + s' : s, s' \in S\}$.
\end{center}

Here, $A = \{ 1, 3, 4\}$, while $B = \{ 2, 3, 4\} \cup C$ with $C \subset [5 , +\infty[$. For every interval $[a , b]$ where $a \in A, b \in B \smallsetminus C$ and such that $[a, b] \cap S = \emptyset$, we find all the quadruples $(u_0 , u_1 , u_2, u_3)$ (up to trivial permutation) such that $I(u_0 , u_1 , u_2, u_3) = [a, b]$. Up to a trivial permutation, this allows to find all the bad quadruples. In the present case, here is the list of all intervals $[a , b]$ where $a \in A, b \in B$ and such that $[a, b] \cap S = \emptyset$, together with the quadruples $(u_0 , u_1 , u_2, u_3)$ such that $I(u_0 , u_1 , u_2, u_3) = [a, b]$. 

\[
\begin{array}{ll}

[3,2] & (2,5,1,1)\\

[3,3] & (2,5,1,2)\\

[3,4] & (2,5,2,2)\\

[4,2] & (1,5,1,1)\\

[4,3] & (1,5,1,2)\\

[4,4] & (1,5,2,2)
\end{array}
\]

Now, let $\tau$ be the transposition of $\{ 0, 1, 2, 3\}$ permuting $1$ and $2$. Let also $T$ be the set of all trivial permutations of $\{ 0, 1, 2, 3\}$. Observe that $T \cup \{ \tau \}$ generates the whole group of permutations of $\{ 0, 1, 2, 3\}$. Thus, we have to check that the set of bad quadruples is closed under all permutations. In practice, however, note that given any permutation $\sigma$ of $\{ 0, 1, 2, 3\}$, $(u_{\sigma (0)} , u_{\sigma (1)} , u_{\sigma (2)}, u_{\sigma (3)})$ is equal to $(u_0, u_1, u_2, u_3)$, to $(u_0, u_1, u_2, u_3) ^* = (u_0, u_2, u_1, u_3)$ or to $(u_0, u_1, u_2, u_3) _* = (u_0, u_3, u_2, u_1)$ up to trivial permutation. Thus, it suffices to show that for every bad quadruple $(u_0, u_1, u_2, u_3)$ above, $(u_0, u_1, u_2, u_3)^*$ and $(u_0, u_1, u_2, u_3)_*$ are also bad. Observe also
that there are some cases where checking only $(u_0, u_1, u_2, u_3)^*$ or $(u_0, u_1, u_2, u_3)_*$ is enough. For example, if $u_0 = u_1$, checking that $(u_0, u_2, u_1, u_3)^*$ is bad is sufficient. There are even cases where there is nothing to check, namely when all but one of the $u_i$'s are equal. Here, if $\approx$ denotes equality modulo a trivial permutation:

\begin{center}

$(2,5,1,1)^* = (2,1,5,1) \approx (1,5,1,2)$

$(2,5,1,2)_* = (2,2,1,5) \approx (1,5,2,2)$ 

$(1,5,1,2)^* = (1,1,5,2) \approx (2,5,1,1)$

$(1,5,2,2)^* = (1,2,5,2) \approx (1,5,1,2)$ 
\end{center}

It follows that $S$ satisfies the $4$-values condition. 

\subsection{$s_0 < 2s_0 < s_1 < s_2 \leqslant s_0 + s_1 < 2s_1 \ \ \{ 1, 3, 4\}$.}

\[
A = \{ 1, 2, 3\}, \ \ B = \{ 2\} \cup C, \ \ C \subset [4, +\infty[.
\]
\[
\begin{array}{lll}

[2,2] & (1, 3, 1, 1) & \\

[3,2] & (1, 4, 1, 1) &   

\end{array}
\]

$\{ 1, 3, 4\}$ satisfies the $4$-values condition. 

\subsection{$s_0 < 2s_0 < s_1 < s_0 + s_1 < s_2 \leqslant 2s_1 \ \ \{ 1, 3, 6\}$.} 

\[
A = \{ 2, 3, 5\}, \ \ B = \{ 2, 4\} \cup C, \ \ C \subset [6, +\infty[. 
\]
\[
\begin{array}{lll}
[2,2] & (1, 3, 1, 1) & \\

[3,2] & (3, 6, 1, 1) & (3,6,1,1)^* = (3,1,6,1) \approx (1,6,1,3)\\

[5,2] & (1, 6, 1, 1) & \\

[5,4] & (1, 6, 1, 3) & (1,6,1,3)^* = (1,1,6,3) \approx (3,6,1,1)

\end{array}
\]

$\{ 1, 3, 6\}$ satisfies the $4$-values condition.

\section{$|S|=4$.}

For $|S|=4$, there are more cases to consider. Recall that for $|S|=3$, the sets we had to check with the $4$-values criterion were provided by the following inequalities:

\begin{center}
(1a) $s_0 < s_1 < s_2 \leqslant 2s_0 < s_0 + s_1 < 2s_1$

(1b) $s_0 < s_1 \leqslant 2s_0 < s_2 \leqslant s_0 + s_1 < 2s_1$

(1d) $s_0 < s_1 \leqslant 2s_0 < s_0 + s_1 < 2s_1 < s_2$ 

\end{center}

\begin{center}

(2a) $s_0 < 2s_0 < s_1 < s_2 \leqslant s_0 + s_1 < 2s_1$

(2b) $s_0 < 2s_0 < s_1 < s_0 + s_1 < s_2 \leqslant 2s_1$

(2c) $s_0 < 2s_0 < s_1 < s_0 + s_1 < 2s_1 < s_2$ 

\end{center}

We look at how $s_0 + s_2$, $s_1 + s_2$ and $2s_2$ may be inserted in these chains:

\

For (1a): 

\begin{center}

$s_0 < s_1 < s_2 < 2s_0 < s_0 + s_1 < s_0 + s_2 < 2s_1 < s_1 + s_2 < 2s_2$ 

$s_0 < s_1 < s_2 < 2s_0 < s_0 + s_1 < 2s_1 < s_0 + s_2 < s_1 + s_2 < 2s_2$ 

\end{center}

\

For (1b):

\begin{center}

$s_0 < s_1 < 2s_0 < s_2 < s_0 + s_1 < s_0 + s_2 < 2s_1 < s_1 + s_2 < 2s_2$ 

$s_0 < s_1 < 2s_0 < s_2 < s_0 + s_1 < 2s_1 < s_0 + s_2 < s_1 + s_2 < 2s_2$ 

\end{center}

\

For (1d):

\begin{center}

$s_0 < s_1 < 2s_0 < s_0 + s_1 < 2s_1 < s_2 < s_0 + s_2 < s_1 + s_2 < 2s_2$ 

\end{center}

\

For (2a):

\begin{center}

$s_0 < 2s_0 < s_1 < s_2 < s_0 + s_1 < s_0 + s_2 < 2s_1 < s_1 + s_2 < 2s_2$ 

$s_0 < 2s_0 < s_1 < s_2 < s_0 + s_1 < 2s_1 < s_0 + s_2 < s_1 + s_2 < 2s_2$

\end{center}

\

For (2b):

\begin{center}

$s_0 < 2s_0 < s_1 < s_0 + s_1 < s_2 < s_0 + s_2 < 2s_1 < s_1 + s_2 < 2s_2$ 

$s_0 < 2s_0 < s_1 < s_0 + s_1 < s_2 < 2s_1 < s_0 + s_2 < s_1 + s_2 < 2s_2$ 

\end{center}

\

For (2c):

\begin{center}

$s_0 < 2s_0 < s_1 < s_0 + s_1 < 2s_1 < s_2 < s_0 + s_2 < s_1 + s_2 < 2s_2$ 

\end{center}

\

We now insert $s_3$ in these chains and check if the $4$-values condition holds for all the corresponding sets.

\subsection{$s_0 < s_1 < s_2 < 2s_0 < s_0 + s_1 < s_0 + s_2 < 2s_1 < s_1 + s_2 < 2s_2 \ \ \{ 5, 7, 8\}$.}

\
\vspace{1em}

\subsubsection{$s_2 < s_3 \leqslant 2s_0 \ \ \{ 5, 7, 8, 11\}$.}
\

No metric restriction. $S$ satisfies the $4$-values condition. 

\

\subsubsection{$2s_0 < s_3 \leqslant s_0 + s_1 \ \ \{ 5, 7, 8, 11\}$.} 

\label{5,7,8,11}

\[
A \subset [0, 6], \ \ B \subset [10, +\infty[.
\]

No bad quadruple. 
$S$ satisfies the $4$-values condition. 

\
\vspace{1em} 

\subsubsection{$s_0 + s_1 < s_3 \leqslant s_0 + s_2 \ \ \{ 5, 7, 8, 13\}$.} 

\label{5,7,8,13}

\[
A \subset [0, 8], \ \ B \subset [10, +\infty[.
\]

No bad quadruple. 
$S$ satisfies the $4$-values condition. 

\
\vspace{1em}  

\subsubsection{$s_0 + s_2 < s_3 \leqslant 2 s_1 \ \ \{ 5, 7, 8, 14\}$.}
\

$(5, 14, 5, 7)$ is a bad quadruple while $(5, 14, 5, 7)^* = (5, 5, 14, 7)$ is not. 
$S$ does not satisfy the $4$-values condition. 

\
\vspace{1em}  

\subsubsection{$2s_1 < s_3 \leqslant s_1 + s_2 \ \ \{ 5, 7, 8, 15\}$.} 
\

\label{5,7,8,15}

$(5, 15, 5, 7)$ is a bad quadruple while $(5, 15, 5, 7)^* = (5, 5, 15, 7)$ is not. 
$S$ does not satisfy the $4$-values condition.

\
\vspace{1em}   

\subsubsection{$s_1 + s_2 < s_3 \leqslant 2 s_2 \ \ \{ 5, 7, 8, 16\}$.} 
\

\label{5,7,8,16}

$(7, 16, 7, 8)$ is a bad quadruple while $(7, 16, 7, 8)^* = (7, 7, 16, 8)$ is not. 
$S$ does not satisfy the $4$-values condition.

\
\vspace{1em}   

\subsubsection{$2 s_2 < s_3 \ \ \{ 5, 7, 8, 17\}$.}
\

\label{5,7,8,17}

$S = S' \cup \{ t \}$ where $S'$ satisfies the $4$-values condition and $2 \max S' < t$. It is easy to check that the $4$-values condition is always satisfied in such a situation.

\subsection{$s_0 < s_1 < s_2 < 2s_0 < s_0 + s_1 < 2s_1 < s_0 + s_2 < s_1 + s_2 < 2s_2 \ \ \{ 5, 6, 9\}$.}

\
\vspace{1em}   

\subsubsection{$s_2 < s_3 \leqslant 2s_0 \ \ \{ 5, 6, 9, 10\}$.}
\

No metric restriction. $S$ satisfies the $4$-values condition.

\
\vspace{1em}   

\subsubsection{$2s_0 < s_3 \leqslant s_0 + s_1 \ \ \{ 5, 6, 9, 11\}$.}
\

$s_2$ does not appear in any non-metric triangle with labels in $S$. 
$4$-values condition is satisfied.

\
\vspace{1em}   

\subsubsection{$s_0 + s_1 < s_3 \leqslant 2s_1 \ \ \{ 5, 6, 9, 12\}$.}
\

Same as previous case. 
$4$-values condition is satisfied.

\
\vspace{1em}   

\subsubsection{$2s_1 < s_3 \leqslant s_0 + s_2 \ \ \{ 5, 6, 9, 14\}$.}
\

Same as previous case. 
$4$-values condition is satisfied.

\
\vspace{1em}   

\subsubsection{$s_0 + s_2 < s_3 \leqslant s_1 + s_2 \ \ \{ 5, 6, 9, 15\}$.}
\

$\{ 5, 6, 9, 15\} \sim \{ 5, 7, 8, 15 \}$. So according to \ref{5,7,8,15}, $S$ does not satisfy the $4$-values condition.

\
\vspace{1em}   

\subsubsection{$s_1 + s_2 < s_3 \leqslant 2s_2 \ \ \{ 5, 6, 9, 18\}$.}
\
$\{ 5, 6, 9, 18\} \sim \{ 5, 7, 8, 16 \}$. So according to \ref{5,7,8,16}, $S$ does not satisfy the $4$-values condition.

\
\vspace{1em}   

\subsubsection{$2s_2 < s_3 \ \ \{ 5, 6, 9, 19\}$.}
\

$\{ 5, 6, 9, 19\} \sim \{ 5, 7, 8, 17 \}$. So according to \ref{5,7,8,17}, $S$ satisfies the $4$-values condition.

\subsection{$s_0 < s_1 < 2s_0 < s_2 < s_0 + s_1 < s_0 + s_2 < 2s_1 < s_1 + s_2 < 2s_2 \ \ \{ 4, 7, 9\}$.}

\
\vspace{1em}   

\subsubsection{$s_2 < s_3 \leqslant s_0 + s_1 \ \ \{ 4, 7, 9, 11\}$.}
\

$s_1$ does not appear in any non-metric triangle with labels in $S$. 
$4$-values condition is satisfied.

\
\vspace{1em}   

\subsubsection{$s_0 + s_1 < s_3 \leqslant s_0 + s_2 \ \ \{ 4, 7, 9, 12\}$.}
\

\label{4,7,9,12}

$\{ 4, 7, 9, 13\} \approx \{ 1, 2, 3, 4\}$, and $4$-values condition is satisfied as $\{ 1, 2, 3, 4\}$ is an initial segment of a set which is closed under sums.  

\newpage

\subsubsection{$s_0 + s_2 < s_3 \leqslant 2s_1 \ \ \{ 4, 7, 9, 14\}$.}
\

$(4, 14, 4, 7)$ is a bad quadruple while $(4, 14, 4, 7)^* = (4, 4, 14, 7)$ is not. 
$S$ does not satisfy the $4$-values condition.

\
\vspace{1em}  

\subsubsection{$2s_1 < s_3 \leqslant s_1 + s_2 \ \ \{ 4, 7, 9, 16\}$.}
\

\label{4,7,9,16}

$(4, 16, 4, 7)$ is a bad quadruple while $(4, 16, 4, 7)^* = (4, 4, 16, 7)$ is not. 
$S$ does not satisfy the $4$-values condition.

\
\vspace{1em}   

\subsubsection{$s_1 + s_2 < s_3 \leqslant 2s_2 \ \ \{ 4, 7, 9, 18\}$.}
\

\label{4,7,9,18}

$(7, 18, 4, 9)$ is a bad quadruple while $(7, 18, 4, 9)^* = (7, 4, 18, 9)$ is not. 
$S$ does not satisfy the $4$-values condition.

\
\vspace{1em}   

\subsubsection{$2s_2 < s_3 \ \ \{ 4, 7, 9, 19\}$.}
\

$4$-values condition is satisfied as $S = S' \cup \{ t \}$ with $S'$ satisfying the $4$-values condition and $2\max S' < t$.

\subsection{$s_0 < s_1 < 2s_0 < s_2 < s_0 + s_1 < 2s_1 < s_0 + s_2 < s_1 + s_2 < 2s_2 \ \ \{ 8, 14, 21 \}$.}

\
\vspace{1em}   

\subsubsection{$s_2 < s_3 \leqslant s_0 + s_1 \ \ \{ 8, 14, 21, 22 \}$.}
\

$s_1$ does not appear in any non-metric triangle with labels in $S$. $4$-values condition is satisfied.

\
\vspace{1em}   

\subsubsection{$s_0 + s_1 < s_3 \leqslant 2s_1 \ \ \{ 8, 14, 21, 28 \}$.} 
\

$\{ 8, 14, 21, 28\} \sim \{ 4, 7, 9, 12 \}$. Thus, according to \ref{4,7,9,12}, $S$ satisfies the $4$-values condition.

\
\vspace{1em}   

\subsubsection{$2s_1 < s_3 \leqslant s_0 + s_2 \ \ \{ 8, 14, 21, 29\}$.}
\

$(14, 29, 8, 8)$ is a bad quadruple while $(14, 29, 8, 8)^* = (14, 8, 29, 8)$ is not. 
$S$ does not satisfy the $4$-values condition.

\
\vspace{1em}   

\subsubsection{$s_0 + s_2 < s_3 \leqslant s_1 + s_2 \ \ \{ 8, 14, 21, 35\}$.}
\

$\{ 8, 14, 21, 35\} \sim \{ 4, 7, 9, 16 \}$. Thus, according to \ref{4,7,9,16}, $S$ does not satisfy the $4$-values condition.

\
\vspace{1em}   

\subsubsection{$s_1 + s_2 < s_3 \leqslant 2s_2 \ \ \{ 8, 14, 21, 42\}$.}
\

$\{ 8, 14, 21, 42\} \sim \{ 4, 7, 9, 18 \}$. According to \ref{4,7,9,18}, $S$ consequently does not satisfy the $4$-values condition.

\
\vspace{1em}   

\subsubsection{$2s_2 < s_3 \ \ \{ 8, 14, 21, 43 \}$.}
\
 
$4$-values condition is satisfied as $S = S' \cup \{ t \}$ with $S'$ satisfying the $4$-values condition and $2\max S' < t$. 

\subsection{$s_0 < s_1 < 2s_0 < s_0 + s_1 < 2s_1 < s_2 < s_0 + s_2 < s_1 + s_2 < 2s_2 \ \ \{ 2, 3, 7\}$.} 

\
\vspace{0.5em}
\subsubsection{$s_2 < s_3 \leqslant s_0 + s_2 \ \ \{ 2, 3, 7, 9 \}$.}

\[ A = \{ 1, 2, 4, 5, 6, 7\}, \ \ B = \{ 4, 5, 6 \} \cup C, \ \ C \subset [9, +\infty[. \]
\[
\begin{array}{lll}
[4, 4] & (3, 7, 2, 2) & (3, 7, 2, 2)_* = (3, 2, 2, 7) \approx (2, 7, 2, 3)\\

[4, 5] & (3, 7, 2, 3) & (3, 7, 2, 3)_* = (3, 3, 2, 7) \approx (2, 7, 3, 3)\\

[4, 6] & (3, 7, 3, 3) & \\

[5, 4] & (2, 7, 2, 2) & \\

[5, 5] & (2, 7, 2, 3) & (2, 7, 2, 3)^* = (2, 2, 7, 3) \approx (3, 7, 2, 2)\\

[5, 6] & (2, 7, 3, 3) & (2, 7, 3, 3)^* = (2, 3, 7, 3) \approx (3, 7, 2, 3)\\

[6, 4] & (3, 9, 2, 2) & (3, 9, 2, 2)^* = (3, 2, 9, 2) \approx (2, 9, 2, 3)\\

[6, 5] & (3, 9, 2, 3) & (3, 9, 2, 3)_* = (3, 3, 2, 9) \approx (2, 9, 3, 3)\\  

[6, 6] & (3, 9, 3, 3) & \\

[7, 4] & (2, 9, 2, 2) & \\

[7, 5] & (2, 9, 2, 3) & (2, 9, 2, 3)^* = (2, 2, 9, 3) \approx (3, 9, 2, 2)\\  

[7, 6] & (2, 9, 3, 3) & (2, 9, 3, 3)^* = (2, 3, 9, 3) \approx (3, 9, 2, 3)

\end{array}
\]

$S = \{ 2, 3, 7, 9\}$  satisfies the $4$-values condition.

\
\vspace{0.5em}   
 
\subsubsection{$s_0 + s_2 < s_3 \leqslant s_1 + s_2 \ \ \{ 2, 3, 7, 10\}$.}
\

$(2, 10, 2, 7)$ is a bad quadruple while $(2, 10, 2, 7)^* = (2, 2, 10, 7)$ is not. 
$S$ does not satisfy the $4$-values condition.

\
\vspace{0.5em}   

\subsubsection{$s_1 + s_2 < s_3 \leqslant 2s_2 \ \ \{ 2, 3, 7, 14\}$.}

\[
A = \{ 1, 4, 5, 7, 11, 12\}, \ \ B = \{ 4, 5, 6, 9, 10 \} \cup C, \ \ C \subset [14, +\infty[. 
\]
\[
\begin{array}{lll}

[4, 4] & (3, 7, 2, 2) & (3, 7, 2, 2)^* = (3, 2, 7, 2) \approx (2, 7, 2, 3)\\ 

[4, 5] & (3, 7, 2, 3) & (3, 7, 2, 3)_* = (3, 3, 2, 7) \approx (2, 7, 3, 3)\\

[4, 6] & (3, 7, 3, 3) & \\

[5, 4] & (2, 7, 2, 2) & \\

[5, 5] & (2, 7, 2, 3) & (2, 7, 2, 3)^* = (2, 2, 7, 3) \approx (3, 7, 2, 2)\\

[5, 6] & (2, 7, 3, 3) & (2, 7, 3, 3)^* = (2, 3, 7, 3) \approx (3, 7, 2, 3)\\ 

[7, 4] & (7, 14, 2, 2) & (7, 14, 2, 2)^* = (7, 2, 14, 2) \approx (2, 14, 2, 7)\\

[7, 5] & (7, 14, 2, 3) & (7, 14, 2, 3)^* = (7, 2, 14, 3) \approx (3, 14, 2, 7)\\

& & (7, 14, 2, 3)_* = (7, 3, 2, 14) \approx (2, 14, 3, 7)\\

[7, 6] & (7, 14, 3, 3) & (7, 14, 3, 3)^* = (7, 3, 14, 3) \approx (3, 14, 3, 7)\\   

[11, 4]& (3, 14, 2, 2) & (3, 14, 2, 2)^* = (3, 2, 14, 2) \approx (2, 14, 2, 3)\\

[11, 5]& (3, 14, 2, 3) & (3, 14, 2, 3)_* = (3, 3, 2, 14) \approx (2, 14, 3, 3)\\

[11, 6]& (3, 14, 3, 3) & \\

[11, 9]& (3, 14, 2, 7) & (3, 14, 2, 7)^* = (3, 2, 14, 7) \approx (7, 14, 2, 3)\\

& & (3, 14, 2, 7)_* = (3, 7, 2, 14) \approx (2, 14, 3, 7) \\

[11, 10]& (3, 14, 3, 7) & (3, 14, 3, 7)^* = (3, 3, 14, 7) \approx (7, 14, 3, 3)\\

[12, 4]& (2, 14, 2, 2) & \\

[12, 5]& (2, 14, 2, 3) & (2, 14, 2, 3)^* = (2, 2, 14, 3) \approx (3, 14, 2, 2)\\

[12, 6]& (2, 14, 3, 3) & (2, 14, 3, 3)^* = (2, 3, 14, 3) \approx (3, 14, 2, 3)\\

[12, 9]& (2, 14, 2, 7) & (2, 14, 2, 7)^* = (2, 2, 14, 7) \approx (7, 14, 2, 2)\\

[12, 10]& (2, 14, 3, 7) & (2, 14, 3, 7)^* = (2, 3, 14, 7) \approx (7, 14, 2, 3)\\

& & (2, 14, 3, 7)_* = (2, 7, 3, 14) \approx (3, 14, 2, 7)

\end{array}
\] 

$S = \{ 2, 3, 7, 14\}$  satisfies the $4$-values condition.

\
\vspace{1em}   

\subsubsection{$2s_2 < s_3 \ \ \{ 2, 3, 7, 15\}$.} 
\
$4$-values condition is satisfied as $S = S' \cup \{ t \}$ with $S'$ satisfying the $4$-values condition and $2\max S' < t$.  

\subsection{$s_0 < 2s_0 < s_1 < s_2 < s_0 + s_1 < s_0 + s_2 < 2s_1 < s_1 + s_2 < 2s_2 \ \ \{ 2, 6, 7\}$.}

\
\vspace{1em}   

\subsubsection{$s_2 < s_3 \leqslant s_0 + s_1 \ \ \{ 2, 6, 7, 8\}$.}
\[
A = \{ 1, 2, 4, 5, 6 \}, B = \{ 4 \} \cup C, C \subset [8, +\infty[.  
\]
\[
\begin{array}{ll}

[4, 4] & (2, 6, 2, 2) \\

[5, 4] & (2, 7, 2, 2) \\

[6, 4] & (2, 8, 2, 2) 

\end{array}
\] 

$S = \{ 2, 6, 7, 8\}$  satisfies the $4$-values condition.

\
\vspace{1em}   

\subsubsection{$s_0 + s_1 < s_3 \leqslant s_0 + s_2 \ \ \{ 2, 6, 7, 9\}$.} 
\

$(6, 9, 2, 2)$ is a bad quadruple while $(6, 9, 2, 2)^* = (6, 2, 9, 2)$ is not. 
$S$ does not satisfy the $4$-values condition.

\
\vspace{1em}   

\subsubsection{$s_0 + s_2 < s_3 \leqslant 2s_1 \ \ \{ 2, 6, 7, 12\}$.}
\[
A = \{ 1, 4, 5, 6, 10\}, \ \ B = \{ 4, 8, 9 \} \cup C, \ \  C \subset [12, +\infty[.  
\]
\[
\begin{array}{lll}

[4, 4] & (2, 6, 2, 2) & \\

[5, 4] & (2, 7, 2, 2) & \\

& (7, 12, 2, 2) & (7, 12, 2, 2)^* = (7, 2, 12, 2) \approx (2, 12, 2, 7) \\

[6, 4] & (2, 8, 2, 2) & \\

& (6, 12, 2, 2) & (6, 12, 2, 2)^* = (6, 2, 12, 2) \approx (2, 12, 2, 6) \\ 

[10, 4] & (2, 12, 2, 2) & \\

[10, 8] & (2, 12, 2, 6) & (2, 12, 2, 6)^* = (2, 2, 12, 6) \approx (6, 12, 2, 2) \\

[10, 9] & (2, 12, 2, 7) & (2, 12, 2, 7)^* = (2, 2, 12, 7) \approx (7, 12, 2, 2)

\end{array}
\] 

$S = \{ 2, 6, 7, 12\}$  satisfies the $4$-values condition.

\
\vspace{1em}   

\subsubsection{$2s_1 < s_3 \leqslant s_1 + s_2 \ \ \{ 2, 6, 7, 13\}$.}
\

$(2, 13, 6, 6)$ is a bad quadruple while $(2, 13, 6, 6)^* = (2, 6, 13, 6)$ is not. 
$S$ does not satisfy the $4$-values condition.

\
\vspace{1em}   

\subsubsection{$s_1 + s_2 < s_3 \leqslant 2s_2 \ \ \{ 2, 6, 7, 14\}$.}
\

$(6, 14, 2, 7)$ is a bad quadruple while $(6, 14, 2, 7)^* = (6, 2, 14, 7)$ is not. 
$S$ does not satisfy the $4$-values condition.

\
\vspace{1em}   

\subsubsection{$2s_2 < s_3 \ \ \{ 2, 6, 7, 15\}$.}
\

$4$-values condition is satisfied as $S = S' \cup \{ t \}$ with $S'$ satisfying the $4$-values condition and $2\max S' < t$.

\subsection{$s_0 < 2s_0 < s_1 < s_2 < s_0 + s_1 < 2s_1 < s_0 + s_2 < s_1 + s_2 < 2s_2$.}
\

This chain of inequalities is not consistent: If $s_2 \leqslant s_0 + s_1$ and $2s_1 \leqslant s_0 + s_2$ then $s_1 \leqslant 2s_0$. 

\subsection{$s_0 < 2s_0 < s_1 < s_0 + s_1 < s_2 < s_0 + s_2 < 2s_1 < s_1 + s_2 < 2s_2 \ \ \{ 1, 4, 6\}$.}

\
\vspace{1em}   

\subsubsection{$s_2 < s_3 \leqslant s_0 + s_2 \ \ \{ 1, 4, 6, 7\}$.}

\label{1,4,6,7}

\[
A = \{ 1, 2, 3, 5, 6 \}, B = \{ 2, 5 \} \cup C, \ \ C \subset [7, +\infty[.  
\]
\[
\begin{array}{lll}

[2, 2] & (4, 6, 1, 1) & (4, 6, 1, 1)^* = (4, 1, 6, 1) \approx (1, 6, 1, 4)\\

[3, 2] & (4, 7, 1, 1) & (4, 7, 1, 1)^* = (4, 1, 7, 1) \approx (1, 7, 1, 4)\\

 & (1, 4, 1, 1) & \\

[5, 2] & (1, 6, 1, 1) & \\

[5, 5] & (1, 6, 1, 4) & (1, 6, 1, 4)^* = (1, 1, 6, 4) \approx (4, 6, 1, 1) \\

[6, 2] & (1, 7, 1, 1) & \\

[6, 5] & (1, 7, 1, 4) & (1, 7, 1, 4)^* = (1, 1, 7, 4) \approx (4, 7, 1, 1)

\end{array}
\] 

$S = \{ 1, 4, 6, 7 \}$  satisfies the $4$-values condition.

\
\vspace{1em}   

\subsubsection{$s_0 + s_2 < s_3 \leqslant 2s_1 \ \ \{ 1, 4, 6, 8\}$.}
\[
A = \{ 2, 3, 4, 5, 7 \}, \ \ B = \{ 2, 5, 7 \} \cup C, \ \ C \subset [8, +\infty[.  
\]
\[
\begin{array}{lll}

[2, 2] & (4, 6, 1, 1) & (4, 6, 1, 1)^* = (4, 1, 6, 1) \approx (1, 6, 1, 4)\\

& (6, 8, 1, 1) & (6, 8, 1, 1)^* = (6, 1, 8, 1) \approx (1, 8, 1, 6) \\

[3, 2] & (1, 4, 1, 1) & \\

[4, 2] & (4, 8, 1, 1) & (4, 8, 1, 1)^* = (4, 1, 8, 1) \approx (1, 8, 1, 4) \\

[5, 2] & (1, 6, 1, 1) & \\

[5, 5] & (1, 6, 1, 4) & (1, 6, 1, 4)^* = (1, 1, 6, 4) \approx (4, 6, 1, 1) \\

[7, 2] & (1, 8, 1, 1) & \\

[7, 5] & (1, 8, 1, 4) & (1, 8, 1, 4)^* = (1, 1, 8, 4) \approx (4, 8, 1, 1) \\

[7, 7] & (1, 8, 1, 6) & (1, 8, 1, 6)^* = (1, 1, 8, 6) \approx (6, 8, 1, 1)

\end{array}
\]  

$S = \{ 1, 4, 6, 8 \}$  satisfies the $4$-values condition.

\
\vspace{1em}   

\subsubsection{$2s_1 < s_3 \leqslant s_1 + s_2 \ \ \{ 1, 4, 6, 10\}$.} 

\label{1,4,6,10}

\[
A = \{ 2, 3, 4, 5, 6, 9 \}, \ \ B = \{ 2, 5, 7, 8 \} \cup C, \ \ C \subset [10, +\infty[.  
\]
\[
\begin{array}{lll}

[2, 2] & (4, 6, 1, 1) & (4, 6, 1, 1)^* = (4, 1, 6, 1) \approx (1, 6, 1, 4)\\

[3, 2] & (1, 4, 1, 1) & \\

[4, 2] & (6, 10, 1, 1) & (6, 10, 1, 1)^* = (6, 1, 10, 1) \approx (1, 10, 1, 6) \\

[5, 2] & (1, 6, 1, 1) & \\

[5, 5] & (1, 6, 1, 4) & (1, 6, 1, 4)^* = (1, 1, 6, 4) \approx (4, 6, 1, 1) \\

[6, 2] & (4, 10, 1, 1) & (4, 10, 1, 1)^* = (4, 1, 10, 1) \approx (1, 10, 1, 4) \\

[6, 5] & (4, 10, 1, 4) & (4, 10, 1, 4)_* = (4, 4, 1, 10) \approx (1, 10, 4, 4) \\

[9, 2] & (1, 10, 1, 1) & \\

[9, 5] & (1, 10, 1, 4) & (1, 10, 1, 4)^* = (1, 1, 10, 4) \approx (4, 10, 1, 1) \\

[9, 7] & (1, 10, 1, 6) & (1, 10, 1, 6)^* = (1, 1, 10, 6) \approx (6, 10, 1, 1) \\

[9, 8] & (1, 10, 4, 4) & (1, 10, 4, 4)^* = (1, 4, 10, 4) \approx (4, 10, 1, 4)

\end{array}
\]  

$S = \{ 1, 4, 6, 10 \}$  satisfies the $4$-values condition.

\
\vspace{1em}   

\subsubsection{$s_1 + s_2 < s_3 \leqslant 2s_2 \ \ \{ 1, 4, 6, 12 \}$.}
\

$(4, 12, 4, 6)$ is a bad quadruple while $(4, 12, 4, 6)^* = (4, 4, 12, 6)$ is not. 
$S$ does not satisfy the $4$-values condition.

\
\vspace{1em}   

\subsubsection{$2s_2 < s_3 \ \ \{ 1, 4, 6, 13\}$.}
\

$4$-values condition is satisfied as $S = S' \cup \{ t \}$ with $S'$ satisfying the $4$-values condition and $2\max S' < t$.

\subsection{$s_0 < 2s_0 < s_1 < s_0 + s_1 < s_2 < 2s_1 < s_0 + s_2 < s_1 + s_2 < 2s_2 \ \ \{ 2, 5, 9\}$.}

\
\vspace{1em}   

\subsubsection{$s_2 < s_3 \leqslant 2s_1 \ \ \{ 2, 5, 9, 10\}$.} 
\

$\{ 2, 5, 9, 10\} \sim \{ 1, 4, 6, 7\}$. Thus, according to \ref{1,4,6,7}, $S$ satisfies the $4$-values condition.











\
\vspace{1em}   

\subsubsection{$2s_1 < s_2 \leqslant s_0 + s_2 \ \ \{ 2, 5, 9, 11\}$.}
\

$(5, 11, 2, 5)$ is a bad quadruple while $(5, 11, 2, 5)_* = (5, 5, 2, 11)$ is not. 
$S$ does not satisfy the $4$-values condition.

\
\vspace{1em}   

\subsubsection{$s_0 + s_2 < s_3 \leqslant s_1 + s_2 \ \ \{ 2, 5, 9, 14\} $.}
\

$\{ 2, 5, 9, 14\} \sim \{ 1, 4, 6, 10 \}$ so according to \ref{1,4,6,10}, $S$ satisfies the $4$-values condition.

\
\vspace{1em}   

\subsubsection{$s_1 + s_2 < s_3 \leqslant 2s_2 \ \ \{ 2, 5, 9, 18\}$.}
\

$(5, 18, 5, 9)$ is a bad quadruple while $(5, 18, 5, 9)^* = (5, 5, 18, 9)$ is not. 
$S$ does not satisfy the $4$-values condition.

\
\vspace{1em}   

\subsubsection{$2s_2 < s_3$.}
\

$4$-values condition is satisfied as $S = S' \cup \{ t \}$ with $S'$ satisfying the $4$-values condition and $2\max S' < t$.      

\subsection{$s_0 < 2s_0 < s_1 < s_0 + s_1 < 2s_1 < s_2 < s_0 + s_2 < s_1 + s_2 < 2s_2 \ \ \{ 1, 3, 7\}$.}

\
\vspace{1em}   

\subsubsection{$s_2 < s_3 \leqslant s_0 + s_2 \ \ \{ 1, 3, 7, 8\}$.}

\[
A = \{ 1, 2, 4, 5, 6, 7 \}, \ \ B = \{ 2, 4, 6\} \cup C, \ \ C \subset [8, +\infty[.  
\]
\[
\begin{array}{lll}

[2, 2] & (1, 3, 1, 1) & \\

[4, 2] & (3, 7, 1, 1) & (3, 7, 1, 1)^* = (3, 1, 7, 1) \approx (1, 7, 1, 3) \\

[4, 4] & (3, 7, 1, 3) & (3, 7, 1, 3)_* = (3, 3, 1, 7) \approx (1, 7, 3, 3) \\

[4, 6] & (3, 7, 3, 3) & \\

[5, 2] & (3, 8, 1, 1) & (3, 8, 1, 1)^* = (3, 1, 8, 1) \approx (1, 8, 1, 3) \\

[5, 4] & (3, 8, 1, 3) & (3, 8, 1, 3)_* = (3, 3, 1, 8) \approx (1, 8, 3, 3) \\

[5, 6] & (3, 8, 3, 3) & \\

[6, 2] & (1, 7, 1, 1) & \\

[6, 4] & (1, 7, 1, 3) & (1, 7, 1, 3)^* = (1, 1, 7, 3) \approx (3, 7, 1, 1) \\

[6, 6] & (1, 7, 3, 3) & (1, 7, 3, 3)^* = (1, 3, 7, 3) \approx (3, 7, 1, 3) \\  

[7, 2] & (1, 8, 1, 1)\\

[7, 4] & (1, 8, 1, 3) & (1, 8, 1, 3)^* = (1, 1, 8, 3) \approx (3, 8, 1, 1) \\

[7, 6] & (1, 8, 3, 3) & (1, 8, 3, 3)^* = (1, 3, 8, 3) \approx (3, 8, 1, 3) \\  

\end{array}
\]  

$S = \{ 1, 3, 7, 8 \}$  satisfies the $4$-values condition.

\
\vspace{1em}   

\subsubsection{$s_0 + s_2 < s_3 \leqslant s_1 + s_2 \ \ \{ 1, 3, 7, 10\}$.}

\[
A = \{ 2, 3, 4, 6, 7, 9 \}, \ \ B = \{ 2, 4, 6, 8 \} \cup C, \ \ C \subset [10, +\infty[.  
\]
\[
\begin{array}{lll}

[2, 2] & (1, 3, 1, 1) & \\

[3, 2] & (7, 10, 1, 1) & (7, 10, 1, 1)^* = (7, 1, 10, 1) \approx (1, 10, 1, 7) \\ 

[4, 2] & (3, 7, 1, 1) & (3, 7, 1, 1)^* = (3, 1, 7, 1) \approx (1, 7, 1, 3) \\

[4, 4] & (3, 7, 1, 3) & (3, 7, 1, 3)_* = (3, 3, 1, 7) \approx (1, 7, 3, 3) \\

[4, 6] & (3, 7, 3, 3) & \\

[6, 2] & (1, 7, 1, 1) &  \\

[6, 4] & (1, 7, 1, 3) & (1, 7, 1, 3)^* = (1, 1, 7, 3) \approx (3, 7, 1, 1) \\

[6, 6] & (1, 7, 3, 3) & (1, 7, 3, 3)^* = (1, 3, 7, 3) \approx (3, 7, 1, 3) \\

[7, 2] & (3, 10, 1, 1) & (3, 10, 1, 1)^* = (3, 1, 10, 1) \approx (1, 10, 1, 3)  \\

[7, 4] & (3, 10, 1, 3) & (3, 10, 1, 3)_* = (3, 3, 1, 10) \approx (1, 10, 3, 3) \\

[7, 6] & (3, 10, 3, 3) & \\  

[9, 2] & (1, 10, 1, 1)\\

[9, 4] & (1, 10, 1, 3) & (1, 10, 1, 3)^* = (1, 1, 10, 3) \approx (3, 10, 1, 1) \\

[9, 6] & (1, 10, 3, 3) & (1, 10, 3, 3)^* = (1, 3, 10, 3) \approx (3, 10, 1, 3) \\  

[9, 8] & (1, 10, 1, 7) & (1, 10, 1, 7)^* = (1, 1, 10, 7) \approx (7, 10, 1, 1)

\end{array}
\]   

$S = \{ 1, 3, 7, 10 \}$  satisfies the $4$-values condition.

\
\vspace{1em}   

\subsubsection{$ s_1 + s_2 < s_3 \leqslant 2s_2 \ \ \{ 1, 3, 7, 14\}$.}
\[
A = \{ 2, 4, 6, 7, 11, 13 \}, \ \ B = \{ 2, 4, 6, 8, 10 \} \cup C, \ \ C \subset [14, +\infty[.  
\]
\[
\begin{array}{lll}

[2, 2] & (1, 3, 1, 1) & \\

[4, 2] & (3, 7, 1, 1) & (3, 7, 1, 1)^* = (3, 1, 7, 1) \approx (1, 7, 1, 3) \\

[4, 4] & (3, 7, 1, 3) & (3, 7, 1, 3)_* = (3, 3, 1, 7) \approx (1, 7, 3, 3) \\

[4, 6] & (3, 7, 3, 3) & \\

[6, 2] & (1, 7, 1, 1) &  \\

[6, 4] & (1, 7, 1, 3) & (1, 7, 1, 3)^* = (1, 1, 7, 3) \approx (3, 7, 1, 1) \\

[6, 6] & (1, 7, 3, 3) & (1, 7, 3, 3)^* = (1, 3, 7, 3) \approx (3, 7, 1, 3) \\

[7, 2] & (7, 14, 1, 1) & (7, 14, 1, 1)^* = (7, 1, 14, 1) \approx (1, 14, 1, 7)  \\

[7, 4] & (7, 14, 1, 3) & (7, 14, 1, 3)^* = (7, 1, 14, 3) \approx (3, 14, 1, 7) \\

& & (7, 14, 1, 3)_* = (7, 3, 1, 14) \approx (1, 14, 3, 7) \\

[7, 6] & (7, 14, 3, 3) & (7, 14, 3, 3)^* = (7, 3, 14, 3) \approx (3, 14, 3, 7) \\  

[11, 2] & (3, 14, 1, 1) & (3, 14, 1, 1)^* = (3, 1, 14, 1) \approx (1, 14, 1, 3)  \\

[11, 4] & (3, 14, 1, 3) & (3, 14, 1, 3)_* = (3, 3, 1, 14) \approx (1, 14, 3, 3) \\

[11, 6] & (3, 14, 3, 3) & \\  

[11, 8] & (3, 14, 1, 7) & (3, 14, 1, 7)^* = (3, 1, 14, 7) \approx (7, 14, 1, 3) \\

& & (3, 14, 1, 7)_* = (3, 7, 1, 14) \approx (1, 14, 3, 7) \\

[11, 10] & (3, 14, 3, 7) & (3, 14, 3, 7)^* = (3, 3, 14, 7) \approx (7, 14, 3, 3) \\ 

[13, 2] & (1, 14, 1, 1) & \\

[13, 4] & (1, 14, 1, 3) & (1, 14, 1, 3)^* = (1, 1, 14, 3) \approx (3, 14, 1, 1) \\   

[13, 6] & (1, 14, 3, 3) & (1, 14, 3, 3)^* = (1, 3, 14, 3) \approx (3, 14, 1, 3) \\  

[13, 8] & (1, 14, 1, 7) & (1, 14, 1, 7)^* = (1, 1, 14, 7) \approx (7, 14, 1, 1) \\    

[13, 10] & (1, 14, 3, 7) & (1, 14, 3, 7)^* = (1, 3, 14, 7) \approx (7, 14, 1, 3) \\ 

 & & (1, 14, 3, 7)_* = (1, 7, 3, 14) \approx (3, 14, 1, 7) 

\end{array}
\]   

$S = \{ 1, 3, 7, 14 \}$  satisfies the $4$-values condition.

\
\vspace{1em}   

\subsubsection{$2s_2 < s_3 \ \ \{ 1, 3, 7, 15\} $.}
\

$4$-values condition is satisfied as $S = S' \cup \{ t \}$ with $S'$ satisfying the $4$-values condition and $2\max S' < t$.      






\chapter*{Appendix B. \ \ Indivisibility of $\Ur _S$ when $|S| \leqslant 4$.}

The purpose of this Appendix is to show that for $|S| = 4$ and satisfying the $4$-values condition, the space $\Ur _S$ is indivisible. The main ingredients of the proofs are indivisibility of $\Ur _S$ when $|S| \leqslant 3$, Milliken's theorem (theorem \ref{thm:Milliken}) and Sauer's theorem (theorem \ref{thm:Sauer ultrahomogeneous}). In what follows, the numbering of the cases corresponds to the sections in Appendix A. 

\

\textbf{2.1.1.} $\{ 5, 7, 8, 10\}$

$\Ur _S$ can be seen as a complete version of the Rado graph with four kinds of edges. An easy variation of the proof working for the Rado graph shows that this space is indivisible. 

\

\textbf{2.1.2.} $\{ 5, 7, 8, 11\}$ 

$8$ does not appear in any non-metric triangle with labels in $S$. Thus, $\Ur _S$ is indivisible thanks to Sauer's theorem. 

\

\textbf{2.1.3.} $\{ 5, 7, 8, 13\}$ 

Same as previous case. 

\

\textbf{2.1.7.} $\{ 5, 7, 8, 17\}$ 

$\Ur _S$ is composed of countably many disjoint copies of $\Ur _{\{ 5, 7, 8\}}$ and the distance between any two points not in the same copy of $\Ur _{\{ 5, 7, 8\}}$ is always $17$. The indivisibility of $\Ur _{\{ 5, 7, 8\}}$ consequently implies that $\Ur _S$ is indivisible. 

\

\textbf{2.2.1.} $\{ 5, 6, 9, 10\}$  

$\{ 5, 6, 9, 10\} \sim \{ 5, 7, 8, 10\}$, so $\Ur _S$ is isomorphic to the space in $2.1.1$ and hence indivisible. 

\

\textbf{2.2.2.} $\{ 5, 6, 9, 11\}$  

$9$ does not appear in any non-metric triangle with labels in $S$. Thus, $\Ur _S$ is indivisible thanks to Sauer's theorem. 

\

\textbf{2.2.3.} $\{ 5, 6, 9, 12\}$

Same as previous case. 

\

\textbf{2.2.4.} $\{ 5, 6, 9, 13\}$    

Same as previous case.  

\

\textbf{2.2.7.} $\{ 5, 6, 9, 19\}$  

$\{ 5, 6, 9, 19\} \sim \{ 5, 7, 8, 17\}$, so $\Ur _S$ is isomorphic to the space in $2.1.7$ and hence indivisible. 

\

\textbf{2.3.1.} $\{ 4, 7, 9, 11\}$  

$7$ does not appear in any non-metric triangle with labels in $S$. Thus, $\Ur _S$ is indivisible thanks to Sauer's theorem. 

\

\textbf{2.3.2.} $\{ 4, 7, 9, 13\}$  

$\{ 4, 7, 9, 13\} \sim \{ 1, 2, 3, 4\}$ so essentially, $\Ur _S$ is $\Ur _4$. Thus, $\Ur _S$ is indivisible. 
\

\textbf{2.3.6.} $\{ 4, 7, 9, 19\}$

$\Ur _S$ is composed of countably many disjoint copies of $\Ur _{\{ 4, 7, 9\}}$ and the distance between any two points not in the same copy of $\Ur _{\{ 4, 7, 9\}}$ is always $19$. The indivisibility of $\Ur _{\{ 4, 7, 9 \}}$ consequently implies that $\Ur _S$ is indivisible. 

\

\textbf{2.4.1.} $\{ 8, 14, 21, 22\}$

$14$ does not appear in any non-metric triangle with labels in $S$. Thus, $\Ur _S$ is indivisible thanks to Sauer's theorem. 

\

\textbf{2.4.2.} $\{ 8, 14, 21, 28\}$

Elements in $\M _S$ are isomorphic to elements in $\M _{S'}$ with $S'$ as in 2.3.2. This case is consequently equivalent to indivisibility of $\Ur _4$ and $\Ur _S$ is indivisible.  

\

\textbf{2.4.6.} $\{ 8, 14, 21, 43\}$

$\Ur _S$ is composed of countably many disjoint copies of $\Ur _{\{ 8, 14, 21 \}}$ and the distance between any two points not in the same copy of $\Ur _{\{ 8, 14, 21 \}}$ is always $43$. The indivisibility of $\Ur _{\{ 8, 14, 21 \}}$ consequently implies that $\Ur _S$ is indivisible.  

\

\textbf{2.5.1.} $\{ 2, 3, 7, 9\}$

The proof of indivisibility for $\Ur _S$ is a simple adaptation of the proof of indivisibility of $\Ur _{\{ 1, 3, 4\}}$: 
Fix an $\omega$-linear ordering $<$ on $2 ^{< \omega}$ extending the tree ordering and consider the following graph structure on $2^{< \omega}$: 

\begin{center}
$\forall s < t \in 2^{< \omega} \ \ \{ s , t \} \in E \leftrightarrow \left( |s| < |t| , t(|s|) = 1 \right)$. 
\end{center}

Now, define $d$ on the set $[2 ^{< \omega}]^2$ of pairs of $2 ^{< \omega}$ as follows: Let $\{ s, t \}_< , \{ s' , t' \}_< $ be in $[2 ^{< \omega}]^2$. Then $d (\{ s, t \}_< , \{ s' , t' \}_< )$ is:

\begin{displaymath}
\left \{ \begin{array}{ll}
 2 & \textrm{if $s=s'$ and $\{ t, t'\} \in E$.}\\
 3 & \textrm{if $s=s'$ and $\{ t, t' \} \notin E$.}\\
 7 & \textrm{if $s \neq s'$ and $\{ t, t'\} \in E$.}\\
 9 & \textrm{if $s \neq s'$ and $\{ t, t' \} \notin E$.}
 \end{array} \right.
\end{displaymath}

One can check that $d$ is a metric. Since $d$ takes its values in $\{ 2, 3, 7, 9\}$, $([2 ^{< \omega}]^2 , d)$ embeds into $\Ur _S$. We now show that $\Ur _S$ embeds into the subspace $\m{X}$ of $([2 ^{< \omega}]^2 , d)$ supported by the set 

\begin{center}
$ X = \{ \{ s , t \}_< \in [2 ^{< \omega}]^2 : |s| < |t|, \ s <_{lex} t, \ t(|s|) = 0 \}$.
\end{center} 

The embedding is constructed inductively. Let $\{ x_n : n \in \omega \}$ be an enumeration of $\Ur _S$. We are going to construct a sequence $(\{ s_n , t_n \})_{n \in \omega}$ of elements in $X$ such that 

\begin{center} 
$\forall m, n \in \omega \ \ d(\{ s, t \}_< , \{ s' , t' \} _<) = d^{\Ur _S } (x_m , x_n)$. 
\end{center}

For $\{ s_0 , t_0 \} _<$, take $s_0 = \emptyset$ and $t_0 = 0$. Assume now that $\{ s_0 , t_0 \} _< ,\ldots , \{ s_n , t_n \} _<$ are constructed such that all the elements of $\{ s_0 ,\ldots , s_n\} \cup \{ t_0 ,\ldots , t_n \}$ have different heights and all the $s_i$'s are strings of $0$'s. Set 

\begin{center}
$M = \{ m \leqslant n : d^{\Ur _S} (x_m , x_{n+1}) \in \{2, 3 \} \}$. 
\end{center}

If $M = \emptyset $, choose $s_{n+1}$ to be a string of $0$'s longer that all the elements constructed so far. Otherwise, there is $s \in 2^{< \omega}$ such that 

\begin{center} 
$\forall m \in M \ \ s_m = s$. 
\end{center}

Set $s_{n+1} = s$. Now, choose $t_{n+1}$ above all the elements constructed so far and such that 

\vspace{0.5em}
\hspace{1em}
i) $\forall m \in M \ \ (t_{n+1} (|t_m|) = 1) \leftrightarrow (d^{\Ur _S }(x_{n+1} , x_m) = 2)$.

\vspace{0.5em}
\hspace{1em}
ii) $\forall m \notin M \ \ (t_{n+1} (|t_m|) = 1) \leftrightarrow (d^{\Ur _S }(x_{n+1} , x_m) = 7)$.

\vspace{0.5em}
\hspace{1em}
iii) $\{ s_{n+1} , t_{n+1}\}_< \in X$. 

\vspace{0.5em} 

i) and ii) are easy to satisfy because all the $t_m$'s have different heights. As for iii), $|s_{n+1}| < |t_{n+1}|$ and $ t_{n+1} (|s_{n+1}|) = 0$ are also easy (again because all heights are different) while $s_{n+1} <_{lex} t_{n+1}$ is satisfied because $s_{n+1}$ being a $0$ string, $|s_{n+1}| < |t_{n+1}|$ implies $s_{n+1} <_{lex} t_{n+1}$. After $\omega$ steps, we are left with $\{ \{ s_n , t_n \} : n \in \omega \} \subset \m{X}$ isometric to $\Ur _S$. Observe that actually, this construction shows that $\Ur _S$ embeds into any subspace of $([2 ^{< \omega}]^2 , d)$ supported by a strong subtree of $2^{< \omega}$. 

Now, to prove that $\Ur _S$ is indivisible, it suffices to prove that given any $\chi : \funct{([2 ^{< \omega}]^2 , d)}{k}$ where $k \in \omega$ is strictly positive, there is a strong subtree $\m{T}$ of $2^{< \omega}$ such that $\chi$ is constant on $[T]^2 \cap X$. But this is guaranteed by Milliken theorem: Indeed, consider the subset $A := \{ 0, 01\}$. Then using the notation introduced for theorem \ref{thm:Milliken}, $[A]_{\mathrm{Em}} = X$. So the restriction $\restrict{\chi}{[A]_{\mathrm{Em}}}$ is really a coloring of $X$, and there is a strong subtree $\m{T}$ of height $\omega$ such that $\restrict{[A]_{\mathrm{Em}}}{T} = [T]^2 \cap X$ is $\chi$-monochromatic.

\

\textbf{2.5.3.} $\{ 2, 3, 7, 14\}$ 

$\Ur _S$ is obtained from $\Ur _2$ by multiplying the distances by $7$ and then blowing up the points to copies of $\Ur _{\{ 2, 3\}}$. $\Ur _2$ and $\Ur _{\{ 2, 3\}}$ being indivisible, so is $\Ur _S$.

\

\textbf{2.5.4.} $\{ 2, 3, 7, 15\}$

$\Ur _S$ is composed of countably many disjoint copies of $\Ur _{\{ 2, 3, 7\}}$ and the distance between any two points not in the same copy of $\Ur _{\{ 2, 3, 7\}}$ is always $15$. The indivisibility of $\Ur _{\{ 2, 3, 7 \}}$ consequently implies that $\Ur _S$ is indivisible. 

\

\textbf{2.6.1.} $\{ 2, 6, 7, 8\}$  

In this case, indivisibility of $\Ur _S$ can be proved thanks to the method of 2.5.1. except that instead of $[2 ^{< \omega}]^2$, one works with $[3 ^{< \omega}]^2$ and $d (\{ s, t \}_< , \{ s' , t' \}_< )$ defined on the set $[3 ^{< \omega}]^2$ of pairs of $3 ^{< \omega}$ by:

\begin{displaymath}
\left \{ \begin{array}{ll}
 2 & \textrm{if $s=s'$}\\
 6 & \textrm{if $s \neq s'$ and $t'(|t|) = 0$.}\\
 7 & \textrm{if $s \neq s'$ and $t'(|t|) = 1$.}\\
 8 & \textrm{if $s \neq s'$ and $t'(|t|) = 2$.}
 \end{array} \right.
\end{displaymath}

\newpage

\textbf{2.6.3.} $\{ 2, 6, 7, 12\}$   

Again, we apply Milliken's theorem. Consider $E$ the standard graph structure on $2^{< \omega}$ and define $d (\{ s, t \}_< , \{ s' , t' \}_< )$ by:

\begin{displaymath}
\left \{ \begin{array}{ll}
 2 & \textrm{if $s=s'$ and $\{ t, t'\} \in E$.}\\
 6 & \textrm{if $s \neq s'$ and $\{ s, s' \} \notin E$ and $\{ t, t' \} \notin E$.}\\
 7 & \textrm{if $s \neq s'$ and $\{ s, s' \} \notin E$ and $\{ t, t'\} \in E$.}\\
 12 & \textrm{if $s \neq s'$ and $\{ s, s' \} \in E$.}
 \end{array} \right.
\end{displaymath}

Then one can check that $d$ is a metric on $[2^{< \omega}]^2$ and that $([2^{< \omega}]^2, d)$ and $\Ur _S$ embed into each other. Milliken's theorem provides indivisibility.

\

\textbf{2.6.6.} $\{ 2, 6, 7, 15\}$  

$\Ur _S$ is composed of countably many disjoint copies of $\Ur _{\{ 2, 6, 7\}}$ and the distance between any two points not in the same copy of $\Ur _{\{ 2, 6, 7\}}$ is always $15$. The indivisibility of $\Ur _{\{ 2, 6, 7\}}$ consequently implies that $\Ur _S$ is indivisible. 

\

\textbf{2.8.1.} $\{ 1, 4, 6, 7\}$  

Let $f : \funct{\{ 1, 4, 6, 7\}}{\{ 2, 6, 7, 12\}}$ be such that $f(1) = 2$, $f(4)=7$, $f(6)=6$ and $f(7) = 12$. Then observe that $f$ establishes an isomorphism between $\Ur _S$ and $\Ur _{\{ 2, 6, 7, 12\}}$ (case 2.6.3). $\Ur _{\{ 2, 6, 7, 12\}}$ being indivisible, so is $\Ur _S$. 

\

\textbf{2.8.2.} $\{ 1, 4, 6, 8\}$  

$\Ur _S$ is obtained from $\Ur _{\{4, 6, 8 \}}$ after having blown the points up to copies of $\Ur _1$. Its indivisibility is a direct consequence of the basic infinite pigeonhole principle and of the indivisibility of $\Ur _{\{4, 6, 8 \}}$.

\

\textbf{2.8.3.} $\{ 1, 4, 6, 10\}$  

$\Ur _S$ is obtained from $\Ur _{\{4, 6, 10 \}}$ after having blown the points up to copies of $\Ur _1$. Its indivisibility is a direct consequence of the basic infinite pigeonhole principle and of the indivisibility of $\Ur _{\{4, 6, 10 \}}$.

\

\textbf{2.8.5.} $\{ 1, 4, 6, 13\}$  

$\Ur _S$ is composed of countably many disjoint copies of $\Ur _{\{ 1, 4, 6\}}$ and the distance between any two points not in the same copy of $\Ur _{\{ 1, 4, 6\}}$ is always $13$. The indivisibility of $\Ur _{\{ 1, 4, 6 \}}$ consequently implies that $\Ur _S$ is indivisible. 

\

\textbf{2.9.1.} $\{ 2, 5, 9, 10\}$  

$\{ 2, 5, 9, 10\} \sim \{ 1, 4, 6, 7\}$, so $\Ur _S$ is isomorphic to the space in $2.8.1$ and is indivisible. 

\

\textbf{2.9.3.} $\{ 2, 5, 9, 14\}$  

$\{ 5, 6, 9, 14\} \sim \{ 1, 4, 6, 10 \}$, so $\Ur _S$ is isomorphic to the space in $2.8.3$ and is indivisible. 

\

\textbf{2.9.5.} $\{ 2, 5, 9, 19\}$  

$\Ur _S$ is composed of countably many disjoint copies of $\Ur _{\{ 2, 5, 9\}}$ and the distance between any two points not in the same copy of $\Ur _{\{ 2, 5, 9\}}$ is always $19$. The indivisibility of $\Ur _{\{ 2, 5, 9\}}$ consequently implies that $\Ur _S$ is indivisible. 

\

\textbf{2.10.1.} $\{ 1, 3, 7, 8\}$  

This case is another instance where Milliken's theorem is useful. Consider $E$ the standard graph structure on $2^{< \omega}$ and define $d (\{ s, t, u \}_< , \{ s' , t', u' \}_< )$ by:

\begin{displaymath}
\left \{ \begin{array}{ll}
 1 & \textrm{if $s=s'$ and $t=t'$.}\\
 3 & \textrm{if $s=s'$ and $t \neq t'$.}\\
 7 & \textrm{if $s \neq s'$ and $\{ u, u' \} \in E$.}\\
 8 & \textrm{if $s \neq s'$ and $\{ u, u' \} \notin E$.}
 \end{array} \right.
\end{displaymath}

Then one can check that $d$ is a metric on $[2^{< \omega}]^3$. $([2^{< \omega}]^3, d)$ embeds into $\Ur _S$ because $d$ takes values in $S$. Conversely, given any strong subtree $T$ of $2^{<\omega}$, $\Ur _S$ embeds into $[T]^3 \cap Y$ where $Y \subset [2 ^{< \omega}]^3$ given by all the triples $\{ s, t, u \}_<$ such that

\begin{displaymath}
\left \{ \begin{array}{l}
|s| < |t| < |u|\\
s <_{lex} t <_{lex} u\\
t(|s|) = u(|s|) = u(|t|) = 0
 \end{array} \right.
\end{displaymath} 

Equivalently, $Y = [B]_{\mathrm{Em}}$ with $B = \{ 0, 10, 110 \}$. These facts allow to apply Milliken's theorem and to deduce indivisibility of $\Ur _S$.

\

\textbf{2.10.2.} $\{ 1, 3, 7, 10\}$  

$\Ur _S$ is obtained from $\Ur _{\{3, 7, 10 \}}$ after having blown the points up to copies of $\Ur _1$. Its indivisibility is a direct consequence of the basic infinite pigeonhole principle and of the indivisibility of $\Ur _{\{3, 7, 10 \}}$.   

\

\textbf{2.10.3.} $\{ 1, 3, 7, 14\}$  

$\Ur _S$ is obtained from $\Ur _{\{3, 7, 14 \}}$ after having blown the points up to copies of $\Ur _1$. Its indivisibility is a direct consequence of the basic infinite pigeonhole principle and of the indivisibility of $\Ur _{\{3, 7, 14\}}$.   

\

\textbf{2.10.4.} $\{ 1, 3, 7, 15\}$ 

$\Ur _S$ is ultrametric with four distances, hence indivisible.





\chapter*{Appendix C. \ \ On the universal Urysohn space $\Ur$.}

The purpose of this appendix is to provide some additional information about the Urysohn space $\Ur$\index{$\Ur$}. As already mentionned, $\Ur$ was originally constructed by P. Urysohn in 1925 in order to show that there is a separable metric space into which every separable metric space embeds isometrically. In the original paper, $\Ur$ was obtained as the completion of $\Ur _{\Q}$ which was constructed by hand and inductively. Here are the main features of $\Ur$ as presented in \cite{U} but using our terminology:

\begin{thm}[Urysohn]

\index{Urysohn!theorem}

\

\begin{enumerate}

\item For every finite subspace $\m{X} \subset \Ur$ and every Kat\v{e}tov map $f$ over $\m{X}$, there is $x \in \Ur$ realizing $f$ over $\m{X}$. 

\item Every separable metric space embeds isometrically into $\Ur$.  

\item $\Ur$ is ultrahomogeneous. 
 
\item $\Ur$ is the unique complete separable metric space satisfying (2) and (3). 

\item $\Ur$ is path connected and locally path connected.

\item $\Ur$ includes two isometric subspaces $\m{X}$ and $\m{Y}$ such that no isometry from $\Ur$ onto itself maps $\m{X}$ onto $\m{Y}$.   

\end{enumerate}

\end{thm}

Some 30 years later, in \cite{Hu}, Huhunai\v{s}vili improved the result (3) about ultrahomogeneity:

\begin{thm}[Huhunai\v{s}vili]

\index{Huhunai\v{s}vili theorem}
Let $\varphi : \funct{\m{X}}{\m{Y}}$ be a bijective isometry between two compact subspaces of $\Ur$. Then $\varphi$ can be extended to an isometry of $\Ur$ onto itself. 
\end{thm}

However, together with an article by Sierpinski \cite{Si}, Huhunai\v{s}vili's contribution represents the only study about $\Ur$ between 1927 and 1986 (There is an article in 1971 by Joiner but the main result is only the rediscovery of a subcase covered by Huhunai\v{s}vili's theorem). In 1986, Kat\v{e}tov provided in \cite{K} the construction of $\Ur _{\Q}$ presented in Chapter 1 (Note that what we called here Kat\v{e}tov functions were undoubtedly introduced and used earlier. For example, they already appear in some work by Isbell in 1964, see \cite{I}, or later in 1974 and 1984 in some articles of Flood, cf \cite{Fl1} and \cite{Fl2}). Thanks to the work of Uspenskij, this new approach became the starting point of a new period of interest for $\Ur$. Today, research about $\Ur$ and the topological group $\iso (\Ur)$ of its surjective isometries (equipped with the pointwise convergence topology) is well alive, as illustrated by the workshop organized recently in Be'er Sheva (May 2006). In what follows, we present a short selection of the main results from the last 20 years. For a more detailed presentation, the reader should refer to \cite{GK}, \cite{Pe1}, \cite{Pe1'}, or to the original papers. Another source of reference is also \cite{Pr}, the proceedings volume of the aforementioned workshop which appeared in \emph{Topology and its applications}.

We start with a result which completes the work carried out by Urysohn and Huhunai\v{s}vili about ultrahomogeneity. It is quite surprising that after having remained unsolved for such a long time, it was obtained recently, independently and simultaneously by two persons. 

\begin{thm}[Ben Ami \cite{BA}, Melleray \cite{Me3}]

\index{Ben Ami theorem}
\index{Melleray!theorem on compact ultrahomogeneity}

Let $\m{X}$ be a Polish metric space. TFAE:

\vspace{0.5em}
\hspace{1em}
i) $\m{X}$ is compact.  

\vspace{0.5em}
\hspace{1em}
ii) If $\m{X} _0$ and $\m{X} _1$ are isometric copies of $\m{X}$ inside $\Ur$ and $\varphi : \funct{\m{X} _0}{\m{X} _1}$, 

\hspace{2em} then $\varphi$ can be extended to an isometry of $\Ur$ onto itself. 

\end{thm}

Here are two other theorems about the intrinsic geometry of $\Ur$: 

\begin{thm}[Melleray, \cite{Me3}]
\index{Melleray!first theorem on fixed points of isometries in $\Ur$}
Let $\varphi \in \iso (\Ur)$ whose orbits have compact closure. Then the set of fixed points of $\varphi$ is either empty or isometric to $\Ur$. 
\end{thm}

\begin{thm}[Melleray, \cite{Me3}]
\index{Melleray!second theorem on fixed points of isometries in $\Ur$}
Let $\m{X}$ be a Polish metric space. Then there is $\varphi$ in $\iso (\Ur)$ whose set of fixed points in $\Ur$ is isometric to $\Ur$.  
\end{thm}

Next, we present the structures which are supported by $\Ur$. We start with the topological characterization of $\Ur$:

\begin{thm}[Uspenskij \cite{Us2}]
\index{Uspenskij!theorem on the topology of $\Ur$}
$\Ur$ is homeomorphic to $\ell _2$. 
\end{thm}

Next, recall that a group is \emph{monothetic} if it contains a dense subgroup isomorphic to the additive group of the integers $\mathbb{Z}$. 

\begin{thm}[Cameron-Vershik \cite{CV}]
\index{Cameron-Vershik theorem}
$\Ur$ admits the structure of a monothetic Polish group. 
\end{thm}

This result has to be compared with the following one, due to Holmes: 

\begin{thm}[Holmes \cite{Ho}]
\index{Holmes theorem}
When $\Ur$ is embedded isometrically into a Banach space with a fixed point $x_0$ sent to the zero element of the Banach space, any finite subset of the copy of $\Ur$ which does not contain $x_0$ is linearly independent and the closed linear span of the copy of $\Ur$ is uniquely determined up to linear isometry.   
\end{thm}

It follows that $\Ur$ does \emph{not} support the structure of Banach space. Indeed, calling $\left\langle \Ur \right\rangle$ the Banach space provided by the previous theorem, $\left\langle \Ur \right\rangle$ cannot have $\Ur$ as underlying set: Otherwise, $\left\langle \Ur \right\rangle$ would be an ultrahomogeneous Banach space but we mentioned in Chapter 1 that the only ultrahomogeneous Banach space is $\ell _2$. $\left\langle \Ur \right\rangle$ is a wild object but is better understood today in the context of so-called \emph{Lipschitz-free spaces}. For example, a recent theorem from Godefroy and Kalton \cite{GoKa} allows to show that every separable Banach space embeds linearly and isometrically into $\left\langle \Ur \right\rangle$. However, many basic questions about $\left\langle \Ur \right\rangle$ remain unanswered. For example, does that space admit a basis? Nevertheless, $\left\langle \Ur \right\rangle$ turned out to be helpful in the resolution of certain problems, as in \cite{Me4} where it allowed to reach a result about the complexity of the isometry relationship between separable Banach spaces. 

We finish our first list of properties related to $\Ur$ by a theorem due to Vershik \cite{Ve1}. We wrote in the introduction that in some cases, Fra\"iss\'e limits can be seen as random objects. $\Ur$ is only the completion of a Fra\"iss\'e limit but a result of very similar flavor seems to hold. We state it following Pestov (\cite{Pe1}, p.143):

\begin{thm}[Vershik]
\index{Vershik!theorem on the randomness of $\Ur$}
Let $M$ be the set of all metrics on $\omega$ and let $\mathbb{P} (M)$ be the set of all probability measures on $M$. Then, for a generic $\mu \in \mathbb{P} (M)$, the completion of $(\omega , d)$ is isometric to $\Ur$ $\mu$-almost surely in $d \in M$. 
\end{thm} 

We now turn to properties related to $\iso (\Ur)$, starting with the following theorem due to Uspenskij:

\begin{thm}[Uspenskij \cite{Us1}]
\index{Uspenskij!theorem on universality of $\iso (\Ur)$}
Every second countable topological group is isomorphic to a topological subgroup of $\iso (\Ur)$. 
\end{thm}

In fact, more can be said: 

\begin{thm}[Melleray \cite{Me1}]
\index{Melleray!theorem on representations of Polish groups}
For every Polish group $G$, there is a closed subspace $\m{X}$ of $\Ur$ such that $G \cong \{ \varphi \in \iso (\Ur) : \varphi '' \m{X} = \m{X} \}$.
\end{thm}

On the other hand, there are also some informations about the actions of $\iso (\Ur)$:

\begin{thm}[Pestov \cite{Pe0}]
\index{Pestov!theorem on $\iso (\Ur)$}
Every continuous action of $\iso (\Ur)$ on a compact space admits a fixed point. 
\end{thm}

As mentioned several times in the body of the present paper, this result is particularly important for our present work because it can be proved via combinatorial methods. However, we should emphasize that in fact, $\iso (\Ur)$ satisfies a stronger property called the \emph{L\'evy property} and which implies the previous theorem, see \cite{Pe1} or \cite{Pe2}. 

Other problems concerning $\iso (\Ur)$ can be attacked via combinatorics. For example, the following result announced by Vershik \cite{Ve2} and proved independently by Solecki \cite{So} can be seen as a metric version of the well-known result about the extension of partial isomorphisms of finite graphs due to Hrushovski \cite{Hr}.  

\begin{thm}[Solecki \cite{So}, Vershik \cite{Ve2}]
\index{Solecki!theorem on extension of isometries}
\index{Vershik!theorem on extension of isometries}

Let $\m{X}$ be a finite metric space. Then there is a finite metric space $\m{Y}$ such that $\m{X} \subset \m{Y}$ and such that every isometry $\varphi$ with $\dom (\varphi) , \ran (\varphi) \subset \m{X}$ of $\m{X}$ extends to an isometry of $\m{Y}$ onto itself. 
\end{thm}

The importance of this result is related to the following concepts. For a Polish group $G$ and $n \in \omega$, the \emph{diagonal action of $G$ on $G^n$} is the action defined by: 

\begin{center}
$g \cdot (h_1 ,\ldots , h_n) = (gh_1 g^{-1},\ldots , gh_n g^{-1})$. 
\end{center}

An element $(h_1 ,\ldots , h_n)$ of $G^n$ is \emph{cyclically dense} if for some $g \in G$, the set $\{ g^k \cdot (h_1 ,\ldots , h_n)  : k \in \omega \}$ is dense in $G^n$.  

\begin{thm}[Solecki \cite{So}]
\index{Solecki!theorem on diagonal actions of $\iso (\Ur)$}
All the diagonal actions of $\iso (\Ur)$ have cyclically dense elements. 
\end{thm}

\begin{thm}[Solecki \cite{So}]
\index{Solecki!theorem on the generation of $\iso (\Ur)$}
There are two elements of $\iso (\Ur)$ generating a dense subgroup. 
\end{thm}

The last result we finish with comes from \cite{KR} and provides a so-called \emph{reconstruction theorem}. The core of the proof is again related to metric combinatorics and extension properties in the Urysohn space. However, it seems to us that this result deserves a particular attention because while most of the previous results deal with isometries, this one concerns a broader class of maps: For metric spaces $\m{X}$ and $\m{Y}$, call a homeomorphism $g : \funct{\m{X}}{\m{Y}}$ \emph{locally bi-Lipschitz} if every $x \in \m{X}$ has a neighborhood $U$ such that $\restrict{g}{U}$ is bi-Lipschitz. Let $L (\m{X})$ denotes the set of all bi-Lipschitz homeomorphisms of $\m{X}$, then:

\begin{thm}[Kubi\'s-Rubin]
\index{Kubi\'s-Rubin theorem}
Let $\m{X}$ and $\m{Y}$ be open subspaces of $\Ur$. Suppose that $\varphi$ is a group isomorphism between $L (\m{X})$ and $L (\m{Y})$. Then there is a locally bi-Lipschitz homeomorphism $\tau$ between $\m{X}$ and $\m{Y}$ such that:

\begin{center}
$\forall g \in L (\m{X}) \ \ \varphi (g) = \tau \circ g \circ \tau ^{-1}$. 
\end{center}

\end{thm}




\printindex

\end{document}